\documentclass{conm-p-l}

\usepackage{amsmath,amsthm,amsfonts,amssymb}
\usepackage{graphicx}

\usepackage{amssymb}
\usepackage{psfig}
\usepackage{amsmath}

\newtheorem{maintheorem}{Theorem}[]
\newtheorem{theorem}{Theorem}[section]
\newtheorem{lemma}[theorem]{Lemma}
\newtheorem{prop}[theorem]{Proposition}
\newtheorem{cy}[theorem]{Corollary}

\theoremstyle{definition}
\newtheorem{definition}[theorem]{Definition}
\newtheorem{example}[theorem]{Example}

\newtheorem{notation}[theorem]{Notation}
\theoremstyle{remark}
\newtheorem{remark}[]{Remark}

\newcommand{\be}{\begin{enumerate}}
\newcommand{\ee}{\end{enumerate}}
\newcommand{\beq}{\begin{equation}}
\newcommand{\eeq}{\end{equation}}
\newcommand{\bi}{\begin{itemize}}
\newcommand{\ei}{\end{itemize}}

\begin{document}
\setcounter{page}{87}

\title{Effective JSJ Decompositions}
\author[O. Kharlampovich]{Olga Kharlampovich}
\address{Department of Mathematics and Statistics, McGill University, Montreal, QC, Canada, H3A2K6}

\email{olga@math.mcgill.ca}
\thanks{The first author was supported
 by NSERC Grant. }

\author[A. G.  Myasnikov]{Alexei G. Myasnikov}
\address{Department of Mathematics and Statistics, McGill University, Montreal, QC, Canada, H3A2K6}
\email{alexeim@att.net}
\thanks{The second author was supported by NSERC Grant and by NSF
Grant DMS-9970618}

\pagestyle{myheadings} \markright{{\sl $\bullet$ Effective JSJ,
version 46 $\bullet$ 06.24.04}}

\date{Version 46, June 24, 2004}



\begin{abstract}
In this paper we describe an elimination process which is a
deterministic rewriting procedure that on each elementary step
transforms one system of equations over free groups into a
finitely many new ones. Infinite branches of this process
correspond to  cyclic splittings of the coordinate group of the
initial system of equations. This allows us to construct
algorithmically Grushko's decompositions of finitely generated
fully residually free groups and cyclic [abelian] JSJ
decompositions of freely indecomposable finitely generated fully
residually free groups. We apply these results to obtain an
effective description of the set of homomorphisms from a given
finitely presented group into a free group, or, more generally,
into an NTQ group.
\end{abstract}

\maketitle

\tableofcontents

\section*{Introduction}

\subsection*{JSJ decompositions: a bit of history and intentions}
\label{se:intro-1}

 A splitting of a group $G$ is a presentation of
$G$ as the fundamental group of a graph of groups.  According to
Bass-Serre theory, a  group  has a non-trivial splitting if and
only if it acts non-trivially on a tree. A splitting is cyclic
[abelian] if all  edge groups are  cyclic [abelian]. In
\cite{shyp} Sela introduced  universal cyclic splittings of a
particular type, so called JSJ decompositions, which  contain (in
an appropriate sense) all other cyclic splittings of the group.
Such splittings are analogues to the characteristic submanifold
constructions for irreducible 3-manifolds described by Jaco and
Shalen \cite{JS} and Johannson \cite{Jo}, hence the name. Sela
showed that JSJ decompositions exist for  freely indecomposable
torsion-free hyperbolic groups. Existence of JSJ decompositions
for arbitrary finitely presented groups was shown by Rips and Sela
in \cite{RS}, their method is based on group actions on ${\mathbb
R}$-trees. A different way to obtain  JSJ splittings of one-ended
hyperbolic groups was proposed by Bowditch \cite{Bo}, he used the
local cut point structures on the boundaries. More general
approaches to JSJ decompositions  of groups have been described by
Dunwoody and Sageev \cite{DS} and by Fujiwara and Papasoglu
\cite{FP}. They use tracks on 2-complexes and  actions on products
of trees, correspondingly. All these results prove the existence
of JSJ decompositions in various classes of groups,  but do not
give an algorithm to construct them.

In this paper we algorithmically construct Grushko's
decompositions of finitely generated fully residually free groups
and cyclic [abelian] JSJ decompositions of freely indecomposable
finitely generated fully residually free groups. We apply these
results to obtain an effective description of the set of
homomorphisms from a given finitely presented group into a free
group, or, more generally, into an NTQ group.

 Non-trivial JSJ decompositions of a given group $G$
are very closely related with automorphisms of $G$. Originally,
they were used mostly to describe the automorphism group ${\rm
Aut}(G)$ of $G$ \cite{shyp}. It turned out recently, that various
types of JSJ splittings are very useful  in studying equations
over groups and, more generally, elementary theories of groups. In
particular,
 we show here that such splittings  provide a group theoretic counterpart to various
pieces of Makanin-Razborov process and give algebraic semantics to
our elimination rewriting process from \cite{KMIrc}. To explain
the main idea of our approach to JSJ we need some notation and
definitions.

\subsection*{Three approaches to fully residually free groups}
\label{se:intro-2}

Let $F=F(A)$ be a free group with basis  $A$. Fully residually
free groups, which are also known as freely discriminated  groups
\cite{BMR}, and $\omega$-residually free groups  \cite{Rem}, play
a crucial role in the theory of equations and model theory of free
groups.  Recall that a group $G$ is discriminated by $F$ (or
freely discriminated)  if for any finite subset $K \subseteq G$ of
non-trivial elements there exists a homomorphism $\phi: G
\rightarrow F$ such that $g^\phi \neq 1$ for each $g \in K$. We
describe  these groups in Subsection \ref{se:2-5} and discuss
their properties in Subsection \ref{se:2-6}. It is worthwhile to
mention here that finitely generated freely discriminated groups
can be viewed in three different ways each of which gives one
essentially different tools to deal with them. First, they appear
as coordinate groups of irreducible systems of equations over free
groups, so methods of algebraic geometry over groups can be used
here, as well as Makanin-Razborov's techniques; secondly, they can
be obtained from free groups by finitely many free products with
amalgamation and HNN extensions of a very particular type, hence
one can easily apply methods of Bass-Serre theory; thirdly, these
groups can be faithfully presented by infinite words over abelian
ordered groups $\mathbb{Z}^n$, which provides one with various
cancellation techniques, Nielsen method, and Stallings' foldings
apparatus to study their finitely generated subgroups. We
frequently use all these techniques below. Notice, that finitely
generated fully residually free groups have some other interesting
characterizations.  Remeslennikov showed  that a finitely
generated group $H$ is fully residually free if and only if $H$
has exactly the same universal theory as $F$ \cite{Rem}. Recently
Sela described  finitely generated  fully residually free groups
precisely as the limit groups \cite{Se}. Champetier and Guirardel
gave another characterization of these groups as limits of free
groups in a compact space of marked groups \cite{CG}.

Now we describe briefly the major components of our method.

\subsubsection*{Fully residually free groups as coordinate groups}

We start with a brief description of a few notions from algebraic
geometry over groups.  For a detailed  discussion on the subject
we refer to \cite{BMR} and \cite{KM9}, see also  Section
\ref{se:2-4}. By $V_F(S)$ we denote the set of all solutions of
the system $S(X,A) =1$  in $F$ (the {\em algebraic set defined by}
$S$). The algebraic set $V_G(S)$ uniquely corresponds to the
normal subgroup ({\it the radical of $S$}):
$$ R(S) = \{ T(X) \in F(A \cup X) \ \mid \ \forall A\in V_F(S)(T(A) = 1) \} $$
of the free group $F(A \cup X)$. The quotient group
$$F_{R(S)}=F(A \cup X)/R(S)$$ is the {\em coordinate group} of the
algebraic set  $V(S).$ Observe,  that if $V_F(S) \neq \emptyset$
then  $F$ is a subgroup of $F_{R(S)}$. It has been shown in
\cite{BMR} that algebraic sets $V_F(S_1)$ and $V_F(S_2)$ are
rationally equivalent if and only if there exists an isomorphism
between their coordinate groups which is identical on $F$. One can
define Zariski topology on $F^m$ ($m \in \mathbb{N}$) taking
algebraic sets as the closed subsets. Guba proved  in \cite{Gub}
that free groups are equationally Noetherian, i.e.,  this Zariski
topology is Noetherian, so every closed set is a finite union of
irreducible components. It turned out \cite{BMR} that  an
algebraic set is irreducible if and only if its coordinate group
is fully residually free.  Furthermore, every finitely generated
fully residually free group can be realized as the coordinate
group of a finite system of equations over $F$. From the group
theoretic view-point the elimination process mentioned above is
all about the coordinate groups of the systems involved. This
allows one to transform  pure combinatorial and algorithmic
results obtained in the process into statements about fully
residually free groups.

\subsubsection*{Fully residually free groups as fundamental groups of graphs of  groups}
In order to describe solutions of equations in free  groups Lyndon
introduced the notion of a group with parametric exponents in an
associative unitary ring $A$ \cite{Ly2}. In particular, he
described and studied   free exponential groups $F^{
\mathbb{Z}[t]}$ over a ring of polynomials $\mathbb{Z}[t]$. One of
the principle outcomes of his study is that the group $F^{
\mathbb{Z}[t]}$ is discriminated by $F$. A modern treatment of
exponential groups is due to
 Myasnikov and  Remeslennikov \cite{MyasExpo2}.
They showed, in particular,  that the group $ F^{Z[t]}$ can be
obtained  as union of an infinite chain, which starts at $F$,  of
HNN-extensions of a very specific type -  extensions of
centralizers (see Subsection \ref{se:2-5}).
 This implies that a  finitely
generated subgroup of $ F^{Z[t]}$ is  a subgroup of a group which
can be obtained from $F$  by finitely many extensions of
centralizers, so  one can apply  Bass-Serre theory to describe the
structure of these subgroups.

We proved in \cite{KMIrc} that finitely generated fully residually
 free groups are embeddable
into $F^{Z[t]}$. This, in view of the Lyndon's result mentioned
above,  gives a complete characterization of such groups. It
follows from Bass-Serre theory now that all these groups, except
for abelian ones, admit an essential cyclic splitting.

\subsubsection*{Fully residually free groups via infinite words}

Another key component  of our approach to  finitely generated
fully residually free groups is based on Lyndon's length functions
and infinite words. To axiomatize Nielsen cancellation argument
Lyndon introduced in \cite{Ly3} abstract length functions on
groups (now called Lyndon's length functions). He showed that a
group $G$ with a free Lyndon's length function with values in
$\mathbb{Z}$ (viewed as ordered abelian group) is embeddable into
a free group (represented by  finite reduced words with natural
length function) and the embedding preserves the length, hence $G$
is free.  It turned out that a similar result holds for finitely
generated fully residually free groups.  Namely, it has been shown
in \cite{MRS} that  elements of the free Lyndon's group
$F^{\mathbb{Z}[t]}$ can be represented by infinite words in the
alphabet $X^{\pm 1}$. These words are functions of the type
$$w: [1, \alpha_w] \rightarrow X^{\pm 1},$$
 where $\alpha_w \in \mathbb{Z}[t]$ and $[1,\alpha_w] = \{g \in \mathbb{Z}[t] \mid 1
\leqslant g \leqslant \alpha_w\}$ is a closed interval in
$\mathbb{Z}[t]$ with respect to the standard lexicographical order
$\leq$. The function $L: w \rightarrow \alpha_w$ gives rise to a
regular free Lyndon length function on $F^{\mathbb{Z}[t]}$ with
values in the additive group of $\mathbb{Z}[t]$, viewed as an
abelian  ordered group. This implies that every finitely generated
fully residually free group has a free length function with values
in  a free abelian group $\mathbb{Z}^n$  of  finite  rank with the
lexicographic order. Once a presentation of elements of
$F^{\mathbb{Z}[t]}$ by infinite words is established, a host of
problems about $F^{\mathbb{Z}[t]}$ can be solved  precisely in the
same way as in the standard free group $F$. In particular, if $H$
is a finitely generated subgroup of a finitely generated fully
residually free group $G$, then one can effectively embed $G$
(hence $H$) into $F^{\mathbb{Z}[t]}$. Now using the representation
above of elements of $G$ by infinite words one can algorithmically
construct a finite labelled graph $\Gamma_H$ (finite automata)
which accepts precisely elements of $H$ (given by their canonical
forms in $F^{\mathbb{Z}[t]}$) (see \cite{MRS2}). The process of
constructing the graph $\Gamma_H$ is a direct generalization of
the Stallings' folding procedure, but on infinite words.
 This allows one to treat finitely generated subgroups of fully
 residually free groups precisely in the same manner as in the
 standard free groups. This approach was realized in \cite{KMRS} where
 a host of algorithmic problems about finitely generated subgroups
 of finitely generated fully residually free groups was solved (see also Section
 \ref{algor2}).

\subsection*{Main ideas: elimination process vs JSJ}

Our main tool in constructing  JSJ decompositions of groups is  a
so-called elimination process  which is a symbolic rewriting
process of a certain type that transforms formal systems of
equations in groups or semigroups. This process can be viewed as a
non-commutative  analog of the classical elimination process  in
algebraic geometry, hence the name. The original version of the
process was introduced by Makanin in \cite{Mak82}. He showed that
the process gives a decision algorithm to verify whether a given
system is consistent or not (decidability of Diophantine problem
over free groups) by estimating the length of the minimal solution
(if it exists).  In \cite{Razborov} Razborov developed the process
 much further, so that the Makanin-Razborov's process
 produces {\em all}  solutions of a system of equations in a free group.
  We refined Razborov's version of the process in \cite{KMIrc} and \cite{Imp}
  to show that up to the rational equivalence irreducible algebraic sets over
  $F$ can be characterized  by non-degenerate triangular quasi-quadratic
  systems,  or NTQ systems (see Section \ref{ntq}), thus introducing NTQ groups
  which play an important role in model theory of free groups.
  Here the "non-degenerate"  part corresponds  to the extension
  theorems in the classical theory of elimination for polynomials.
In fact, our rewriting process from \cite{KMIrc} or \cite{Imp} by
no means is unique, every time it can be easily adjusted to some
particular needs. However, there exist several fundamental common
features that unify all variations of our's and Razborov's
processes. It is worthwhile to mention here two of them.  First,
only three precisely defined infinite branches (subprocesses) can
occur in the process: linear case (Cases 7--10), quadratic case
(Cases 11--12), and general JSJ case (Cases 13--15) which includes
periodic structures and abelian vertex groups (Case 2). Secondly,
groups of automorphisms of the coordinate groups are used in
encoding the infinite branches of the process. In what follows we
refer to processes of this  type as to elimination processes.
Observe, that the original Makanin's process lacks some of these
features (for example, there are no infinite branches
corresponding to the periodic structures).

    Notice that the principal ideas from
the  Makanin's process (elementary and entire transformations,
complexity, etc.)   were used in \cite{BF,Gab} to prove  the
classification theorem for finitely generated groups acting freely
on $\mathbb R$-trees and to describe  stable actions on $\mathbb
R$-trees, via so-called Rip's machine. Later, these results played
a key part in the proof of existence of  JSJ decompositions of
finitely presented groups with a single end \cite{shyp,RS}.

 In this paper we show that a modification of the  elimination process
 from \cite{KMIrc} can be used directly to construct effectively
  JSJ decompositions of finitely generated fully residually free
groups, thus avoiding algorithmically difficult detour into the
limit actions. Theorem \ref{th:spl} establishes a correspondence
between infinite branches  of an elimination process and
 splittings of the coordinate group of the system under
consideration:  the linear case corresponds to thin (or Levitt)
type of a subcomplex in \cite {BF} and \cite {Gab}, the quadratic
case corresponds to the surface type (or interval exchange),
periodic structures correspond to toral (or axial) type. Moreover,
the automorphism associated with infinite branches of the process
are precisely the canonical automorphism of the JSJ decomposition
associated with the splittings. Now we are in a position to
explain the main idea of the proof.

Let $G$ be a non-abelian finitely generated fully residually free
group given as the coordinate group of a finite system of
equations $S(A,X) = 1$ over $F$.  Every solution in $F$   of the
system  $S(X,A)=1$ gives rise to a homomorphism $\phi
:G\rightarrow F$. Composition of $\phi$ with an automorphism
$\sigma$ of $G$ provides  a new solution of the system $S(X,A)=1$
in $F$. Thus, different types of canonical automorphisms
associated with a JSJ decomposition of $G$ induce solutions of the
system $S(X,A) = 1$ of a particular type, that can be recognized
in  the elimination process as infinite branches of the
corresponding type. Conversely, infinite branches in the
elimination process for $S(X,A) = 1$ provide infinite families of
canonical automorphisms of $G$, that determine the corresponding
elementary splitting of $G$. This, in combination with Bass-Serre
methods and length functions techniques, allows one  to
reconstruct the JSJ decompositions of $G$.

\subsection*{Main results and organization of the paper}

 In the next section we collect some known results on
  fully residually free groups, algebraic
geometry over groups, and   quadratic equations.

In Section~\ref{se:split} we discuss splittings, JSJ
decompositions, and canonical automorphisms.   Here we prove an
important technical lemma on induced $QH$-subgroups (Lemma
\ref{le:1.6}).

Section \ref{algor} contains various results on decidability of
algorithmic problems in finitely generated fully residually free
groups. Some of them are known (complete proofs are in
\cite{KMIrc,KMRS,MRS}) and some are new.

In Sections \ref{se:ge}--\ref{subs:per}  we introduce and study an
elimination process.  Initially, we describe this process  as an
infinite branching deterministic  procedure which starts at a
generalized equation $\Omega$ and on each elementary step rewrites
a given generalized equation into a finitely many  new ones. The
procedure gives rise to an infinite directed tree $T(\Omega)$
which is equipped with some extra objects (enriched tree)
associated with vertices and edges. Namely,  the coordinate groups
of the generalized equations (obtained in the process) are
associated with vertices, while  homomorphisms of the
corresponding coordinate groups  (related to the elementary
rewriting steps) are associated with edges. As was mentioned above
infinite branches of this tree are the focal points of these
sections. The periodic structures provide the most difficult case,
they occupy the whole Section \ref{subs:per}.

In Sections \ref{se:5.3} and  \ref{se:5.5} we study splittings of
the coordinate group and the structure of solutions of a given
finite system of equations over free groups. To this end we cut
off the infinite branches of the tree $T(\Omega)$ and encode them
by groups of canonical automorphisms of the corresponding
coordinate groups. This results into  two finite enriched trees.
In one of them, $T_{{\rm dec}}(\Omega)$, every branch describes a
splitting of the coordinate group $F_{R(\Omega)}$ of the root
equation $\Omega$. In the other, $T_{{\rm sol}}(\Omega)$, every
vertex $v$ is associated with the coordinate group
$F_{R(\Omega_v)}$ of the corresponding system $\Omega_v$ and also
with a group of the canonical automorphisms $A_v$ of this
coordinate group, as before edges $u\rightarrow v$ of $T_{{\rm
sol}}(\Omega)$ are marked by homomorphisms $ F_{R(\Omega_u)}
\rightarrow F_{R(\Omega_v)}$. We refer to the  latter  graphs as
to Hom-diagrams.  The tree $T_{{\rm dec}}(\Omega)$ is the main
instrument to recognize all splittings essential to build JSJ
decompositions of the coordinate group $F_{R(\Omega)}$ (Theorem
\ref{th:onelevel}). The tree $T_{{\rm sol}}(\Omega)$ gives, in
particular,  a description of all solutions in $F$ of the initial
equation (Theorem \ref{th:5.3.1}).

In Section \ref{se:canonical} we prove

\medskip
\textsc{Theorem~\ref{th:suff}.} {\em Let $G$ be a finitely
generated freely indecomposable {\rm [}freely indecomposable
modulo $F${\rm ]}
 fully residually free group $G$. Then a non-degenerate JSJ
 $\mathbb Z$-decomposition {\rm [}modulo $F${\rm ]} of $G$  is a sufficient
splitting of
 $G$ {\rm (}see the definition in Subsection {\rm \ref{se:suffs}}{\rm )}.}

\medskip
In the same section we prove Theorem \ref{qq} about canonical
embeddings of finitely generated fully residually free groups into
NTQ groups. This result explains how the top quadratic system of
equations in an NTQ system assigned to the NTQ group corresponds
to a cyclic JSJ decomposition of the group that we embed.

In Section \ref{efd} we give an algorithm to find Grushko's
decompositions of finitely generated fully residually free groups.
 Recall that a free decomposition
 $$G=G_1 \ast \cdots \ast G_k \ast F_r$$
is a Grushko's decomposition of a group $G$ if the factors
$G_1,\dots ,G_k$ are freely indecomposable non-cyclic groups and
the factor $F_r$ is a free group of rank $r$. It is known that
such decompositions of $G$ are essentially unique, i.e., any other
decomposition of $G$ of this type has the same numbers $k$ and $r$
and the indecomposable non-cyclic factors are conjugated in $G$ to
one of the factors $G_1, \dots ,G_k$.  The next result shows that
Gruschko's decompositions (even modulo a given finitely generated
subgroup) of $G$ are effective.

\medskip
\textsc{Theorem~\ref{th:1.7}.} {\em There is an algorithm which
for every finitely generated fully residually free group $G$ and
its finitely generated subgroup $H$ determines whether or not $G$
has a nontrivial  free decomposition such that $H$ belongs to one
of the factors. Moreover, if $G$ does have such a decomposition,
the algorithm finds one {\rm(}by giving finite generating sets of
the factors{\rm)}.}

\medskip
This theorem allows one to reduce algorithmic problems for
finitely generated fully residually free groups to the  freely
indecomposable ones. For such groups we prove the following
theorem which is one of the major results of  this paper.

\begin{maintheorem} \label{th1} There exists  an algorithm to obtain a
 cyclic [abelian]
JSJ decomposition of a freely indecomposable finitely generated
fully residually free group. The algorithm constructs a
presentation of this group as the fundamental group of a JSJ graph
of groups.\end{maintheorem}

Actually,   we prove (in Section \ref{se:eff-JSJ}) the following
much more general result which implies Theorem \ref{th1}.

\medskip
\textsc{Theorem~\ref{modulo}.} {\em There exists an algorithm to
obtain a cyclic [abelian] $JSJ$ decomposition for a finitely
generated fully residually free group modulo a given finite family
of finitely generated subgroups.}
\medskip

 Combining Theorems~\ref{th:1.7}  and \ref{modulo} one  can
 find effectively  cyclic [abelian] JSJ
decompositions of the freely indecomposable non-abelian
factors\linebreak $G_1,\dots ,G_k$ of the Grushko's decomposition
of $G$. This gives a cyclic [abelian] splitting of $G$ which
completely describes the structure of $G$ as the fundamental group
of a graph of groups.

A variation  of Theorem \ref{th1} first appeared  in \cite{KM7}
(in terms of the elimination processes), where it was used in the
proof of the second Tarski's conjecture on the decidability of the
elementary theory of free non-abelian groups.

Theorem \ref{th1} has numerous applications.  Using Theorem
\ref{th1} we refine in Theorem \ref{th:hom-machine} the effective
description of a solution set of a  system of equations over $F$
from \cite{KMIrc}. In Section~\ref{sec:mathcal F} we give a
description of solution sets of systems of equations over finitely
generated fully residually free groups. To do this we need the
notion of a Hom-diagram described in Section~\ref{se:homs}.

\medskip
 \textsc{Theorem~\ref{th:hom-machine}.} {\em Let $G$ be a finitely generated group and $F$ a free group.
Then:
 \begin{enumerate}
  \item   there exists a complete canonical
 $Hom(G,F)$-diagram ${\mathcal C}$.
 \item  if $F$ is a fixed subgroup of $G$ then there exists an $F$-complete canonical
 $Hom_F(G,F)$-diagram ${\mathcal C}$.
 \end{enumerate}
Moreover,  if the group $G$ is finitely presented then  the
Hom-diagrams from {\rm (1)} and {\rm (2)} can be found
effectively.}

\medskip
This theorem implies the following.

\medskip
\textsc {Theorem~\ref{th:12.3}.} {\em Let $S(X)=1$ be a finite
system of equations and \[ G=(F*F[X])/ ncl(S).\] Then one can
effectively construct an $F$-complete canonical $Hom(G,F)$-diagram
${\mathcal C}$ such that all solutions of $S(X)=1$ in $F$ are
exactly all $F$-homomorphisms in ${\mathcal C}$.}

\medskip

 \textsc {Theorem~\ref{th:hom-machineNTQ}.} {\em Let $G$ be a
finitely generated group and $N$ an NTQ group. Then:
 \begin{enumerate}
  \item   there exists a complete canonical
 $Hom(G,N)$-diagram ${\mathcal C}$.
 \item  if $N$ is a fixed subgroup of $G$ then there exists an $N$-complete canonical
 $Hom_N(G,N)$-diagram ${\mathcal C}$.
 \end{enumerate}
Moreover,  if the group $G$ is finitely presented then  the
diagrams from (1) and (2) can be found effectively. }

\medskip
 This theorem implies the following result

\medskip
\textsc {Corollary~\ref{cor:15.2}.} {\em Let $G$ be a finitely
generated group and $H$ a finitely generated fully residually free
group. Then:
 \begin{enumerate}
  \item   there exists a complete canonical
 $Hom(G,H)$-diagram.
 \item  if $H$ is a fixed subgroup of $G$ then there exists an $H$-complete canonical
 $Hom_H(G,H)$-diagram.
 \end{enumerate}
}

\medskip
 In Section~\ref{lef} 
 we construct free regular Lyndon length
functions on NTQ groups.

\medskip
\textsc{Theorem~\ref{th:13.6}.} {\em If $G$ is the coordinate
group of a regular quadratic equation, then there exists a free
regular Lyndon length function $\ell :G\rightarrow {\mathbb Z}^2,$
where ${\mathbb Z}^2$ is ordered lexicographically.}

\medskip
\textsc{Theorem~\ref{th:13.7}.} {\em If $G$ is the coordinate
group of an NTQ system, then there exists a free regular Lyndon
length function $\ell :G\rightarrow {\mathbb Z}^t,$ where
${\mathbb Z}^t$ is ordered lexicographically.}

In Section~\ref{sec:NTQ} 
we prove some results about
monomorphisms into fully residually free groups.

 Let  $G$ and $H$ be
 finitely generated fully residually free groups.
 Consider monomorphisms from
  $G$ to $H$. One can define an
equivalence relation on the set of all such monomorphisms: two
monomorphisms $\phi$ and $\psi$ are equivalent if $\psi$ is a
composition of $\phi$ and conjugation by an element from $H$.

\medskip
\textsc {Theorem~\ref{monom}.} {\em Let $G$ {\rm[} $H${\rm]} be a
finitely generated fully residually free group, and let $\mathcal
A=\{A_1,\dots ,A_n\}$ {\rm[}respectively, $\mathcal B=\{B_1,\dots
,B_n\}${\rm]} be a finite set of non-conjugated maximal abelian
subgroups of $G$ {\rm[}respectively, $H${\rm]} such that the
abelian decomposition of $G$ modulo $\mathcal A$ is trivial. The
number of equivalence classes of monomorphisms from $G$ to $H$
that map subgroups from $\mathcal A$ onto conjugates of the
corresponding subgroups from $\mathcal B$ is finite. A set of
representatives of the equivalence classes can be effectively
found.}

\medskip
\textsc {Corollary~\ref{lem:BPproperty}.} {\em Let $G$ be a
finitely generated fully residually free group, and let $\mathcal
A=\{A_1,\dots ,A_n\}$ be a finite set of maximal abelian subgroups
of $G$. Denote by $Out(G;{\mathcal A})$ the set of those outer
automorphisms of $G$ which map each $A_i\in {\mathcal A}$ onto a
conjugate of it. If $Out(G;{\mathcal A})$ is infinite, then $G$
has a non-trivial abelian splitting, where each subgroup in
$\mathcal A$ is elliptic. There is an algorithm to decide whether
$Out(G;{\mathcal A})$ is finite or infinite. If $Out(G;{\mathcal
A})$ is infinite, the algorithm finds the splitting. If
$Out(G;{\mathcal A})$ is finite, the algorithm finds all its
elements. }

\medskip
 This  can be compared to the following Paulin's result
\cite{Paulin}. If $G$ is a hyperbolic group that has infinite
$Out(G)$, then $G$ admits a non-trivial small action on a real
tree. In this case $G$ splits over an elementary group \cite{BF}.

\section{Preliminaries}

\subsection{Free monoids and free groups}
\label{se:2-1}

Let $A= \{a_1, \dots, a_m\}$ be a set. By $F_{mon}(A)$ we denote
the free monoid generated by  $A$ which is defined as the set of
all words (including the empty word $1$)  over the alphabet $A$
with concatenation as multiplication. For a word $w = b_1 \dots
b_n,$ where $ b_i \in A$,    by $|w|$ or $d(w)$ we denote the
length $n$ of $w$.

To each $a\in A$ we associate a symbol $a^{-1}$. Put $A^{-1} =
\{a^{-1} \mid a \in A\},$ and suppose that $A\cap
A^{-1}=\emptyset.$   We assume that $a^{1} = a$, $(a^{-1})^{-1}=
a$ and $A^1 = A$. Denote  $A^{\pm 1} = A \cup A^{-1}$. If $w =
b_1^{\varepsilon_1} \cdots b_n^{\varepsilon_n} \in F_{mon}(A^{\pm
1})$, where  $ (\varepsilon_i \in \{1, -1\}) $, then we put
$w^{-1} = b_n^{-\varepsilon_n} \cdots b_1^{-\varepsilon_1};$ we
see that $w^{-1}\in M(A^{\pm 1})$ and say that $w^{-1}$ is an
inverse of $w$. Furthermore, we put $1^{-1}=1$.

A word $w \in F_{mon}(A^{\pm 1})$ is called {\it reduced}  if it
does not contain subwords  $bb^{-1}$ for $b \in A^{\pm 1}$. If $w
= w_1bb^{-1}w_2,$ $w  \in F_{mon}(A^{\pm 1})$ then we say that
$w_1w_2$ is obtained from $w$ by an elementary reduction $ bb^{-1}
\rightarrow 1$. A reduction process for $w$ consists of finitely
many reductions which bring $w$ to a reduced word $\bar w$. This
$\bar w$ does not depend on a particular reduction process and is
called the {\em reduced form} of $w$.

Consider a congruence relation on $F_{mon}(A^{\pm 1})$, defined
the following way:  two words are congruent if a reduction process
brings them to the same reduced word. The set of congruence
classes with respect to this relation forms a free group $F(A)$
with basis $A$. If not said otherwise, we assume that $F(A)$ is
given as the set of all reduced words in $A^{\pm 1}$.
Multiplication in $F(A)$ of two words $u, w$ is given by the
reduced form of their concatenation, i.e.,  $u\cdot v = {\bar
{uv}}$. A word $w \in F_{mon}(A^{\pm 1})$ determines the element
$\bar w \in F(A)$, in this event we sometimes say that $w$ is an
element of $F(A)$ (even though $w$ may not be reduced).

 Words $u, w \in F_{mon}(A^{\pm 1})$ are {\it graphically equal} if they are
equal in the monoid $F_{mon}(A^{\pm 1})$ (for example,
$a_1a_2a_2^{-1}$ is not graphically equal to $a_1$).

Let $X=\{x_1,\dots ,x_n\}$be a finite set of elements disjoint
with $A$. Let $w(X)= w(x_1,\dots,x_n)$ be a word in the alphabet
$(X\cup A)^{\pm 1}$ and $U = (u_1(A), \dots, u_n(A))$ be a tuple
of words in the alphabet $A^{\pm 1}$. By $w(U)$ we denote the word
which is obtained from $w$ by replacing each $x_i$ by $u_i$.
Similarly, if $W = (w_1(X), \dots, w_m(X))$ is an $m$-tuple of
words in variables $X$ then by $W(U)$ we denote the tuple
$(w_1(U),\dots,w_m(U))$. For any set $S$ we denote by $S^n$ the
set of all $n$-tuples of elements from $S$. Every word $w(X)$
gives rise to a map $p_w : (F_{mon}(A^{\pm 1}))^n \rightarrow
F_{mon}(A^{\pm 1})$ defined by $p_w(U) = w(U)$ for $U \in
F_{mon}(A^{\pm 1})^n$.  We call $p_w$ the word map defined by
$w(X)$. If $W(X) = (w_1(X), \dots, w_m(X))$ is an $m$-tuple of
words in variables $X$ then we define a word map $P_W:
(F_{mon}(A^{\pm 1}))^n \rightarrow F_{mon}(A^{\pm 1})^m$ by the
rule $P_W(U) = W(U)$.

\subsection{ $G$-groups}
\label{se:2-2}

 For the purpose of algebraic geometry over a given fixed group
 $G$,
one has to consider the category of $G$-groups, i.e., groups which
contain the group $G$ as a distinguished subgroup. If $H$ and $K$
are $G$-groups then a homomorphism $\phi: H \rightarrow K$ is a
$G$- homomorphism if $g^\phi = g$ for every $g \in G$, in this
event we write $\phi: H \rightarrow_G K$. In this category
morphisms are $G$-homomorphisms; subgroups are $G$-subgroups, etc.
By $Hom_G(H,K)$ we denote the set of all $G$-homomorphisms from
$H$ into $K$. It is not hard to see that the free product $G \ast
F(X)$ is a free object in the category of $G$-groups. This group
is called  a free $G$-group with basis $X$,  and we denote it  by
$G[X]$.  A $G$-group $H$ is termed {\em finitely generated
$G$-group} if there exists a finite subset $A \subset H$ such that
the set $G \cup A$ generates $H$. We refer to \cite{BMR} for a
general discussion on $G$-groups.

To deal with cancellation in the group $G[X]$ we need the
following notation. Let  $u = u_1 \cdots u_n \in G[X] = G \ast
F(X)$. We say that $u$  is {\em reduced} (as written) if $u_i \neq
1$, $u_i$ and $u_{i+1}$ are in
 different factors of the free product,  and if $u_i \in F(X)$ then it is reduced in the free group
 $F(X)$. By $red(u)$ we denote the reduced form of $u$.
  If  $red(u)  = u_1 \cdots u_n \in G[X]$, then we define  $|u| = n$, so $|u|$ is
  the syllable length  of $u$ in the free product $G[X]$.
 For reduced $u,v \in G[X]$, we write $u\circ v$ if the product $uv$ is reduced as written.
 If $u = u_1 \cdots u_n $ is  reduced  and $u_1, u_n$ are in
 different factors, then  we say that $u$ is {\em cyclically reduced}.

If $u = r \circ s$, $ v = s^{-1} \circ  t$, and $rt = r\circ t$
then we say that the word $s$ {\it cancels out in reducing}  $uv$,
or, simply, $s$ cancels out in $uv$. Therefore $s$ corresponds to
the {\it maximal} cancellation in $uv$.

\subsection{Elements of algebraic geometry over groups}
\label{se:2-4}

Here  we  introduce some basic notions of algebraic geometry over
groups. We refer to \cite{BMR} and \cite{KM9} for details.

Let $G$ be a group generated by a finite set $A$, $F(X)$ be a free
group with basis $X = \{x_1, x_2, \ldots  x_n\}$, $G[X] = G \ast
F(X)$ be a free product of $G$ and $F(X)$. If $S \subset G[X]$
then the expression $S = 1$ is called {\em a system of equations}
over $G$. As an element of the free product, the left side of
every equation in $S = 1$ can be written as a product of some
elements from $X \cup X^{-1}$ (which are called {\em variables})
and some elements from $A$ ({\em constants}). To emphasize this we
sometimes write $S(X,A) = 1$.

A {\em solution} of the system $ S(X) = 1$ over a group $G$ is a
tuple of elements $g_1, \ldots, g_n \in G$ such that after
replacement of each $x_i$ by $g_i$ the left hand side of every
equation in $S = 1$ turns into the trivial element of $G$.
Equivalently, a solution of the system $S = 1$ over $G$ can be
described as a $G$-homomorphism $\phi : G[X] \longrightarrow G$
such that $\phi(S) = 1$. Denote by $ncl(S)$ the normal closure of
$S$ in $G[X]$, and by $G_S$ the quotient group $G[X]/ncl(S)$. Then
every solution of $S(X) = 1$ in $G$ gives rise to a
$G$-homomorphism $G_S \rightarrow G$, and vice versa. By $V_G(S)$
we denote the set of all solutions in $G$ of the system $ S = 1$,
it is called the {\em algebraic set defined by} $S$. This
algebraic set $V_G(S)$ uniquely corresponds to the normal subgroup
$$ R(S) = \{ T(x) \in G[X] \ \mid \ \forall A\in G^n (S(A) = 1
\rightarrow T(A) = 1 \} $$ of the group $G[X]$. Notice that if
$V_G(S) = \emptyset$, then $R(S) = G[X]$. The subgroup $R(S)$
contains $S$, and it is called the {\it radical of $S$}. The
quotient group
$$G_{R(S)}=G[X]/R(S)$$ is the {\em coordinate group} of the
algebraic set  $V(S).$ Again, every solution of $S(X) = 1$ in $G$
can be described as a $G$-homomorphism $G_{R(S)} \rightarrow G$.

We recall from \cite{MyasExpo2} that a group $G$ is called a {\em
CSA group} if every maximal abelian subgroup $M$ of $G$ is
malnormal, i.e., $M^g \cap M = 1$ for any $g \in G \smallsetminus
M.$ The class of CSA-groups is quite substantial. It includes all
abelian groups, all torsion-free hyperbolic groups
\cite{MyasExpo2}, all groups acting freely on $\Lambda$-trees
\cite{bass}, and  many one-relator groups  \cite{gkm}.

We define a Zariski topology on $G^n$ by taking algebraic sets in
$G^n$ as a sub-basis for the closed sets of this topology. If $G$
is a  non-abelian CSA group, in particular, a non-abelian fully
residually free group (see definition in Subsection \ref{se:2-5}),
then the union of two algebraic sets is again algebraic. Therefore
the closed sets in the Zariski topology over $G$ are precisely the
algebraic sets.

A $G$-group $H$ is called {\it equationally Noetherian} if every
system $S(X) = 1$ with coefficients from $G$ is equivalent over
$G$ to a finite subsystem $S_0 = 1$, where $S_0 \subset S$, i.e.,
$V_G(S) = V_G(S_0)$. If $G$ is $G$-equationally Noetherian, then
we say that $G$ is equationally Noetherian. It is known that
linear groups (in particular, fully residually free groups) are
equationally Noetherian \cite{Gub,Br,BMR}. If $G$ is equationally
Noetherian then the Zariski topology over $G^n$ is {\em
Noetherian} for every $n$, i.e., every proper descending chain of
closed sets in $G^n$ is finite. This implies that every algebraic
set $V$ in $G^n$ is a finite union of irreducible subsets (they
are called {\it irreducible components} of $V$), and such
decomposition of $V$ is unique. Recall that a closed subset $V$ is
{\it irreducible} if it is not a union of two proper closed (in
the induced topology) subsets. Denote by $F_{\omega R(S)}$ a
coordinate group of an irreducible subvariety of $V(S)$.

\subsection{Discrimination and big powers}
\label{se:2-5}

Let $H$ and $K$ be $G$-groups.  We say that a family of $G$-
homomorphisms $\Phi \subset Hom_G(H,K)$ {\it separates} ({\it
discriminates}) $H$ into $K$ if for every non-trivial element $h
\in H$ (every finite set of non-trivial elements $H_0 \subset H$)
there exists $\phi \in \Phi$ such that $h^\phi \neq 1$ ($h^\phi
\neq 1$ for every $h \in H_0$). In this case we say that $H$ is
$G$-{\it separated} ($G$-{\it discriminated}) by $K$. Sometimes we
do not mention $G$ and simply say, $H$ is separated
(discriminated) by $K$. In the event when $K$ is a free group we
say that $H$ is {\it freely separated} ({\it freely
discriminated}).

Below we describe a  method of discrimination which is called a {\
big powers} method.  We refer to \cite{MyasExpo2} and \cite{KvM}
for details about BP-groups.

Let $G$ be a group. We say that a tuple $u =
(u_{1},\dots,u_{k})\in G^{k}$ has {\it commutation} if $[u_{i},
u_{i+1}] = 1$ for some $i=1,\dots,k-1.$ Otherwise we call $u$ {\it
commutation-free}.

\begin{definition}  A group $G$ satisfies the {\em big powers condition} (BP) if
for any commutation-free tuple $ u = (u_{1},\dots,u_{k})$ of
elements from $G$ there exists an integer $n(u)$ ({\it a boundary
of separation} for $u$) such that
$${u_{1}}^ {{\alpha}_{1}}\cdots u_{k}^{{\alpha}_{k}} \neq 1$$
 for
any integers $\alpha_1, \ldots, \alpha_k \geqslant n(u)$. Such
groups are called {\em BP-groups}. \end{definition}

The following  provides a host of examples of BP-groups.
Obviously, a subgroup of a BP-group is a BP-group; a group
discriminated by a BP-group is a BP- group \cite{MyasExpo2}; every
torsion-free hyperbolic group is a BP-group \cite{Ol'sh1}. It
follows that every freely discriminated group is a BP-group.

Let $G$ be a non-abelian CSA group and $u \in G$   be a non-proper
power. The following HNN-extension $$G(u,t) =  \langle G, t \mid
g^t = g (g \in C_G(u))\rangle$$ is called a free extension of the
centralizer $C_G(u)$ by a letter $t$.  It is not hard to see that
for any integer $k$ the map $t \rightarrow u^k$ extends uniquely
to a $G$-homomorphism $\xi_k : G(u,t) \rightarrow G$.

The result below is the essence of the big powers method of
discrimination.

\begin{theorem}
 \cite{MyasExpo2} Let $G$ be a non-abelian CSA BP-group
and $u\in G$ a non-proper power. If $G(u,t)$ is a free extension
of the centralizer of $u$ by $t$, then  the family of
$G$-homomorphisms $\{\xi_k \mid k \ is \ an \ integer \}$
discriminates $G(u,v)$ into $G$. More precisely, for every $w \in
G(u,t)$ there exists an integer  $N_w$ such that for every $k
\geqslant N_w$ $w^{\xi_k} \neq 1$.\end{theorem}

\smallskip
If $G$ is  a non-abelian CSA BP-group and $X$ is a finite set,
then   the group $G[X]$ is $G$-embeddable into $G(u,t)$ for any
non-proper power $u \in G$. It follows from the theorem above that
$G[X]$ is $G$-discriminated by $G$.

Unions of chains of extensions of centralizers play an important
part in this paper. Let   $G$ be a non-abelian CSA BP- group and
 $$G = G_0 <  G_1 <  \ldots < G_n $$
  be a chain of extensions of centralizers $G_{i+1} = G_i(u_i,t_i)$. Then every
  $n$-tuple of integers $p = (p_1, \ldots, p_n)$ gives rise to a $G$-homomorphism
  $\xi_p: G_n \rightarrow G$ which is  composition of homomorphisms
  $\xi_{p_i}: G_i \rightarrow G_{i-1}$ described above.

  A set $P$  of
  $n$-tuples  of integers is called {\it unbounded} if  for every integer $d$
  there exists a tuple $p  = (p_1, \ldots, p_n) \in P$ with $p_i \geqslant d$ for each $i$.
   The following result is a consequence of  the theorem above.

\begin{cy}
   Let $G_n$ be as above. Then for every unbounded set of tuples $P$
   the set of $G$-homomorphisms $\Xi_P = \{\xi_p \mid p \in P\}$ $G$-discriminates
   $G_n$ into $G$.
\end{cy}


   Similar results hold for  infinite chains
   of extensions of centralizers (see \cite{MyasExpo2}) and \cite{BMR3}).
   For example,   Lyndon's  free
   $Z[x]$-group $F^{Z[x]}$ can be realized as  union of a  countable chain of extensions of
   centralizers which starts with the free group $F$ (see  \cite{MyasExpo2}),
   hence  there exists a family of $F$-homomorphisms which discriminates
   $F^{Z[x]}$ into $F$.

\subsection{Fully residually free groups}
\label{se:2-6}

 Denote by ${\mathcal F}$
the class of all finitely generated fully residually free groups.
 In this section we describe some
properties of  groups from ${\mathcal F}$,  which are crucial for
our considerations.

It is not hard to see that every freely discriminated group is a
torsion-free CSA group \cite{BMR}.

Notice that every CSA group is commutation transitive
\cite{MyasExpo2}. A group $G$ is called {\em commutation
transitive} if commutation is transitive on the set of all
non-trivial elements of $G$, i.e., if $a,b,c \in G \smallsetminus
\{1\}$ and $[a,b] = 1, [b,c] = 1$, then $[a,c] = 1.$ Clearly,
commutation transitive groups are precisely the groups in which
centralizers of non-trivial elements are commutative.

\medskip
\begin{theorem}\cite{Rem}    Let $F$ be a free non-abelian group.
Then a finitely generated $F$-group $G$ is freely
$F$-discriminated by $F$ if and only if $G$ is $F$-universally
equivalent to $F$ {\rm(}i.e., $G$ and $F$ satisfy precisely the
same universal sentences in the language $L_A${\rm)}.\end{theorem}

\smallskip

\medskip
\begin{theorem}\cite{BMR,KM9}   Let $F$ be a free non-abelian group. Then a
finitely generated $F$-group $G$ is the coordinate group of a
non-empty irreducible algebraic set over $F$ if and only if $G$ is
freely $F$-discriminated by $F$.\end{theorem}

\smallskip

\medskip
\begin{theorem}\cite{KMNull} Let $F$ be a non-abelian free group. Then a
finitely generated $F$-group is  the coordinate group $F_{R(S)}$
of an irreducible non-empty algebraic set $V(S)$ over $F$ if and
only if $G$ is $F$-embeddable into the  free Lyndon's $Z[t]$-group
$F^{Z[t]}$.  \end{theorem}

\smallskip

Here we state one corollary of the results mentioned above.

\begin{theorem}\label{nt}\cite{KMIrc}
Every group $G \in {\mathcal F}$  can be obtained from free groups
by finitely many operations of the following type:
\begin{enumerate}
\item free products;

\item free products with amalgamation along  cyclic subgroups with
at least one of them being maximal;

\item separated HNN extensions along
 cyclic subgroups with at least one of them being maximal;

 \item free extensions of centralizers {\rm(}and for each element its centralizer extends
  at most once{\rm )}.

\end{enumerate}\end{theorem}

\begin{cy} \label{cy:fp} \cite{KMIrc,Se}
Every finitely generated fully residually free group is finitely
 presented.
 \end{cy}

\subsection{Quadratic equations over freely discriminated groups}
\label{se:2-7}

In this section we collect some known results about quadratic
equations over fully residually free groups, which will be in use
throughout this paper.

\smallskip
 Let $S \subset G[X]$. Denote by $var(S)$ the set of variables that
occur in $S$.
\begin{definition} A set $S \subset G[X]$ is called quadratic  if every variable
from
 $var(S)$ occurs in $S$ not more then twice. The set $S$ is strictly
quadratic if every letter from $var(S)$ occurs in $S$ exactly
twice.

A system $S = 1$ over $G$ is {\em quadratic [strictly quadratic]}
if the corresponding set $S$ is quadratic [strictly quadratic].
\end{definition}

\begin{definition}
A standard quadratic equation over the group $G$ is an equation of
the one of the following forms (below $d,c_i$ are nontrivial
elements from $G$):
\begin{equation}\label{eq:st1}
\prod_{i=1}^{n}[x_i,y_i] = 1, \ \ \ n > 0;
\end{equation}
\begin{equation}\label{eq:st2}
\prod_{i=1}^{n}[x_i,y_i] \prod_{i=1}^{m}z_i^{-1}c_iz_i d = 1,\ \ \
n,m \geqslant 0, m+n \geqslant 1 ;
\end{equation}
\begin{equation}\label{eq:st3}
\prod_{i=1}^{n}x_i^2 = 1, \ \ \ n > 0;
\end{equation}
\begin{equation}\label{eq:st4}
\prod_{i=1}^{n}x_i^2 \prod_{i=1}^{m}z_i^{-1}c_iz_i d = 1, \ \ \
n,m \geqslant 0, n+m \geqslant 1.
\end{equation}

Equations (\ref{eq:st1}), (\ref{eq:st2}) are called {\em
orientable} of genus $n$, equations (\ref{eq:st3}), (\ref{eq:st4})
are called {\em non-orientable} of genus $n$.
\end{definition}

\begin{lemma} \label{EC1}
Let $W$ be a strictly quadratic word over a group $G$. Then there
is a $G$-automorphism $f \in Aut_G(G[X])$ such that  $W^f$ is a
standard quadratic word over $G.$
\end{lemma}
\begin{proof} See \cite{LS}.\end{proof}

\begin{definition}
Strictly quadratic words  of the type $ [x,y], \ x^2, \ z^{-1}cz,
$ where $c \in G$, are called {\em atomic quadratic words} or
simply {\em atoms}.
\end {definition}

By definition a standard quadratic equation $S = 1$  over $G$  has
the form $$ r_1 \ r_2 \cdots r_k d = 1,$$ where $r_i$ are atoms,
$d \in G$.  This number $k$ is called the {\it atomic rank } of
this equation,  we denote it by $r(S)$.
 The {\em size} of $S=1$ is a pair $size (S)=(g(S), r(S))$,
 where $g(S)$ is the genus of $S=1$. We compare sizes lexicographically from the left.

 In
Section \ref{se:2-4} we defined the notion of the coordinate group
$G_{R(S)}.$ Every solution of the system $S=1$ is a homomorphism
$\phi : G_{R(S)} \rightarrow G$.
\begin{definition}  Let $S = 1$ be a standard quadratic equation written in the
atomic form
 $r_1r_2\cdots r_kd = 1 $ with $k \geqslant 2$.  A solution $\phi : G_{R(S)}
\rightarrow G$
 of $S = 1$  is called:
 \begin{enumerate}
 \item [1)] degenerate, if $r_i^\phi = 1$ for some $i$, and
 non-degenerate otherwise;
 \item [2)] commutative, if $[r_i^{\phi},r_{i+1}^{\phi}]=1$ for all
$i=1,\ldots ,k- 1,$  and  non-commutative otherwise;
 \item [3)] in a general position, if $[r_i^{\phi},r_{i+1}^{\phi}] \neq 1$ for all
$i=1,\ldots ,k-1,$.
 \end{enumerate}
 \end{definition}

 Observe that if a standard quadratic equation $S(X) = 1$ has a
 degenerate non-commutative solution then it has a non-degenerate
 non-commutative solution \cite{KMNull}.

\begin{theorem} \cite{KMNull} Let $G$ be a freely discriminated
group and $S = 1$ a standard quadratic equation over $G$ which has
a solution in $G$.
 In the following cases $S=1$
 always has a solution in $G$ in a general position:
\begin{enumerate}
 \item $S=1$ is of the form {\rm (\ref{eq:st1})}, $n>2$;
 \item $S=1$ is of the form {\rm (\ref{eq:st2})},  $n>0, \ n+m>1$;
 \item $S=1$ is of the form {\rm (\ref{eq:st3})}, $n>3;$
 \item $S=1$ is of the form {\rm (\ref{eq:st4})}, $n>2;$
 \item  $r(S) \geqslant 2$  and $S=1$ has a non-commutative solution.
\end{enumerate}
\end{theorem}
The following theorem describes the radical $R(S)$ of a standard
quadratic equation $S=1$ which has at least one  solution in a
freely discriminated group $G$.

\begin{theorem} \cite{KMNull}
\label{th:Nul} Let $G$ be a freely discriminated  group and let
$S=1$ be a standard quadratic equation over $G$ which has a
solution in $G$. Then
 \begin{enumerate}
\item If $S = [x,y]d$ or $S = [x_1,y_1][x_2,y_2]$, then $R(S) =
ncl(S)$; \item If $S = x^2d$, then $R(S) = ncl(xb)$ where $b^2 =
d$; \item If $S=c^zd$, then $R(S) = ncl([zb^{-1},c])$ where
$d^{-1} = c^b$; \item If $S = x_1^2x_2^2$, then $R(S) =
ncl([x_1,x_2])$; \item If $S = x_1^2x_2^2x_3^2$, then $R(S) =
ncl([x_1,x_2], [x_1,x_3], [x_2,x_3])$; \item If  $r(S) \geqslant
2$ and $S=1$
 has a non-commutative solution,  then $R
(S)=ncl (S)$; \item If  $S = 1$ is of the type (\ref{eq:st4}) and
all solutions of $S=1$
 are commutative, then $R(S)$ is the normal closure of the following
system:
  $$  \{x_1\cdots x_n=s_1\cdots s_n, [x_k,x_l]=1, [a_i^{-1}z_i,x_k]=1,
[x_k,C] = 1, [a_i^{-1}z_i,C] = 1,  $$ $$ [a_i^{-1}z_i,a_j^{-1}z_j]
= 1 \ \ (k,l = 1,\ldots,n; \ i,j = 1, \ldots , m) \} , $$ where
$x_k \rightarrow s_k, z_i \rightarrow a_i$ is a solution of $S=1$
and
$$C=C_G(c_1^{a_1},\dots ,c_m^{a_m},s_1,\ldots ,s_n)$$ is the
corresponding centralizer.  The group $G_{R(S)}$ is an extension
of the centralizer $C$.
  \end{enumerate}
  \end{theorem}

Put
 $$\kappa(S) = |X| + \varepsilon(S),$$
  where $\varepsilon(S) = 1$  if
the coefficient $d$ occurs in $S$, and  $\varepsilon(S) = 0$
otherwise.

\begin{definition}\label{regular}
A standard quadratic equation $S(X) = 1$ is  {\em regular} if
$\kappa(S) \geqslant 4$ and there is a non-commutative solution of
$S(X) = 1$ in $G$, or it is an equation of the type $[x,y]d = 1$.
\end{definition}

Notice, that if $S(X) = 1$ has a solution in $G$,  $\kappa(S)
\geqslant 4$, and $n
> 0$ in the orientable case ($n > 2$ in the non-orientable case),
then the equation $S = 1$ has a non-commutative solution, hence
regular.


\begin{cy} \ \begin{enumerate}
\item Every consistent orientable quadratic equation $S(X) = 1$ of
positive genus is regular, unless it is the equation $[x,y] = 1$;

\item Every consistent non-orientable equation of positive genus
is regular, unless  it is an equation of the type $x^2c^z = a^2c$,
$x_1^2x_2^2 = a_1^2a_2^2, x_1^2x_2^2x_3^2 = 1$, or $ S(X) =1$ can
be transformed to the form $[\bar z_i, \bar z_j]=[\bar z_i, a]=1,
\ i,j=1,\ldots ,m$ by changing variables.

\item Every standard quadratic equation $S(X) = 1$ of genus 0 is
regular unless either it is an equation of the type $c_1^{z_1} =
d, c_1^{z_1}c_2^{z_2} = c_1c_2$, or $ S(X) =1$ can be transformed
to the form $[\bar z_i, \bar z_j]=[\bar z_i, a]=1, \ i,j=1,\ldots
,m$ by changing variables.\end{enumerate}\end{cy}

\section{Splittings}
\label{se:split}

\subsection{Graphs}

A {\em directed graph} $X$ consists of a set of vertices $V(X)$
and a set of edges $E(X)$ together with two functions $\sigma :
E(X) \rightarrow V(X)$, $\tau : E(X) \rightarrow V(X)$. For an
edge $e\in E(X)$ the vertices $\sigma (e)$ and $\tau (e)$ are
called the {\it origin} and the {\it terminus} of $e$. A {\em
non-oriented graph} is a directed graph $X$ with involution $-:
E(X) \rightarrow E(X)$ which satisfies the following conditions:
$${\overline{\overline e}}=e,\,\, e\not = {\bar{e}},\,\, \sigma (\bar e )=\tau (e).$$
We refer to a pair $\{ e,\, \bar e \}$ as a non-oriented edge.

A {\it path} $p$ in a graph $X$ is a sequence if edges
$e_{1},\ldots , e_{n}$ such that $\tau (e_{i}) = \sigma
(e_{i+1}),\,\, i\in \{ 1,\ldots , n-1 \}$. Put $\sigma (p) =
\sigma (e_{1})$, $\tau (p) = \tau (e_{n})$. A path $p= e_{1}\cdots
e_{n}$ is {\it reduced} if $e_{i+1}\not = {\bar e}_{i}$ for each
$i$. A path $p$ is {\it closed} (or a {\it loop}) if $\sigma (p) =
\tau (p)$.

\subsection{Graphs of groups}

A {\em graph of groups} $\Gamma = {\mathcal G}(X)$ is defined by
the following data:
\begin{enumerate} \item [1)] a connected graph $X$;
\item [2)]  a function $\mathcal G$ which for every vertex $v\in
V(X)$ assigns a group $G_{v}$, and for each edge $e\in E(X)$
assigns a group $G_{e}$ such that $G_{\bar e} = G_{e}$, $\sigma
(G_e)=\tau (G_{\bar e})$;
 \item [3)] For each edge $e\in E(X)$ there are monomorphisms  $\sigma :
G_{e}\rightarrow G_{e\sigma }$ and  $\tau : G_e\rightarrow
G_{e\tau}$.
\end{enumerate}

Let ${\mathcal G}(X)$ be a graph of groups  and $T$ a maximal
subtree of $X$. The fundamental group $\pi ({\mathcal G}(X),T)$
 of the graph of groups ${\mathcal G}(X)$ with respect to the tree $T$
 is defined by generators and relations as follows: $\pi ({\mathcal
 G}(X),T)=$
$$\langle (*_{v\in V(X)}G_v), t_e (e\in
E(X))\mid t_e=1 (e\in T), t_e^{-1}gt_e=g^{\tau} (g\in G_e),
t_et_{\bar e}=1\rangle  .$$ It is known that $\pi ({\mathcal
G}(X),T)$ is independent (up to isomorphism) of $T$. Therefore, we
will omit sometimes the tree $T$ from the notations and write
simply $\pi ({\mathcal G}(X))$.

If some presentation is fixed, the non-trivial generators $t_e$
will be called {\em stable letters}. The group $\pi ({\mathcal
G}(X))$ can be obtained from the vertex groups by a tree product
with amalgamation and then by HNN-extensions. The following lemma
shows that subgroups of $\pi ({\mathcal G}(X))$ are again
fundamental groups of some special graphs of groups related to
${\mathcal G}(X).$

\begin{lemma}\label{Coh}\cite{Cohen} Let ${\mathcal G}(X)$ be a graph of groups, and let
$H\leqslant \pi ({\mathcal G}(X)).$ Then $H=\pi ({\mathcal G}(Y))$
where the vertex groups of ${\mathcal G}(Y))$ are $H\cap
gG_vg^{-1}$ for all vertices $v\in X$, and $g$ runs over a
suitable set of $(H,G_v)$ double coset representatives, and the
edge groups are $H\cap gG_eg^{-1}$ for all edges $e\in X$, where
$g$ runs over a suitable set of $(H,G_e)$ double coset
representatives.\end{lemma}

\subsection{Definitions and elementary properties of splittings}
Let $\pi ({\mathcal G}(X);T)$ be the fundamental group of graph of
groups ${\mathcal G}(X)$ with respect to a maximal subtree $T$.
Let $\phi: G\rightarrow \pi ({\mathcal G}(X);T)$ be an isomorphism
of groups. In this event the triple $D=(G,({\mathcal G}(X);T),\phi
)$ is called a {\em splitting} of $G$. A splitting $D$ is a
${\mathbb Z}$-splitting [abelian splitting] if every edge group is
infinite cyclic [abelian]. Splittings  of the type $G=A*_{C}B$ or
$G=A*_{C},$ are called {\em elementary} ${\mathbb Z}$-splittings.
An elementary abelian splitting $D$ is called {\it essential} if
the images of the edge group under the boundary monomorphisms do
not have finite index in the corresponding vertex groups.   A
splitting is {\em reduced} if all vertex groups of valency one and
two properly contain the images of groups of adjacent edges. An
abelian splitting is called {\em non-degenerate} if its graph of
groups is reduced and has an edge. A splitting is non-trivial if
its graph of groups has an edge.

Recall, that a splitting of a group $G$ is called 2-acylindrical
if for every non-trivial element $g\in G$, the fixed set of $g$
when acting on the Bass-Serre tree corresponding to the splitting
has diameter at most 2.
 A splitting
of a group is a \emph{star of groups}, if its underlying graph is
a tree $T$ which has diameter $2$.

\subsection{Elementary transformations of graphs of groups and splittings}

A {\em conjugation of a splitting} of $G$ is a  conjugation of
$\pi ({\mathcal G}(X);T)$.

A {\em sliding} is a modification of  a graph of groups according
to the relation
$$(A_1*_{C_1}A_2)*_{C_2}A_3\cong (A_1*_{C_1}A_3)*_{C_2}A_2 $$
in the case when $C_1\leqslant C_2$. More precisely, suppose a
graph of groups $\Gamma ={\mathcal G}(X)$ contains vertices
$v_1,v_2,v_3$ with vertex groups $A_1,A_2,A_3$ respectively, and
edges $e_1=(v_1,v_2), e_2=(v_2,v_3)$ with edge groups $C_1$ and
$C_2$ and $\tau _{e_1}(C_1)\leqslant\sigma _{e_2}(C_2)$. We
replace $e_1$ by the edge $\bar e_1=(v_1,v_3)$ with edge group
$C_1$ and $\sigma _{\bar e_1}(C_1)=\sigma _{e_1}(C_1), \tau _{\bar
e_1}(C_1)=\tau _{e_2}\sigma _{e_2}^{-1}\tau _{e_1}(C_1).$ Notice
that two of the vertices may coincide.

If $\Gamma$ is a graph of groups, and $G_e$ is an edge group in
$\Gamma$ such that $\sigma (G_e)=G_{e\sigma}$, then {\em
conjugation of the boundary monomorphism} $\sigma$  is a
replacement of the monomorphism $\sigma$ by $\sigma _h$ such that
$\sigma _h(g)=h^{-1}\sigma (g)h$ for some $h\in G_{e\sigma}$ and
any $g\in G_e$.

A splitting $G=A*_{C_1}B_1$ is obtained by {\em folding} from the
splitting $G=A*_CB$ if $C$ is a proper subgroup of $C_1$ and
$B_1=C_1*_CB$. An {\em unfolding} is the inverse operation to
folding.  A splitting is {\em unfolded} if one cannot apply an
unfolding to it.

\begin{lemma} Elementary transformations preserve the fundamental
groups of graphs of groups up to isomorphism.
\end{lemma}
\begin{proof} The statement is obvious for a conjugation.
Suppose the graph of groups $\Gamma _1={\mathcal G}_1(X_1)$ is
obtained from $\Gamma={\mathcal G}(X)$ by an elementary
transformation. Different choice of $T$ in the graph $X$ defines
an isomorphism of the fundamental group  $\pi ({\mathcal
G}(X),T)$. Therefore it is enough to prove the isomorphism of the
fundamental groups for a suitable choice of $T$ in $X$ and $ T_1$
in $X_1$. For sliding we consider only the case when vertices
$v_1,v_2,v_3$ are different. In this case we choose the tree $T$
such a way that edges  $e_1$ and $e_2$ belong to $T$, and obtain
$T_1$ by the sliding defined above from $T$. Then the isomorphism
becomes obvious.

Suppose $\Gamma _1={\mathcal G}_1(X)$ is obtained from $\Gamma $
by the conjugation of the boundary monomorphism $\sigma$. Consider
first the case when  it is possible to choose  $T$ such a way that
$e\not\in T$. In this case $G=\pi ({\mathcal G}(X);T)$ has a
stable letter $t$ corresponding to $e$ and a relation
$t^{-1}\sigma (g)t=\tau (g)$ for any $g\in G_e$. Let $t_h$ be a
stable letter in the presentation $\pi ({\mathcal G}_1(X);T)$
which has the corresponding relation $t_h^{-1}\sigma _h
(g)t_h=\tau (g)$. These two presentations define the same group,
because replacing the generator $t_h$ in the second presentation
by $h^{-1}t$, we obtain the first presentation. Now consider the
case when it is not possible to choose $T$ such a way that
$e\not\in T$. Then removing $e$ from $X$ we obtain two connected
components $X_1$ and $X_2$. Let $G_{e\sigma}\in X_1$, $G_1$ and
$G_2$ are fundamental groups of the graphs of groups corresponding
to $X_1$ and $X_2$. Then
$$\pi ({\mathcal G}(X);T)=G_1*_{\sigma ( G_e)=\tau (G_e)}G_2$$ and
$$\pi ({\mathcal G}_1(X);T)=G_1*_{h^{-1}\sigma ( G_e)h=\tau
(G_e)}(\bar G_2).$$ These presentations define the same group in
different generators, and $\bar x=x^h,$ for $x\in G_2$ and
corresponding $\bar x\in\bar G_2$.

Similarly, in the case of foldings we choose $T$ so that the edge
corresponding to $C$ belongs to $T$.

\end{proof}

\subsection{Freely decomposable groups}

Recall that  a group $G$ is {\em freely decomposable} if it is
isomorphic to a non-trivial free product (in which there are at
least two non-trivial factors). Otherwise, $G$ is called {\em
freely indecomposable}. A free decomposition $$G=G_1*\cdots
*G_n*F(Y)$$ is called  {\em Grushko's decomposition} if all the
factors $G_1,\ldots ,G_{n}$ are non-cyclic freely indecomposable
groups, and $F(Y)$ is a free group with basis $Y$ (pehaps empty).
If there is another Grushko's decomposition $$G=H_1*\cdots *H_m
*F(Z),$$ then $n=m$, corresponding factors $G_i$ and $H_i$ (after
reordering) are conjugated, and $|Z|=|Y|$.

\begin{definition} Let $G$ be a group and $H$ be a subgroup
of $G$. We say that $G$ has a {\em non-trivial free decomposition
modulo $H$}, if
 $$ G \simeq G_1 \ast G_2,$$
 $H \leqslant G_1,$  and $G_2 \neq 1.$
 \end{definition}

 Now we generalize the definition above.
 \begin{definition}
Let $G$ be a group and $K_1, \ldots, K_n$ be  subgroups of $G$. We
say that $G$ has a non-trivial free decomposition modulo $K_1,
\ldots, K_n$, if
 $$ G \simeq G_1 \ast G_2,$$
 $K_1 \leqslant G_1,$  and for each $i$ there exists $g_i\in G$ such that $K_i^{g_i}\leqslant G_1$.
 \end{definition}

\begin{lemma}
\label{le:1.5} Let $G$ be a group with  a non-trivial splitting
$D$ of $G$. If $D$ contains an edge with the trivial associated
group then $G$ is decomposable into a non-trivial free product.
\end{lemma}
\begin{proof}  Let $D = {\mathcal G}(X)$ be a non-trivial splitting of $G$
with an edge with the trivial associated subgroup.
 Denote by $X^*$ the subgraph of $X$
 formed by all vertices in $X$ and all edges in $X$  whose associated groups are non-trivial.
Assume, first, that  $X*$ is not connected and $Y_1, \ldots, Y_n$
($n \geqslant 2$) are the connected components of $X^*$.  By
collapsing in $D$ every  subgraph $Y_i$ into a single vertex, say
$v_i$, with the associated group $\pi_1(Y_i)$ one gets a graph of
groups $Z$ in which all the edge groups are trivial. Since all the
vertex groups in $D$ are nontrivial  the groups $\pi_1(Y_i)$ are
non-trivial. By Lemma \ref{Coh} the group $\pi_1(Z)$ is the  free
product of  its vertex groups and a free group, so it is freely
decomposable. Observe, that $G \simeq \pi_1(Z)$ since collapses
preserve the fundamental groups. Therefore $G$ is  freely
decomposable.

Suppose now that  $X^*$ is connected. In this case a given maximal
subtree $T$ of $X^*$ is also a maximal subtree of $X$. If $e$ is
an edge in $X$ with the trivial associated subgroup $G_e$ then $e
\not \in T$, so $t_e \neq 1$ in $G = \pi_1(X)$. Since $G_e = 1$
the infinite cyclic group $\langle t_e\rangle$
  is  a free factor of $G$. Hence $G$ is freely decomposable.
  This completes the proof of the lemma.
\end{proof}

\begin{definition}
Let $G$ be a group and ${\mathcal K} = \{K_1,\ldots, K_n\}$  be a
set of subgroups of $G$. We say that a free decomposition of $G$
$$G \simeq G_1*\cdots *G_k$$ is compatible with ${\mathcal K}$  if
each subgroup in ${\mathcal K}$  is a conjugate of a subgroup  of
one of the factors $G_j$. Denote by ${\mathcal K}_j$ the set of
all conjugates of subgroups in ${\mathcal K}$ which belong to
$G_j$. This decomposition is called {\em reduced} if none of the
$G_j$ has a non-trivial compatible free decomposition modulo
${\mathcal K}_j$.
\end{definition}
\begin{prop}Let $$G=G_1*\cdots
*G_n*F(Y)= G=H_1*\cdots *H_m *F(Z),$$ be two compatible with
$\mathcal K$ reduced free decompositions, then $n=m$,
corresponding factors $G_i$ and $H_i$ (after reordering) are
conjugated, and $|Z|=|Y|$.\end{prop}
\begin{proof} Consider a Bass-Serre tree $T$ corresponding to the
first decomposition. Subgroups from $\mathcal K$ fix some vertices
of this tree. Each subgroup $H_i$ acts on $T$, and, since it does
not have a free decomposition compatible with $\mathcal K$, fixes
a vertex of this tree. Therefore it is conjugated into some
factor, say $G_i$, of the first decomposition. Conversely, each
$G_i$ is conjugated into some $H_j$. Each factor $H_i, G_i$ is
malnormal, therefore $H_i$ and $G_i$ are conjugated. The normal
closure generated by $H_1*\cdots *H_m $ and by $G_1*\cdots *G_m$
is the same, therefore $|Y|=|Z|.$\end{proof}

\begin{remark}\label{rk1}
In the notations above, if $$G \simeq G_1*\cdots *G_k$$ is  a
compatible with ${\mathcal K}$  decomposition of $G$,  then we
will always assume (taking a conjugation of $G$ if necessary and
renaming subgroups $G_i$) that $K_1 \leqslant G_1$. \end{remark}

\begin{definition} Let $G$ be a group and ${\mathcal K} = \{K_1,\ldots, K_n\}$ be a
set of subgroups of $G$. An abelian splitting $D$ is called a {\em
splitting modulo $\mathcal K$} if all subgroups from ${\mathcal
K}$ are conjugated into vertex groups in $D$.\end{definition}

\subsection{Splittings of finitely generated fully residually free groups}
\label{frf}

The following result follows immediately from Theorem \ref{nt}.
\begin{theorem}
Every freely indecomposable  non-abelian group from ${\mathcal F}$
 has an essential ${\mathbb Z}$-splitting.\end{theorem}

\subsection{Elliptic and hyperbolic subgroups}\label{ehs}

If $H$ and $K$ are subgroups of a group $G$, we say that $H$ can
be {\em conjugated into} $K$ if it is a conjugate of a subgroup of
$K$.

An element $g\in G$ [a subgroup $H \leqslant G$] is called {\em
elliptic}  in a given splitting of $G$ if $g$ [correspondingly,
$H$] can be conjugated into a vertex group, and {\em hyperbolic}
otherwise.

\begin{lemma}\cite{RS} Let $G$ be a freely indecomposable group. If  $D_i$ is an
 elementary ${\mathbb Z}$-splitting of $G$ with the edge
group $C_i, (i = 1,2),$  then $C_1$ is hyperbolic {\rm
(}elliptic{\rm )} in $D_2$ if and only if $C_2$ is hyperbolic {\rm
(}elliptic{\rm )} in $D_1$.\end{lemma}

A pair of elementary ${\mathbb Z}$-splittings $D_i (i = 1.2)$  is
called {\em intersecting} if they form a hyperbolic-hyperbolic
pair, namely $C_1$ is hyperbolic with respect to $D_2$ and $C_2$
is hyperbolic with respect to $D_1$.

\begin{lemma}\label{ae}
Let $G\in {\mathcal F}$ and $N$ a maximal abelian non-cyclic
subgroup of $G$. Then the following holds:
 \begin{enumerate}
 \item If $G = A \ast _C B$ is an abelian splitting of $G$ then
 $N$ is elliptic in this splitting;
 \item  If $G=A \ast _C$ is an abelian splitting of $G$ then one of the following
 holds:
  \begin{itemize}
   \item [a)] $N$ is elliptic in this splitting;
    \item [b)] $C$ is the centralizer $C_A(v)$ of some element $ v \in A$, $C \leqslant N^g$ for some $g \in
 G$, and   $G = A \ast _C N^g$ is an extension of the centralizer $C$.
 \end{itemize}
 \end{enumerate}
\end{lemma}
 \begin{proof} The first statement follows from the description of
commuting elements in a free product with amalgamation (see, for
example, \cite{MKS}) . The second one is a direct corollary of
Lemma 2 and Theorem 5 from \cite{gkm}.
\end{proof}

\begin{definition} An abelian splitting $D$ of a group $G$ is
called {\em normal} if all maximal abelian non-cyclic subgroups of
$G$ are elliptic in $D$. By ${\mathcal D}(G)$ we denote the set of
all normal splittings of $G$.
  Denote by ${\mathcal D}_F(G)$ the class of all normal
splittings of $G$ such that $F$ belongs to a vertex
group.\end{definition}

\subsection{Quadratically hanging subgroups}
\label{subsec:QH-subgroups}

 Let $D=(G,({\mathcal G}(X);T),\phi )$  be an  abelian splitting
of $G$. A vertex group $Q = G_v, v \in V$ is called {\em
quadratically hanging} in $D$  (in short, QH-subgroup in $D$), if
the following conditions hold:
\begin{enumerate}

\item [1)] $Q$  admits one of the  following presentations

\begin{equation}\label{5}
\left\langle p  _1,\ldots ,p_m,a_1,\ldots ,a_g,b_1,\ldots ,b_g
\mid\prod _{j=1}^g[a_i,b_i]\prod _{k=1}^mp_k\right\rangle  ,
g\geqslant 0, m\geqslant 1,
\end{equation}
and if $g=0,$ then $m\geqslant 4,$

\begin{equation}\label{6}
\left\langle p  _1,\ldots ,p_m,a_1,\ldots ,a_g\mid\prod
_{j=1}^ga_i^2\prod _{k=1}^mp_k\right\rangle   , g\geqslant 1,
m\geqslant 1,
\end{equation}
(in particular, $Q$ is a free group);

\item [2)] for every edge $e \in E$ outgoing from $v$,
 the edge group $G_e$ is conjugate to one of the subgroups
 $\langle p_i\rangle , i=1,\ldots ,m$.

\item [3)] for each $p_i$ there is an edge $e_i\in E$ outgoing
from $v$ such that $G_{e_i}$ is a conjugate of $\langle p_i\rangle
$.\end{enumerate}

Notice that for a freely indecomposable group $G$ property 3) is
automatically satisfied. Moreover, using slidings one  can modify
the splitting $D$ in such a way that there are exactly $m$
outgoing edges $e_1,\ldots ,e_m$ from the QH-vertex $Q$  and
$\sigma(G_{e_i}) = \langle p_i\rangle $ for each edge $e_i$.

A QH-subgroup $Q$ is called a {\em maximal QH-subgroup} (in short,
MQH-subgroup) if    for every elementary  abelian splitting $D$ of
$G$ ( where $G =A*_CB $ or $G=A*_C$)  either $Q$ is elliptic in
$D$, or the edge group $C$ can be conjugated into $Q$, in which
case $D$ is inherited from the ${\mathbb Z}$-splitting of $Q$
along $C$.

Non-QH non-abelian vertex groups of $D$ are called {\em rigid}.

\subsection{QH-subgroups and quadratic equations}\label{QHquad}
 Let $D$ be a splitting of $G$.  One can choose the maximal subtree
 $T$ of the graph of groups $\Gamma$ in $D$ as follows.
Let $T_1$ be  a maximal subforest  of  the subgraph of $\Gamma$
spanned by non-QH non-abelian vertex groups.  Then one can extend
$T_1$ to a maximal subforest $T_2$ of the subgraph of $\Gamma$
spanned by all non-QH vertex groups, and then extend $T_2$  to a
maximal tree $T$.

{\bf Orientable case.} Let an MQH subgroup $Q$ in $D$ be given by
a presentation
$$\prod _{i=1}^{n}[x_i, y_i]p_1\cdots p_{m}=1.$$
 We may assume (see Section \ref{subsec:QH-subgroups}) that there are exactly $m$
outgoing edges $e_1,\ldots ,e_m$ from the QH-vertex $Q$  and
$\sigma(G_{e_i}) = \langle p_i\rangle $ for each edge $e_i$.
Denote also $\tau(G_{e_i}) = \langle c_i\rangle $.

In the standard presentation of $G$,  by generators and relations
as the fundamental group of $D$, the relations corresponding to
QH-vertex $Q$ have the following form
$$\prod _{i=1}^{n}[ x_i, y_i]p_1\cdots p_{m}=1, \ c_i^{z_i}=p_i, \ i=1,\ldots ,m, $$
where  $z_i=1$ if $e_i\in T,$  and $e_m \in T$. Excluding $p_i$
from the relations above one gets the relation
\begin{equation}\label{QQQ}\prod _{i=1}^{n}[ x_i,
y_i]c_1^{z_1}\cdots c_{m-1}^{z_{m-1}}c_m=1\end{equation}
 in the presentation of $G$. We refer to this relation as to the
 {\it quadratic relation} in $G$ corresponding to $Q$.

 Conversely, if $H \in {\mathcal F}$ and $c_1, \ldots, c_m \in H$
 then the quotient group
  $$ G = H[x_1,y_1,\ldots,x_n,y_n, z_1,\ldots ,z_{m-1}]/ncl\left(\prod _{i=1}^{n}[ x_i,
y_i]c_1^{z_1}\cdots c_{m-1}^{z_{m-1}}c_m\right)$$
 can be represented as the  fundamental group of a graph of groups  with
  a QH vertex  $Q$ given by the presentation
   $$\prod _{i=1}^{n}[ x_i, y_i]p_1\cdots p_{m}=1,$$
    with outgoing edges $e_i,\ i=1,\ldots ,m,$ such that
    $\sigma(G_{e_i}) = \langle p_i\rangle $ and $\tau(G_{e_i}) = \langle c_i\rangle $ for each edge $e_i$.

Every homomorphism $\phi :G\rightarrow F$ gives rise to a solution
in $F$ of the quadratic equation (corresponding to $Q$):
\begin{equation}\prod _{i=1}^{n}[ x_i,
y_i]c_1^{\phi z_1}\cdots c_{m-1}^{\phi z_{m-1}}c^\phi
_m=1.\end{equation}

{\bf Non-orientable case.} If an MQH subgroup $Q$ in $D$ is given
by
$$x_1^2\cdots x_n^2p_1\cdots p_{m}=1,  \ c_i^{z_i}=p_i, \ i=1,\ldots ,m, $$
where  $z_i=1$ if $e_i\in T,$  and $e_m \in T$. Then, similarly to
the orientable case, in the presentation of $G$ (relative to $D$)
the subgroup $Q$ inputs a  quadratic relation
$$x_1^2\cdots x_n^2c_1^{z_1} \cdots c_{m-1}^{z_{m-1}}c_m = 1. $$
 The converse is also true, so the quotient group
 $$G = H[x_1,\ldots,x_n, z_1,\ldots ,z_{m-1}]/ncl\left(x_1^2\cdots x_n^2 c_1^{z_1}\cdots c_{m-1}^{z_{m-1}}c_m\right)$$
 can be represented as the  fundamental group of a graph of groups  with
  a  QH vertex  $Q$ corresponding to the quadratic relation above.
 Again, homomorphisms $\phi: G \rightarrow F$ give rise to solutions
 of the equation
$$x_1^2\cdots x_n^2c_1^{\phi z_1} \cdots c_{m-1}^{\phi z_{m-1}}c_m = 1. $$

\subsection{Induced $QH$-subgroups}

\begin{lemma}
\label{le:1.6} Let $G\in {\mathcal F}$, $\Gamma = \Gamma(X) $ a
cyclic {\rm[}abelian{\rm]} splitting of $G$, and $Q$ a QH-subgroup
in $\Gamma$ associated with a vertex $v \in X$ with outgoing edges
$e_1, \ldots, e_m$. Denote by $Y_1, \ldots, Y_k$ the connected
components of the graph $X \smallsetminus \{e_1, \ldots, e_m\}$
and by $P_1, \ldots, P_k$ - the fundamental groups of the graphs
of groups induced from $\Gamma$ on $Y_1, \ldots, Y_k$. If $H$ is a
finitely generated non-cyclic subgroup of $G$ then one of the
following conditions holds:
 \be
 \item  $H$ is a nontrivial  free product;
 \item  $H \leqslant P_i^g$ for some $g \in G$ and $1 \leqslant i \leqslant k$;
 \item  $H$ is freely indecomposable, and for some $g\in G$
 the subgroup $H\cap Q^g$ has finite index in $Q^g$.
In this event $H \cap Q^g$ is a QH-vertex group in $H$. \ee

 If $H_Q=H\cap Q$ is non-trivial and has infinite index in $Q$,
 then $H_Q$ is a free product of some conjugates of $p_1^{\alpha
 _1},\ldots ,p_m^{\alpha _m}, p^{\alpha}$ and a free group $F_1$ {\rm(}maybe trivial{\rm)} which
 does not intersect any conjugate of $\langle p_i\rangle $ for
 $i=1,\ldots ,m.$

 \end{lemma}
\begin{proof} We first prove the lemma for a freely indecomposable group $G$.
Denote by $D$ the  splitting of $G$ obtained from the
splitting $\Gamma$ by collapsing all the  edges in $X$ not
adjacent to the vertex $v$.  Let  $D_H$ be the splitting of $H$
induced  by $D$. Suppose some conjugate of $Q$ has a non-trivial
intersection with $H$. Without loss of generality we can assume it
is $Q$ itself. Put $Q_H = Q \cap H$. The group $Q_H$ is finitely
generated as the intersection of two finitely generated subgroups
of a fully residually free group (see Lemma \ref{KMRS1} in Section
\ref{algor})  and, therefore, finitely presented.

Suppose that $Q$ has a presentation
\begin{equation}
\label{3}Q=\langle x_1,y_1,\ldots ,x_n,y_n,p_1,\ldots ,p_m, p \mid
\prod _{i=1}^{n}[x_{i},y_{i}]\prod _{i=1}^m p_i p=1\rangle  .
\end{equation}
Let $$\hat Q=\langle\hat x_1,\hat y_1,\ldots ,\hat x_n,\hat
y_n,\hat p_1,\ldots ,\hat p_m, \hat p \mid \prod _{i=1}^{n}[\hat
x_{i},\hat y_{i}]\prod _{i=1}^m \hat p_i \hat p=1\rangle  .$$

Denote by $\bar Q$ a  group  generated by $x_i,y_i,\hat x_i,\hat
y_i, p_i,\hat p_i,p,\hat p, t_i$ and given by  defining relations
$$\prod _{i=1}^{n}[x_{i},y_{i}]\prod _{i=1}^m p_{i}p=1,$$
$$\prod _{i=1}^{n}[\hat x_{i},\hat y_{i}]\prod _{i=1}^m \hat p_{i}\hat p=1,$$
$$p_i^{t_i}=\hat p_i, \ \ p=\hat p.$$
The group $\bar Q$ is a closed surface group. It is also a
fundamental group of a graph of groups with two vertices $w$ and
$\hat w$ and $m+1$ edges $e_0,\ldots ,e_{m}$ joining these
vertices. The vertex groups are $G_w=Q,\ G_{\hat w}=\hat Q$ and
the edge groups are $G_{e_0}=\langle p  \rangle  ,\
G_{e_i}=\langle p  _i\rangle  .$

There are two possibilities.

1. $Q_H$ does not intersect non trivially a conjugate of any
subgroup $\langle p_i\rangle  , i=1,\ldots m$ and $ \langle
p\rangle  .$ In this case $H$ is a nontrivial free product,
because $Q_H$ is a nontrivial vertex group in $D_H$ with all the
adjacent edge groups trivial.

2. $Q_H$ intersects  non trivially $\langle p_i\rangle  ^g$.
Without loss of generality we can suppose that $Q_H\cap \langle
p\rangle  ^g$ is non trivial. Taking $Q_H^{g^{-1}}$ instead of
$Q_H$ we can suppose that $Q_H$ intersects non-trivially the
subgroup $\langle p  \rangle  $. Consider a graph of groups
$\Gamma$ with two vertices $v$ and $\hat v$, and vertex groups
$Q_H$, and the double of this group $\hat Q_H\leqslant \hat Q$.
The outgoing edges from $v$ correspond to all possible
intersections of the conjugates of $\langle
 p_1\rangle  ,\ldots ,\langle p_m\rangle  , \langle p\rangle  $
 in $Q$ with $Q_H$ which are not
conjugated in $Q_H$. The outgoing edges from $\hat v$ are their
doubles. To obtain $\Gamma$, we identify each edge with its
double. Let $\bar Q_H$ be the fundamental group of $\Gamma $ with
stable letters $s_1,\ldots ,s_n$. Suppose the edge groups of $\bar
Q_H$ conjugated into $\langle p_i\rangle  $ are $\langle
p_i^{\alpha _jg_j}\rangle  ,\ j=1,\ldots ,k_i,$ where $g_j$ are
representatives of $(\langle p_i\rangle  , Q_H)$ double cosets in
$Q$, and the edge group conjugated into $\langle p\rangle  $ are
$\langle p^{\alpha _lg_l}\rangle  ,\ l=1,\ldots ,k_0$ ($\alpha$'s
are natural numbers). Relations of $\bar Q_H$ consist of relations
of $Q_H$ and $\hat Q_H$ and relations
$$p^{\alpha _lg_ls_{0l}}=\hat p^{\alpha _l\hat g_l},$$ where
$g_1=s_{01}=1,$  and
$$p_i^{\alpha _jg_js_{ij}}=\hat p_i^{\alpha _j\hat g_j}.$$

We will show now that $\bar Q_H<\bar Q$. Notice first, that
$p_i^{g_j}$ is conjugated to $\hat p_i^{\hat g_j}$ by
$g_j^{-1}t_i\hat g_j$. Consider a homomorphism $\phi :\bar
Q_H\rightarrow \bar Q$ that extends the natural embedding of $Q_H$
into $Q$ and $\hat Q_H$ into $\hat Q$ and such that
$s_{0l}^{\phi}=g_l^{-1}\hat g_l,\ s_{ij}^{\phi}=g_j^{-1}t_i\hat
g_j.$ It is a homomorphism, because it preserves the relations of
$\bar Q_H$.

We have to show that it is a monomorphism. Take some element $r\in
\bar Q_H$  in the HNN reduced form $r=h_0 s_{i_1}^{\epsilon_1} h_1
\cdots s_{i_n}^{\epsilon_n} h_n, (n \geqslant 0)$, where $h_i \in
Q_H*_{\langle p^{\alpha _0}=\hat p^{\alpha _0}\rangle  }\hat Q_H,
\epsilon_i \in \{1,-1\}$. Since $r$ is in the reduced form, there
is no consecutive subsequence $s_{ij}^{-1}, h_{k}, s_{ij}$ with
$h_{k} \in \langle p_i^{\alpha _{ij}g_{ij}}\rangle  $ or $s_{ij},
h_k, s_{ij}^{-1}$ with $h_k \in \langle \hat p_i^{\alpha _{ij}\hat
g_{ij}}\rangle  $. Suppose, first, that $n=0,\ r=h_0$. Then
$h_0=q_1\hat q_2\cdots q_{n-1}\hat q_n,$ where $q_i\not\in \langle
p^{\alpha _{01}}\rangle  $ and $\hat q_i\not\in \langle \hat
p^{\alpha _{01}}\rangle  $. This is a reduced form in the
amalgamated product $Q_H*_{\langle p^{\alpha _0}=\hat p^{\alpha
_0}\rangle  }\hat Q_H$. This element considered as an element of
$Q*_{\langle p=\hat p\rangle  }\hat Q$ is also in the reduced
form, therefore is nontrivial. Since $Q*_{\langle p=\hat p\rangle
}\hat Q\leqslant \bar Q$, $r^{\phi}\neq 1.$ Suppose now that
$n>0$. To obtain $r^{\phi}$ we substitute in $r$ elements
$s_{0l}^{\phi}=g_l^{-1}\hat g_l,\ s_{ij}^{\phi}=g_j^{-1}t_i\hat
g_j$ instead of $s_{0l}, s_{ij}$. We have to show that $r^{\phi}$
is in the reduced form in $\bar Q$. Consider first the case, when
$r$ does not contain any $s_{ij}$, $i\neq 0$, therefore,
$r^{\phi}$ does not contain any $t_i$. In this case
$$r^{\phi}=h_0 (g_{i_1}^{-1}\hat g_{i_1})^{\epsilon_1} h_1
\cdots (g_{i_n}^{-1}\hat g_{i_n})^{\epsilon_n} h_n,$$ $(n
\geqslant 0)$, where $h_i \in Q_H*_{\langle p^{\alpha _0}=\hat
p^{\alpha _0}\rangle  }\hat Q_H, \epsilon_i \in \{1,-1\}$. This is
an element from $Q*_{\langle p=\hat p\rangle  }\hat Q$. If we
write $h_1,\ldots ,h_n$ in reduced form, we can see how to write
$r^{\phi }$ in reduced form. Without loss of generality we can
suppose that $\epsilon _1=-1$. Since $g_{ij}\in Q\smallsetminus
Q_H$, the syllable length of $(g_{i_1}^{-1}\hat
g_{i_1})^{\epsilon_1} h_1  (g_{i_2}^{-1}\hat
g_{i_2})^{\epsilon_2}$ can become  less than three only in the
case $\epsilon _2=1$, $h_1\in Q_H$ and $ g_{i_1}h_1
g_{i_2}^{-1}\in \langle  p\rangle  .$  Suppose this is the case.
Then $Q_H$ contains $h_1=g_{i_1}^{-1}p^kg_{i_2},\ p^{\alpha}, \
g_{i_1}^{-1}p^{\alpha _{i_1}}g_{i_1},\ g_{i_2}^{-1}p^{\alpha
_{i_2}}g_{i_2}.$ Therefore it contains $h_1g_{i_2}^{-1}p^{\alpha
_{i_1}-k}g_{i_1}.$ Then $g_{i_1}$ and $g_{i_2}$ represent the same
double coset and thus equal, which implies that corresponding
stable letters $s_{0l}$ are the same and
$r=h_0s_{0l}^{-1}h_1s_{0l}\cdots h_n$ is not in reduced form. This
contradicts to the assumption. Suppose now $r$ does contain
$s_{ij}$ for some $i>1$. Suppose that we have a
 reducible sequence $t_i^{-1}p_i^kt_i$ in $r^{\phi}$.
This can only appear from the sequence $s_{ij_1}^{-1}hs_{ij_2}$ in
reduced form of $r$. We can take for definiteness $j_1=1, j_2=2.$
Then $(s_{i1}^{-1}hs_{i2})^{\phi}=\hat
g_{i_1}^{-1}t_i^{-1}g_{i_1}hg_{i_2}^{-1}t_i\hat g_{i_2}.$ As in
the previous case we have $g_{i_1}hg_{i_2}^{-1}=p_i^k$, and  $Q_H$
contains $h=g_{i_1}^{-1}p_i^kg_{i_2}.$ This implies $s_1=s_2$ and
$r$ is not in  reduced form, that contradicts the assumption.

 Either $\bar Q_H$ is a closed surface group (in this case
 $\bar Q_H$ has finite index in $\bar Q$
 and $Q_H$ is of finite index in $Q$), or
 $\bar Q_H$ is free.

Suppose $\bar Q_H$ is free. Denote the elements $p_i^{\alpha
_{ij}g_{ij}}$ by $q_0,\ldots ,q_n,$ where $q_0$ is a power of $p$,
and rewrite the corresponding relations in the form $q_0=\hat
q_0,$ $q_i^{s_i}=\hat q_i.$ We will prove that for some
$h_1,\ldots ,h_n\in Q_H$ the product $\langle q_0\rangle  *\langle
q_1\rangle  ^{h_1}*\cdots *\langle q_n\rangle  ^{h_n}$ is a free
factor of $Q_H$. Consider the following free groups.
$$R_n=Q_H*\langle r_1\rangle  *\cdots *\langle r_n\rangle
,\ \hat R_n=\hat Q_H*\langle \hat r_1\rangle  *\cdots *\langle
\hat r_n\rangle  .$$ Let $\bar R_n$ be the amalgamated product of
$R_n$ and $\hat R_n$ amalgamated along  subgroups $\langle
q_0,q_1^{r_1},\ldots ,q_n^{r_n}\rangle  $ and $\langle \hat
q_0,\hat q_1^{\hat r_1},\ldots ,\hat q_n^{\hat r_n}\rangle  $ such
that $q_0=\hat q_0,$ $q_i^{r_i}=\hat q_i^{\hat r_i}.$ The group
$\bar Q_H$ is embedded into $\bar R_n$. More precise, $\bar
Q_H*\langle r_1\rangle  *\cdots *\langle r_n\rangle  $ is
isomorphic to $\bar R_n$. Indeed, a mapping $\theta$ that maps
$Q_H$ and $\hat Q_H$ identically onto corresponding subgroups of
$\bar R_n$, maps $s_i$ to $r_i^{-1}\hat r_i$ and $r_i$ to $r_i$
can be extended to an isomorphism.

The group $\bar R _n$ is an amalgamated product of two copies of
the free group $Q_H*\langle r_1\rangle  *\cdots *\langle
r_n\rangle $ along the same subgroup. By a result of Swan (see
\cite{Swan}, Lemma 7.1) such a group is free if and only if the
amalgamated subgroup is a free factor.

Therefore $Q_H*\langle r_1\rangle  *\cdots *\langle r_n\rangle
=\langle q_0\rangle  *\langle q_1^{r_1}\rangle *  \cdots *\langle
q_n^{r_n}\rangle  *F_0$ for some free group $F_0$. Each subgroup
$\langle q_i\rangle  $, in particular $\langle q_0\rangle  $, is a
free factor in $Q_H$ by the same result of Swan. Let
$Q_H=Q'_H*\langle q_0\rangle  .$ Considering $\bar R_n$ for $n=1$
we obtain that $\langle q_0\rangle  *Q'_H*\langle r_1\rangle
=\langle q_0\rangle  *\langle q_1^{r_1}\rangle  *F_0.$ Denote this
group by $T$. We want to show that $\langle q_0\rangle  *\langle
q_1\rangle  ^h$ is a free factor of $Q_H$ for some $h\in Q_H$. Fix
a basis of $T$ that includes $q_0,$ some basis of $Q'_H$ and
$r_1=r$, and represent elements of $\langle q_1^{r}\rangle  $ and
$F_0$ in this basis. We will now be considering bases of $T$
containing $q_0$, $q_1^r$. We order bases of $T$ containing $q_0$,
$q_1^r$ and some elements $\bar z=\{z_1,\ldots ,z_k\}$ the
following way: let $n_{\bar z}=(n_0,\ldots ,n_i,\ldots ),$ where
$n_i$ is the number of elements in $\bar z$ containing $i$
occurrences of $r^{\pm 1}$. One basis is smaller than the other if
the tuple $n_{\bar z}$ for it is larger in the left lexicographic
order. Let $\ell (z_i)$ be the number of occurrences of $r^{\pm
1}$ in $z_i$. Applying Nielsen transformations that do not change
$q_0$ and $q_1^r$ to a minimal basis of $T$ we turn the set of
basis elements that do not contain $r^{\pm 1}$ into a Nielsen set.
We call such a basis a strongly minimal. If, for some $z_i$, $\ell
(q_1^rz_i^{\pm 1})< 2$, then either $q_1^r$ is completely
cancelled, and this contradicts the minimality, or $\ell (z_i)=1$
and $q_1^{rz_{i}^{\pm 1}}$ is a free factor in $Q_H'$, therefore
$\langle q_0\rangle  *\langle q_1\rangle  ^h$ is a free factor in
$Q_H$ and we do not have to prove anything.

We claim now that by applying Nielsen transformations that do not
change $q_0$ and $q_1^r$ to a strongly minimal basis of $T$ we can
turn this basis into a basis with the following properties: (1)
for any two basis elements $y_1,y_2$, $\ell (y_1y_2^{\pm
1})\geqslant \ell (y_1);$  (2) if for a basis element $z$ the
number $\ell (z)$ is odd, then the middle occurrence of $r^{\pm
1}$ in $z$ is not cancelled in any product of basis elements; (3)
in any product of basis elements, elements  with an even number of
occurrences of $r^{\pm 1}$ greater than two have either one of the
middle occurrences of $r^{\pm 1}$ uncancelled or, if all
occurrences of $r^{\pm 1}$ are cancelled, one of the middle
letters not cancelled. (4) in any product of $q_1^r$ and elements
$z_i$, such that each $z_i$ either doesn't contain $r^{\pm 1}$ or
contains two occurrences of $r^{\pm 1}$, the first and last letter
$r^{\pm 1}$ cannot cancel.  The first two properties we already
have for a strongly minimal basis. The proof of (3) is just a
repetition of the proof of the same property in Nielsen's theorem.
We will now prove (4). Suppose that the first and last letters
$r^{\pm 1}$ are cancelled in the product of the form
$q_1^{r}z_{i_1}^{\epsilon _1}\cdots z_{i_j}^{\epsilon _j}$.
 If $z_{i_j}\neq z_i,\ i=i_1,\ldots ,i_{j-1}$ we can
replace $z_{i_j}$ by $q_1^{r}z_{i_1}^{\epsilon _1}\cdots
z_{i_j}^{\epsilon _j}$ and obtain a smaller basis which
contradicts  the minimality of the basis $q_0,
q_1^{r},z_{1},\ldots z_{k}.$ Suppose $z_{i_j}$ appears more than
once in the product $q_1^{r}z_{i_1}^{\epsilon _1}\ldots
z_{i_j}^{\epsilon _j}$. Suppose $z_{i_j}^{\epsilon _j}$ ends with
$rw$, where $w\in Q_H$. Then $w$ cannot be expressed in terms of
basis elements of length zero  and $w\neq 1$ because the
expression is nontrivial. Let $z_{i_t}^{\epsilon _t}$ be the first
element from the right in the expression $q_1^{r}z_{i_1}^{\epsilon
_1}\ldots z_{i_j}^{\epsilon _j}$ that ends with $rw$, where $w\in
Q_H$ and does not begin with $w^{-1}r^{-1}$ (if $z_{i_j}^{\epsilon
_j}$ does not begin with $w^{-1}r^{-1}$, then $t=j$). Multiplying
all the other basis elements of length 2 that end with $rw$ by
$z_{i_t}^{\pm 1}$ we will obtain a new basis such that the unique
element of length two that ends with $rw$ or begins with
$w^{-1}r^{-1}$ (but not both)is $z_{i_t}$. Since the expression
$q_1^{r}z_{i_1}^{\epsilon _1}\ldots z_{i_j}^{\epsilon _j}$ equals
 $w$, in the new basis it ends with $z_{i_t}$, and
$z_{i_t}$ appears only once. Replacing $z_{i_t}$ by this
expression, we obtain a basis smaller than the one we began with,
contradicting minimality. Let $q_0, q_1^r,z_1,\ldots ,z_k$ be a
basis with properties (1)-(4).
 Denote by $F'_0$ the free group
with  basis $\{ z_1,\ldots ,z_k\}$. Then $r\in \langle
q_1^r,z_1,\ldots ,z_k\rangle  $ and $T=\langle q_0\rangle *\langle
q_1\rangle  *(F'_0)^{r^{-1}}.$ Therefore $Q_H=\langle q_0\rangle
*\langle q_1\rangle  *F''_0.$

 We can now apply induction on $n$ and prove that
 $Q_H$ is a free product of conjugates of $p_1^{\alpha _1},\ldots
,p_m^{\alpha _m}, p^{\alpha}$ belonging to $Q_H$ and a group $F_1$
which does not intersect any conjugate of $\langle p_i\rangle  $
for $i=1,\ldots ,m$. Each conjugate of $p_j^{\alpha}$ in $H$ can
be conjugated into one of the groups $\langle p_i^{\alpha
_i}\rangle  $ in $Q_H$. This proves the second statement of the
lemma.

Suppose that for any $g$ the subgroup $Q^g\cap H$ is either
trivial or has the structure described above. Consider now the
decomposition $D_H$. If the group $F_1$ is nontrivial, then $H$ is
freely decomposable, because the vertex group $Q_H$ in $D_H$ is a
free product, and all the edge groups belong to the other factor.
If at least for one subgroup $Q^g$, such a group $F_1$ is
non-trivial, then $H$ is a non-trivial free product. Suppose  each
non-trivial subgroup $H\cap Q^g$ is a free product of conjugates
of some elements $p_i^{\alpha _i},\ \alpha _i\in Z,$ in $Q^g$.
According to the Bass-Serre theory, for the group $G$ and its
decomposition $D$ one can construct a tree such that $G$ acts on
this tree, and stabilizers correspond to vertex and edge groups of
$D$. Denote this Bass-Serre tree by $T_D$. The subgroup $H$ also
acts on $T_D$. Let $T_1$ be a fundamental transversal for this
action. Since  $H$ is not conjugated into any of the subgroups
$P_i$, the amalgamated product of the stabilizers of the vertices
of $T_1$ is a free product of subgroups $H\cap P_i^g$. Therefore
$H$ is either such a free product or is obtained from such a free
product by a sequence of HNN extensions with associated subgroups
belonging to distinct factors of the free product. In both cases
$H$ is freely decomposable.

In the case when $Q_H$ has finite index in $Q$, $\bar Q_H$ has
finite index in $\bar Q$, therefore the closed surface with the
fundamental group $\bar Q_H$ is a finite cover of the closed
surface with the fundamental group $\bar Q$.
Therefore $Q_H$ is a QH subgroup for $H$.

 Suppose now that $G=G_1*\cdots *G_k*F_r$ is a
Grushko's decomposition of $G$, where $G_1,\ldots ,G_k$ are freely
indecomposable non-cyclic groups. Then any other such
decomposition has the same numbers $k$ and $r$ and indecomposable
factors conjugated to $G_1,\ldots ,G_k$. Since $H$ is freely
indecomposable, we can suppose that $H$ is contained in $G_1$. Let
$D$ be a decomposition of $G$ as in the lemma, and $H$ does not
belong to a conjugate of $P_i$. If the intersection $Q^g\cap H$
has infinite index in $Q^g$ for any $g\in G$, then the proof that
$H$ should be freely decomposable does not require the property of
$G$ being freely indecomposable. Suppose there exists $g\in G$
such that this intersection has  finite index. We can suppose
$g=1$, then $Q\leqslant G_1$. Then $Q$ is a vertex group in the
induced decomposition of $G_1$, therefore it is  a QH subgroup of
$G_1$. Now we can apply Lemma \ref{le:1.6} to $H$ and $G_1$ and
conclude that $Q\cap H$ is a QH subgroup of $H$.\end{proof}

\subsection{Quadratic decomposition}

By Theorem 5.6 from \cite{RS} for every f. g. freely
indecomposable torsion free non-surface group $H$ there exists a
reduced (may be trivial) $\mathbb Z$-splitting $D_{quadr}$ (a
quadratic decomposition of $H$ ) with the following properties:
\begin{enumerate}
\item every  $MQH$ subgroup of $H$ can be conjugated to a vertex
group in $D_{quadr}$; every $QH$ subgroup of $H$ can be conjugated
into one of the $MQH$ subgroups of $H$; every vertex with a
non-MQH vertex group is adjacent only to vertices with MQH vertex
groups;

\item if an elementary $\mathbb Z$-splitting $H=A*_CB$ or $H=A*_C$
is hyperbolic in another elementary $\mathbb Z$-splitting of $H$,
then $C$ can be conjugated into some MQH subgroup;

\item  for every elementary $\mathbb Z$- [abelian] splitting
$H=A*_CB$ or $H=A*_C$ from ${\mathcal D}(H)$ which is elliptic in
each elementary $\mathbb Z$- [abelian] splitting from ${\mathcal
D}(H)$, the edge group $C$ can be conjugated into a non-MQH
subgroup  of $D_{quadr}$;

\item  if $D_{quadr}'$ is another splitting that has properties
(1)--(4), then it can be obtained from $D_{quadr}$ by slidings,
conjugations, and modifying boundary monomorphisms by conjugation.
\end{enumerate}
\subsection{JSJ-decompositions}

\label{jsj} All elementary cyclic [abelian] splittings of a
finitely presented (f.p.) torsion free freely indecomposable group
are encoded in a splitting called  a JSJ decomposition.

\begin{prop}\label{Rips}{\rm(} \cite{RS}, part of Theorem
7.1{\rm)}
Let $H$ be a f.p. torsion-free freely
indecomposable group. There exists a  reduced, unfolded ${\mathbb
Z}$-splitting of $H$, called a JSJ decomposition of $H$, with the
following properties.

\bi \item[(1)] Every  $MQH$ subgroup of $H$ can be conjugated to a
vertex group in the JSJ decomposition. Every $QH$ subgroup of $H$
can be conjugated into one of the $MQH$ subgroups of $H$. Every
non-MQH vertex group in the JSJ decomposition is elliptic in every
${\mathbb Z}$-splitting of $H$.

\item[(2)] If an elementary ${\mathbb Z}$-splitting $H=A*_CB$ or
$H=A*_C$ is hyperbolic in another elementary ${\mathbb Z}$
splitting, then $C$ can be conjugated into some MQH subgroup. \ei
\end{prop}

We call a splitting of a finitely presented group {\em almost
reduced} if vertices of valency one and two properly contain the
images of the edge groups except vertices between two MQH
subgroups that may coincide with one of the edge groups.
\begin{prop}\label{existence}
Let $H$ be  a freely indecomposable f.g. fully residually free
group. There exists an almost  reduced unfolded cyclic [abelian]
splitting $D\in{\mathcal D}(H)$ of $H$, where the class $\mathcal
D(H)$ was defined in Section \ref{ehs}, with the following
properties:
\begin{enumerate}

\item Every  $MQH$ subgroup of $H$ can be conjugated to a vertex
group in D; every $QH$ subgroup of $H$ can be conjugated into one
of the $MQH$ subgroups of $H$; non-MQH subgroups in $D$ are of two
types: maximal abelian and non-abelian, every non-MQH vertex group
in $D$ is elliptic in every cyclic [abelian] splitting in
${\mathcal D}(H)$.

\item If an elementary cyclic [abelian] splitting $H=A*_CB$ or
$H=A*_C$ is hyperbolic in another elementary cyclic [abelian]
splitting, then $C$ can be conjugated into some MQH subgroup.

\item Every elementary cyclic [abelian] splitting $H=A*_CB$ or
$H=A*_C$ from ${\mathcal D}(H)$ which is elliptic with respect to
any other elementary cyclic [abelian] splitting from ${\mathcal
D}(H)$ can be obtained from $D$ by a sequence of collapsings,
foldings, conjugations and modifying boundary monomorphisms by
conjugation.

\item If $D_1\in{\mathcal D}(H)$ is another splitting that has
properties {\rm (1)--(2)}, then it can be obtained from $D$ by
slidings, conjugations, and modifying boundary monomorphisms by
conjugation.
\end{enumerate}\end{prop}

We will call such a splitting a {\em cyclic [abelian] JSJ
decomposition of} $H$. Similar result holds for the class of
splittings $\mathcal D _F(H).$ Such a splitting will be called  an
{\em cyclic} ({\em or abelian}) {\em JSJ decomposition} of $H$
modulo $F$.

\begin{proof} By Lemma \ref{ae} an elementary abelian non-cyclic splitting
from ${\mathcal D}(H)$ is elliptic in all
splittings from ${\mathcal D}(H)$. The proof of Lemma 2.1 from
\cite{RS} can be repeated to show that the case when a cyclic
splitting $H=A_1*_{C_1}B_1$ ($H=A_1*_{C_1}$) is hyperbolic in a
non-cyclic abelian splitting  $H=A_2*_{C_2}B_2$ ($H=A_2*_{C_2}$)
is impossible.

To construct a decomposition $D$ with properties 1-4, we first
construct a cyclic decomposition $D_{cyclic}$. We begin with
$D_{quadr}$ and refine it using property (c). We take an
elementary  cyclic splitting from  ${\mathcal D}(H)$, and if there
is a non-MQH vertex group that is not elliptic in this splitting,
we split this vertex group. This process stops, because $H$ is
finitely presented \cite{BeF}. As a result we obtain a cyclic
splitting of $H$. All non-MQH vertex groups of this splitting are
elliptic in all cyclic splittings from ${\mathcal D}(H)$.

 Every subgraph corresponding
to a decomposition $H_1*_{C_1}H_2*_{C_2}\cdots *_{C_{n-1}}H_n$,
such that $C_1,\dots ,C_{n-1}$ are subgroups of a maximal abelian
subgroup $M$ which is a subgroup of $H_1$ can be modified (using
slidings) into a star of groups with $H_1$ at the center. Since
all maximal abelian subgroups of $G$ are elliptic in $D$ and they
are malnormal, using slidings we can transform an abelian
decomposition into 2-acylindrical decomposition. Now we can refine
 $D_{cyclic}$ by splitting non-MQH non-abelian subgroups using
 non-cyclic abelian splittings from ${\mathcal D}(H)$. This
 procedure stops by \cite{We}, and we obtain a decomposition $D$
 with properties (1)--(2). These properties and 2-acylindricity imply
properties (3) and (4).
\end{proof}

\begin{cy}{\rm(}from Theorem {\rm\ref{nt})} Every freely indecomposable
non-abelian non-surface
 group from $\mathcal F$ admits a non-degenerate cyclic
 {\rm[}abelian{\rm]}
 $JSJ$ decomposition. \end{cy}

 Let $G$ be a group and ${\mathcal K} = \{K_1,\ldots, K_n\}$  be a
set of subgroups of $G$. Consider a free decomposition
$G=G_1*\cdots *G_{\ell}$ compatible with $\mathcal K$ and such
that each factor $G_i$ that contains  conjugates of some $K_j,
j=1,\ldots ,n,$ is freely indecomposable modulo these subgroups
and does not have a non-trivial compatible free decomposition
modulo them, and each factor $G_i$ that does not contain any
conjugate of $K_j, j=1,\ldots ,n,$ is freely indecomposable. For
each $G_i$ that does not contain conjugates of any $K_j$ we can
consider an abelian JSJ decomposition. For each $G_i$ containing
conjugates of some $K_{j_1},\ldots ,K_{j_s}$ one can consider
splittings modulo $K_{j_1},\ldots ,K_{j_s}$ and introduce the
notion of a {\em JSJ decomposition modulo} $K_{j_1},\ldots
,K_{j_s}$. This decomposition may be degenerate. Similarly to
Proposition \ref{existence} one can prove.

\begin{prop}
Let $H$ be  a freely indecomposable f.g. fully residually free
group modulo  $\mathcal K=\{K_{j_1},\ldots ,K_{j_s}\}$ such that
there is no non-trivial compatible free decomposition of $H$
modulo $\mathcal K$. There exists an  abelian splitting
$D\in{\mathcal D}(H)$ of $H$ modulo $\mathcal K$ (possibly
degenerate), where the class $\mathcal D(H)$ was defined in
Section \ref{ehs}, with the following properties:
\begin{enumerate}

\item Every  $MQH$ subgroup of $H$ modulo $\mathcal K$ can be
conjugated to a vertex group in D; every $QH$ subgroup of $H$
modulo $\mathcal K$ can be conjugated into one of the $MQH$
subgroups of $H$ modulo $\mathcal K$; non-MQH subgroups in $D$ are
of two types: maximal abelian and non-abelian, every non-MQH
vertex group in $D$ is elliptic in every abelian splitting in
${\mathcal D}(H)$ modulo $\mathcal K$.

\item If an elementary abelian splitting $H=A*_CB$ or $H=A*_C$
modulo $\mathcal K$ is hyperbolic in another elementary abelian
splitting modulo $\mathcal K$, then $C$ can be conjugated into
some MQH subgroup.

\item Every elementary abelian splitting modulo $\mathcal K$,
$H=A*_CB$ or $H=A*_C$ from ${\mathcal D}(H)$ which is elliptic
with respect to any other elementary abelian splitting modulo
$\mathcal K$ from ${\mathcal D}(H)$ can be obtained from $D$ by a
sequence of collapsings, foldings, conjugations and modifying
boundary monomorphisms by conjugation.

\item If $D_1\in{\mathcal D}(H)$ is another splitting modulo
$\mathcal K$ that has properties 1-2, then it can be obtained from
$D$ by slidings, conjugations, and modifying boundary
monomorphisms by conjugation. \end{enumerate}
\end{prop}

To obtain an abelian JSJ decomposition of $G$ modulo $\mathcal K$
we first take a free decomposition as described above
$G=G_1*\cdots *G_{\ell}$  and then an abelian JSJ decomposition of
each factor $G_i$ modulo $\{K_{j_1},\ldots ,K_{j_s}\}.$

If $G$ is an $F$-group, we always suppose that $F\leqslant K_1$.

\subsection{NTQ systems and NTQ groups}\label{ntq} We recall now
the definition of a NTQ group from \cite{Imp} and \cite{KMIrc}.

Let $G$ be a group with a generating set $A$. A system of
equations $S = 1$  is called {\em triangular quasiquadratic}
(shortly, TQ) if it can be partitioned into the following
subsystems

\medskip
$S_1(X_1, X_2, \ldots, X_n,A) = 1,$

\medskip
$\ \ \ \ \ S_2(X_2, \ldots, X_n,A) = 1,$

$\ \ \ \ \ \ \ \ \ \  \ldots$

\medskip
$\ \ \ \ \ \ \ \ \ \ \ \ \ \ \ \ S_n(X_n,A) = 1$

\medskip \noindent
 where for each
$i$ one of the following holds:
\begin{enumerate}
\item [1)] $S_i$ is quadratic  in variables $X_i$;
 \item [2)] $S_i= \{[y,z]=1, [y,u]=1 \mid y, z \in X_i\}$ where $u$ is a
group word in $X_{i+1} \cup  \cdots \cup X_n \cup A$ such that its
canonical image  in $G_{i+1}$ is not a proper power. In this case
we say that $S_i=1$ corresponds to an extension of a centralizer;
 \item [3)] $S_i= \{[y,z]=1 \mid y, z \in X_i\}$;
 \item [4)] $S_i$ is the empty equation.
  \end{enumerate}

Define $G_{i}=G_{R(S_{i}, \ldots, S_n)}$ for $i = 1, \ldots, n$
and put $G_{n+1}=G.$ The  TQ system $S = 1$ is called {\em
non-degenerate} (shortly, NTQ) if each system  $S_i=1$, where
$X_{i+1}, \ldots, X_n$ are viewed as the corresponding constants
from $G_{i+1}$ (under the canonical maps $X_j \rightarrow
G_{i+1}$, $j = i+1, \ldots, n$, has a solution in $G_{i+1}$. The
coordinate group of an NTQ system is called an {\em NTQ group}.

 An NTQ system $S = 1$ is called {\em regular} if each non-empty quadratic
equation in $S_i$ is regular (see Definition \ref{regular}).

We say that an NTQ system   $S(X)=1$ is in {\em standard form} if
all quadratic equations in item 1) of the definition are
 in standard form. Clearly, every NTQ system $S(X)=1$ is rationally
 equivalent to a unique NTQ system
 in standard form, namely, there exists an automorphism
 $\phi\in Aut F[X]$ such that $S^{\phi}=1$
 is a standard NTQ system.
Let $S=1$ be a standard NTQ system. Denote by $S_R$ the set of all
regular quadratic equations in $S=1$

 Let $n_1\geqslant n_2\geqslant \ldots \geqslant n_k$ be the sequence of sizes of  equations from
$S_R$ in the decreasing order. Then the tuple $rsize
(S)=(n_1,\ldots ,n_k)$ is called the {\em regular size} of the
system $S=1$. We compare sizes of systems lexicographically from
the left.

\subsection{Rational equivalence}

Recall, that two systems of equations $S(X) = 1$ and $T(Y) = 1$
with coefficients from $F$ are {\it rationally equivalent} if
there are polynomial maps $Y = P_X(X)$ and $X = P_Y(Y)$ such that
the restriction $P_S$ of the map $P_X:F^{|X|} \rightarrow F^{|Y|}$
onto the algebraic set $V_F(S)$ gives a bijection $P_S:V_F(S)
\rightarrow V_F(T)$, and the restriction $P_T$ of the map  $P_Y$
onto the algebraic set $V_F(T)$ gives the inverse of $P_S$.
Sometimes we refer to the maps $P_X$ and $P_Y$ as to change of
coordinates.  It was shown in (see \cite{BMR}) that systems  $S(X)
= 1$ and $T(Y) = 1$ are rationally equivalent if and only if their
coordinate groups $F_{R(S)}$ and  $F_{R(T)}$ are isomorphic as
$F$-groups. Observe that if $f: F_{R(S)} \rightarrow F_{R(T)}$ is
an $F$-isomorphism, then $X = f(X)$ and $Y = f^{-1}(Y)$ are
corresponding change of coordinates for the systems $S(X) = 1$ and
$T(Y) = 1$.

\begin{definition}  A system $S(X) = 1$ with coefficients in a
group $G$ {\em splits}  if there exist a nontrivial  partition of
$X$ into $k >  1$ disjoint subsets $X = X_1 \cup \cdots \cup X_k$,
and elements $S_i(X_i) \in G[X_i] $ such that $S(X) = 1$ is the
union of the systems $S_i(X_i) = 1$, $i = 1, \dots, k.$
\end{definition}

We say that a system $S(X) = 1$ with coefficients in $G$ {\em
splits up to the rational equivalence} or  {\em rationally splits}
if some system $T(Y) = 1$ which is rationally  equivalent to $S(X)
= 1$ splits.

Notice that $S(X) = 1$ rationally splits
 if and only if its coordinate  group $G_{R(S)}$ is a nontrivial  free product
 of the coordinate groups $G_{R(S_i)}$ with the group of constants $G$
amalgamated. In the case when the systems $S_i(X_i) = 1$  are
coefficient free for $i \neq 1$, one has $G_{R(S)}\simeq G_1 \ast
\cdots \ast G_n.$

\subsection{Canonical automorphisms}

Let $G=A*_{C}B$  be an elementary abelian splitting of $G$.
  For $c\in C$ we  define an automorphism $\phi_c :G\rightarrow G$
 such that $\phi_c(a)=a$ for $a\in A$ and $\phi_c(b)=b^{c}=c^{-1}bc$ for  $b\in B$.

If $G=A*_{C}=\langle A,t\mid c^{t}=c', c\in C\rangle$ then for $c
\in C$  define $\phi_c :G\rightarrow G$ such that  $\phi_c (a)=a$
for $a\in A$ and  $\phi_c (t)=ct$.

We call $\phi_c$ a {\em Dehn twist} obtained from the
corresponding elementary abelian splitting of $G$. If $G$ is an
$F$-group, where $F$ is a subgroup of one of the factors  $A$ or
$B$, then Dehn twists that fix elements of the free group
$F\leqslant A$ are
 called  {\em canonical Dehn twists}. Similarly, one can define
 canonical Dehn twists with respect to an arbitrary  fixed subgroup
 $K$ of $G$.

\begin{definition}
Let  $D \in {\mathcal D}(G)$ [$D \in {\mathcal D}_F(G]$ be an
abelian splitting of a group $G$ and $G_v$ be either a $QH$ or an
abelian vertex of $D$. Then an automorphism $\psi \in {\rm
Aut}(G)$ is called a canonical automorphism corresponding to the
vertex $G_v$ if $\psi$ satisfies the following conditions:
\begin{enumerate}
\item [1)] $\psi$ fixes element-wise all other vertex groups in
$D$ (hence fixes all the edge groups);

\item [2)]  if $G_v$ is a $QH$-vertex in $D$, then $\psi$  is  a
Dehn twist [canonical Dehn twists] corresponding to  some
essential ${\mathbb Z}$-splitting of $G$ along a cyclic subgroup
of $G_v$;

\item [3)] if $G_v$ is an abelian subgroup then $\psi$ acts as an
automorphism on $G_v$ which fixes  all the edge subgroups of
$G_v$.\end{enumerate}
\end{definition}

\begin{definition}
Let  $D \in {\mathcal D}(G)$ [$D \in {\mathcal D}_F(G)]$ be an
abelian splitting of a group $G$ and $e$ an edge in $D$. Then an
automorphism $\psi \in {\rm Aut}(G)$ is called a canonical
automorphism corresponding to the edge $e$ if $\psi$ is a Dehn
twist [canonical Dehn twist] of $G$ with respect to the elementary
splitting of $G$ along the edge $e$ which is induced from $D$.
\end{definition}

 \begin{definition}
Let  $D  \in {\mathcal D}(G)$ [$D \in {\mathcal D}_F(G]$ be an
abelian splitting of a  group [$F$-group] $G$. Then the  canonical
group of automorphisms $A_D = A_D(G)$ of $G$ with respect to $D$
is the subgroup of ${\rm Aut}(G)$ generated by all canonical
automorphisms of $G$ corresponding to all edges, all $QH$
vertices, and all abelian vertices of $D$.
\end{definition}

\subsection{Canonical automorphisms of QH-subgroups (orientable
case)}

In this section we discuss  some canonical automorphisms of
QH-subgroups. Let $G\in {\mathcal F}.$

Let the relation corresponding to a QH subgroup $Q$ in some
decomposition of $G$ be
$$S=\prod _{i=1}^n[ x_i, y_i]\prod _{j=1}^m c_j^{z_j}d^{-1}=1.$$
We define the basic sequence
$$\Gamma = (\gamma_1, \gamma_2, \ldots, \gamma_{K(m,n)})$$
of canonical automorphisms of the group $G$ corresponding to $Q$.
Let $$X=\{x_1,y_1,\dots ,x_n,y_n, z_1,\dots ,z_m\}.$$ We assume
that each $\gamma \in \Gamma$ acts identically on all the
generators from $X$ that are not mentioned in the description of
$\gamma$.

\medskip \noindent
Let $m \geqslant 1, n = 0$. In this case $K(m,0) = m-1.$ Put

\smallskip
$\gamma _{i} \ \ \ :\ z_i\rightarrow
z_i(c_i^{z_i}c_{i+1}^{z_{i+1}}), \ \ \ z_{i+1}\rightarrow
z_{i+1}(c_i^{z_i}c_{i+1}^{z_{i+1}}),$ \ \ \ for $i=1,\ldots ,m-1$.

Notice that this is a Dehn twists corresponding to the splitting
of $G$ as an HNN extension with the cyclic edge group $\langle
p_ip_{i+1}\rangle$ belonging to $Q$.

\medskip \noindent
Let $m = 0, n\geqslant 1$. In this case $K(0,n) = 4n-1.$ Put

\medskip
$\gamma _{4i-3}:\ y_i\rightarrow x_iy_i, $ \ \ \ for $i=1,\ldots
,n$ (this is a Dehn twist corresponding to the splitting of $G$ as
an HNN extension with the edge group $\langle x_i\rangle$ and
stable letter $y_i$);

\smallskip
$\gamma _{4i-2}:\ x_i\rightarrow y_ix_i, $ \ \ \ for $i=1,\ldots
,n$;

\smallskip
$\gamma _{4i-1}:\ y_i\rightarrow x_iy_i, $ \ \ \  for $i=1,\ldots
,n$;

\smallskip
$\gamma _{4i} \ \ \ :\ x_i\rightarrow (y_ix_{i+1}^{-1})^{-1}x_i,\
\ \ y_i\rightarrow y_i^{y_ix_{i+1}^{-1}},\ \ \  x_{i+1}\rightarrow
x_{i+1}^{y_ix_{i+1}^{-1}}, \ \ \ y_{i+1}\rightarrow
(y_ix_{i+1}^{-1})^{-1}y_{i+1},$

\smallskip
  for  $i=1,\ldots ,n-1$.

 \medskip \noindent
 Let $m \geqslant 1, n\geqslant 1$. In this case $K(m,n) = m + 4n-1.$ Put

\smallskip
$\gamma _{i} \ \ \ :\ z_i\rightarrow
z_i(c_i^{z_i}c_{i+1}^{z_{i+1}}), \ \ \ z_{i+1}\rightarrow
z_{i+1}(c_i^{z_i}c_{i+1}^{z_{i+1}}),$ \ \ \ for $i=1,\ldots ,m-1$;

\smallskip
 $\gamma _{m}\ \ \ :\ \ \ z_m\rightarrow z_m(c_m^{z_m}x_1^{-1}),\ \ \
x_1\rightarrow x_1^{c_m^{z_m}x_1^{-1}},\ \ \  y_1\rightarrow
(c_m^{z_m}x_1^{-1})^{-1}y_1;$

\medskip
$\gamma _{m + 4i-3}:\ y_i\rightarrow x_iy_i, $ \ \ \ for
$i=1,\ldots ,n$;

\smallskip
$\gamma _{m + 4i-2}:\ x_i\rightarrow y_ix_i, $ \ \ \ for
$i=1,\ldots ,n$;

\smallskip
$\gamma _{m + 4i-1}:\ y_i\rightarrow x_iy_i, $ \ \ \  for
$i=1,\ldots ,n$;

\smallskip
$\gamma _{m + 4i} \ \ \ :\ x_i\rightarrow
(y_ix_{i+1}^{-1})^{-1}x_i,\ \ \ y_i\rightarrow
y_i^{y_ix_{i+1}^{-1}},\ \ \  x_{i+1}\rightarrow
x_{i+1}^{y_ix_{i+1}^{-1}}, \ \ \ y_{i+1}\rightarrow
(y_ix_{i+1}^{-1})^{-1}y_{i+1},$

\smallskip
  for  $i=1,\ldots ,n-1$ (this is a Dehn twist corresponding to the splitting of $G$ with an edge
group  $\langle y_ix_{i+1}^{-1}\rangle$ belonging to $Q$).

\medskip
Observe, that in the case $m \neq 0, n\neq 0$  the basic sequence
of automorphisms $\Gamma$
 contains the basic automorphisms from the other two cases. This allows
us, without loss of
 generality, to formulate the  results below only for the case $K(m,n) =
m + 4n - 1$.
 Obvious adjustments provide the proper argument in the other cases.

The following lemma describes the action of powers of basic
automorphisms from
 $\Gamma$ on $X$. The proof is obvious, and we omit it.

\begin{lemma}
\label{le:7.1.29} Let $\Gamma = (\gamma_1, \ldots,
\gamma_{m+4n-1})$ be the basic sequence of automorphisms and $p$
be a  positive integer. Then the following hold:

\medskip
$\gamma _{i}^p \ \ \ :\ z_i\rightarrow
z_i(c_i^{z_i}c_{i+1}^{z_{i+1}})^p, \ \ \ z_{i+1}\rightarrow
z_{i+1}(c_i^{z_i}c_{i+1}^{z_{i+1}})^p,$ \ \ \ for $i=1,\ldots
,m-1$;

\medskip
 $\gamma _{m}^p \ \  :\ z_m\rightarrow z_m(c_m^{z_m}x_1^{-1})^p,\ \ \
x_1\rightarrow x_1^{(c_m^{z_m}x_1^{-1})^p},\ \ \ y_1\rightarrow
(c_m^{z_m}x_1^{-1})^{-p}y_1;$

\medskip
 $\gamma _{m+4i-3}^p :\ y_i\rightarrow x_i^py_i, $ \ \ \ for
$i=1,\ldots ,n$;

\medskip
$\gamma _{m+4i-2}^p :\ x_i\rightarrow y_i^px_i, $ \ \ \ for
$i=1,\ldots ,n$;

\medskip
$ \gamma _{m+4i-1}^p :\ y_i\rightarrow x_i^py_i, $ \ \ \ for
$i=1,\ldots ,n$;

\medskip
$\gamma _{m+4i}^p \ \ \ :\ x_i\rightarrow
(y_ix_{i+1}^{-1})^{-p}x_i,\ y_i\rightarrow
y_i^{(y_ix_{i+1}^{-1})^p},\
 x_{i+1}\rightarrow x_{i+1}^{(y_ix_{i+1}^{-1})^p},\ \ \
y_{i+1}\rightarrow
 (y_ix_{i+1}^{-1})^{-p}y_{i+1},$

\medskip
for $i=1,\ldots ,n-1$.
 \end{lemma}

Now we introduce  vector notations for automorphisms of particular
type.

 Let ${\mathbb N}$ be the set of all positive integers and  ${\mathbb N}^k$
 the set of all $k$-tuples of elements from ${\mathbb N}$.
 For $s \in {\mathbb N}$ and $p \in {\mathbb N}^k$ we say that the tuple $p$ is
{\it $s$-large} if every coordinate of $p$  is greater then  $s$.
Similarly, a subset $P \subset  {\mathbb N}^k$ is {\em $s$-large}
if every tuple in $P$ is $s$-large.   We say that the set $P$ is
{\em unbounded} if for any $s \in  {\mathbb N}$ there exists an
$s$-large tuple in $P$.

Let $\delta = (\delta_1, \ldots, \delta_k)$ be a sequence of
automorphisms of the group $F_{R(S)}$, and $p = (p_1, \ldots,p_k)
\in {\mathbb N}^k$.  Then by $\delta^p$ we denote the following
automorphism of $F_{R(S)}$: $$\delta^p = \delta_1^{p_1} \cdots
\delta_k^{p_k}.$$

\begin{notation}
  Let $\Gamma = (\gamma_1, \dots, \gamma_K)$ be the
 basic sequence of automorphisms for $S = 1$. Denote by
$\Gamma_{\infty}$ the
 infinite periodic sequence with period $\Gamma$, i.e.,
 $\Gamma_{\infty} = \{\gamma_i\}_{i \geqslant 1}$ with $\gamma_{i+K} =
\gamma_i$.
 For $j \in {\mathbb N}$ denote by $\Gamma_j$ the initial segment of
 $\Gamma_{\infty}$ of the length $j$. Then for a given $j$ and  $p \in
{\mathbb N}^j$ put
 $$ \phi_{j,p} =
\stackrel{\leftarrow}{\Gamma}_j^{\stackrel{\leftarrow}{p}} =
\gamma_j^{p_j}\gamma_{j-1}^{p_{j-1}} \cdots \gamma _1^{p_1}.$$
Sometimes we omit $p$ from $ \phi_{j,p}$ and write simply
$\phi_j$.
\end{notation}

\subsection{Canonical automorphisms of  QH subgroups (non-orientable case).}

Similarly we introduce the notion of the  basic sequence of
automorphisms for a non-orientable QH subgroup. It is more
convenient to consider a non-orientable relation in the form

\begin{equation}\label{8'}
S=\prod_{i=1}^{m}z_i^{-1}c_iz_i\prod_{i=1}^{n}[x_i,y_i]x_{n+1}^2 =
c_1\cdots c_m\prod_{i=1}^{n}[a_i,b_i]a_{n+1}^2,
\end{equation}

or
\begin{equation}\label{8''}
S=\prod_{i=1}^{m}z_i^{-1}c_iz_i\prod_{i=1}^{n}[x_i,y_i]x_{n+1}^2x_{n+2}^2
= c_1\cdots c_m\prod_{i=1}^{n}[a_i,b_i]a_{n+1}^2a_{n+2}^2.
\end{equation}

Without loss of generality we consider equation (\ref{8''}).
 We define a basic sequence
$$\Gamma =(\gamma_1, \gamma_2, \dots, \gamma_{K(m,n)})$$
 of automorphisms of $F_{R(S)}.$
We assume that
 each $\gamma \in \Gamma$    acts identically on all the generators
from $X$ that are not mentioned in the description of $\gamma$.
Automorphisms $\gamma _i,\ i=1,\ldots ,m+4n-1$ are the same as in
the orientable case.

\medskip \noindent
Let $ n = 0$. In this case $K=K(m,0) = m+2.$ Put

\medskip
 $\gamma _{m}: \
z_m\rightarrow  z_m(c_m^{z_m}x_1^2),\ \ \
 x_1\rightarrow x_1^{(c_m^{z_m}x_1^2)};$

\smallskip
$\gamma _{m+1}:\ x_1\rightarrow x_1(x_1x_{2}),\ \ \
x_{2}\rightarrow (x_1x_{2})^{-1}x_{2};$

\smallskip
$\gamma _{m+2}:\ x_1\rightarrow x_1^{(x_1^2x_{2}^2)},\ \ \
x_{2}\rightarrow x_{2}^{(x_1^2x_{2}^2)}.$

\medskip \noindent
Let $ n\geqslant 1$. In this case $K=K(m,n) = m+4n+2.$ Put

\medskip
 $\gamma _{m+4n}: \
x_n\rightarrow (y_nx_{n+1}^2)^{-1}x_n,\ \ \
 y_n\rightarrow y_n^{(y_nx_{n+1}^2)},\ \ \
 x_{n+1}\rightarrow x_{n+1}^{(y_nx_{n+1}^2)};$

\smallskip
$\gamma _{m+4n+1}:\ x_{n+1}\rightarrow x_{n+1}(x_{n+1}x_{n+2}),\ \
\ x_{n+2}\rightarrow (x_{n+1}x_{n+2})^{-1}x_{n+2};$

\smallskip
$\gamma _{m+4n+2}:\ x_{n+1}\rightarrow
x_{n+1}^{(x_{n+1}^2x_{n+2}^2)},\ \ \ x_{n+2}\rightarrow
x_{n+2}^{(x_{n+1}^2x_{n+2}^2)}.$

A family of automorphisms
$$\Gamma _P=\{\phi _{j,p}, j\in N_0, p\in P\}$$ is called {\em
positive unbounded} if $P$ and $N_0\subseteq {\mathbb N}$ are
unbounded. A family of homomorphisms $\Gamma _P\beta$ from
$F_{R(S)}$ onto $F$, where $\beta$ is a solution in general
position, and $\Gamma _P$ is positive unbounded, is a generic
family (see Definition \ref{generic} below).

\subsection{Minimal solutions and maximal standard quotients}
 Let $G$ and $K$  be  $H$-groups and  $A \leqslant Aut_H(G)$ a group
 of $H$-automorphisms of $G$. Two
$H$-homomorphisms $\phi$ and $\psi$ from $G$ into $K$ are called
{\em $A$-equivalent} (symbolically, $\phi \sim_A \psi$) if  there
exists $\sigma\in A$ such that $\phi= \sigma\psi$ (i.e., $g^\phi =
(g^\sigma)^\psi$ for $g \in G$). Obviously, $\sim_A$ is an
equivalence relation on $Hom_H(G,K)$.

 Let  $D$ be a fixed  abelian splitting of $G$ and $A_D = A_D(G)$.

\begin{definition} Let $S(X,A)=1$ be a system of  equations over $F=F(A)$,
$G=F_{R(S)},$ and $D\in{\mathcal D}(G)$. Suppose that $\phi
_1,\phi _2:F_{R(S)}\rightarrow F(A\cup Y)$ are solutions of $S=1$
in  a free group $\bar F = F(A\cup Y)$.
 We write $\phi _1 <_D \phi _2$ if there exists an
automorphism $\sigma\in A_{D}$ and an endomorphism $\pi \in
Hom_F(\bar F, \bar F)$   such that $\phi _2=\sigma^{-1} \phi
_1\pi$ (see Fig. 1) and
 $$\sum_{x \in X} | x^{\phi _1}| < \sum_{x \in X} | x^{\phi _2}|.$$
  \end{definition}

\begin{figure}[here]
\centering{\mbox{\psfig{figure=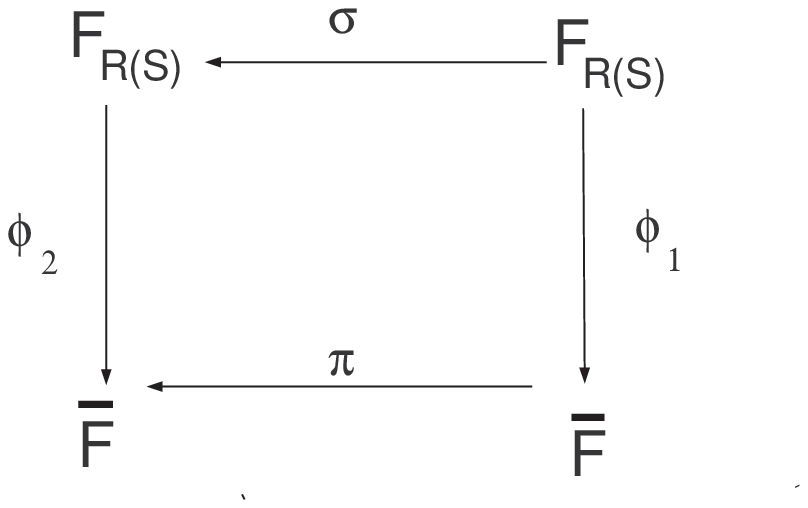,height=2in}}} \caption{$\phi
_1 <_D \phi _2$} \label{ET0.}
\end{figure}
  Notice that if $Y = \emptyset$ then $\phi_1 < \phi_2$ if and only if
  $\phi_1 \sim_{A_D}\phi _2$ and
  $\sum_{x \in X} |x^{\phi _1}| < \sum_{x \in X} |x^{\phi _2}|.$
  This provides a way to chose representatives in  $\sim_{A_D}$-equivalence
  classes of automorphisms. A solution $\phi :G\rightarrow F$ is
  called {\em minimal} if $\phi$ is $<_{D}$-minimal in its
  ${\sim}_{A_D}$-equivalence class.

 Let $R_D$ be the intersection of
the kernels of all minimal  homomorphisms from $Hom_F(G,F)$. Then
$G/R_D$ is called the  {\em maximal standard quotient} of $G$
relative to $D$.

\begin{lemma} Assume that $D$ is a fixed abelian splitting of $G$. Then
every $F$-homomorphism from $G$ onto $F$ can be presented as
composition of a canonical {\rm(}relative to $D${\rm)}
automorphism of $G$, the canonical epimorphism  $G \rightarrow
G/R_D$, and an $F$-homomorphism from $G/R_D$ onto $F$.
\end{lemma}

\begin{lemma}\label{le:dir} Let $G$ be a freely indecomposable $F$-group in
${\mathcal F}$, and $D$ be an abelian JSJ decomposition of $G$.
Then $A_D$ is a direct product of abelian groups generated by
canonical Dehn twists corresponding to edges of $D$ between non-QH
non-abelian subgroups, and groups of canonical automorphisms
corresponding to MQH and abelian vertex groups.\end{lemma}

\begin{proof} Consider two edges $e_1, e_2$ of $D$ such that $e_1\in T$ and $e_2\not\in T$.
Collapsing all the other edges we obtain $G=\langle H*_{\langle
C\rangle }K, t \mid A^t=B\rangle ,$ where $F\leqslant H$,
$A\leqslant H$, $B\leqslant K.$ We write down the action of
canonical Dehn twists corresponding to the edges of this
decomposition on the generators of $G$. Let $c\in C, b\in B.$ Then
$$x^{\phi _c}=x,\ x\in H,\ \ y^{\phi _c}=y^c,\ y\in K,\ \ t^{\phi _c}=tc,$$
$$x^{\phi _b}=x,\ x\in H,\ \ y^{\phi _b}=y,\ y\in K,\ \ t^{\phi _b}=tb.$$
In this case, $\phi _c\phi _b=\phi _b\phi _c.$


All the other possible cases can be verified similarly.\end{proof}

\subsection{Sufficient splittings}\label{se:suffs}
Let $H\leqslant G$. A family of $H$-homomorphisms
$$\Psi
=\{\psi :G\rightarrow H\}$$ is called {\em separating} if for any
nontrivial $g\in G$ there exists $\psi\in\Psi$ such that $\psi
(g)\neq 1$ in $H$.

\begin{definition}
Let $G$ be a group, ${\mathcal K} = \{ K_1,\ldots, K_n\}$  be  a
set of subgroups of $G$ and $F$ is a fixed subgroup of $K_1$. We
say that there is a sufficient splitting
 of
$G$ modulo $\mathcal K$ if one of the following holds:
\begin{enumerate}
\item [1)] $G$ is freely decomposable modulo $\mathcal K$, \item
[2)] there exists a reduced $\mathcal K$-compatible free
decomposition of $G$ in which at least one  factor $G_j$ has   an
abelian splitting $D$ such that all subgroups from  ${\mathcal
K}_j$ are elliptic in $D$ and such that $Rep_{{\mathcal
K}_j}(G_j,D)$ is not a separating family of homomorphisms
($F$-homomorphisms, if $j = 1$) from $G_j$ into $F$.
\end{enumerate}
\end{definition}

In the case of a freely indecomposable group $G$, $D$ is a
sufficient splitting of $G$ if and only if the standard maximal
quotient $G/R_D$ is a proper quotient of $G$.

\subsection{Maximal standard fully residually free quotients}

Let $G$ be a group with a finite generating set $X$ and $D$  an
abelian splitting of $G$. Denote by  $G/R_D =\langle X \mid
S\rangle$ a presentation of the maximal standard quotient $G/R_D$
relative to $X$. Then the coordinate group $F_{R(S)}$ is a
quotient of $G/R_D$. If $S_1 = 1, \ldots, S_k = 1$  are finite
systems that determine the irreducible components of the algebraic
set $V_F(S)$ then (see \cite{BMR})
 $$R(S) = R(S_1) \cap \cdots \cap R(S_k)$$
 and the radicals $R(S_i)$ are uniquely defined (up to
 reordering). Notice that the coordinate groups $F_{R(S_i)}$ are fully residually free,
 they are called the
 {\it standard maximal fully residually free quotients} of $G$
 relative to $D$.

\begin{lemma} The restriction of the canonical epimorphism  $\pi :G \rightarrow G/R_D$
onto a rigid subgroup of $D$ and onto a subgroup of the abelian
vertex group generated by the images of edge groups, is a
monomorphism.
\end{lemma}
\begin{proof} Let $H$ be a rigid subgroups of $D$. Canonical automorphisms corresponding to $D$ act on
$H$ as conjugation. The set of all representatives of the
equivalence classes of solutions $H\rightarrow F$ with respect to
conjugation is a discriminating set for $H$. Similar argument
works for a subgroup of the abelian vertex group generated by the
images of edge groups.\end{proof}
\begin{cy}\label{cy:maxstq} There exists a  standard maximal fully residually free quotient $K$ of
$G$ such that all the restrictions of the canonical epimorphism
$G\rightarrow K$ onto  rigid subgroups of $D$, onto edge subgroups
in $D$,  and onto the subgroups of abelian vertex groups $A$
generated by the images of all the edge groups of edges adjacent
to $A$, are monomorphisms.
\end{cy}
\begin{proof} Let $K_1,\ldots ,K_s$ be all maximal standard fully
residually free quotients of $G$. Let $\phi _i:G/R_D\rightarrow
K_i$ be the canonical epimorphism. Suppose that for each $i$ there
exists an element $u_i$ in some rigid subgroup of $D$ or in a
subgroup of the abelian vertex group generated by the images of
edge groups such that $\phi _i(u_i)=1$. Since $G$ is fully
residually free, and for any homomorphism from $G$ to $F$ the
image of $u_i$ is a conjugate of $\phi _i(u_i)=1$,  there is some
$i$ such that each homomorphism from $G$ to $F$ satisfies the
equation $u_i=1$.\end{proof}

\section{Algorithms over  fully residually free
groups}\label{algor}

\subsection{Algorithms for equations and coordinate groups}
\label{algor1}

 In this section we collect some results on
algorithmic problems concerning equations over free groups and
their coordinate groups. We assume below that a coordinate group
$G$ is given by a finite system of equations $S(X) = 1$ over $F$
in such a way that $G = F_{R(S)}$.

 \begin{theorem}
 \label{le:0.1} The following statements are true:
 \begin{enumerate}
 \item  There is an algorithm which for a given finite system of equations
$S(X) = 1$  over $F$ and a given group word $w(X)$ in $X\cup A$
determines whether $w(X)$ is equal to $1$ in $F_{R(S)}$ or not;
 \item  There is an algorithm which for a given finite systems of
 equations $S(X) = 1$, $T(X) = 1$ over $F$ decides whether or not
 $R(S) = R(T)$.
 \end{enumerate}
\end{theorem}
\begin{proof} (1) Let $S(X) = 1$ be a finite system of equations over $F$
and $w$ a group word in the alphabet $X \cup A$. Then the
universal sentence
 $$\Phi_{S,w} = \forall X (S(X) = 1 \rightarrow w(X) = 1)$$
 is true in the free group $F$ if and only if $w \in R(S)$.
 Since the universal theory $Th_{\forall}(F)$ is decidable
 \cite{Mak84} one can effectively check whether or not $\Phi_{S,w} \in
 Th_{\forall}(F)$. This proves (1). Now (2) follows immediately from
 (1).
\end{proof}

\begin{theorem}
\label{th:0.4.} {\rm(} \cite{KMIrc}, see also Theorem {\rm
\ref{Ase}} below.{\rm)}
 There is an algorithm which for a given   system of equations $S(X) = 1$ over
 $F$ finds finitely many NTQ-systems $Q_1 = 1, \ldots, Q_n = 1$ over $F$ and
  $F$-homomorphisms $\phi_i : F_{R(S)} \rightarrow F_{R(Q_i)}$
   such that any $F$-homomorphism  $\lambda: F_{R(S)} \rightarrow F$ factors
   through one of the homomorphisms $\phi_i, i = 1, \ldots, n.$
 \end{theorem}

\begin{cy}
\label{co:effH}
 There is an algorithms which for a given finite
system of equations $S(X) = 1$ over $F$  finds finitely many
groups $G_1, \ldots, G_n \in {\mathcal F}$ {\rm(}given by  finite
presentations in generators $X\cup A${\rm)} and epimorphisms
$\phi_i: F_{R(S)} \rightarrow G_i$ such that any homomorphism
$\phi: F_{R(S)} \rightarrow F$ factors through one of the
epimorphisms $\phi_1, \ldots, \phi_n$.
\end{cy}

\begin{proof} The result follows from Theorem \ref{th:0.4.} and
Theorem \ref{le:0.2}.
\end{proof}

\begin{theorem} \cite{KMNull}
\label{th:effLyndon}  Given an NTQ system $Q$ over $F$ one can
effectively find an embedding $\phi:F_{R(Q)} \rightarrow
F^{\mathbb{Z}[t]}$.
\end{theorem}

\begin{theorem}
\label{th:effmonic} Given a finite system of equations $S = 1$
over $F$ and an $F$-homomorphism $\phi:F_{R(S)} \rightarrow H$ of
$F_{R(S)}$ into an NTQ $F$-group $H$ one can effectively decide
whether $\phi$ is a monomorphism  or not.
\end{theorem}

\begin{proof} Let $G = F_{R(S)}$ for a finite system $S = 1$ over
$F$ and $H$ an NTQ $F$-group. Suppose $\phi:F_{R(S)} \rightarrow
H$  is a homomorphism given by the set of images $X^\phi \subseteq
H.$ By Theorem \ref{le:0.2} one can effectively find a finite set,
say $T$, of defining relations of the subgroup $G^\phi = \langle
X^\phi \rangle \leqslant H$ with respect to the generating set
$X^\phi$. It follows that $R(S) \leqslant R(T)$. Clearly, $\phi$
is a monomorphism if and only if $R(S) = R(T)$. The latter can be
checked effectively by Theorem \ref{le:0.1}.
\end{proof}

\begin{theorem} \label{th:effNTQ}
Given a finite irreducible system of equations $S = 1$ over $F$
one can effectively find an NTQ system $Q$ over $F$ and an
embedding $\phi:F_{R(S)} \rightarrow F_{R(Q)}$.
\end{theorem}

\begin{proof} By Theorem \ref{th:0.4.} one can effectively find
finitely many NTQ systems $Q_1, \ldots, Q_n$ over $F$ and
$F$-homomorphisms $\phi_i : F_{R(S)} \rightarrow F_{R(Q_i)}$
   such that any $F$-homomorphism  $\lambda: F_{R(S)} \rightarrow F$ factors
   through one of the homomorphisms $\phi_i, i = 1, \ldots, n.$ It
   is known (see \cite{KMNull,KMIrc, BMR}) that in this case at
   least one of the homomorphisms $\phi_i$ is monic. Now by
   Theorem \ref{th:effmonic} one can effectively check which
   homomorphisms among $\phi_1, \ldots, \phi_n$ are monic. This
   proves the theorem.
   \end{proof}

\begin{theorem}
There is an algorithm which for a given finite irreducible system
of equations $S(X) = 1$ over $F$ finds a finite representation of
the coordinate group $F_{R(S)}$ with respect to the generating set
$X \cup A$.
\end{theorem}

\begin{proof}
By Theorem \ref{th:effNTQ} one can effectively find an NTQ system
$Q$ over $F$ and an embedding $\phi:F_{R(S)} \rightarrow
F_{R(Q)}$. It follows that one can effectively find the  finite
generating set $X^\phi  \cup A$ of the subgroup $F_{R(S)}^\phi$
 of the group $F_{R(Q)}$. Now, since the group $F_{R(Q)}$ belongs
 to ${\mathcal F}$ the result follows immediately from Theorem
 \ref{le:0.2}. \end{proof}

 \begin{cy}\label{cy:rem} For every finite irreducible system of
equations $S=1$ one can  effectively find the radical $R(S)$ by
specifying a finite set of generators of $R(S)$ as a normal
subgroup.\end{cy}

\begin{theorem}[\cite{KMNull}]
   \label{th:0.2}
   There is an algorithm which for a given finite system of equations $S(X) = 1$
   over $F$ finds its irreducible components.
   \end{theorem}
   \begin{proof}  We give here another, more direct, proof of this result.
   By Corollary \ref{co:effH}  one can effectively find finitely many groups
$G_1, \ldots, G_n \in {\mathcal F}$, given by  finite
presentations $\langle X \cup A \mid T_i \rangle$,  and
epimorphisms $\phi_i: F_{R(S)} \rightarrow G_i$ such that any
homomorphism $\phi: GF_{R(S)} \rightarrow F$ factors through one
of the epimorphisms $\phi_1, \ldots, phi_n$. It is not hard to see
that the systems $T_i = 1$ are irreducible over $F$ and
 \begin{equation}
 \label{eq:effirr}
 V_F(S) =
V_F(T_1) \cup \cdots\cup V_F(T_n).
 \end{equation}
  Now by Theorem \ref{le:0.1} one can check effectively whether or
  not $V_F(T_i) = V_F(T_j)$ or $V_F(T_i) = V_F(S)$, thus producing
  all irreducible components of $V_F(S)$.
 \end{proof}

\subsection{Algorithms for finitely generated fully residually free
groups} \label{algor2}

 In this section we collect
some  results on algorithmic problems for finitely generated fully
residually free groups from \cite{KMIrc,Imp,KMRS,MRS2}. We assume
below that groups from $\mathcal F$ are given by finite
presentations. Notice that if $\langle X \mid S \rangle $ is a
finite presentation of a group $G \in {\mathcal F}$ then $S(X) =
1$ can be viewed as a finite irreducible coefficient-free system
of equations over $F$ and $R(S) = ncl(S)$ where the radical $R(S)$
is taken in the free group $F(X)$ (without coefficients from $F$).
This allows one to apply the algorithmic results from the Section
\ref{algor1} to groups from ${\mathcal F}$.

 \begin{theorem}
The word problem is decidable in groups from $\mathcal{F}$.
 \end{theorem}
\begin{proof}
Let $G \in {\mathcal F}$ and $g \in G$. Since $G$ is finitely
presented one can effectively enumerate all consequences of
relators of $G$, so if $g = 1$ then $g$ will occur in this
enumeration. On the other hand, one can effectively enumerate all
homomorphisms $\phi_1, \phi_2, \ldots, $ from $G$ into a given
free group $F$ (say of rank 2). If $g \neq 1$ then, since $G$ is
residually free,  there exists $\phi_i$ such that $\phi_i(g) \neq
1$, which can be  verified effectively (by trying one by one all
the images $\phi_1(g), \phi_2(g), \ldots $). This shows that the
word problem is decidable in $G$, as required.
\end{proof}

\begin{theorem} \cite{KMRS} \label{conj}
The conjugacy problem is decidable in groups from $\mathcal{F}$.
\end{theorem}
Notice, that Theorem \ref{conj} also follows from \cite{Bum},
because finitely generated fully residually free groups are
relatively hyperbolic  \cite{Dah}.

\begin{theorem} \cite{MRS2}
  \label{th:0.4}
  The membership problem is decidable in groups from
  $\mathcal{F}$. Namely, there exists an algorithm which for a group
  $G \in {\mathcal F}$ given
   by a finite presentation $\langle X \mid R\rangle$  and a finite tuple
  of words $h_1(X), \ldots, h_k(X), w(X)$ in the alphabet
  $X^{\pm 1}$ decides whether or not the element $w(X)$
  belongs to the subgroup $\langle h_1(X), \ldots, h_k(X) \rangle$.
  \end{theorem}

Using an analog of Stallings' foldings introduced in \cite{MRS}
for finitely generated subgroups of $F^{{\mathbb Z}[t]},$ one can
obtain  the following results.

\begin{theorem} \cite{KMRS} \label{KMRS1} Let $G \in {\mathcal F}$ and
 $H$ and $K$ finitely
generated subgroups of  $G$ given by finite generating sets. Then
$H\cap K$ is finitely generated, and one can effectively find a
finite set of generators of $H\cap K$.
\end{theorem}

\begin{theorem} \cite{KMRS}  \label{int}
\label{th:1} Let $G \in {\mathcal F}$ and
 $H$ and $K$ finitely
generated subgroups of  $G$ given by finite generating sets.  Then
one can effectively find a finite family ${\mathcal J}_G(H,K)$ of
non-trivial finitely generated subgroups of $G$ (given by finite
generating sets), such that
\begin{enumerate}
\item every $J \in {\mathcal J}_G(H,K)$ is of one of the following
types
$$H^{g_1} \cap K,\ \ \ H^{g_1} \cap C_K(g_2),$$
where $g_1 \in G \smallsetminus H$, $g_2 \in K$, moreover $g_1,\
g_2$ can be found effectively; \item for any non-trivial
intersection $H^g \cap K,\ g \in G \smallsetminus H$, there exists
$J \in {\mathcal J}_G(H,K)$ and $f \in K$ such that
$$H^g \cap K = J^f,$$
moreover $J$ and $f$ can be found effectively.
\end{enumerate}
\end{theorem}

\begin{cy}
\label{cy:conjugate-into} Let $H,K$ be finitely generated
subgroups of finitely generated fully residually free group $G$ .
Then  one can effectively verify whether or not $K$ is conjugate
into $H$, and if it is, then find a conjugator.
\end{cy}

\begin{cy}\label{int1} Let $H,K$ be finitely generated subgroups of
finitely generated fully residually free group $G$, let $H$ be
abelian. Then one can effectively find a finite family ${\mathcal
J}$ of non-trivial intersections $J=H^g\cap K\neq 1$ such that any
non-trivial intersection $H^{g_1}\cap K$ has form $J^k$ for some
$k\in K$ and $J\in {\mathcal J}.$ One can effectively find the
generators of the subgroups from ${\mathcal J}.$
\end{cy}

 \begin{theorem}
  \label{le:0.2.new}
  Given  a group $G \in {\mathcal F}$,   a splitting $D$ of $G$, and a  finitely generated
  freely indecomposable subgroup $H$ of $G$ {\rm(}given by a finite generating set $Y${\rm)} one can effectively   find the
  splitting $D_H$
  of $H$ induced from $D$. Moreover, one can describe  all the vertex and edge groups,
  and homomorphisms,
  which occur in $D_H$,  explicitly as words in generators  $Y$.
   \end{theorem}
 \begin{proof} Following Lemma \ref{Coh} one can construct
effectively the graph of groups for the subgroup $H$ using Theorem
\ref{int}. Indeed, for every vertex $v \in X$ by Theorem \ref{int}
one can find effectively the complete finite family   of
non-trivial subgroups such that each  intersection $H \cap G_v^g$
were $g$ runs over $G$  contains a conjugate of one of them  in
$H$. If all these intersections are trivial then by Lemma
\ref{Coh} the subgroup $H$ is cyclic. Otherwise there exists a
vertex $v \in X$ and an element $g \in G$ such that the subgroup
$H_v = H \cap G_v^g$ is non-trivial. We start building the  graph
of groups $\Gamma_H$ for $H$ induced from $G$ with the graph of
groups $\Gamma_1$ consisting of the vertex $v$ and the subgroup
$H_v$ associated with it.  Now for every edge $e$ outgoing from
$v$ in the graph $X$ we find by Corollary \ref{int1} the complete
finite set (up to conjugation in $H_v$)  of non-trivial
intersections $H_e = H_v \cap G_e^g$ were $g$ runs over $G$.  If
there are no non-trivial intersections of this type then either
$H_v = H$ or $H_v$ is a free factor of $H.$ Since, by Theorem
\ref{th:0.4}  the membership problem is decidable in $G$ one can
effectively check whether $H = H_v$ or not.
Suppose  that $H_e = H_v \cap G_e^g \neq 1$  for some $g \in G$
and an edge $e$ outgoing from $v$.  If $u$ is the terminal vertex
of $e$ then $H_u \neq 1$ and we have reconstructed an edge in the
graph $\Gamma_H$.  Denote by $\Gamma_2$ the graph of groups
obtained from $\Gamma_1$ by adding the edge $e$ and the vertex $u$
with the associated subgroups $H_e$ and $H_u$. If the fundamental
group $\pi(\Gamma_2)$ is equal to $H$ then we are finished,
otherwise we continue as above. In finitely many steps we will get
a graph of groups $\Gamma_k$ such that $\pi(\Gamma_k) = H$. This
proves the theorem.
 \end{proof}

The Theorem \ref{le:0.2.new} allows one to prove an effective
version of Theorem \ref{nt}.
\begin{theorem} \label{th:effextcent}
Let $G$ be a finitely generated group which is given as a finite
sequence
  \begin{equation}
  \label{eq:chain}
  F = G_0 \leqslant G_1 \leqslant \ldots \leqslant G_n = G
   \end{equation}
of extensions of centralizers $G_{i+1} = G_i(u_i,t_i)$. Then given
a finite set of elements $Y \subseteq G$ one can effectively
construct the subgroup $\langle Y \rangle $ generated by $Y$ in
$G$ from free groups by finitely many operations of the following
type:
\begin{enumerate}
\item free products;

\item free products with amalgamation along  cyclic subgroups with
at least one of them being maximal;

\item separated HNN extensions along
 cyclic subgroups with at least one of them being maximal;

\item free extensions of centralizers;
\end{enumerate}
  in such a way that all groups and homomorphisms which occur during this process
  are given explicitly as words in generators  $Y$.
\end{theorem}
\begin{proof}
The result follows from Theorem \ref{le:0.2.new} by induction on
the length of the sequence (\ref{eq:chain}).
\end{proof}

\begin{cy}
There is an algorithm which  for a given finitely generated fully
residually free group $G$ determines whether $G$ is hyperbolic or
not.
\end{cy}
\begin{proof} Let $G$ be a finitely generated fully residually
free group. By Theorems \ref{th:effLyndon} and \ref{th:effNTQ} one
can effectively embed $G$ into a group $G_n$ which is obtained
from a free group $F$ by a finite sequence of extensions of
centralizers as in \ref{eq:chain}. Now, by Theorem
\ref{th:effextcent} one can effectively construct $G$ from $F$ by
finitely many operations of the type (1)--(4). It is known (see,
for example, \cite {KM0}) that in this event  $G$ is hyperbolic if
and only if no operations of the type (4) occurred in the
construction of  $G$. The latter can be checked algorithmically
when the finite sequence of operations is given.
\end{proof}

As  a corollary of the Theorem \ref{le:0.2.new} one can
immediately obtain the following result.

 \begin{theorem}
  \label{le:0.2}
  There is an algorithm which for a given  NTQ system $Q$ over $F$ and
   a finitely generated subgroup $H  \leqslant F_{R(Q)}$, given by a finite
   generating set $Y$,  finds
   a finite presentation for $H$ in the generators $Y$.
   \end{theorem}

\begin{proof} By Theorem \ref{th:effLyndon}  one can effectively
embed $F_{R(Q)}$ into $F^{\mathbb{Z}[t]}$. It follows that one can
effectively embed $F_{R(Q)}$ into a finitely generated group $H$
which is obtained from $F$ by a finite sequence of extensions of
centralizer. Now the result follows from Theorem
\ref{th:effextcent}.
\end{proof}

\begin{theorem}
\label{le:0.2.2}
  There is an algorithm which for a given group $G \in {\mathcal F}$
   and  a finitely generated subgroup $H  \leqslant G$ {\rm(}given by a finite
   generating set $Y${\rm)},  finds
   a finite presentation for $H$ in the generators $Y$.
   \end{theorem}

\begin{proof}
Let $G = \langle X \mid S\rangle$ be a finite presentation of $G$.
As we have mentioned above one can view the relations $S(X) = 1$
as a  coefficient-free finite irreducible system of equations over
$F$ with the coordinate group $G$. By Theorem \ref{th:effNTQ} one
can effectively embed the group $G$  into the coordinate group
$F_{R(Q)}$ of an  NTQ system $Q$ over $F$. Now the result follows
from Theorem \ref{le:0.2}.
 \end{proof}
Observe also, that Theorem \ref{le:0.2} can be derived from
 \cite{KMW}, Theorem 5.8.

 \begin{theorem}
 \label{le:0.2.}
 There is an algorithm which for a given homomorphism $\phi: G  \rightarrow
 H$ between two groups from ${\mathcal F}$  decides
 whether or not:
\begin{enumerate}
 \item $\phi$ is an epimorphism;

 \item $\phi$ is a monomorphism;

 \item $\phi$ is an isomorphism.\end{enumerate}
 \end{theorem}
 \begin{proof}  Let $G, H \in {\mathcal F}$ and $\langle X \mid S \rangle$  a given finite presentation of $G$.

 (1) Obviously, $\phi$ is onto if and only if $G^\phi =
 H$. Observe that $G^\phi$ is generated by a finite set $X^\phi$, so by Theorem \ref{th:0.4} one
  can effectively verify whether $G^\phi = H$ or not, as required.

  (2)  By Theorem \ref{le:0.2} one can
 effectively find a finite presentation, say $\langle X \mid T\rangle $,  of the subgroup $G^\phi$ of $H$
 with respect to the generating set $X^\phi$. Clearly, $\phi$ is monic
 if and only if $R(S)=  R(T)$ which can be effectively verified by Theorem \ref{le:0.1}.
 This proves (2) and the theorem since (3) follows from (1) and (2).

\end{proof}

\begin{theorem} \cite{KMS} The Diophantine problem is decidable in groups from ${\mathcal F}$.
 Namely, there is an algorithm which for  a given group $G \in {\mathcal F}$ and an equation $S = 1$ over $G$
 decides whether or not the equation $S = 1$ has a solution in $G$ {\rm(}and finds a solution if it exists\/ {\rm)}.
\end{theorem}

\section{Generalized equations over free groups}\label{se:ge}
Makanin \cite{Mak82} introduced the concept of a generalized
equation constructed for a finite system of equations in a free
group $F=F(A)$.   Geometrically a generalized equation consists of
three kinds of objects: bases, boundaries and items. Roughly it is
 a long interval with marked division points. The marked division
points are the boundaries. Subintervals between division points
are items (we assign a  variable to each item). Line segments
below certain subintervals, beginning at some boundary and ending
at some other boundary, are bases.  Each base either corresponds
to a letter from $A$ or has a double.

\begin{figure}[here]
\centering{\mbox{\psfig{figure=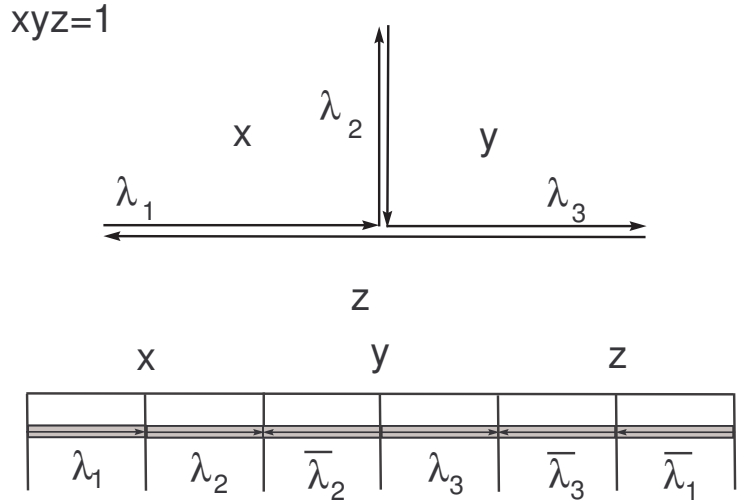,height=3in}}} \caption{From
the cancellation tree for the equation $xyz=1$ to the generalized
equation ($x=\lambda _1\circ\lambda _2,\ y=\lambda
_2^{-1}\circ\lambda _3,\ z=\lambda _3^{-1}\circ\lambda _1^{-1}).$}
\label{ET0}
\end{figure}
This concept becomes crucial to our subsequent work and is
difficult to understand. This is one of the main tools used to
describe solution sets of systems of equations. In subsequent
papers we will use it also to obtain effectively different
splittings of groups. Before we give a formal definition we will
try to motivate it with a simple example.

Suppose we have the simple equation $xyz = 1$ in a free group.
Suppose that we have a solution to this equation denoted by
$x^\phi,y^\phi,z^\phi$ where is $\phi$ is a given homomorphism
into a free group $F(A)$.  Since $x^\phi,y^\phi,z^\phi$ are
reduced words in the generators $A$ there must be complete
cancellation. If we take a concatenation of the geodesic subpaths
corresponding to $x^{\phi}, y^{\phi}$ and $z^{\phi}$ we obtain a
path in the Cayley graph corresponding to this complete
cancellation. This is called a cancellation tree (see Fig. 2). In
the simplest situation $x = \lambda_1\circ \lambda_2$, $y =
\lambda_2^{-1}\circ \lambda_3$ and $z = \lambda_3^{-1}\circ
\lambda_1^{-1}$.  The generalized equation would then be the
following interval.

 The
boundaries would be the division points, the bases are the
$\lambda's$ and the items in this simple case are also the
$\lambda's$.  In a more complicated equation where the variables
$X,Y,Z$ appear more than one time this basic interval would be
extended, Since the solution of any equation in a free group must
involve complete cancellation this drawing of the interval is
essentially the way one would solve such an equation. Our picture
above depended on one fixed solution $\phi$. However for any
equation there are only finitely many such cancellation trees and
hence only finitely many generalized equations.

\subsection{Generalized equations}
\label{se:4-1}

 Let $A= \{a_1, \ldots, a_m\}$ be a set of constants and $X = \{x_1,
\ldots, x_n\}$
  be a set of variables. Put   $G = F(A)$ and  $G[X] = G \ast F(X).$

\begin{definition}
A combinatorial generalized equation $\Omega$ (with constants from
$A^{\pm 1}$) consists  of the following objects:

1. A finite set of \emph{bases} $BS = BS(\Omega)$.  Every base is
either a constant base or a variable base. Each constant base is
associated with exactly one letter from $A^{\pm 1}$. The set of
variable bases ${\mathcal M}$ consists of $2n$ elements ${\mathcal
M} = \{\mu_1, \ldots, \mu_{2n}\}$. The set ${\mathcal M}$ comes
equipped with two functions: a function $\varepsilon: {\mathcal M}
\rightarrow \{1,-1\}$ and an involution $\Delta: {\mathcal M}
\rightarrow {\mathcal M}$ (i.e., $\Delta$ is a bijection such that
$\Delta^2$ is an identity on  ${\mathcal M}$). Bases $\mu$ and
$\Delta(\mu)$ (or $\bar\mu$) are called {\it dual bases}.  We
denote variable bases by $\mu, \lambda, \ldots.$

2.  A set of \emph{boundaries} $BD = BD(\Omega)$. $BD$ is  a
finite initial segment of the set of positive integers  $BD = \{1,
2, \ldots, \rho+1\}$. We use letters $i,j, \ldots$ for boundaries.

3. Two functions $\alpha : BS \rightarrow BD$ and $\beta : BS
\rightarrow BD$. We call $\alpha(\mu)$ and $\beta(\mu)$ the
initial and terminal boundaries of the base $\mu$ (or endpoints of
$\mu$). These functions satisfy the following conditions:
$\alpha(b) <  \beta(b)$  for every base $b \in BS$; if $b$ is a
constant base then $\beta(b) = \alpha(b) + 1$.

4. A finite set of \emph{boundary connections} $BC = BC(\Omega)$.
A boundary connection is a triple $(i,\mu,j)$ where $i, j \in BD$,
$\mu \in {\mathcal M}$ such that $\alpha(\mu) <  i < \beta(\mu)$
and $\alpha(\Delta(\mu)) <  j < \beta(\Delta(\mu))$. We will
assume for simplicity,  that if $(i,\mu,j) \in BC$ then
$(j,\Delta(\mu), i) \in BC$. This allows one to identify
connections $ (i,\mu,j)$ and  $(j,\Delta(\mu), i)$.
\end{definition}

   For a combinatorial generalized equation $\Omega$, one can canonically
associate a system of equations in \emph{variables} $h_1, \ldots,
h_\rho$ over $F(A)$ (variables $h_i$ are sometimes  called {\it
items}). This system is called a \emph{generalized equation}, and
(slightly abusing the language) we  denote it by the same symbol
$\Omega$. The generalized equation  $\Omega$  consists of the
following three types of equations.

1. Each pair of dual variable bases $(\lambda, \Delta(\lambda))$
provides an equation
 $$[h_{\alpha
(\lambda )}h_{\alpha (\lambda )+1}\cdots h_{\beta (\lambda )-1}]^
{\varepsilon (\lambda)}= [h_{\alpha (\Delta (\lambda ))}h_{\alpha
(\Delta (\lambda ))+1} \cdots h_{\beta (\Delta (\lambda ))-1}]^
{\varepsilon (\Delta (\lambda))}.$$ These equations are called
\emph{basic equations}.

2. For each constant base $b$ we write down a \emph{coefficient
equation}
 $$ h_{\alpha(b)} = a,$$
 where $a \in A^{\pm 1}$ is the constant associated with $b$.

3.  Every boundary connection $(p,\lambda,q)$ gives rise to a
\emph{ boundary equation}
 $$[h_{\alpha (\lambda )}h_{\alpha (\lambda
)+1}\cdots h_{p-1}]= [h_{\alpha (\Delta (\lambda ))}h_{\alpha
(\Delta (\lambda ))+1} \cdots h_{q-1}],$$ if $\varepsilon
(\lambda)= \varepsilon (\Delta(\lambda))$ and $$[h_{\alpha
(\lambda )}h_{\alpha (\lambda )+1}\cdots h_{p-1}]= [h_{q}h_{q+1}
\cdots h_{\beta (\Delta (\lambda ))-1}]^{-1} ,$$ if $\varepsilon
(\lambda)= -\varepsilon (\Delta(\lambda)).$

\begin{remark}
We  assume that every generalized equation comes associated with a
combinatorial one.
\end{remark}

\begin{example} Consider as an example the Malcev equation $[x,y][b,a]=1$,
where $a,b\in A.$ Consider the following solution of this
equation:
$$  x^{\phi}=(( b^{n_1}a)^{n_2}b)^{n_3}b^{n_1}a, \ \ y^{\phi}=(
b^{n_1}a)^{n_2}b.$$ Fig. 3 shows the cancellation tree and the
generalized equation for this solution.
\begin{figure}[here]
\centering{\mbox{\psfig{figure=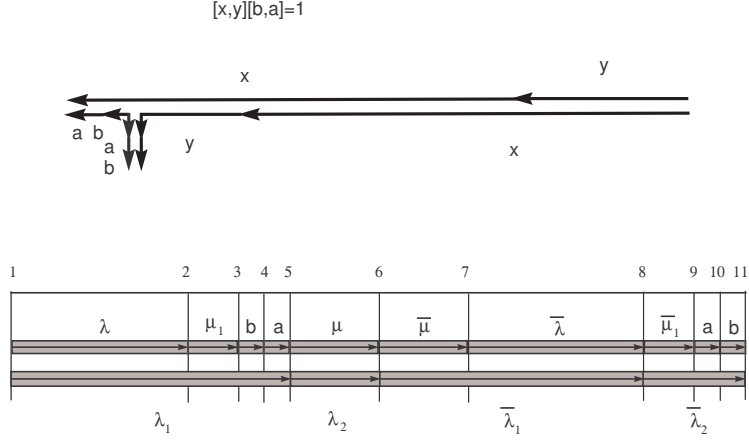,height=3.3in}}} \caption{A
cancellation tree and the generalized equation corresponding to
this tree for the equation $[x,y][b,a]=1.$} \label{ET00}
\end{figure}
This generalized equation has ten variables $h_1,\ldots ,h_{10}$
and eleven boundaries. The system of basic equations for this
generalized equation is the following
$$h_1=h_7, \ h_2=h_8,\ h_5=h_6,\ h_1h_2h_3h_4=h_6h_7,\
h_5=h_8h_9h_{10}.$$ The system of coefficient equations is
$$h_3=b,\ h_4=a,\ h_9=a,\ h_{10}=b.$$
\end{example}

\begin{definition}
Let $\Omega(h) = \{L_1(h)=R_1(h), \ldots, L_s(h) = R_s(h)\}$ be a
generalized equation in variables $h = (h_1, \ldots,h_{\rho})$
with constants from $A^{\pm 1 }$. A sequence of reduced   nonempty
words $U = (U_1(A), \ldots, U_{\rho}(A))$  in the alphabet $A^{\pm
1}$ is a {\em solution} of $\Omega $ if:
\begin{enumerate}
\item [1)] all words $L_i(U)$, $R_i(U)$  are reduced as written;

\item [2)] $L_i(U) =  R_i(U),  \ \ i = 1, \dots, s.$
\end{enumerate}\end{definition} The notation $(\Omega, U)$ means that $U$ is a
solution of the generalized equation $\Omega $.

\begin{remark}
Notice that a solution $U$ of a generalized equation $\Omega$ can
be viewed as  a solution of $\Omega$ in the free monoid
$F_{mon}(A^{\pm 1})$ (i.e., the equalities $L_i(U) = R_i(U)$ are
graphical) which satisfies an additional  condition  $ U \in F(A)
\leqslant F_{\rm mon}(A^{\pm 1})$.
\end{remark}

 Obviously, each solution  $U$ of $\Omega$ gives rise
to  a solution of $\Omega$ in the free group $F(A)$. The converse
does not hold in general, i.e., it might happen that  $U$ is a
solution of $\Omega$ in $F(A)$ but not in $F_{mon}(A^{\pm 1})$,
i.e., all equalities $L_i(U) = R_i(U)$ hold only after a free
reduction but not graphically. We introduce the following notation
which will allow us  to distinguish in which structure
($F_{mon}(A^{\pm 1})$ or $F(A)$) we are looking for solutions  for
$\Omega$.

 If
 $$S =  \{L_1(h)=R_1(h), \ldots, L_s(h) =
R_s(h)\}$$ is an arbitrary system of equations with constants from
$A^{\pm 1}$, then by $S^*$ we denote the system of equations $$S^*
= \{L_1(h)R_1(h)^{-1} = 1, \ldots, L_s(h)R_s(h)^{-1} = 1\}$$  over
the free group $F(A)$.

\begin{definition}
A generalized equations $\Omega$ is called {\em formally
consistent} if it satisfies the following conditions. \be \item
[1)] If $\varepsilon (\mu)=-\varepsilon (\Delta (\mu))$, then the
bases $\mu$ and $\Delta (\mu )$ do not intersect, i.e. non of the
the items $h_{\alpha(\mu)}, h_{\beta(\mu)-1}$ is contained in
$\Delta (\mu )$. \item [2)] If two boundary equations have
respective parameters $(p,\lambda ,q)$ and $(p_1,\lambda ,q_1)$
with $p\leqslant p_1,$ then $q\leqslant q_1$ in the case when
$\varepsilon (\lambda)\varepsilon (\Delta (\lambda))=1,$ and
$q\geqslant q_1$ in the case $\varepsilon (\lambda)\varepsilon
(\Delta (\lambda))=-1,$ in particular, if $p=p_1$ then $q = q_1$.
 \item [3)] Let $\mu$ be
a base such that  $\alpha (\mu)=\alpha (\Delta (\mu))$ (in this
case we say that bases $\mu$ and $\Delta(\mu)$ form a matched pair
of dual bases). If $(p,\mu ,q)$ is  a boundary connection related
to $\mu$ then  $p=q.$

\item  [4)]A variable cannot occur in two distinct coefficient
equations, i.e., any two constant bases with the same left
end-point are labelled by the same letter from $A^{\pm 1}$.

\item [5)]If $h_i$ is a variable from some coefficient equation,
and if $(i,\mu,q_1), (i+1,\mu,q_2)$ are boundary connections, then
$|q_1- q_2|=1.$ \ee
\end{definition}

\begin{lemma}
\label{le:4.1}\ \begin{enumerate} \item If a generalized equation
$\Omega$ has a solution then $\Omega$ is formally consistent;

\item There is an algorithm which for every generalized equation
checks whether it is formally consistent or not.\end{enumerate}
\end{lemma}
The proof is easy and we omit it.
\smallskip

\begin{remark} In the sequel we  consider only formally consistent
generalized equations.
\end{remark}

It is convenient to  visualize  a generalized equation $\Omega$ as
follows.

\begin {center}
\begin{picture}(200,100)(0,0)
\put(10,85){\line(0,-1){80}} \put(30,85){\line(0,-1){80}}
\put(50,85){\line(0,-1){80}} \put(70,85){\line(0,-1){80}}
\put(90,85){\line(0,-1){80}} \put(110,85){\line(0,-1){80}}
\put(130,85){\line(0,-1){80}} \put(150,85){\line(0,-1){80}}
\put(170,85){\line(0,-1){80}} \put(10,80){\line(1,0){160}} \put
(10,87){$1$} \put (30,87){$2$} \put (50,87){$3$} \put
(140,87){$\rho -1$} \put (170,87){$\rho$} \put
(10,65){\line(1,0){60}} \put (70,45){\line(1,0){20}} \put
(90,55){\line(1,0){40}} \put (110,30){\line(1,0){40}} \put
(20,67){$\lambda$} \put (92,57){$\Delta(\lambda )$} \put
(80,47){$\mu$} \put (112,37){$\Delta(\mu)$}
\end{picture}
\end{center}

\subsection{Reduction to generalized equations}
\label{se:4-2}

In this section, following Makanin \cite{Mak82},   we show how,
for a given finite system of equations $S(X,A) = 1$ over a free
group $F(A)$, one can canonically associate a finite collection
of generalized equations ${\mathcal GE}(S)$ with constants from
$A^{\pm 1}$, which to some extent describes all solutions of the
system $S(X,A) = 1$.

 Let $S(X,A) = 1$ be a finite system of equations $S_ 1 = 1, \ldots, S_m
= 1$
  over a free group $F(A)$.  We write $S(X,A) = 1$  in the form
\begin{equation}\label{*}\begin{array}{c}
 r_{11}r_{12}\cdots r_{1l_1}=1,\\
 r_{21}r_{22}\cdots r_{2l_2}=1,\\
 \vdots \\
 r_{m1}r_{m2}\cdots r_{ml_m}=1,\\ \end{array} \end{equation} where $r_{ij}$ are letters in
the alphabet $ X^{\pm 1}\cup A^{\pm 1}.$

A {\it partition table} $T$  of the system above  is  a set of
reduced words
$$T = \{V_{ij}(z_1, \ldots ,z_p)\} \ \ (1\leqslant i\leqslant m, 1\leqslant j\leqslant l_i)$$ from a free
group $F[Z] = F(A \cup Z)$, where $Z = \{z_1,\ldots ,z_p\}$, which
satisfies the following conditions:
\begin{enumerate}
\item [1)] The equality $V_{i1}V_{i2} \cdots V_{il_i}=1,
1\leqslant i\leqslant m,$ holds in  $F[Z]$; \item [2)]
$|V_{ij}|\leqslant l_i - 1$; \item [3)] if $r_{ij} = a \in A^{\pm
1}$, then $V_{ij}= a.$\end{enumerate}

Since $|V_{ij}|\leqslant l_i - 1$ then  at most $|S| = \sum_{i =
1}^m (l_i - 1)l_i$ different  letters $z_i$ can occur in a
partition table of the equation $S(X,A) = 1$. Therefore we will
always assume  that $p \leqslant d(S)$.

Each partition table encodes a particular type of cancellation
that happens when one substitutes  a particular solution $W(A) \in
F(A)$ into  $S(X,A) = 1$ and then freely reduces the words in
$S(W(A),A)$ into the empty word.

\begin{lemma}
\label{le:4.2} Let $S(X,A) = 1$ be a finite system of equations
over $F(A)$. Then \be

\item the set $PT(S)$ of all partition tables  of $S(X,A) = 1$ is
finite, and its cardinality is bounded by a number which depends
only on  $|S(X,A)|$;

\item one can effectively enumerate the set $PT(S)$.
\ee\end{lemma}
\begin{proof} Since the words $V_{ij}$ have  bounded length,  one
can effectively enumerate the finite set  of all collections of
words $\{V_{ij}\}$ in $F[Z]$ which satisfy the conditions 2), 3)
above. Now for each such collection $\{V_{ij}\}$,  one can
effectively check  whether the  equalities $V_{i1}V_{i2} \cdots
V_{il_i}=1, 1\leqslant i\leqslant m$ hold in the free group $F[Z]$
or not. This allows one to list effectively all partition tables
for $S(X,A) = 1$. \end{proof}

To each partition table $T=\{V_{ij}\}$ one can assign a
generalized equation $\Omega _T$ in the following way (below we
use the notation $\doteq $ for graphical equality). Consider the
following word $V$ in $M(A^{\pm 1} \cup Z^{\pm 1}):$
 $$V\doteq V_{11}V_{12}\cdots V_{1l_1}\cdots V_{m1}V_{m2} \cdots  V_{ml_m} = y_1
 \cdots y_\rho, $$
where $y_i \in A^{\pm 1} \cup Z^{\pm 1}$ and $\rho = d(V)$ is the
length of $V$. Then the generalized equation $\Omega_T =
\Omega_T(h)$ has $\rho + 1$ boundaries and $\rho$ variables
$h_1,\ldots ,h_{\rho}$ which are denoted by   $h = (h_1,\ldots
,h_{\rho})$.

Now we define bases of $\Omega_T$ and the functions $\alpha,
\beta, \varepsilon$.

 Let $z \in Z$.  For any two distinct occurrences of  $z$ in $V$ as
 $$y_i = z^{\varepsilon _i}, \ \ \ y_j = z^{\varepsilon _j} \ \ \
 (\varepsilon _i, \varepsilon _j \in \{1,-1\})$$
 we introduce a pair of dual variable bases $\mu_{z,i}, \mu_{z,j}$
such that $\Delta(\mu_{z,i}) = \mu_{z,j}$ (say, if $i < j$). Put
$$\alpha(\mu_{z,i}) = i, \ \ \ \beta(\mu_{z,i}) = i+1, \ \ \
\epsilon(\mu_{z,i}) = \varepsilon _i.$$
 The basic equation that
corresponds to this pair of dual bases is
 $h_{i}^{\varepsilon_i}=h_{j}^{\varepsilon _j} .$

Let $x \in X$.  For any two distinct occurrences of  $x$ in
$S(X,A) = 1$  as
 $$r_{i,j} = x^{\varepsilon_{ij}}, \ \ \ r_{s,t} =
x^{\varepsilon_{st}} \ \ \ (\varepsilon _{ij}, \varepsilon _{st}
\in \{1,-1\})$$
 we introduce a pair of dual bases $\mu_{x,i,j}$ and $\mu_{x,s,t}$
 such that $\Delta(\mu_{x,i,j}) = \mu_{x,s,t}$ (say, if $(i,j) <  (s,t)$
 in the left lexicographic order).  Now let $V_{ij}$ occurs in the
 word $V$ as a subword
 $$V_{ij} = y_c \cdots y_d.$$
 Then we put
$$\alpha(\mu_{x,i,j}) = c, \ \ \ \beta(\mu_{x,i,j}) = d+1, \ \ \
\epsilon(\mu_{x,i,j}) = \varepsilon_{ij}.$$
 The basic equation which
corresponds to these dual bases can be written in the form
 $$[h_{\alpha(\mu_{x,i,j})}\cdots h_{\beta(\mu_{x,i,j})-1}]^{\varepsilon _{ij}}=
 [h_{\alpha(\mu_{x,s,t})}\cdots h_{\beta(\mu_{x,s,t})-
1}]^{\varepsilon _{st}}.$$

Let $r_{ij} = a \in A^{\pm 1}$.  In this case we introduce a
constant base $\mu_{ij}$ with the label $a$. If $V_{ij}$ occurs in
$V$ as $V_{ij} = y_c$, then we put
 $$ \alpha(\mu_{ij}) = c, \beta(\mu_{ij}) = c+1.$$
 The corresponding coefficient equation is written as  $h_{c}=a$.

The list of boundary connections here (and hence the boundary
equations) is empty. This defines the generalized equation
$\Omega_T$. Put
 $${\mathcal GE}(S) = \{\Omega_T \mid T \ \mbox{is   a  partition  table
 for}\ S(X,A)= 1 \}.$$
  Then ${\mathcal GE}(S)$ is a finite collection of generalized
  equations which can be effectively constructed for a given
  $S(X,A) = 1$.

  For a generalized equation $\Omega $ we can also  consider the same system of equations in a free group.
    We denote this system by $\Omega ^*$. By $F_{R(\Omega)}$ we denote the coordinate group of
    $\Omega ^*.$
    Now we explain relations between  the coordinate groups of
  $S(X,A) = 1$ and $\Omega_T^*$.

For a letter $x$ in $X$ we choose an arbitrary   occurrence of $x$
in $S(X,A) = 1$ as
 $$r_{ij} = x^{\varepsilon_{ij}}.$$
  Let $\mu = \mu_{x,i,j}$ be the base that corresponds to this occurrence
  of $x$.  Then $V_{ij}$ occurs  in $V$ as the subword
  $$V_{ij} = y_{\alpha(\mu)} \cdots y_{\beta(\mu) -1}.$$
 Define a word $P_x(h) \in F[h]$ (where $h = \{h_1, \ldots,h_\rho\}$) as
 $$ P_x(h,A) = h_{\alpha(\mu)} \cdots h_{\beta(\mu)-1}^{\varepsilon_{ij}},$$
 and put
 $$P(h) = (P_{x_1}, \ldots, P_{x_n}).$$
 The tuple of words $P(h)$ depends on a choice of occurrences of letters from
  $X$ in $V$. It follows from the construction above that  the map $ X
 \rightarrow F[h]$ defined by $x \rightarrow  P_x(h,A)$
 gives rise to an $F$-homomorphism
 $$\pi : F_{R(S)}\rightarrow F_{R(\Omega _T)}.$$
Observe that the image  $\pi (x)$ in $F_{R(\Omega _T)}$ does not
depend on a particular choice of the occurrence of $x$ in $S(X,A)$
(the basic equations of $\Omega_T$ make these images equal). Hence
$\pi$ depends only on $\Omega_T$.

  Now we relate solutions of $S(X,A) = 1$  with solutions of generalized
  equations  from ${\mathcal GE}(S)$.
  Let $W(A)$ be a solution of $S(X,A) = 1$ in $F(A)$. If in the system
  (\ref{*}) we make the substitution  $\sigma : X \rightarrow W(A)$,
   then
   $$(r_{i1}r_{i2}\cdots r_{il_i})^{\sigma} =
   r_{i1}^{\sigma}r_{i2}^{\sigma}\cdots r_{il_i}^{\sigma} = 1$$
    in $F(A)$ for  every $i = 1, \dots, m$.
Hence every product $R_i = r_{i1}^{\sigma}r_{i2}^{\sigma}\cdots
r_{il_i}^{\sigma}$   can be reduced to the empty word by a
sequence of free   reductions. Let us fix a particular reduction
process for  each $R_i$.   Denote by ${\tilde z}_1, \ldots,
{\tilde z}_p$ all the (maximal) non-trivial subwords of
$r_{ij}^{\sigma}$  that cancel out in some $R_i$ ($i = 1, \ldots
,m$) during the chosen reduction process. Since every word
$r_{ij}^\sigma$ in this process  cancels out completely, that
implies that
$$r_{ij}^\sigma = V_{ij}({\tilde z}_1, \ldots,  {\tilde z}_p)$$
 for some reduced words $V_{ij}(Z)$ in variables $Z = \{z_1, \ldots, z_p\}$.
  Moreover, the equality above is graphical. Observe also that if
  $r_{ij} = a \in A^{\pm 1}$ then $r_{ij}^{\sigma} = a$ and we have
 $V_{ij} = a$. Since every word $r_{ij}^\sigma$
 in $R_i$ has at most one cancellation with any other word
 $r_{ik}^\sigma$
 and does not have cancellation with itself, we have $|V_{ij}| \leqslant l_i - 1$.
  This shows that
   the set $T = \{V_{ij}\}$ is  a partition
table for $S(X,A) = 1$. Obviously,
$$U(A) = ({\tilde z}_1, \ldots, {\tilde z}_p)$$
is the  solution of the generalized equation $\Omega_T$, which is
induced by  $W(A)$. From the construction of the map $P(H)$ we
deduce that $W(A) = P(U(A))$.

The reverse is also true: if $U(A)$ is an arbitrary  solution of
the generalized equation $\Omega_T$, then $P(U(A))$ is a solution
of $S(X,A) = 1$.

We  summarize the discussion above in the following lemma, which
is essentially due to Makanin \cite{Mak82}.

\begin{lemma}
\label{le:R1} For a given  system of equations $S(X,A)=1$ over  a
free group $F = F(A)$,  one can  effectively  construct a finite
set
$${\mathcal GE}(S) = \{\Omega_T \mid T \ is \  a \ partition \ table \
 for\ S(X,A)= 1 \}$$
 of generalized equations  such  that
\begin{enumerate}
\item If the set ${\mathcal GE}(S)$ is empty, then  $S(X,A)= 1$
has no solutions in $F(A)$; \item   for each $\Omega (H) \in
{\mathcal GE}(S)$ and  for each $x \in X$  one can effectively
find  a word $P_x(H,A) \in F[H]$   of length at most $|H|$ such
that the map $x: \rightarrow P_x(H,A)$ ($x \in X$) gives rise to
an $F$-homomorphism $\pi_\Omega: F_{R(S)}\rightarrow
F_{R(\Omega)}$; \item for any solution $W(A)  \in F(A)^n$ of the
system $S(X,A)=1$ there exists $\Omega (H) \in {\mathcal GE}(S)$
and a solution $ U(A)$ of $\Omega(H)$ such that $W(A) = P(U(A))$,
where $P(H) = (P_{x_1}, \dots, P_{x_n})$,
 and this equality is graphical;
\item for any $F$-group $\tilde F$, if a  generalized equation
$\Omega (H) \in {\mathcal GE}(S)$
 has a solution $\tilde U$ in $\tilde F$, then $P(\tilde U)$
is a solution of $S(X,A) = 1$ in $\tilde F$.
\end{enumerate}
\end{lemma}
\begin{cy}
\label{co:R1} In the notations of Lemma {\rm \ref{le:R1}}  for any
solution $W(A)  \in F(A)^n$ of the system $S(X,A)=1$ there exists
$\Omega (H) \in {\mathcal GE}(S)$  and a solution $ U(A)$ of
$\Omega(H)$ such that the following diagram commutes.
\end{cy}
\medskip

\begin{center}

\begin{picture}(100,100)(0,0)
\put(0,100){$F_{R(S)}$} \put(100,100){$F_{R(\Omega)}$}
\put(0,0){$F$} \put(15,103){\vector(1,0){80}}
\put(5,93){\vector(0,-1){78}} \put(95,95){\vector(-1,-1){80}}
\put(-10,50){$\pi_{W}$} \put(55,45){$\pi _{ U}$}
\put(50,108){$\pi$}
\end{picture}

\end{center}

\medskip

\subsection{Generalized equations with parameters}
\label{se:parametric}

In this section, following \cite{Razborov} and \cite{KMIrc},  we
consider generalized equations with {\it parameters}. These kinds
of equations appear naturally in Makanin's type rewriting
processes and provide a convenient tool to organize induction
properly.

Let $\Omega$ be a generalized equation. An  item $h_i$ {\it
belongs} to a base $\mu$ (and, in this event, $\mu$ {\it contains}
$h_i$) if $\alpha (\mu)\leqslant i\leqslant \beta (\mu)-1.$ An
item $h_i$ is {\it constant} if it belongs to a constant base,
$h_i$ is {\it free } if it does not belong to any base. By
$\gamma(h_i) = \gamma_i$ we denote the number of bases which
contain $h_i$. We call $\gamma_i$ the {\it degree } of $h_i$.

A boundary $i$ {\it crosses} (or {\it intersects\/})  the base
$\mu$ if $\alpha (\mu)< i < \beta (\mu).$ A boundary $i$ {\it
touches} the base $\mu$ (or $i$ is an end-point of $\mu$)  if
$i=\alpha (\mu)$ or $i=\beta (\mu)$. A boundary is said to be {\em
open} if it crosses at least one base, otherwise it is called {\it
closed}. We say that a boundary $i$ is {\it tied} (or {\it bound})
by a base $\mu$ (or {\it $\mu$-tied})  if there exists a boundary
connection $(p,\mu,q)$ such that $i = p$ or $i = q$. A boundary is
{\em free} if it does not touch any base and it is not tied by a
boundary connection.

 A set of consecutive  items $[i,j] =
\{h_i,\ldots, h_{i+j-1}\}$   is called a {\em section}.
 A section
is said to be {\em closed} if the boundaries $i$ and $i+j$ are
closed and all the boundaries between them are open. A base $\mu$
is {\it contained } in a base $\lambda$ if $\alpha(\lambda)
\leqslant \alpha(\mu) <  \beta(\mu) \leqslant \beta(\lambda)$. If
$\mu$ is a base then by $\sigma(\mu)$ we denote the section
$[\alpha(\mu),\beta(\mu)]$ and by $h(\mu)$ we denote the product
of items $h_{\alpha(\mu)}\cdots h_{\beta(\mu)-1}$. In general for
a section $[i,j]$ by $h[i,j]$ we denote the product  $h_i \cdots
h_{j-1}$.

\begin{definition}
\label{de:gepar} Let $\Omega$ be a generalized equation. If the
set $\Sigma = \Sigma_\Omega$ of all closed sections  of $\Omega$
is partitioned into a disjoint union of subsets
  \begin{equation}
  \label{eq:ge1-0}
\Sigma_\Omega = V\Sigma\cup  P\Sigma\cup C\Sigma ,
\end{equation}
then $\Omega $ is called a {\em generalized equation with
parameters} or a {\it parametric} generalized equation. Sections
from $V\Sigma, P\Sigma$, and $C\Sigma$ are called correspondingly,
{\em variable, parametric}, and {\em constant} sections. To
organize the branching process properly, we usually divide
variable sections into two disjoint parts:
 \begin{equation}
 \label{eq:geq1-1}
 V\Sigma = A\Sigma\cup NA\Sigma
\end{equation}
Sections from $A\Sigma$ are called {\em active}, and sections from
$NA\Sigma$ are {\em non-active}. In the case when partition
(\ref{eq:geq1-1}) is not specified we assume that $A\Sigma =
V\Sigma$. Thus, in general, we have a partition
  \begin{equation}
  \label{eq:ge1}
\Sigma_\Omega = A\Sigma\cup NA\Sigma\cup P\Sigma\cup C\Sigma
\end{equation}
If $\sigma \in \Sigma$, then every base or item from $\sigma$ is
called active, non-active, parametric, or constant, with respect
to the type of $\sigma$.
\end{definition}

We will see later that every parametric generalized equation can
be written in a particular {\it standard} form.

\begin{definition}
 We say that a parametric generalized equation $\Omega$ is {\em in a standard
 form}
 if the following conditions hold:
\begin{enumerate}
\item [1)]  all non-active sections from $NA\Sigma_\Omega$ are
located to the right of all active sections from $A\Sigma$, all
parametric sections from $P\Sigma_\Omega$ are located to
 the right  of all non-active sections, and all constant sections from $C\Sigma$
  are located to the right of all parametric sections; namely, there are
  numbers $1 \leqslant \rho_A  \leqslant  \rho_{NA} \leqslant \rho_P \leqslant \rho_C \leqslant \rho_ = \rho_
  \Omega$ such that $[1,\rho_A +1]$, $[\rho_A
  +1, \rho_{NA}+1]$, $[\rho_{NA}+1, \rho_P+1]$, and
  $[\rho_P +1, \rho_\Omega +1]$ are, correspondingly, unions of
  all active,  all non-active, all parametric, and
  all constant sections;
 \item [2)] for every letter $a \in A^{\pm 1}$ there is at most one
 constant base in $\Omega$ labelled by $a$, and  all such bases are located in the
 $C\Sigma$;
\item [3)] every free variable (item) $h_i$ of $\Omega$ is located
in $C\Sigma$.
 \end{enumerate}
 \end{definition}
Sometimes, we label the item associated with some letter by this
letter.

Now we describe a typical method for  constructing generalized
equations with parameters starting with a system of ordinary group
equations with constants from $A$.

 \subsection{Parametric generalized equations
corresponding to group equations} \label{se:4-4}

 Let
\begin{equation}
\label{eq:systS} S(X,Y_1,Y_2, \ldots, Y_k,A) = 1 \end{equation}
 be a finite system of equations  with constants from $A^{\pm 1}$ and with the
 set of variables partitioned into a  disjoint union
  \begin{equation}
  \label{eq:ge3}
   X \cup Y_1 \cup  \cdots \cup Y_k
   \end{equation}
 Denote by  ${\mathcal GE}(S)$ the  set
of generalized equations corresponding to $S = 1$ from Lemma
\ref{le:R1}. Put $Y = Y_1 \cup \cdots \cup Y_k$.   Let $\Omega \in
{\mathcal GE}(S)$.  Recall that every  base $\mu$ occurs in
$\Omega$
 either related to some occurrence of a variable  from $X \cup Y$ in the
 system $S(X,Y,A) = 1$, or
related to an occurrence of a letter $z \in Z$ in the word $V$
(see Lemma \ref{le:4.2}), or is a constant base. If $\mu$
corresponds to a variable $x \in X$ ($y \in Y_i$)  then we say
that $\mu$ is an {\it $X$-base} ({\it $Y_i$-base}). Sometimes we
refer to $Y_i$-bases as to $Y$ bases.  For a base $\mu$ of
$\Omega$ denote by $\sigma_\mu$ the section $ \sigma_\mu =
[\alpha(\mu),\beta(\mu)]$. Observe that the section $\sigma_\mu$
is closed in $\Omega$ for every $X$-base, or $Y$-base. If $\mu$ is
an $X$-base ($Y$-base or $Y_i$-base), then the section
$\sigma_\mu$ is called an {\it $X$-section} ({\it $Y$-section} or
{\it $Y_i$-section}). If $\mu$ is a constant base and the section
$\sigma_\mu$ is closed then we call $\sigma_\mu$ a {\it constant }
section. Using the derived transformation D2 we transport all
closed $Y_1$-sections to the right end of the generalized
equations behind all the sections of the equation (in an arbitrary
order), then we transport all $Y_2$-sections and put them behind
all $Y_1$-sections, and so on. Eventually, we transport  all
$Y$-sections to the very end of the interval and they appear there
with respect to the partition (\ref{eq:ge3}). After that we take
all the constant sections and put them behind all the parametric
sections. Now, let  $A\Sigma$ be  the set of all $X$-sections,
$NA\Sigma = \emptyset$, $P\Sigma$ be the set of all $Y$-sections,
and $C\Sigma$ be the set of all constant sections. This defines a
parametric generalized equation $\Omega = \Omega_Y$ with
parameters corresponding to the set of variables $Y$. If the
partition of variables (\ref{eq:ge3})
 is fixed we will omit $Y$ in the notation above and call $\Omega$ the
{\it parameterized} equation obtained from $\Omega$.
 Denote by
 $$ {\mathcal GE}_{par}(\Omega) = \{\,\Omega_{Y} \mid \Omega \in  {\mathcal
 GE}(\Omega)\,\}$$
the set of all parameterized equations of the system
(\ref{eq:systS}).

\section{Elimination process: construction of $T(\Omega )$}
\label{se:5}

In  \cite{KMIrc} and \cite{Imp}, Section 5, we described a general
process of transformations of a generalized equation and
construction of fundamental sequences of solutions of equations.
In those papers we called this process ``Makanin-Razborov
process'', and sometimes it caused a confusion because the process
that we described is quite different. Since we use different
variations of the process, we prefer to call any such variation
here an {\em Elimination process}.

\subsection{ Elementary transformations}
\label{se:5.1}

 In this section we describe {\em elementary transformations} of generalized
 equations which were introduced by Makanin in \cite{Mak82}.
 Recall that we consider only formally consistent equations. In
general, an elementary transformation $ET$ associates   to  a
generalized equation $\Omega$ a finite set of generalized
equations $ET(\Omega) = \{\Omega_1, \ldots, \Omega_r\}$  and a
collection of surjective homomorphisms $\theta _i:G_{R(\Omega
)}\rightarrow G_{R(\Omega _i)}$ such that for every pair $(\Omega
, U)$ there exists a unique pair of the type $(\Omega _i, U_i)$
for which the following diagram commutes.

\vspace{.3cm}
\medskip

\begin{center}

\begin{picture}(100,100)(0,0)
\put(0,100){$F_{R(\Omega )}$} \put(100,100){$F_{R(\Omega _i)}$}
\put(0,0){$F(A)$} \put(25,103){\vector(1,0){70}}
\put(5,93){\vector(0,- 1){78}} \put(95,95){\vector(-1,-1){80}}
\put(-10,50){$\pi _{U_i}$} \put(55,45){$\pi _{U_i}$}
\put(50,108){$\theta _i$}
\end{picture}

\end{center}

\medskip
Here $\pi _{U}(X)=U.$  Since the pair $(\Omega_i,U_i)$ is defined
uniquely, we have a well-defined map  \newline $ET : (\Omega,U)
\rightarrow (\Omega_i,U_i).$

 ET1  ({\em Cutting a base}). Suppose  $\Omega$ contains a boundary
connection $\langle p,\lambda ,q\rangle  $. Then we replace (cut
in $p$) the base $\lambda$ by two new bases  $\lambda _1$ and
$\lambda _2$ and also replace (cut in $q$) $\Delta(\lambda)$ by
two new bases $\Delta (\lambda _1)$ and $\Delta (\lambda _2)$ such
that the following conditions hold.

If $\varepsilon(\lambda) =  \varepsilon(\Delta(\lambda))$, then
$$\alpha(\lambda_1) = \alpha(\lambda), \ \  \beta(\lambda_1) = p, \ \ \ \
\alpha(\lambda_2) = p, \ \  \beta(\lambda_2) = \beta(\lambda); $$
 $$
\alpha(\Delta(\lambda_1)) = \alpha(\Delta(\lambda)), \ \
\beta(\Delta(\lambda_1)) = q, \ \ \ \ \alpha(\Delta(\lambda_2)) =
q, \ \ \beta(\Delta(\lambda_2)) = \beta(\Delta(\lambda)); $$

If  $\varepsilon(\lambda) = - \varepsilon(\Delta(\lambda))$, then
$$\alpha(\lambda_1) = \alpha(\lambda), \ \  \beta(\lambda_1) = p, \ \ \ \
\alpha(\lambda_2) = p, \ \  \beta(\lambda_2) = \beta(\lambda); $$
 $$
\alpha(\Delta(\lambda_1)) = q,  \ \  \beta(\Delta(\lambda_1)) =
\beta(\Delta(\lambda)), \ \ \ \ \alpha(\Delta(\lambda_2)) =
\alpha(\Delta(\lambda)),  \ \  \beta(\Delta(\lambda_2)) = q; $$

Put $\varepsilon(\lambda_i) = \varepsilon(\lambda), \ \
\varepsilon(\Delta(\lambda_i)) = \varepsilon(\Delta(\lambda)), \ i
= 1,2.$

Let $(p', \lambda, q')$ be  a boundary connection in $\Omega$.

If    $p' <  p$,  then replace $(p', \lambda, q')$ by  $(p',
\lambda_1, q')$.

If    $p' >   p$, then replace   $(p', \lambda, q')$ by $(p',
\lambda_2, q')$.

Notice, since the equation $\Omega$ is formally consistent, then
the conditions above define boundary connections in the new
generalized equation. The resulting generalized equation
$\Omega^\prime$ is formally consistent.   Put $ET(\Omega) =
\{\Omega^\prime \}$.  Fig. 4 below explains the name of the
transformation ET1.

\begin{figure}[here]
\centering{\mbox{\psfig{figure=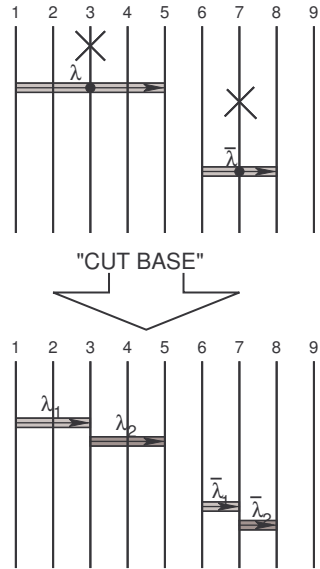}}} \caption{Elementary
transformation ET1.} \label{ET1}
\end{figure}

ET2  ({\em Transfer of a base}).  Let    a base $\mu $ of a
generalized equation $\Omega$  be contained in the base $\lambda$,
i.e., $\alpha (\lambda)\leqslant \alpha (\mu)< \beta
(\mu)\leqslant\beta (\lambda))$. Suppose that the boundaries
$\alpha(\mu)$ and $\beta(\mu))$ are $\lambda$-tied, i.e., there
are boundary connections of the type $\langle \alpha
(\mu),\lambda, \gamma_1\rangle  $ and $\langle \beta
(\mu),\lambda,\gamma_2\rangle $. Suppose also that every
$\mu$-tied boundary is $\lambda$-tied. Then we transfer $\mu$ from
its location  on the base $\lambda$ to the corresponding location
on the base $\Delta (\lambda)$ and adjust all the basic and
boundary equations (see Fig. 5). More formally, we replace $\mu$
by a new base $\mu^\prime$ such that $\alpha(\mu^\prime) =
\gamma_1, \beta(\mu^\prime) = \gamma_2$ and replace each
$\mu$-boundary connection $(p,\mu,q)$ with a new one
$(p^\prime,\mu^\prime,q)$ where $p$ and $p^\prime$ come from the
$\lambda$-boundary connection $(p,\lambda, p^\prime)$. The
resulting equation is denoted by $\Omega^\prime = ET2(\Omega)$.

\begin{figure}[here]
\centering{\mbox{\psfig{figure=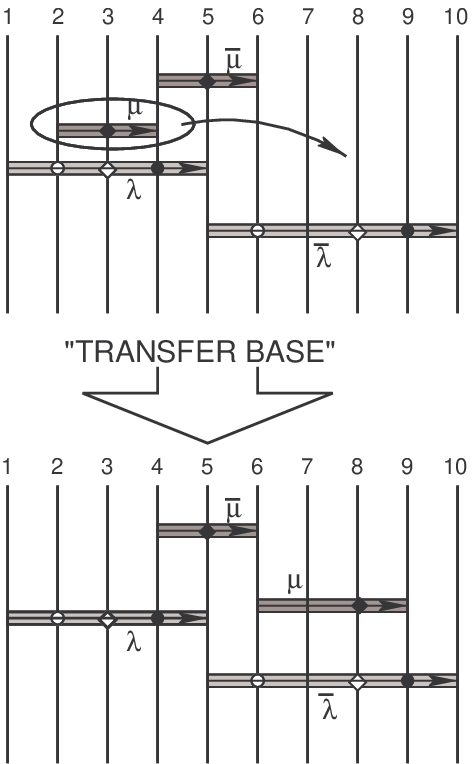}}} \caption{Elementary
transformation ET2.} \label{ET2}
\end{figure}

 ET3 ({\em Removal of a pair of  matched bases} (see Fig. 6)). Let
$\lambda$ and $\Delta(\lambda)$ be a pair of matched bases in
$\Omega$. Since $\Omega$ is formally consistent one  has
$\varepsilon(\lambda) = \varepsilon(\Delta(\lambda))$,
$\beta(\lambda) = \beta(\Delta(\lambda))$ and every
$\lambda$-boundary connection is of the type $(p,\lambda,p)$.
Remove the pair of bases $\lambda, \Delta(\lambda)$ with all
boundary connections related to $\lambda$. Denote the new
generalized equation by $\Omega^\prime$.

\smallskip
 {\bf Remark.} Observe that, for $i = 1,2,3$  $ETi(\Omega)$ consists of
a single equation $\Omega^\prime$,  such that  $\Omega$ and
$\Omega^\prime$ have the same set of variables $H$,   and the
identity map $F[H] \rightarrow F[H]$ induces an $F$-isomorphism
$F_{R(\Omega )} \rightarrow F_{R(\Omega^{\prime })}$. Moreover,
 $U$ is a solution of $\Omega$ if and only if $U$ is  a solution of
$\Omega^\prime$.

\begin{figure}[here]
\centering{\mbox{\psfig{figure=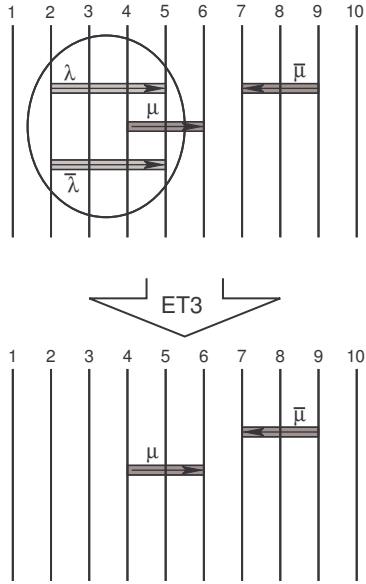}}} \caption{Elementary
transformation ET3.} \label{ET3}
\end{figure}

ET4 ({\em Removal of a lonely base} (see Fig. 7)).  Suppose in
$\Omega$ a variable base $\lambda$ does not intersect any other
variable base, i.e., the items $h_{\alpha(\lambda)}, \ldots,
h_{\beta(\lambda)-1}$ are contained in only one variable base
$\lambda$. Suppose also that all boundaries in $\lambda$ are
$\lambda$-tied, i.e., for every $i$ ($\alpha(\lambda)+1 \leqslant
i\leqslant \beta -1$) there exists a boundary $b(i)$ such that
$(i,\lambda ,b(i))$ is a boundary connection in $\Omega$.  For
convenience  we define: $b(\alpha(\lambda)) =
\alpha(\Delta(\lambda))$ and $b(\beta(\lambda)) =
\beta(\Delta(\lambda))$ if
$\varepsilon(\lambda)\varepsilon(\Delta(\lambda)) = 1$, and
$b(\alpha(\lambda)) = \beta(\Delta(\lambda))$ and
$b(\beta(\lambda)) = \alpha(\Delta(\lambda))$ if
$\varepsilon(\lambda)\varepsilon(\Delta(\lambda)) = -1$.

 The transformation ET4 carries $\Omega $ into a  unique
generalized equation $\Omega _1$ which is obtained from $\Omega $
by deleting the pair of bases $\lambda$ and $\Delta(\lambda)$;
deleting all the boundaries $\alpha(\lambda)+1, \ldots,
\beta(\lambda)-1$ ( and renaming  the rest $\beta(\lambda) -
\alpha(\lambda) - 1$ boundaries) together with all
$\lambda$-boundary connections; replacing every constant base
$\theta$ which is contained in $\lambda$ by a constant base
$\theta^\prime$ with the same label as $\theta$ and such that
$\alpha(\theta^\prime) = b(\alpha(\theta)), \beta(\theta^\prime) =
b(\beta(\theta))$.

 We define the homomorphism $\pi: F_{R(\Omega )} \rightarrow
F_{R(\Omega^{\prime })}$
 as follows:
$\pi(h_j)=h_j$ if $j< \alpha (\lambda)$ or $j\geqslant \beta
(\lambda);$

$$\pi(h_{i})=\{\begin{array}{ll} h_{b(i)}\cdots h_{b(i)-1},& if \varepsilon
(\lambda)=\varepsilon (\Delta\lambda),\\ h_{b(i)}\cdots
h_{b(i-1)-1},& if \varepsilon (\lambda)=-\varepsilon
(\Delta\lambda)
\end{array}$$
for $\alpha +1 \leqslant i\leqslant\beta (\lambda)-1.$ It is not
hard to see that $\pi$ is  an $F$-isomorphism.

\begin{figure}[here]
\centering{\mbox{\psfig{figure=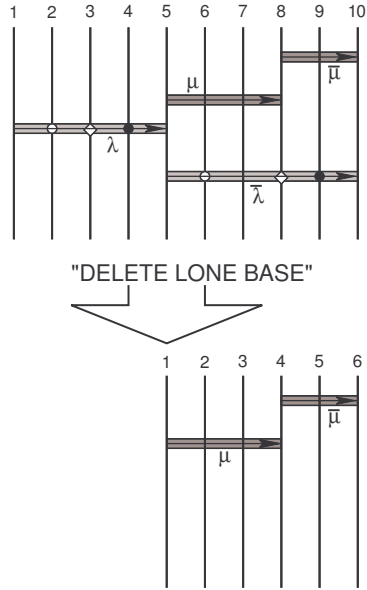}}} \caption{Elementary
transformation ET4.} \label{ET4}
\end{figure}

 ET5 ({\em Introduction of a boundary} (see Fig. 8)).  Suppose  a
point $p$ in a base $\lambda$ is not
 $\lambda$-tied. The transformation ET5 $\lambda$-ties it in all possible ways,
producing
 finitely many different generalized equations. To this end, let $q$ be
a boundary
 on $\Delta (\lambda)$. Then we perform one of the following two
transformations:

1. Introduce the boundary connection $\langle p,\lambda ,q\rangle
$ if the resulting equation $\Omega_q$ is formally consistent. In
this case the corresponding $F$-homomorphism
 $\pi_q: F_{R(\Omega )}$ into $F_{R(\Omega _q)}$ is induced by
the identity isomorphism on  $F[H]$. Observe that $\pi_q$  is not
necessary an isomorphism.

2. Introduce a new boundary $q^\prime$ between $q$ and $q+1$ (and
rename all the boundaries); introduce  a new boundary connection
$(p,\lambda,q^\prime)$. Denote the resulting equation by
$\Omega_q^\prime$. In this case the corresponding $F$-homomorphism
 $\pi_{q^\prime}: F_{R(\Omega )}$ into $F_{R(\Omega _{q^{\prime}})}$
is induced by the map $\pi(h) = h,$ if $h  \neq h_q$,  and
$\pi(h_q) = h_{q^\prime}h_{q^\prime +1}$.  Observe that  $\pi
_{q^\prime}$ is an $F$-isomorphism.

\begin{figure}[here]
\centering{\mbox{\psfig{figure=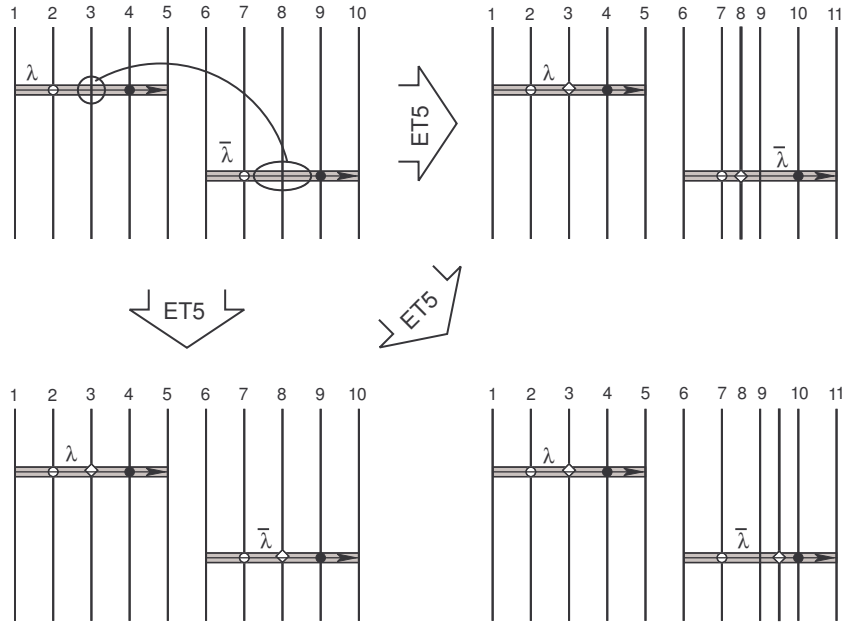}}} \caption{Elementary
transformation ET5.} \label{ET5}
\end{figure}

Let $\Omega$ be a generalized equation and $E$ be an elementary
transformation. By $E(\Omega)$ we denote a generalized equation
obtained from $\Omega$ by the elementary transformation $E$
(perhaps several such equations) if $E$ is applicable to $\Omega$,
otherwise we put $E(\Omega) = \Omega$.   By $\phi_E :
F_{R(\Omega)} \rightarrow F_{R(E(\Omega))}$ we denote the
canonical homomorphism of the coordinate groups (which has been
described above in the case $E(\Omega) \neq \Omega)$, otherwise,
the identical  isomorphism.

\begin{lemma}
\label{le:hom-check} There exists an algorithm which for every
generalized equation $\Omega$ and every elementary transformation
$E$ determines whether the canonical homomorphism $\phi_E:
F_{R(\Omega)} \rightarrow F_{R(E(\Omega))}$  is an isomorphism or
not.
\end{lemma}
\begin{proof} The only non-trivial case is when $E = E5$ and no new
boundaries were introduced. In this case $E(\Omega)$ is obtained
from $\Omega$ by adding a new particular equation, say $s = 1$,
which is effectively determined by $\Omega$ and $E(\Omega)$. In
this event, the coordinate group
$$F_{R(E(\Omega))} = F_{R(\Omega \cup \{s\})}$$
is a quotient group of $F_{R(\Omega)}$. Now $\phi_E$ is an
isomorphism if and only if $R(\Omega) = R(\Omega \cup \{s\})$, or,
equivalently, $s \in R(\Omega)$. The latter condition holds if and
only if $s$ vanishes on all solutions of the system of
(group-theoretic) equations $\Omega = 1$ in $F$, i.e., if the
following formula holds in $F$:
 $$ \forall x_1 \cdots \forall x_\rho (\Omega(x_1, \ldots, x_\rho) = 1
 \rightarrow s(x_1, \ldots, x_\rho) = 1).$$
 This can be checked effectively, since the universal theory of a free
 group $F$ is decidable \cite{Mak84}.\end{proof}

\subsection{Derived transformations and auxiliary transformations}
\label{se:5.2half}

 In this section we describe several useful transformations of
generalized equations. Some of them can be realized as  finite
sequences of elementary transformations, we call them {\it derived
} transformations. Other transformations result in equivalent
generalized equations but cannot be realized by finite sequences
of elementary moves.

 D1  ({\em Closing a section}).

Let $\sigma$ be a section of $\Omega $. The transformation D1
makes the section $\sigma$ closed. To perform D1 we  introduce
boundary connections (transformations ET5) through the end-points
of $\sigma$ until these end-points are tied by every base
containing them, and then cut through the end-points all the bases
containing them (transformations ET1) (see Fig. 9).
\begin{figure}[here]
\centering{\mbox{\psfig{figure=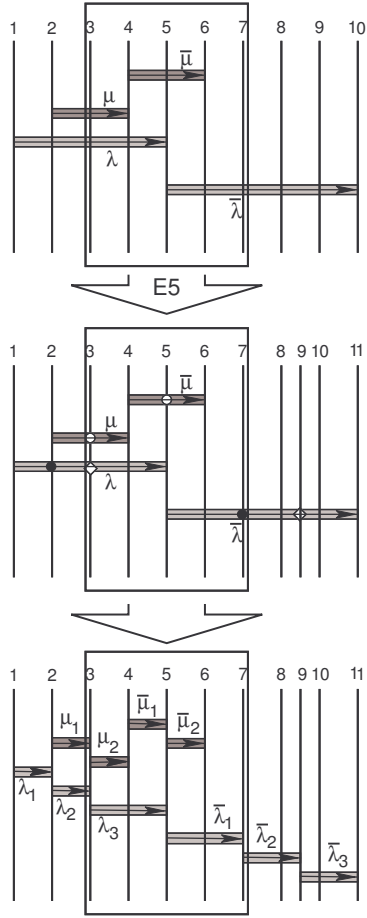}}} \caption{Derived
transformation D1.}\label{ D1}
\end{figure}

\newpage
D2 ({\it Transporting a closed section}).

Let $\sigma$ be a closed section of a generalized equation
$\Omega$. We cut $\sigma$ out of the interval $[1,\rho_\Omega]$
together with all the bases and boundary connections on $\sigma$
and put $\sigma$ at the end of the interval or between any two
consecutive   closed sections of $\Omega$. After that we
correspondingly re-enumerate all the items and boundaries of the
latter equation to bring it to the proper form. Clearly, the
original equation $\Omega$ and the new one $\Omega^\prime$ have
the same solution  sets and their coordinate groups are isomorphic
(see Fig. 10).

\begin{figure}[here]
\centering{\mbox{\psfig{figure=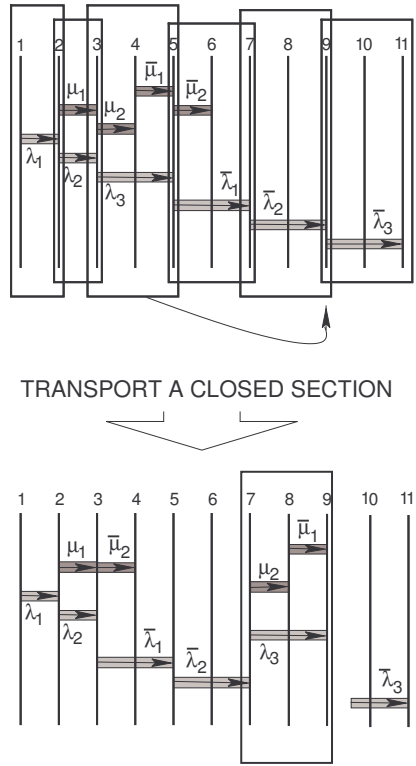}}} \caption{Derived
transformation D2.} \label{D2}
\end{figure}

D3 ({\it  Complete cut}).

 Let $\Omega $ be  a generalized
equation. For every boundary connection $(p,\mu,q)$ in $\Omega$ we
cut the base $\mu$ at $p$  applying ET1. The resulting generalized
equation $\tilde\Omega$ is obtained from $\Omega$ by a consequent
application of all possible ET1 transformations. Clearly,
$\tilde\Omega$ does not depend on a particular choice of the
sequence of transformations ET1. Since ET1 preserves  isomorphism
between the coordinate groups, equations $\Omega$ and
$\tilde\Omega$ have isomorphic coordinate groups, and the
isomorphism arises from the identity map $F[H] \rightarrow F[H]$.

D4 ({\it Kernel of a generalized equation}).

Suppose  that a generalized equation $\Omega$  does not contain
boundary connections. An active base $\mu  \in A\Sigma_\Omega$ is
called {\em eliminable} if at least one of the following holds:

 a) $\mu$ contains an  item $h_i$ with $\gamma(h_i)=1$;

b) at least one of the boundaries $\alpha (\mu),\beta (\mu)$ is
different from $1,\rho +1 $ and it does not touch any other base
(except  $\mu$).

A {\it cleaning  process} for $\Omega$ consists of consequent
removals  of  eliminable bases until no eliminable bases left in
the equation. The resulting generalized equation is called a {\it
kernel} of $\Omega$ and we denote it by ${\rm Ker}(\Omega )$. It
is easy to see that ${\rm Ker}(\Omega )$ does not depend on a
particular cleaning process. Indeed, if $\Omega$ has two different
eliminable bases $\mu_1$, $\mu_2$, and deletion of $\mu_i$ results
in an equation $\Omega_i$ then  by induction (on the number of
eliminations) ${\rm Ker}(\Omega_i)$ is uniquely defined for $i =
1,2$. Obviously, $\mu_1$ is still eliminable in $\Omega_2$, as
well as $\mu_2$ is eliminable in $\Omega_1$. Now eliminating
$\mu_1$ and $\mu_2$ from $\Omega_2$ and $\Omega_1$ we get one and
the same equation $\Omega_0$. By induction ${\rm Ker}(\Omega_1) =
{\rm Ker}(\Omega_0) = {\rm Ker}(\Omega_2)$ hence the result.
 We say that a variable $h_i$ {\it belongs to
the kernel} ($h_i \in {\rm Ker}(\Omega)$), if either $h_i$ belongs
to at least one base in the kernel, or it is parametric, or it is
constant.

Also, for an equation $\Omega$ by $\overline{\Omega}$ we denote
the equation which is obtained from $\Omega$ by deleting all free
variables. Obviously,
$$F_{R(\Omega)} = F_{R(\overline {\Omega})} \ast F(Y)$$
where $Y$ is the set of free variables in $\Omega$.

Let us consider what happens on the group level in the cleaning
process.

 We start with the case  when just one base is eliminated.
 Let $\mu$ be an eliminable base in $\Omega = \Omega(h_1, \ldots,
h_\rho)$. Denote by $\Omega_1$ the equation resulting from
$\Omega$ by eliminating $\mu$.

1) Suppose  $h_i \in \mu$ and $\gamma(h_i) = 1$. Then the variable
$h_i$ occurs only once in $\Omega$ -- precisely in the equation
$s_\mu = 1$ corresponding to the base $\mu$.  Therefore, in the
coordinate group $F_{R(\Omega )}$ the relation $s_\mu = 1$ can be
written as $h_i = w$, where $w$ does not contain $h_i$. Using
Tietze transformations we can rewrite the presentation of
$F_{R(\Omega )}$ as $F_{R(\Omega^\prime)}$, where $\Omega^\prime$
is obtained from $\Omega$ by deleting $s_\mu$ and the item $h_i$.
It follows immediately that
$$F_{R(\Omega_1)} \simeq  F_{R(\Omega^\prime)} \ast \langle h_i \rangle$$
and
 \begin{equation}
 \label{eq:ker1}
 F_{R(\Omega)} \simeq F_{R(\Omega^\prime)} \simeq F_{R(\overline
{\Omega_1})} \ast F(Z) \end{equation}
 for some free group $F(Z)$.
Notice that all the groups and equations which occur above can be
found effectively.

2) Suppose now that  $\mu$ satisfies  case b) above with respect
to a boundary $i$. Then in the equation $s_\mu = 1$ the variable
$h_{i-1}$ either occurs only once or it occurs precisely twice and
in this event the second occurrence of $h_{i-1}$ (in
$\Delta(\mu)$) is a part of the subword $(h_{i-1}h_i)^{\pm 1}$. In
both cases it is easy to see that the tuple
$$(h_1, \ldots, h_{i-2}, s_\mu, h_{i-1}h_i, h_{i+1}, \ldots,h_\rho)$$
forms a basis of the ambient free group generated by $(h_1,
\ldots, h_\rho)$  and constants from $A$. Therefore, eliminating
the relation $s_\mu = 1$,  we can rewrite the presentation of
$F_{R(\Omega)}$ in generators $Y = (h_1, \ldots, h_{i-2},
h_{i-1}h_i, h_{i+1}, \ldots,h_\rho)$. Observe also  that  any
other equation $s_\lambda = 1$ ($\lambda \neq \mu$) of $\Omega$
either does not contain variables $h_{i-1}, h_i$ or it contains
them as parts of the subword $(h_{i-1}h_i)^{\pm 1}$, i.e., any
such a word  $s_\lambda$ can be expressed as a word $w_\lambda(Y)$
in terms of generators $Y$ and constants from $A$. This shows that
$$F_{R(\Omega)} \simeq F(Y \cup A)_{R(w_\lambda(Y) \mid \lambda \neq
\mu)} \simeq F_{R(\Omega^\prime)},$$ where $\Omega^\prime$ is a
generalized equation obtained from $\Omega_1$ by deleting the
boundary $i$. Denote by $\Omega^{\prime \prime}$ an equation
obtained from $\Omega^\prime$ by adding a free variable $z$ to the
right end of $\Omega^\prime$.
 It follows now that
 $$F_{R(\Omega_1)} \simeq  F_{R(\Omega^{\prime \prime})} \simeq
 F_{R(\Omega)} \ast \langle z \rangle$$
and
\begin{equation}
\label{eq:ker2} F_{R(\Omega)} \simeq F_{R(\overline
{\Omega^\prime})} \ast F(Z) \end{equation}
 for some free group $F(Z)$. Notice that all the
groups and equations which occur above can be found effectively.

By induction on the number of steps in a cleaning process we
obtain the following lemma.
\begin{lemma}\label{7-10}
 $$F_{R(\Omega)} \simeq F_{R(\overline {{\rm Ker}\, \Omega})} \ast F(Z)$$
 where $F(Z)$ is a free group on ${ Z}$. Moreover, all the groups and equations
which occur above can be found effectively.\end{lemma}
 \begin{proof} Let
 $$\Omega = \Omega_0 \rightarrow \Omega_1 \rightarrow \ldots
 \rightarrow \Omega_l = {\rm Ker}\, \Omega$$
 be a cleaning process for $\Omega$. It is easy to see (by induction on $l$)
 that for every $j = 0, \ldots,l-1$
 $$\overline{{\rm Ker}\, \Omega_j} = \overline{Ker\, \overline{\Omega_j}}.$$
Moreover, if $\Omega_{j+1}$ is obtained from $\Omega_j$ as in the
case 2) above, then (in the notation above)
 $$\overline{{\rm Ker} (\Omega_j)_1} = \overline{{\rm Ker}\, \Omega_j^\prime} .$$
Now the statement of the lemma follows from the remarks above and
equalities (\ref{eq:ker1}) and (\ref{eq:ker2}).
\end{proof}

D5 ({\it  Entire transformation}).

We need a few further definitions.
 A base $\mu $ of the equation $\Omega$ is called a {\it leading} base if
$\alpha(\mu)=1$. A leading base is said to be {\it maximal} (or a
{\it carrier}) if $\beta (\lambda)\leqslant \beta (\mu),$ for any
other leading base $\lambda $.  Let $\mu $ be a carrier base of
 $\Omega.$ Any active base $\lambda \neq \mu$ with $\beta(\lambda )\leqslant \beta (\mu )$
 is called a {\it transfer} base (with respect to $\mu$).

  Suppose now that $\Omega$ is a generalized equation with
$\gamma(h_i)\geqslant 2$ for each $h_i$ in the active part of
$\Omega$. {\em An entire transformation} is a sequence of
elementary transformations which are performed as follows. We fix
a carrier base $\mu$ of $\Omega$. For any transfer base $\lambda $
we
 $\mu$-tie (applying $ET5$) all boundaries in $\lambda $. Using ET2 we transfer
 all transfer
bases from $\mu$ onto $\Delta (\mu)$. Now, there exists some $i <
\beta (\mu)$ such that $h_1,\ldots ,h_i$ belong to only one base
$\mu,$ while
 $h_{i+1}$ belongs to at least two bases. Applying ET1 we cut $\mu$
along the boundary $i+1$. Finally, applying ET4 we delete the
section $[1,i+1]$.

D6 ({\it Identifying closed constant sections}).

Let $\lambda$ and $\mu$ be two constant bases in $\Omega$ with
labels $a^{\varepsilon_\lambda}$ and $a^{\varepsilon_\mu}$, where
$a \in A$ and $\varepsilon_\lambda, \varepsilon_\mu \in \{1,-1\}$.
Suppose that the sections $\sigma(\lambda) = [i,i+1]$ and
$\sigma(\mu) = [j,j+1]$ are closed. Then we introduce a new
variable base $\delta$ with its dual $\Delta(\delta)$ such that
$\sigma(\delta) = [i,i+1]$, $\sigma(\Delta(\delta)) = [j,j+1]$,
$\varepsilon(\delta) = \varepsilon_\lambda$,
$\varepsilon(\Delta(\delta)) = \varepsilon_\mu$. After that we
transfer all bases from $\delta$ onto $\Delta(\delta)$ using ET2,
remove the bases $\delta$ and $\Delta(\delta)$, remove the item
$h_i$, and enumerate the items in a proper order. Obviously, the
coordinate group of the resulting equation is   isomorphic to the
coordinate group of the original equation.

\subsection{ Construction of the tree  $T(\Omega)$}
\label{se:5.2}

In this section we describe a branching rewrite  process for a
generalized equation $\Omega$. This process results in an
(infinite) tree $T(\Omega)$. At the end of the section we describe
infinite paths in $T(\Omega)$.

{\bf Complexity of a parametric generalized equation.}

 Denote by $\rho _A$ the number of variables $h_i$ in all active
sections of $\Omega ,$  by $n_A=n_A(\Omega )$ the number of bases
in active sections of $\Omega$,   by $\nu '$ - the number of open
boundaries in the active sections, by $\sigma '$ - the number of
closed boundaries in the active sections.

The number of closed active sections containing no bases,
precisely one base, or more than one base are denoted by $t_{A0},
t_{A1}, t_{A2}$ respectively. For a closed section $\sigma \in
\Sigma_\Omega$ denote by $n(\sigma)$, $\rho(\sigma)$  the number
of bases and, respectively, variables  in $\sigma$.

$$\rho _A= \rho_A(\Omega) = \sum_{\sigma \in
A\Sigma_\Omega}\rho(\sigma)$$
$$n_A= n_A(\Omega) = \sum_{\sigma \in
A\Sigma_\Omega} n(\sigma)$$

The {\em complexity} of a parametric  equation $\Omega $ is the
number
$$\tau = \tau (\Omega) = \sum_{\sigma \in
A\Sigma_\Omega} \max\{0, n(\sigma)-2\}.$$

Notice that the entire transformation (D5) as well as  the
cleaning process (D4) do not increase complexity of equations.

Let  $\Omega $ be a parametric generalized equation. We construct
a tree $T(\Omega)$ (with associated structures), as a directed
tree oriented from a root $v_0$, starting at $v_0$ and proceeding
by induction from vertices at distance $n$ from the root to
vertices at distance $n+1$ from the root.

We start with a general description of the tree $T(\Omega)$.  For
each vertex $v$ in $T(\Omega)$ there exists a unique generalized
equation $\Omega_v$ associated with $v$. The initial equation
$\Omega$ is associated with the root $v_0$, $\Omega_{v_0} =
\Omega$. For each edge $v\rightarrow v'$ (here $v$ and $v'$ are
the origin and the terminus of the edge) there exists a unique
surjective homomorphism $\pi(v,v'):F_{R(\Omega _v )}\rightarrow
F_{R(\Omega _v' )}$ associated with $v\rightarrow v'$.

 If
 $$v\rightarrow v_1\rightarrow\ldots\rightarrow
v_s\rightarrow u$$ is a path in $T(\Omega )$, then by $\pi (v,u)$
we denote  composition of corresponding homomorphisms
$$\pi (v,u) = \pi (v,v_1)   \cdots
\pi (v_s,u).$$

The set of edges of $T(\Omega)$ is subdivided into two classes:
{\it principal} and {\it auxiliary}. Every newly constructed edge
is principle, if not said otherwise. If $v \rightarrow v^\prime$
is a principle edge then there exists a finite sequence of
elementary  or derived transformations from $\Omega_v$ to
$\Omega_{v^\prime}$ and the homomorphism $\pi(v,v')$ is
composition of the homomorphisms corresponding to these
transformations. We also assume that active [non-active] sections
in $\Omega_{v^\prime}$ are naturally inherited from $\Omega_v$, if
not said otherwise.

Suppose the tree $T(\Omega)$ is constructed by induction up to a
level $n$, and suppose  $v$ is a vertex at distance $n$ from the
root $v_0$. We describe now how to extend the tree from $v$.  The
construction of the  outgoing edges at  $v$ depends on which case
described below takes place at the vertex $v$. We always assume
that if we have Case $i,$ then all Cases $j$, with $j \leqslant
i-1$, do not take place at $v$. We will see from the description
below that there is an effective procedure to check whether or not
a given case takes place at a given vertex. It will be obvious for
all cases, except Case 1. We treat this case below.

{\bf Preprocessing.}

Case 0.  In $\Omega_v$ we transport closed sections using D2 in
such a way that all active sections are at the left end of the
interval (the active part of the equation), then come all
non-active sections (the non-active part of the equation), then
come parametric sections (the parametric part of the equation),
and behind them  all constant sections are located (the constant
part of the equation).

 {\bf Termination conditions.}

Case 1. The homomorphism $\pi (v_0,v)$ is not an isomorphism (or
equivalently, the homomorphism $\pi(v_1,v)$, where $v_1$ is the
parent of $v$, is not an isomorphism). The vertex $v$ is called a
{\it leaf} or an {\it end} vertex. There are no outgoing edges
from $v$.

\begin{lemma}
There is an algorithm  to verify whether the homomorphism
$\pi(v,u)$, associated with an edge $v \rightarrow u$ in
$T(\Omega)$  is an isomorphism or not.
\end{lemma} \begin{proof} We will see below (by a straightforward inspection of
Cases 1-15 below) that every homomorphism of the type $\pi(v,u)$
is a composition of the canonical homomorphisms corresponding to
the elementary (derived) transformations. Moreover, this
composition is effectively given. Now the result follows from
Lemma \ref{le:hom-check}.
\end{proof}

Case 2. $\Omega _v$ does not contain active sections. The vertex
$v$ is called a {\it leaf} or an {\it end} vertex. There are no
outgoing edges from $v$.

\begin{figure}[h!]
\centering{\mbox{\psfig{figure=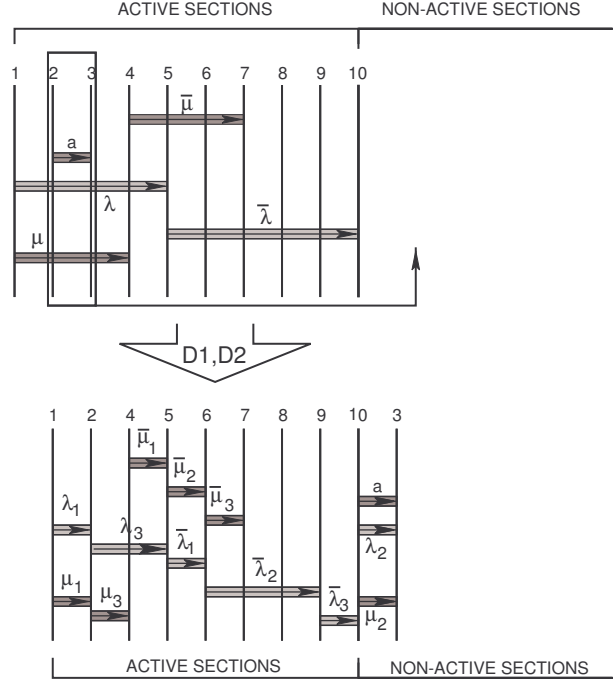}}} \caption{Case 3-4:
Moving constant bases.} \label{constants}
\end{figure}

{\bf Moving constants to the right.}

Case 3. $\Omega _v$ contains a constant base $\delta$ in an active
section such that  the section $\sigma(\delta)$ is not closed.

Here we  close the section $\sigma(\delta)$ using the derived
transformation $D1$.

Case 4. $\Omega _v$ contains a constant base $\delta$ with a label
$a \in A^{\pm 1}$ such that the section $\sigma(\delta)$ is
closed.

Here we transport the section $\sigma(\delta)$ to the location
right after all variables and parametric sections in $\Omega_v$
using the derived transformation D2. Then we identify  all closed
sections of the type $[i,i+1]$, which contain a constant base with
the label $a^{\pm 1}$, with the transported section
$\sigma(\delta)$, using the derived transformation D6. In the
resulting generalized equation $\Omega_{v^\prime}$ the section
$\sigma(\delta)$ becomes a constant section, and the corresponding
edge $(v,v^\prime)$ is auxiliary. See Fig. 11.

{\bf Moving free variables to the right.}

 Case 5. $\Omega _v$ contains a free variable $h_q$ in an  active section.

Here we close the section $[q,q+1]$ using D1, transport it to the
very end of the interval behind all items in $\Omega_v$ using D2.
In the resulting generalized equation $\Omega_{v^\prime}$ the
transported section  becomes a constant section, and the
corresponding edge $(v,v^\prime)$ is auxiliary.

\begin{remark} If Cases 0--5 are not possible at $v$ then the parametric
generalized equation $\Omega_v$ is in standard form.
\end{remark}

 Case 6. $\Omega _v$ contains a pair of matched bases in an
active section.

Here we perform ET3 and delete it. See Fig. 12.

\begin{figure}[here]
\centering{\mbox{\psfig{figure=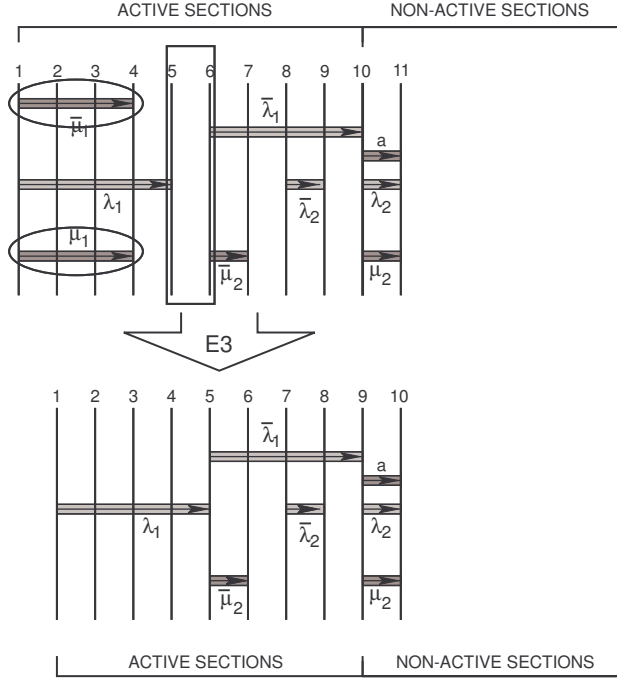}}} \caption{Case 5-6:
Trivial equations and useless variables.} \label{useless}
\end{figure}

{\bf Eliminating linear variables.}

 Case 7. In $\Omega_v$ there is $h_i$ in an active section with
$\gamma _i=1$ and  such that both boundaries $i$ and $i+1$ are
closed.

Here we remove the closed section $[i,i+1]$ together with the lone
base using ET4.

Case 8. In $\Omega_v$ there is $h_i$ in an active section with
$\gamma _i=1$ and  such that  one of the boundaries $i,i+1$ is
open, say $i+1$, and the other is closed.

Here we perform ET5 and $\tau$-tie $i+1$ through the only base
$\tau$ it intersects; using ET1 we cut $\tau$ in $i+1$; and then
we delete the closed section $[i,i+1]$ by ET4.  See Fig. 13.

\begin{figure}[here]
\centering{\mbox{\psfig{figure=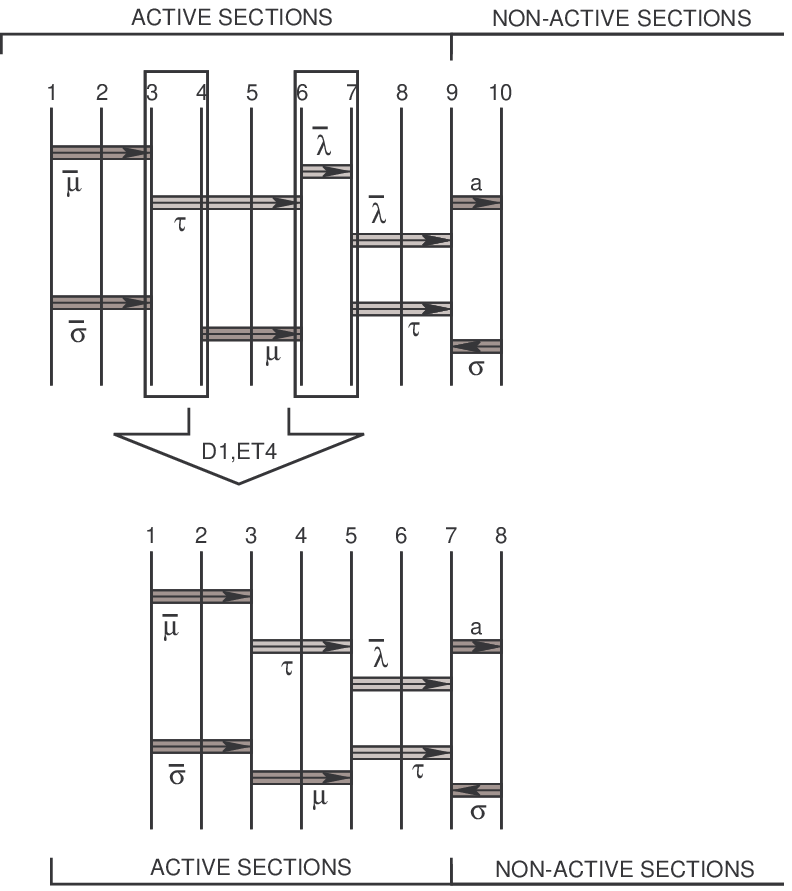}}} \caption{Case 7--10:
Linear variables.} \label{linear'}
\end{figure}

Case 9. In $\Omega_v$ there is $h_i$ in an active section with
$\gamma _i=1$ and  such that both boundaries $i$ and $i+1$ are
open. In addition, assume that there is a closed section $\sigma$
containing exactly two (not matched)  bases $\mu _1$ and $\mu _2$,
such that $\sigma = \sigma(\mu_1) = \sigma(\mu_2)$ and in the
generalized equation $\tilde {\Omega}_v$ (see the derived
transformation D3) all the bases obtained from $\mu _1,\mu _2$  by
ET1 in constructing $\tilde {\Omega}_v$ from $\Omega_v$, do not
belong to the kernel of $\tilde{\Omega}_v.$

Here, using ET5, we  $\mu _1$-tie  all the boundaries inside
$\mu_1$;
 using ET2, we transfer $\mu _2$ onto $\Delta (\mu_1)$; and
 remove $\mu _1$ together with the closed section $\sigma$  using
 ET4.

Case 10. $\Omega_v$  satisfies the first assumption of Case 9 and
does not satisfy the second one.

In this event we close the section $[i,i+1]$ using D1 and remove
it using ET4.

{\bf Tying a free boundary.}

Case 11. Some boundary $i$ in  the active part of $\Omega_v$ is
free. Since we do not have Case 5 the boundary $i$ intersects at
least one base, say, $\mu$.

Here we $\mu$-tie $i$  using ET5.

{\bf Quadratic case.}

Case 12. $\Omega_v$ satisfies the condition  $\gamma _i = 2$ for
each $h_i$ in the active part.

We apply the entire transformation D5.

Let us analyze the structure of $F_{R(\Omega _v)}$ in this case.
Let $\sigma$ be a closed section in the active part $\sigma
=[h_1,\ldots ,h_i]$. Let $F_1$ be a free group with basis
$\{h_1,\ldots ,h_i\}$. If $\sigma$ contains an open boundary $j,$
then we can consider a new generalized equation $\Omega  _v'$
obtained from $\Omega _v$ by replacing the product $h_{j-1}h_j$ by
a new variable $h_{j-1}'$ and represent $F_{R(\Omega _v)}=
F_{R(\Omega _v')}*\langle h_j\rangle $. Therefore we can suppose
that $\sigma$ does not contain open boundaries. We say that $\mu$
is a {\em quadratic} base if  $\sigma$ contains $\mu$ and $\Delta
(\mu )$, otherwise $\mu$ is a {\em quadratic-coefficient base}.
Denote the set of quadratic-coefficient bases by $C$. Suppose that
$\sigma$ contains  quadratic bases. Let $F_1/\sigma$ be the
quotient of $F_1$ over the normal closure of elements $ h[\alpha
(\mu ),\beta (\mu )] h[\alpha (\Delta (\mu ) ),\beta (\Delta (\mu
))]^{-1}$. Let  $\mathcal M$ be a set that contains exactly one
representative of each pair of double bases on $\sigma$ and
contains also each base $\mu$ such that $\mu\in\sigma$ and $\Delta
(\mu )\not\in\sigma$. If we identify each base on $ \sigma$ with
its double, then  the product $h_1\cdots h_i$ can be written as a
product of bases from $\mathcal M$ in exactly two different ways
$\mu _{i_1}\cdots \mu _{i_k},$ and $\mu_{j_1}\cdots \mu _{j_t}$.
Then $F_1/\sigma$ is isomorphic to the quotient of the free group
$F(\mathcal M )$ over the relation $\mu _{i_1}\cdots \mu
_{i_k}=\mu_{j_1}\cdots \mu _{j_t}.$ Every element of $\mathcal M$
occurs in this relation at most twice. Applying an automorphism of
$F(\mathcal M )$ identical on $C$, we can obtain another basis
$X\cup T\cup C$ of this group such that in this basis the relation
has form of a standard quadratic equation in variables from $X$
with coefficients in $F(C)$ \cite{LS}, Section 1.7, (variables
from $T$ do not participate in any relations).
 The
quadratic equation corresponding to $\sigma$ can be written in the
standard form with coefficients expressed in terms of non-active
variables. If the equation is regular, then $F_{R(\Omega _v)}$ has
a QH subgroup corresponding to this equation as described in
Subsection \ref{QHquad}.

If $\sigma$ does not have quadratic-coefficient bases, then
$F_{R(\Omega _v)}$ splits as a free product with one factor being
a closed surface group.

Case 13. $\Omega_v$ satisfies the condition $\gamma _i\geqslant 2$
for each $h_i$ in the active part,  and $\gamma _i>2$ for at least
one such $h_i$. In addition,  for some active base $\mu$ section
$\sigma(\mu) = [\alpha (\mu ),\beta (\mu )]$ is closed.

In this case using $ET5$, we $\mu$-tie every boundary inside
$\mu$; using $ET2$, we transfer all bases from $\mu$ to $\Delta
(\mu )$; using $ET4$, we remove the lone base $\mu$ together with
the section $\sigma(\mu)$.

Case 14.  $\Omega_v$ satisfies the condition $\gamma _i\geqslant
2$ for each $h_i$ in the active part, and $\gamma _i>   2$ for at
least one such $h_i$. In addition,  some boundary $j$ in the
active part touches some base $\lambda$, intersects some base
$\mu$,  and $j$ is not $\mu$-tied.

Here we $\mu$-tie $j$.

{\bf General JSJ-case.}

 Case 15. $\Omega_v$ satisfies the condition $\gamma _i\geqslant 2$ for
each $h_i$ in the active part, and $\gamma _i>   2$ for at least
one such $h_i$.  We apply, first,  the entire transformation D5.

Here for every boundary $j$ in the active part that touches at
least one base, we $\mu$-tie $j$ by every base $\mu$ containing
$j$. This results in finitely many new vertices
$\Omega_{v^\prime}$ with principle edges $(v,v^\prime)$.

If, in addition,  $\Omega_v$ satisfies the following condition (we
called it  Case 15.1 in \cite{KMIrc}) then we construct the
principle edges as was described above, and also  construct a few
more auxiliary edges outgoing from the vertex $v$:

Case 15.1. The carrier base $\mu$ of the equation $\Omega _v$
intersects with its dual $\Delta (\mu)$.

Here we  construct an auxiliary    equation ${\hat \Omega _{ v}}$
(which does not occur in $T(\Omega)$)  as follows. Firstly, we add
a new constant section $[\rho _v+1,\rho _v+2]$ to the right of all
sections in $\Omega_v$ (in particular, $h_{\rho _v+1}$ is a new
free variable). Secondly, we introduce a new pair of bases
$(\lambda ,\Delta (\lambda))$ such that
$$\alpha(\lambda)=1, \beta (\lambda)= \beta
(\Delta (\mu)), \alpha (\Delta(\lambda))=\rho _v+1, \beta (\Delta
(\lambda))=\rho _v+2.$$
 Notice that $\Omega _v$ can be obtained from ${\hat \Omega
_{v}}$ by ET4:  deleting $\delta(\lambda)$ together with the
closed section $[\rho _v+1,\rho _v+2]$. Let
$${\hat \pi}_{v}: F_{R(\Omega _v)}\rightarrow F_{R({\hat \Omega}_{v})}$$
 be the isomorphism induced by ET4.  Case 15 still holds for ${\hat
 \Omega _{ v}}$,  but now $\lambda$ is the carrier base. Applying to
 ${\hat \Omega _{v}}$ transformations described in  Case 15, we obtain a
 list of new vertices $\Omega_{v^\prime}$ together with isomorphisms
$${\eta}_{v'}: F_{R({\hat \Omega}_{ v})} \rightarrow F_{R(\Omega _{v'})}.$$

Now for each such $v'$ we add to $T(\Omega)$ an auxiliary edge
$(v,v')$ equipped  with composition of homomorphisms $\pi(v,v') =
\eta_{v'} \circ {\hat \pi}_{v}$ and assign $\Omega_{v'}$ to the
vertex $v'$.

If none of the Cases 0-15 is possible, then we stop, and the tree
$T(\Omega)$ is constructed. In any case, the tree $T(\Omega)$ is
constructed by induction. Observe that, in general, $T(\Omega)$ is
an infinite locally finite tree.

If Case $i\ (0\leqslant i\leqslant 15)$ takes place at  a vertex
$v$ then we say that $v$ has type $i$ and write $tp(v)=i$.

{\bf Example.} Consider the generalized equation in Fig. 3. For
this generalized equation we first have to move constants to the
right (Cases 3, 4). After this we have  Case 11 (tieing free
boundaries). Then we have Case 12 (the entire transformation). At
the beginning $\lambda _1$ is the leading base, then $\mu$ is the
leading base. After these transformations we obtain the
generalized equation in Fig. 14.

\begin{figure}[here]
\vspace{-10ex}
\centering{\mbox{\psfig{figure=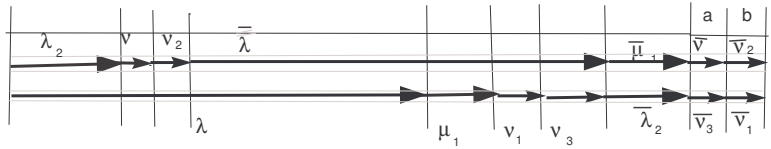,height=3.1in}}}
\vspace{-20ex} \caption{Case 12 for the generalized equation
$\Omega$ for the equation $[x,y][b,a]=1.$} \label{ET001}

\end{figure}

For this equation we still have Case 12. Beginning at this vertex
of the tree $T(\Omega)$ we have an infinite branch such that for
all the vertices $v_i$ of this branch $\Omega _{v_i}=\Omega
_{v_{i+1}}$.

\begin{lemma} \cite[Lemma 3.1]{Razborov} \label{3.1} If $u\rightarrow
v$ is a principal edge of the tree $T(\Omega)$, then
\begin{enumerate}
\item $n_{A}(\Omega_v) \leqslant n_{A}(\Omega_u)$, if $tp(u)\not =
3,10,$ this inequality is proper if $tp(u)=6,7,9,13;$ \item If
$tp(u)=10,$ then $n_{A}(\Omega_v) \leqslant n_{A}(\Omega_u) + 2;$
\item $\nu^\prime(\Omega_v)
 \leqslant \nu^\prime(\Omega_u)$ if $tp(u)\leqslant 13$ and $tp (u)\not = 3,11$; \item
$\tau(\Omega_v)  \leqslant \tau(\Omega_u)$, if $tp (u)\not = 3.$
\end{enumerate}\end{lemma}
\begin{proof} Straightforward verification.
\end{proof}

\begin{lemma} \label{3.2} Let
 $$v_1\rightarrow v_2\rightarrow\ldots \rightarrow v_r \rightarrow \ldots $$
  be an infinite path in the tree $T(\Omega ).$
Then there exists a natural number $N$ such that all the edges
$v_n \rightarrow v_{n+1}$ of this path with $n \geqslant N$ are
principal edges, and one of the following situations holds:
\begin{enumerate}
\item [1)] {\rm (\emph{linear case})}\  $7\leqslant
tp(v_n)\leqslant 10 $ for all $n\geqslant N;$ \item [2)] {\rm
(\emph{quadratic case})} \ $tp(v_n)=12$ for all $n\geqslant N;$
\item [3)] {\rm (}\emph{general JSJ case}{\rm )}\ $tp(v_n)=15$ for
all $n\geqslant N$.
\end{enumerate}
\end{lemma}

\begin{proof} Observe that starting with a generalized equation
$\Omega$ we can have Case 0 only once, afterward in all other
equations the active part is at the left, then comes the
non-active part, then - the parametric part,  and at the end - the
constant part. Obviously, Cases 1 and 2 do not occur on an
infinite path. Notice also that Cases 3 and 4 can only occur
finitely many times, namely, not more then $2t$ times where $t$ is
the number of constant bases in the original equation $\Omega$.
Therefore, there exists a natural number $N_1$ such that
$tp(v_i)\geqslant 5$ for all $i \geqslant N_1$.

Now we show that the number of vertices $v_i$ ($i \geqslant N$)
for which $tp (v_i)=5$ is not more than the minimal number of
generators of the group $F_{R(\Omega )}$, in particular, it cannot
be greater than $\rho +1 +|A|$, where $\rho =  \rho (\Omega).$
Indeed, if a path from the root $v_0$ to a vertex $v$ contains $k$
vertices of type 5, then $\Omega _v$ has at least $k$ free
variables in the constant part. This implies that the coordinate
group $F_{R(\Omega _v)}$ has a free group of rank $k$ as a free
factor, hence it cannot be generated by less than $k$ elements.
Since $\pi(v_0,v): F_{R(\Omega )} \rightarrow F_{R(\Omega _v)}$ is
a surjective homomorphism, the group $F_{R(\Omega )}$ cannot be
generated by less then $k$ elements. This shows that $k \leqslant
\rho +1 +|A|$.  It follows that there exists  a number $N_2
\geqslant N_1$ such that  $tp(v_i) >   5$ for every $i \geqslant
N_2$.

Suppose $i >   N_2$. If $tp(v_i)=12,$ then it is easy to see that
$tp(v_{i+1})=6$ or $tp(v_{i+1})=12$. But if $tp(v_{i+1})=6,$ then
$tp(v_{i+2})=5$ -- contradiction with $i >   N_2$. Therefore,
$tp(v_{i+1})=tp(v_{i+2})= \ldots = tp(v_{i+j}) = 12$ for every $j
>   0$ and we have  situation 2) of the lemma.

Suppose now $tp(v_{i})\neq 12$ for all $i \geqslant N_2$. By Lemma
\ref{3.1} $\tau(\Omega_{v_{j+1}}) \leqslant \tau(\Omega_{v_{j}})$
for every principle edge $v_j \rightarrow v_{j+1}$ where $j
\geqslant N_2$. If $v_j \rightarrow v_{j+1}$, where $j \geqslant
N_2$, is an auxiliary edge then  $tp(v_j) = 15$ and, in fact, Case
15.1 takes place at $v_j$.  In the notation of Case 15.1
$\Omega_{v_{j+1}}$ is obtained from ${\hat \Omega}_{v_j}$ by
transformations from Case 15. In this event, both bases $\mu$ and
$\Delta (\mu)$ will be transferred from the new carrier base
$\lambda $ to the constant part, so the complexity will be
decreased at least by two: $\tau(\Omega_{v_{j+1}}) \leqslant
\tau({\hat \Omega}_{v_j}) - 2$. Observe also that $\tau ({\hat
\Omega_{v_j}}) = \tau(\Omega_{v_j}) + 1$. Hence
$\tau(\Omega_{v_{j+1}})
 <  \tau(\Omega_{v_j}).$

It follows that there exists a number $N_3 \geqslant N_2$ such
that $\tau(\Omega_{v_j}) = \tau(\Omega_{v_{N_3}})$ for every $j
\geqslant N_3$, i.e., complexity stabilizes. Since every auxiliary
edge gives a decrease of complexity, this implies that for every
$j \geqslant N_3$ the edge $v_j \rightarrow v_{j+1}$ is principle.

Suppose now that $i \geqslant N_3$. We claim that $tp(v_i) \neq
6$. Indeed, if $tp(v_i)=6$, then the closed section, containing
the matched bases $\mu, \Delta (\mu)$, does not contain any other
bases (otherwise the complexity of $\Omega_{v_{i+1}}$ would
decrease).  But in this event $tp(v_{i+1})=5$ which is impossible.

So  $tp (v_i)\geqslant 7$ for every $i \geqslant N_3$. Observe
that ET3 (deleting match bases) is the only elementary
transformation that can produce new free boundaries. Observe also
that ET3 can be applied only in Case 6.  Since Case 6 does not
occur anymore along the path  for $i \geqslant N_3$, one can see
that no new free boundaries occur in equations $\Omega_{v_j}$ for
$j \geqslant N_3$. It follows that there exists a number $N_4
\geqslant N_3$ such that $tp (v_i)\neq 11$ for every $j \geqslant
N_4$.

Suppose now that for some $i \geqslant N_4$,  $13\leqslant
tp(v_i)\leqslant 15.$ It is easy to see from the description of
these cases  that in this event $tp(v_{i+1})\in\{6,13,14,15\}.$
Since $tp(v_{i+1})\neq 6$, this implies that $13\leqslant
tp(v_j)\leqslant 15$ for every $j \geqslant i$. In this case the
sequence $n_A(\Omega_{v_j})$ stabilizes by lemma \ref{3.1}. In
addition, if $tp (v_j)=13,$ then $ n_A(\Omega_{v_{j+1}}) <
n_A(\Omega_{v_j}).$ Hence there exists a number $N_5 \geqslant
N_4$ such that $tp (v_j)\not = 13$ for all $j \geqslant N_5.$

Suppose $i \geqslant N_5$.  There cannot be more than
$8(n_A(\Omega_{v_i})^2$ vertices of type 14 in a row starting at a
vertex $v_i$; hence there exists $j\geqslant i$ such that
$tp(v_j)=15$. The series of transformations ET5 in Case 15
guarantees the inequality $tp (v_{j+1})\not = 14$; hence
$tp(v_{j+1})=15,$ and we have situation 3) of the lemma.

So we can suppose $tp (v_i)\leqslant 10$ for all the vertices of
our path. Then we have situation 1) of the lemma. \end{proof}

\section{Elimination process: periodic structures}\label{subs:per}
Recall that a reduced word $P$ in a free group $F$ is called a
{\it period} if it is cyclically reduced and not a proper power. A
word $w \in F$ is called $P$-{\em periodic}  if $d(w) \geqslant
d(P)$ and it is a subword of some (integer) power of $P$. Every
$P$-periodic word $w$ can be presented in the form
\begin{equation}
\label{2.50} w = A^rA_1
\end{equation}
where $A$  is a cyclic permutation of $P^{\pm 1}$,  $r \geqslant
1$, $A = A_1 \circ A_2$, and $A_2 \neq 1$. This representation is
unique if $r \geqslant 2$. The number $r$ is called the {\em
exponent} of $w$. A maximal exponent of $P$-periodic subword in
$u$ is called the {\em exponent of $P$-periodicity in $u$}. We
denote it $e_P(u).$

\begin{definition}\label{11'}
Let $\Omega$ be a  standard generalized equation. A solution
$H:h_i \rightarrow H_i$ of $\Omega$ is called {\em periodic with
respect to a period $P$}, if for every  variable section $\sigma$
of $\Omega$ one of the following conditions hold:
\begin{enumerate}
 \item [1)] $H(\sigma)$ is $P$-periodic with exponent $r \geqslant 2$;
 \item [2)] $d(H(\sigma)) \leqslant d(P)$;
 \item [3)]  $H(\sigma)$ is $A$-periodic and $d(A) \leqslant d(P)$;
\end{enumerate}
Moreover, condition 1) holds at least for one such $\sigma$.
\end{definition}

The section $\sigma$ is called {\em $P$-periodic}. In this section
we will prove the following result.
\begin{prop}\label{pr:per} \

\bi \item[1.] Let $\Omega$ be a generalized equation. One can
effectively find a number $N$ such that every periodic solution of
$\Omega$ with some period $P$ {\rm (}see Definition {\rm
\ref{11'})} is a composition of a canonical automorphism of the
coordinate group $F_{R(\Omega)}$ with  either a  solution with
exponent of periodicity for $P$ less than $N$ or a solution of a
fixed equation not in $R(\Omega )$. These canonical automorphisms
correspond to Dehn twists of $F_{R(\Omega)}$ which are related to
the splitting of $F_{R(\Omega)}$ {\rm (}which comes from a
so-called periodic structure{\rm )} over an abelian edge group.

\item[2.] Let $\Omega$ be a generalized equation with parameters
$Y$ that has  a periodic solution $H$ of $\Omega$ with period $P$.
Then one can effectively find a number $N$ such that there are two
possibilities.
 \begin{enumerate} \item $F_{R(\Omega )}$
has a splitting modulo $\langle Y\rangle$ and $H$ is a composition
of a canonical automorphism of the coordinate group
$F_{R(\Omega)}$ corresponding to this splitting (fixing the
subgroup of parameters) with  a solution with  exponent of $P$
bounded by $Ne(H(Y))$  or a solution of an equation not in
$R(\Omega )$.\item $e(H)\leqslant Ne(H(Y))$.
\end{enumerate}
\ei

 One can effectively find $N$ and the splitting.\end{prop}

We will prove this proposition in Subsection \ref{proofper}. We
will later formulate a more detailed statement about splittings of
$F_{R(\Omega )}$ (see Theorem \ref{th:spl}).
\subsection{Definition of a periodic structure}
In this subsection we introduce a notion of a periodic structure
which allows one to describe periodic solutions of generalized
equations.

We fix till the end of the section a generalized equation $\Omega$
in standard form. Recall that in $\Omega$ all closed sections
$\sigma$, bases $\mu$, and variables $h_i$ belong to either the
variable part $V\Sigma$, or the parametric part $P\Sigma$, or the
constant part $C\Sigma$ of $\Omega$.

\begin{definition} \label{above}
Let  $\Omega$ be a generalized equation in  standard form with no
boundary connections. A {\em periodic structure} on $\Omega$ is a
pair $\langle {\mathcal P}, R \rangle$, where
\begin{enumerate}
\item [1)]  ${\mathcal P}$ is a set consisting of some variables
$h_i$, some bases $\mu$, and some closed sections $\sigma$ from
$V\Sigma$  and such that the following conditions are satisfied:
\begin{itemize}
\item[a)] if $h_i \in {\mathcal P}$ and $h_i \in \mu$,  and
$\Delta(\mu) \in V\Sigma$, then $\mu \in {\mathcal P}$;

\item[b)] if $\mu \in {\mathcal P}$, then $\Delta(\mu) \in
{\mathcal{P}}$;

\item[c)] if $\mu \in {\mathcal P}$ and $\mu \in \sigma$, then
$\sigma \in {\mathcal P}$;

\item [d)] there exists a function ${\mathcal X}$ mapping the set
of closed sections from ${\mathcal P}$ into $\{-1, +1\}$ such that
for every $\mu, \sigma_1, \sigma_2 \in {\mathcal P}$, the
condition that $\mu \in \sigma_1$ and $\Delta(\mu) \in \sigma_2$
implies

$\varepsilon(\mu) \cdot \varepsilon(\Delta(\mu)) = {\mathcal
X}(\sigma_1) \cdot {\mathcal X}(\sigma_2)$;
\end{itemize}
\item [2)]  $R$ is an equivalence relation on a certain set
${\mathcal B}$ (defined below) such that the following conditions
are satisfied:
\begin{itemize}

 \item[e)] Notice, that for every boundary $l$ either
there exists a unique closed section $\sigma(l)$ in ${\mathcal P}$
containing $l$, or there exist precisely two closed section
$\sigma_{\rm left}(l) = [i,l], \sigma_{\rm right} =  [l,j]$ in
${\mathcal P}$ containing $l$. The set of boundaries of the first
type we denote
 by ${\mathcal B}_1$, and of the second  type - by  ${\mathcal B}_2$. Put
$${\mathcal B} = {\mathcal B}_1  \cup \{l_{\rm left}, l_{\rm right}  \mid l \in {\mathcal B}_2
\}$$
 here $l_{\rm left}, l_{\rm right}$ are two  "formal copies" of $l$. We
 will use the following agreement: for any
base $\mu$ if $\alpha(\mu) \in {\mathcal B}_2$ then by
$\alpha(\mu)$ we mean $\alpha(\mu)_{\rm right}$ and, similarly, if
$\beta(\mu) \in {\mathcal B}_2$ then by $\beta(\mu)$ we mean
$\beta(\mu)_{\rm left}$.

\item[f)]  Now, we define $R$ as follows. If $\mu \in {\mathcal
P}$ then
$$\alpha(\mu) \sim_R \alpha(\Delta(\mu)), \ \ \beta(\mu)
\sim_R \beta(\Delta(\mu)) \ \ if \ \varepsilon(\mu) =
\varepsilon(\Delta(\mu))$$
$$\alpha(\mu)) \sim_R \beta(\Delta(\mu)),\
 \beta(\mu)\sim_R  \alpha(\Delta(\mu)) \ \ if \ \varepsilon(\mu) = -
 \varepsilon(\Delta(\mu)).$$
\end{itemize}
\end{enumerate}
\end{definition}
\begin{remark}
For a given $\Omega$ one can effectively find all periodic
structures on $\Omega$.
\end{remark}

 Let $\langle {\mathcal P}, R \rangle $ be a periodic structure of
$\Omega$. Put
$$ {\mathcal NP} = \{\mu \in B\Omega \mid \exists h_i \in {\mathcal P} \mbox{ such
that } h_i \in \mu \mbox{ and } \Delta(\mu) \ \mbox{is parametric
or
 constant} \}$$

Now we will show how one can associate with a $P$-periodic
solution $H$ of  $\Omega$  a periodic structure ${\mathcal P}(H,
P) = \langle {\mathcal P}, R \rangle$. We define ${\mathcal P}$ as
follows. A closed section $\sigma$ is  in ${\mathcal P}$ if and
only if $\sigma \in V\Sigma$  and  $H(\sigma)$ is $P$-periodic
with exponent $\geqslant 2$. A variable $h_i$ is in ${\mathcal P}$
if and only if $h_i \in \sigma$ for some $\sigma \in {\mathcal P}$
and $d(H_i) \geqslant 2 d(P)$. A base $\mu$ is in ${\mathcal P}$
if and only if both $\mu$ and $\Delta(\mu)$  are in $v\Sigma$  and
one of them contains $h_i$ from ${\mathcal P}$.

\begin{lemma}
\label{le:PP} Let $H$ be a periodic solution of $\Omega$. Then
${\mathcal P}(H, P)$ is a periodic structure on $\Omega$.
\end{lemma}
\begin{proof} Let ${\mathcal P}(H, P) = \langle {\mathcal P}, R
\rangle$.
 Obviously,  ${\mathcal P}$ satisfies a) and b) from the definition \ref{above}.

Let $\mu \in {\mathcal P}$ and $\mu \in [i,j]$. There exists an
unknown $h_k \in {\mathcal P}$ such that $h_k \in \mu$ or $h_k \in
\Delta({\mu})$. If $h_k \in \mu$, then, obviously, $[i,j] \in
{\mathcal P}$. If $h_k \in \Delta(\mu)$ and $\Delta(\mu) \in [i',
j']$, then $[i', j'] \in {\mathcal P}$, and hence, the word
$H[\alpha(\Delta(\mu)), \beta(\Delta(\mu))]$ can be written in the
form $Q^{r'} Q_1$, where $Q=Q_1 Q_2$; $Q$ is a cyclic shift of the
word $P^{\pm 1}$ and $r' \geqslant 2$. Now let (\ref{2.50}) be a
presentation for the section $[i,j]$. Then $H[\alpha(\mu),
\beta(\mu)] = B^s B_1$, where $B$ is a cyclic shift of the word
$A^{\pm 1}$, $d(B) \leqslant d(P)$, $B = B_1 B_2$, and $s
\geqslant 0$. From the equality $H[\alpha(\mu),
\beta(\mu)]^{\varepsilon(\mu)} = H[\alpha(\Delta(\mu)),
\beta(\Delta(\mu)))]^{\varepsilon(\Delta(\mu))}$ and Lemma 1.2.9
\cite{1} it follows that $B$ is a cyclic shift of the word
$Q^{\pm1}$. Consequently, $A$ is a cyclic shift of the word
$P^{\pm 1}$, and $r \geqslant 2$ in (\ref{2.50}), since $d(H[i,j])
\geqslant d(H[\alpha(\mu), \beta(\mu)]) \geqslant 2 d(P)$.
Therefore, $[i,j] \in {\mathcal P}$; i.e, part c) of the
definition \ref{above} holds.

Put ${\mathcal X}([i,j]) = \pm 1$ depending on whether in
(\ref{2.50}) the word $A$ is conjugate to $P$ or to $P^{-1}$. If
$\mu \in [i_1, j_1]$, $\Delta(\mu) \in [i_2, j_2]$, and $ \mu \in
{\mathcal P}$, then the equality $\varepsilon(\mu) \cdot
\varepsilon(\Delta(\mu))$ = ${\mathcal X}([i_1, j_1]) \cdot
{\mathcal X}([i_2, j_2])$ follows from the fact that given $A^r
A_1 = B^s B_1$ and $r,s \geqslant 2$, the word $A$ cannot be a
cyclic shift of the word $B^{-1}$. Hence part d) also holds.

Now let $[i,j]\in {\mathcal P}$ and $ i \leqslant l \leqslant j$.
Then there exists a subdivision $P = P_1P_2$ such that if
${\mathcal X} ([i,j]) =1$, then the word $H[i,l]$ is the end of
the word $(P^\infty)P_1$, where $P^\infty$ is the infinite word
obtained by concatenations of powers of $P$, and $H[l,j]$ is the
beginning of the word $P_2(P^\infty)$, and if ${\mathcal
X}([i,j])= -1$, then the word $H[i,l]$ is the end of the word
$(P^{-1})^\infty P_2^{-1}$ and $H[l,j]$ is the beginning of
$P_1^{-1}(P^{-1})^\infty$. Again, Lemma 1.2.9 \cite{1} implies
that the subdivision $P=P_1P_2$ with the indicated properties is
unique; denote it by $\delta(l)$. Let us define a relation $R$ in
the following way: $R(l_1, l_2) \rightleftharpoons \delta(l_1) =
\delta(l_2)$. Condition e) of the definition of a periodic
structure obviously holds.

Condition f) follows from the graphic equality $H[\alpha(\mu),
\beta(\mu)]^{\varepsilon(\mu)}=$ $H[\alpha(\Delta(\mu)),$ $
\beta(\Delta(\mu))]^{\varepsilon(\Delta(\mu))}$ and Lemma 1.2.9
\cite{1}.

This proves the lemma.\end{proof}

Now let us fix a nonempty periodic structure $\langle {\mathcal
P}, R \rangle$. Item d) allows us to assume (after replacing the
variables $h_i, \ldots, h_{j-1}$ by $h_{j-1}^{-1}, \ldots,
h_i^{-1}$ on those sections $[i,j] \in {\mathcal P}$ for which
${\mathcal X}([i,j])=-1$) that $\varepsilon(\mu)=1$ for all $\mu
\in {\mathcal P}$. For a boundary $k$, we will denote by $(k)$ the
equivalence class of the relation $R$ to which it belongs.

 Let us construct
an oriented  graph $\Gamma$ whose set of vertices is the set of
$R$--equivalence classes. For each unknown $h_k$ lying on a
certain closed section from ${\mathcal P}$, we introduce an
oriented edge $e$ leading from $(k)$ to $(k+1)$ and an inverse
edge $e^{-1}$ leading from $(k+1)$ to $(k)$. This edge $e$ is
assigned the label $h(e) \rightleftharpoons h_k$ (respectively,
$h(e^{-1}) \rightleftharpoons h_k^{-1}$.) For every path
$r=e_1^{\pm 1} \ldots e_s^{\pm 1}$ in the graph $\Gamma$, we
denote by $h(r)$ its label $h(e_1^{\pm 1}) \ldots h(e_j^{\pm 1})$.
The periodic structure $\langle {\mathcal P}, R \rangle$ is called
{\em connected}, if the graph $\Gamma$ is connected. Suppose first
that $\langle {\mathcal P}, R \rangle$ is connected.  We can
always suppose that each boundary of $\Omega $ is a boundary
between two bases.

\begin{lemma} \label{2.9}
Let $H$ be a $P$-periodic solution of a generalized equation
$\Omega$,    $\langle {\mathcal P}, R \rangle = {\mathcal P}(H,
P)$; $c$ be a cycle in the graph $\Gamma$ at the vertex $(l)$;
$\delta(l)=P_1P_2$. Then there exists $n \in {\mathbb Z}$ such
that $H(c) = (P_2P_1)^n.$
\end{lemma}

\noindent \begin{proof} \hspace{2mm} If $e$ is an edge in the
graph $\Gamma$ with initial vertex $V'$ and terminal vertex $V''$,
and $P = P_1'P_2', $ $P = P_1 '' P_2''$ are two subdivisions
corresponding to the boundaries from $V'$, $V''$ respectively,
then, obviously, $H(e) = P_2' P^{n_k}P_1''$ $(n_k \in {\mathbb
Z})$. The claim is easily  proven by multiplying together the
values $H(E)$  for all the edges $e$ taking part in the cycle $c$.
\end{proof}

\begin{definition}
\label{2.51}
 A generalized equation $\Omega$ is called {\em
periodized} with respect to a given  periodic structure $\langle
{\mathcal P}, R \rangle$ of $\Omega$ , if for every  two cycles
$c_1$ and $c_2$ with the same initial vertex  in the graph
$\Gamma$ , there is a relation
 $[h(c_1), h(c_2)]=1$  in
 $F_{R(\Omega)}$.
\end{definition}

\subsection{ Case 1. Set $N\mathcal P$ is empty.}

Let $\Gamma_0$ be the subgraph of the graph $\Gamma$ having  the
same set of vertices and consisting of the edges $e$ whose labels
do not belong to ${\mathcal P}$. Choose a maximal subforest  $T_0$
in the graph $\Gamma_0$ and extend it to a maximal subforest $T$
of the graph $\Gamma$. Since $\langle {\mathcal P}, R \rangle$ is
connected by assumption, it follows that  $T$ is a tree. Let $v_0$
be an arbitrary vertex of the graph $\Gamma$ and $r(v_0, v)$ the
(unique) path from $v_0$ to $v$ all of whose vertices belong to
$T$. For every edge $e: v \rightarrow v'$ not lying in $T$, we
introduce a cycle $c_e = r(v_0, v) e (r(v_0, v'))^{-1}$. Then  the
fundamental group $\pi_1(\Gamma, v_0)$ is generated by the cycles
$c_e$ (see, for example, the proof of Proposition 3.2.1
\cite{LS}). This and decidability of the universal theory of a
free group imply that the property of a generalized equation ``to
be periodized with respect to a given periodic structure'' is
algorithmically decidable.

Furthermore, the set of elements

\begin{equation} \label{2.52}
\{h(e) \mid e \in T \} \cup \{h(c_e) \mid e \not \in T \}
\end{equation}

\noindent forms a basis of the free group with the set of
generators $\{h_k \mid h_k$ is an unknown lying on a closed
section from ${\mathcal P} \}$. If $\mu \in {\mathcal P}$, then
$(\beta(\mu)) = (\beta(\Delta(\mu)))$, $(\alpha(\mu)) =
(\alpha(\Delta(\mu)))$ by part f) from Definition \ref{above} and,
consequently, the word $h[\alpha(\mu), \beta(\mu)]
h[\alpha(\Delta(\mu)), \beta(\Delta(\mu))]^{-1}$ is the label of a
cycle $c'(\mu)$ from $\pi_1 (\Gamma, (\alpha(\mu)))$. Let $c(\mu)
= r(v_0, (\alpha(\mu)))c'(\mu) r(v_0, (\alpha(\mu)))^{-1}$. Then

\begin{equation} \label{2.53}
h(c(\mu)) = uh[\alpha(\mu), \beta(\mu)] h[\alpha(\Delta(\mu)),
\beta(\Delta(\mu))]^{-1} u^{-1},
\end{equation}

\noindent where $u$ is a certain word. Since $c(\mu) \in
\pi_1(\Gamma, v_0)$, it follows that $c(\mu) = b_\mu ( \{c_e \mid
e \not \in T \})$, where $b_\mu$ is a certain word in the
indicated generators  which can be effectively constructed (see
Proposition 3.2.1 \cite{LS}).

Let $\tilde{b}_\mu$ denote the image of the word $b_\mu$ in the
abelianization of $\pi (\Gamma ,v_0).$ Denote by $\widetilde{Z}$
the free abelian group consisting of formal linear combinations
$\sum_{e \not \in T} n_e \tilde{c}_e$ $(n_e \in {\mathbb Z})$, and
by $\widetilde{B}$ its subgroup generated by the elements
$\tilde{b}_\mu$ $(\mu \in {\mathcal P})$ and the elements
$\tilde{c}_e$ $(e \not \in T, \ h(e) \not \in {\mathcal P}).$ Let
$\widetilde{A} = \widetilde{Z} / \widetilde{B}$,
$T(\widetilde{A})$ the  torsion subgroups of the group
$\widetilde{A}$, and $\widetilde{Z}_1$ the preimage of
$T(\widetilde{A})$ in $\widetilde{Z}$. The group $\widetilde{Z} /
\widetilde{Z}_1$ is free; therefore, there exists a decomposition
of the form

\begin{equation} \label{2.54}
\widetilde{Z} = \widetilde{Z}_1 \oplus \widetilde{Z}_2, \
\widetilde{B} \subseteq \widetilde{Z}_1, \ (\widetilde{Z}_1 :
\widetilde{B}) < \infty .
\end{equation}

Note that it is possible to express effectively a certain basis
$\widetilde{\bar{c}}^{(1)}$,  $\widetilde{\bar{c}}^{(2)}$ of the
group $\widetilde{Z}$ in terms of the generators $\widetilde{c}_e$
so that for the subgroups $\widetilde{Z}_1$, $\widetilde{Z}_2$
generated by the sets $\widetilde{\bar{c}}^{(1)}$,
$\widetilde{\bar{c}}^{(2)}$ respectively, relation (\ref{2.54})
holds. For this it suffices, for instance, to look through the
bases one by one, using the fact that under the condition
$\widetilde{Z} = \widetilde{Z}_1 \oplus \widetilde{Z}_2$ the
relations $\widetilde{B} \subseteq \widetilde{Z}_1$,
$(\widetilde{Z}_1 : \widetilde{B}) <  \infty$ hold if and only if
the generators of the groups $\widetilde{B}$ and $\widetilde{Z}_1$
generate the same linear subspace over ${\bf Q}$, and the latter
is easily verified algorithmically. Notice, that a more economical
algorithm can be constructed by analyzing the proof of the
classification theorem for finitely generated abelian groups. By
Proposition 1.4.4  \cite{LS}, one can effectively construct a
basis $\bar{c}^{(1)}$, $\bar{c}^{(2)}$ of the free (non-abelian)
group $\pi_1(\Gamma, v_0)$ so that $\widetilde{\bar{c}}^{(1)}$,
$\widetilde{\bar{c}}^{(2)}$ are the natural images of the elements
$\bar{c}^{(1)}$, $\bar{c}^{(2)}$ in $\widetilde{Z}$.

Now assume that $\langle {\mathcal P}, R \rangle$ is an arbitrary
periodic structure of a periodized generalized equation $\Omega$,
not necessarily connected. Let $\Gamma_1, \ldots, \Gamma_r$ be the
connected components of the graph $\Gamma$. The labels of edges of
the component $\Gamma_i$ form in the equation $\Omega$ a union of
closed sections from ${\mathcal P}$; moreover, if a base $\mu \in
{\mathcal P}$ belongs to such a section, then its dual
$\Delta(\mu)$, by condition f) of Definition \ref{above}, also
possesses this property. Therefore, by taking for ${\mathcal P}_i$
the set of labels of edges from $\Gamma_i$ belonging to ${\mathcal
P}$, sections to which these labels belong, and bases $\mu \in
{\mathcal P}$ belonging to these sections, and restricting in the
corresponding way the relation $R$, we obtain a periodic connected
structure $\langle {\mathcal P}_i, R_i \rangle$ with the graph
$\Gamma_i$.

The notation $\langle {\mathcal P}', R'  \rangle$ $\subseteq$
$\langle {\mathcal P}, R \rangle$ means that  ${\mathcal P}'
\subseteq {\mathcal P},$ and the relation $R'$ is a restriction of
the relation $R$. In particular, $\langle {\mathcal P}_i, R_i
\rangle$ $\subseteq$ $\langle {\mathcal P}, R \rangle$ in the
situation described in the previous paragraph. Since $\Omega$ is
periodized, the periodic structure must be connected.

Let $e_1, \ldots, e_m$ be all the edges of the graph $\Gamma$ from
$T \smallsetminus T_0$. Since $T_0$ is the spanning forest of the
graph $\Gamma_0$, it follows that $h(e_1), \ldots, h(e_m) \in
{\mathcal P}$. Let $F(\Omega)$ be a free group generated by the
variables of $\Omega$. Consider in the  group $F(\Omega )$ a new
basis $A\cup\bar x$ consisting of $A,$ variables not belonging to
the closed sections from $\mathcal P$ (we denote by $\bar t$ the
family of these variables), variables $\{\,h(e) \mid e\in T\,\}$
and words $h(\bar c^{(1)})$, $h(\bar c^{(2)})$. Let $v_i$ be the
initial vertex of the edge $e_i$. We introduce new variables $\bar
u^{(i)}=\{\, u_{ie}\mid e\not\in T,\ e\not\in {\mathcal P}\,\},$
$\bar z^{(i)}=\{\,z_{ie}\mid e\not\in T,\, e\not\in {\mathcal
P}\,\}$ for $1\leqslant i\leqslant m,$ as follows
\begin{equation}\label{2.59}
u_{ie}=h(r(v_0,v_i)^{-1}h(c_e)h(r(v_0,v_i)),\end{equation}
\begin{equation}
\label{2.60} h(e_i)^{-1}u_{ie}h(e_i)=z_{ie}.\end{equation}

Notice, that without loss of generality we can assume that $ v_0$
corresponds to the beginning of the period $P$.

\begin{lemma} \label{2.10''} Let $\Omega $ be a consistent generalized
equation periodized with respect to a periodic structure $\langle
{\mathcal P},R\rangle  $ with empty set $N(\mathcal P )$.   Then
the following is true.
\begin{enumerate}\item [(1)]
One can choose the basis $\bar c^{(1)}$ so that  for any solution
$H$ of $\Omega$ periodic with respect to a period $P$ and
${\mathcal P}(H,P)=\langle {\mathcal P},R\rangle $ and any $c\in
\bar c^{(1)}$, $H(c)=P^n$, where $|n|< 2\rho$.

\item [(2)] In a fully residually free quotient of $F_{R(\Omega)}$
discriminated by  solutions from (1) the image of $\langle h(\bar
c^{(1)})\rangle $ is either trivial or a cyclic subgroup.

\item [(3)] Let $K$ be the subgroup of $F_{R(\Omega )}$ generated
by $\bar t$, $h(e), e\in T_0$, $h(\bar c^{(1)})$, $\bar u^{(i)}$
and $\bar z^{(i)}$, $i=1,\dots ,m.$ If $|\bar c^{(2)}|=s\geqslant
1$, then the group $F_{R(\Omega )}$ splits as a fundamental group
of a graph of groups with two vertices, where   one vertex group
is $K$ and the other is a free abelian group generated  by $h(\bar
c^{(2)})$ and $h(\bar c^{(1)})$. The corresponding edge group is
generated by $h(\bar c^{(1)})$. The other edges are loops at the
vertex with vertex group $K$, have stable letters $h(e_i),\
i=1,\ldots ,m$, and associated subgroups $\langle \bar
u^{i}\rangle ,$ $\langle \bar z^{i}\rangle .$ If $\bar
c^{(2)}=\emptyset,$ then there is no vertex with abelian vertex
group.

\item [(4)]  Let $A\cup \bar x$ be the generators of the group
$F_{R(\Omega )}$ constructed above. If $e_i\in {\mathcal P}\cap
T$, then the mapping defined as $h(e_i)\rightarrow u_{ie}^kh(e_i)
$ ($k$ is any integer) on the generator $h(e_i)$ and fixing all
the other generators can be extended to an automorphism of
$F_{R(\Omega )}.$

\item [(5)]  If $c\in\bar c^{(2)}$ and   $c'$  is a cycle with
initial vertex $v_0$, then the mapping defined by $h(c)\rightarrow
h(c')^kh(c)$ and fixing all the other generators can be extended
to an automorphism of $F_{R(\Omega )}.$

\end{enumerate}

\end{lemma}

\begin{proof} To prove assertion (1) we have to show that each
simple cycle in the graph $\Gamma _0$ has length less than
$2\rho$. This is obvious, because the total number of edges in
$\Gamma _0$ is not more than $\rho$ and corresponding variables do
not  belong to $\mathcal P$.

(2) The image of the group $\langle h(\bar c^{(1)})\rangle $ in
$F$ is cyclic, therefore
 one of the finite number of equalities $h(c_1)^{n}=h(c_2)^m$,
 where $c_1,c_2\in c^{(1)},\ n,m<2\rho$
must hold for any solution. Therefore in a fully residually free
quotient the group generated by the image of $\langle\ h(\bar
c^{(1)})\rangle $ is a cyclic subgroup.

To prove (3) we are to study in more detail how the unknowns
$h(e_i)$ ($1 \leqslant i \leqslant m$) can participate in the
equations from $\Omega^\ast$ rewritten in the set of variables
$\bar x\cup A.$

If $h_k$ does not lie on a closed section from ${\mathcal P}$, or
$h_k \not\in {\mathcal P}$, but $e \in T$ (where $h(e)=h_k$), then
$h_k$ belongs to the basis $\bar{x}\cup A$ and is distinct from
each of $h(e_1), \ldots, h(e_m)$. Now let $h(e) = h_k$, $h_k \not
\in {\mathcal P}$ and $e \not \in T$. Then $e=r_1c_er_2,$ where
$r_1,r_2$ are paths in $T$. Since $e \in \Gamma_0$, $h(c_e)$
belongs to $\langle c^{(1)}\rangle$ modulo commutation of cycles.
The vertices $(k)$ and $(k+1)$ lie in the same connected component
of the graph $\Gamma_0$, and hence they are connected by a path
$s$ in the forest $T_0$. Furthermore, $r_1$ and $sr_2^{-1}$ are
paths in the tree $T$ connecting the vertices $(k)$ and $v_0$;
consequently, $r_1 = s r_2^{-1}$. Thus, $e=sr_2^{-1}c_er_2$ and
$h_k = h(s) h(r_2)^{-1} h(c_e)h(r_2)$. The unknown $h(e_i)$ ($1
\leqslant i \leqslant m$) can occur in the right-hand side of the
expression obtained (written in the basis $\bar{x}\cup A$) only in
$h(r_2)$ and at most once. Moreover, the sign of this occurrence
(if it exists) depends only on the orientation of the edge $e_i$
with respect to the root $v_0$ of the tree $T$. If $r_2 = r_2'
e_i^{\pm 1}r_2 ''$, then
 all the occurrences of the unknown $h(e_i)$ in the words
$h_k$ written in the basis $\bar{x}\cup A$, with $h_k \not\in
{\mathcal P}$, are contained in the occurrences of words of the
form $h(e_i)^{\mp 1} h((r_2')^{-1}c_er_2')h(e_i)^{\pm 1}$, i.e.,
in occurrences of the form $h(e_i)^{\mp 1} h(c) h(e_i)^{\pm 1}$,
where $c$ is a certain cycle of the graph $\Gamma$ starting at the
initial vertex of the edge $e_i^{\pm 1}$.

Therefore all the occurrences of $h(e_i),\ i=1,\ldots ,m$ in the
equations corresponding to $\mu\not\in\mathcal P$ are of the form
$h(e_i^{-1})h(c)h(e_i)$. Also, $h(e_i)$ does not occur in the
equations corresponding to $\mu\in\mathcal P$ in the basis
$A\cup\bar x.$ The system $\Omega ^*$ is equivalent to the
following system in the variables $ \bar x, \bar z^{(i)},\bar
u^{(i)},A, i=1,\ldots ,m$
 : equations (\ref{2.59}), (\ref{2.60}),

\begin{equation}\label{2.61}
[u_{ie_1},u_{ie_2}]=1,\end{equation}
\begin{equation}
[h(c_1),h(c_2)]=1, \ c_1,c_2\in c^{(1)}, c^{(2)},\end{equation}
and a system $\bar \psi (h(e),e\in T\smallsetminus {\mathcal P},
h(\bar c^{(1)}), \bar t,\bar z^{(i)},\bar u^{(i)},A)=1$, such that
either $h(e_i)$ or $\bar c^{(2)}$ do not occur in $\bar\psi $. Let
$K=F_{R(\bar\psi )}.$ Then to obtain $F_{R(\Omega )}$ we fist take
an HNN extension of the group $K$ with abelian associated
subgroups generated by $\bar u ^{(i)}$ and $\bar z^{(i)}$ and
stable letters $h(e_i)$, and then extend the centralizer of the
image of $\langle \bar c^{(1)}\rangle$ by the free abelian
subgroup generated by the images of $\bar c^{(2)}.$

 Statements (4) and (5) follow from (3).

\end{proof}

We now introduce the notion of a {\em canonical group of
automorphisms corresponding to a connected periodic structure}.
\begin{definition} \label{canon} In the case when the family of bases $N{\mathcal P}$ is empty
automorphisms described in Lemma \ref{2.10''} for $e_1,\ldots
,e_m\in T\smallsetminus T_0$ and all $c_e$ for $e\in {\mathcal
P}\smallsetminus T$ generate the {\em canonical group of
automorphisms} $P_0$ corresponding to a connected periodic
structure.
\end{definition}

\begin{lemma}\label{2.10'}
Let $\Omega $ be a nondegenerate generalized equation with no
boundary connections, periodized with respect to the periodic
structure $\langle \mathcal P, R\rangle  $. Suppose that the set
$N\mathcal P$ is empty. Let  $H$ be a solution of $\Omega$
periodic with respect to a period $P$ and ${\mathcal
P}(H,P)=<{\mathcal P},R>$.
 Combining canonical automorphisms of
$F_{R(\Omega )}$
  one can get a solution $H^+$ of $\Omega$ with the property that
for any $h_k\in{\mathcal P}$ such that $H_k=P_2P^{n_k}P_1$ ($P_2$
and $P_1$ are an end and a beginning of $P$),
$H_k^+=P_2P^{n^+_k}P_1$, where $n_k, n_k^+ >  0$ and the numbers
$n_k^+$'s are bounded by a certain computable function $f_2(\Omega
,{\mathcal P},R)$. For all $h_k\not\in{\mathcal P}\ H_k=H_k^+.$
\end{lemma}

\begin{proof}   Let $\delta((k)) = P_1^{(k)} P_2^{(k)}$. Denote by
$t(\mu, h_k)$ the number of occurrences of the edge with label
$h_k$ in the cycle $c_{\mu}$, calculated taking into account the
orientation. Let

\begin{equation} \label{2.65}
H_k = P_2^{(k)} P^{n_k} P_1^{(k+1)}
\end{equation}

\noindent ($h_k$ lies on a closed section from ${\mathcal P}$),
where the equality in (\ref{2.65}) is graphic whenever $h_k \in
{\mathcal P}$. Direct calculations show that

\begin{equation} \label{2.66}
H(b_\mu) = P^{\sum_k t(\mu, h_k)(n_k+1)}.
\end{equation}
This equation implies that the vector $\{n_k\}$ is a solution to
the following system of Diophantine equations in variables
$\{z_k|h_k\in{\mathcal P}\}$:

\begin{equation}\label{2..}\sum _{h_k\in{\mathcal P}}t(\mu ,h_k)z_k+
\sum _{h_k\not\in{\mathcal P}}t(\mu ,h_k)n_k=0,
\end{equation}
$\mu\in {\mathcal P}.$ Note that the number of unknowns is
bounded, and coefficients of this system are bounded from above
($|n_k|\leqslant 2$ for $h_k\not\in {\mathcal P}$) by a certain
computable function of $\Omega, {\mathcal P},$ and $R$. Obviously,
$(P_2^{(k)})^{-1}H^+_kH^{-1}_kP_2^{(k)}=P^{n_k^+-n_k}$ commutes
with $H(c)$, where $c$ is a cycle such that $H(c)=P^{n_0},\ n_0<
2\rho$.

If system (\ref{2..}) has only one solution, then it is bounded.
Suppose it has infinitely many solutions. Then $(z_1,\ldots
,z_k,\ldots )$  is a composition of a bounded solution of
(\ref{2..}) and a linear combination of independent solutions of
the corresponding homogeneous system. Applying canonical
automorphisms from Lemma \ref{2.10''} we can decrease the
coefficients of this linear combination to obtain a bounded
solution $H^+$. Hence for $h_k=h(e_i),\ e_i\in{\mathcal P}$, the
value $H_k$ can be obtained by a composition of a canonical
automorphism (Lemma \ref{2.10''}) and a suitable bounded solution
$H^+$ of $\Omega$. \end{proof}

\subsection{ Case 2. Set $N\mathcal P $ is non-empty.}

We construct an oriented graph $B\Gamma$ with the same set of
vertices as $\Gamma$. For each item $h_k\not\in\mathcal P$ such
that $h_k$ lie on a certain closed section from $\mathcal P$
introduce an edge $e$ leading from $(k)$ to $(k+1)$ and $e^{-1}$
leading from $(k+1)$ to $(k)$. For each pair of bases $\mu ,\Delta
(\mu )\in \mathcal P$ introduce an edge $e$ leading from $(\alpha
(\mu))=(\alpha (\Delta (\mu )))$ to $(\beta (\mu))=(\beta (\delta
(\mu )))$ and $e^{-1}$ leading from $(\beta (\mu))$ to $(\alpha
(\mu ))$. For each base $\mu\in N\mathcal P$ introduce an edge $e$
leading from $(\alpha (\mu)$ to $(\beta (\mu))$ and $e^{-1}$
leading from $(\beta (\mu))$ to $(\alpha (\mu ))$. denote by
$B\Gamma _0$ the subgraph with the same set of vertices and edges
corresponding to items not from ${\mathcal P}$ and bases from
$\mu\in N\mathcal P$. Choose a maximal subforest $BT_0$ in the
graph $B\Gamma _0$ and extend it to a
 maximal subforest
$BT$ of the graph $B\Gamma $. Since $\mathcal P$ is connected,
$BT$ is a tree. The proof of the following lemma is similar to the
proof of Lemma \ref{2.9}.

\begin{lemma} \label{2.9b}
Let $H$ be a solution of a generalized equation $\Omega$ periodic
with respect to a period $P$, $\langle {\mathcal P}, R \rangle =
{\mathcal P}(H, P)$;  $c$ be a cycle in the graph $B\Gamma$ at the
vertex $(l)$; $\delta(l)=P_1P_2$. Then there exists $n \in
{\mathbb Z}$ such that $H(c) = (P_2P_1)^n.$
\end{lemma}

As we did in the graph $\Gamma$, we choose a vertex $v_0$. Let
$r(v_0,v)$ be the unique path in $BT$ from $v_0$ to $v$. For every
edge $e=e(\mu ):v\rightarrow v'$  not lying in $BT$, introduce a
cycle $c_{\mu}= r(v_0,v)e(\mu)r(v_0,v')^{-1}$. For every edge
$e=e(h_k ):V\rightarrow V'$  not lying in $BT$, introduce a cycle
$c_{h_k}= r(v_0,v)e(h_k)r(v_0,v')^{-1}$.

It suffices to restrict ourselves to the case of a connected
periodic structure. If $e=e(h_k)$, we denote $h(e)=h_k$; if
$e=e(\mu )$, then $h(e)=\mu $. Let $e_1, \ldots, e_m$ be all the
edges of the graph $B\Gamma$ from $BT \smallsetminus BT_0$. Since
$BT_0$ is the spanning forest of the graph $B{\Gamma} _0$, it
follows that $h(e_1), \ldots, h(e_m) \in {\mathcal P}$. Consider
in the free group $F(\Omega )$ a new basis $A\cup\bar x$
consisting of $A$, items  $h_k$ such that $h_k$ does not belong to
closed sections from $\mathcal P$ (denote this set by $\bar t$),
variables $\{h(e)|e\in T\}$ and words from $h( C^{(1)}),h(
C^{(2)})$, where the set $ C^{(1)}, C^{(2)}$ form a basis of the
free group $\pi (B\Gamma ,v_0)$, $C^{(1)}$ correspond to the
cycles that represent the identity in $F_{R(\Omega  )}$ (if $v$
and $v'$ are initial and terminal vertices of some closed section
in $\mathcal P$ and  $r$ and $r_1$ are different paths from $v$ to
$v'$, then $r(v_0,v)rr_1^{-1}r(v_0,v)^{-1}$ represents the
identity),  cycles $c_{\mu},\mu\in N\mathcal P$ and $c_{h_k}$,
$h_k\not\in \mathcal P$;
 and $C^{(2)}$ contains the rest of the basis of $\pi (B{\Gamma },v_0)$.

We study in more detail how the unknowns $h(e_i)$ ($1 \leqslant i
\leqslant m$) can participate in the equations from $\Omega^\ast$
rewritten in this basis.

If $h_k$ does not lie on a closed section from ${\mathcal P}$, or
$h_k=h(e), h(\mu)=h(e) \not\in {\mathcal P}$, but $e \in T$, then
$h(\mu)$ or $h_k$  belongs to the basis $\bar{x}\cup {A}$ and is
distinct from each of $h(e_1), \ldots, h(e_m)$. Now let $h(e) =
h(\mu)$, $h(\mu) \not \in {\mathcal P}$ and $e \not \in T$. Then
$e=r_1c_er_2$, where $r_1,r_2$ are path in $BT$ from $(\alpha (\mu
))$ to $v_0$ and from $(\beta (\mu ))$ to $v_0$. Since $e \in
B\Gamma_0$, the vertices $(\alpha (\mu))$ and $(\beta (\mu))$ lie
in the same connected component of the graph $B\Gamma _0$, and
hence are connected by a path $s$ in the forest $BT_0$.
Furthermore, $r_1$ and $sr_2^{-1}$ are paths in the tree $BT$
connecting the vertices $(\alpha (\mu ))$ and $v_0$; consequently,
$r_1 = s r_2^{-1}$. Thus, $e=sr_2^{-1}c_er_2$ and $h(\mu ) = h(s)
h(r_2)^{-1} h(c_e)h(r_2)$. The unknown $h(e_i)$ ($1 \leqslant i
\leqslant m$) can occur in the right-hand side of the expression
obtained (written in the basis $\bar{x}\cup A$) only in $h(r_2)$
and at most once. Moreover, the sign of this occurrence (if it
exists) depends only on the orientation of the edge $e_i$ with
respect to the root $v_0$ of the tree $T$. If $r_2 = r_2' e_i^{\pm
1}r_2 ''$, then
 all the occurrences of the unknown $h(e_i)$ in the words
$h (\mu)$ written in the basis $\bar{x}\cup A$, with $h(\mu)
\not\in {\mathcal P}$, are contained in the occurrences of words
of the form $h(e_i)^{\mp 1} h((r_2')^{-1}c_er_2')h(e_i)^{\pm 1}$,
i.e., in occurrences of the form $h(e_i)^{\mp 1} h(c) h(e_i)^{\pm
1}$, where $c$ is a certain cycle of the graph $B\Gamma$ starting
at the initial vertex of the edge $e_i^{\pm 1}$. Similarly, all
the occurences of the unknown $h(e_i)$ in the words $h_k$ written
in the basis $\bar{x}, A$, with $h_k \not\in {\mathcal P}$, are
contained in  occurrences of words of the form $h(e_i)^{\mp 1}
h(c) h(e_i)^{\pm 1}$.

Therefore all the occurences of $h(e_i),\ i=1,\ldots ,m$ in the
equations corresponding to $\mu\not\in {\mathcal P}$ are of the
form $h(e_i^{-1})h(c)h(e_i)$. Also,
 cycles from $C^{(1)}$ that represent the identity and not in $B\Gamma _0$
 are basis elements themselves.
This implies
\begin{lemma} \label{2n} \
\begin{enumerate}
\item [(1)] Let $K$ be the subgroup of $F_{R(\Omega )}$ generated
by $\bar t$, $h(e), e\in BT_0$, $h(C^{(1)})$ and $\bar u^{(i)},
\bar z^{(i)}, i=1,\ldots ,m,$ where elements $\bar z^{(i)}$ are
defined similarly to the case of empty $\mathcal NP$.

If $|C^{(2)}|=s\geqslant 1$, then the group $F_{R(\Omega )}$
splits as a fundamental group of a graph of groups with two
vertices, where one vertex group is $K$ and the other is a free
abelian group generated  by $h(C^{(2)})$ and  $ h(C^{(1)})$. The
edge group is generated by $h(C^{(1)})$. The other edges are loops
at the vertex with vertex group $K$ and have stable letters $h(e),
e\in BT\smallsetminus BT_0$. If $C^{(2)}=\emptyset,$ then there is
no vertex with abelian vertex group.

\item [(2)] Let $H$ be a solution of $\Omega$ periodic with
respect to a period $P$ and $\langle {\mathcal P},R\rangle =
{\mathcal P}(H,P).$ Let $P_1P_2$ be a partition of $P$
corresponding to the initial vertex of $e_i$. A transformation
$H(e_i)\rightarrow P_2P_1H(e_i)$, $i\in\{1,\ldots ,m\}$,  which is
identical on all the other elements from $A, H(\bar x)$, can be
extended to another solution of $\Omega ^*$. If $c$ is a cycle
beginning at the initial vertex of $e_i$, then the transformation
$h(e_i)\rightarrow h(c)h(e_i)$ which is identical on all other
elements from $A\cup\bar x$, is an automorphism of
$F_{R(\Omega)}.$

\item [(3)]  If $c(e)\in C^{(2)}$, then the transformation
$H(c(e))\rightarrow PH(c(e))$ which is identical on all other
elements from $A, H(\bar x)$ , can be extended to another solution
of $\Omega ^*$. A transformation $h(c(e))\rightarrow h(c)h(c(e))$
which is identical on all other elements from $A\cup\bar x$, is an
automorphism of $F_{R(\Omega )}.$\end{enumerate}
\end{lemma}

\begin {definition} If $\Omega $ is a nondegenerate   generalized equation
periodic with respect to a connected periodic structure $\langle
\mathcal P, R\rangle  $ and the set $N\mathcal P$ is non-empty, we
consider the group ${\bar A}(\Omega )$  of transformations  of
solutions of $\Omega ^*$, where $\bar A(\Omega )$ is generated by
the transformations  defined in Lemma \ref{2n}. If these
transformations are automorphisms, the group will be denoted
$A(\Omega ).$
\end{definition}

\begin{definition} In the case when for a connected periodic structure $\langle \mathcal P,R\rangle  ,$ the set
$C^{(2)}$ has more than one element or $C^{(2)}$ has one element,
and $C^{(1)}$ contains a cycle formed by edges $e$ such that
variables $h_k=h(e)$ are not from $\mathcal P$,  the periodic
structure will be called {\em singular}.
\end{definition} This definition coincides with the definition of singular
periodic structure given in \cite{KMIrc}) in the case of empty set
$\Lambda$.

Lemma \ref{2n} implies the following
\begin{lemma}\label{2.11}
Let $\Omega $ be a nondegenerate generalized equation with no
boundary connections, periodized with respect to a singular
periodic structure $\langle \mathcal P, R\rangle  $. Let $H$ be a
solution of $\Omega$ periodic with respect to a period $P$ and
$\langle {\mathcal P},R\rangle = {\mathcal P}(H,P).$ Combining
canonical automorphisms from $A(\Omega )$ one can get a solution
$H^+$ of $\Omega ^*$ with the following properties:
\bi
 \item[1)]
For any $h_k\in {\mathcal P}$ such that $H_k=P_2P^{n_k}P_1$ {\rm
(}$P_2$ and $P_1$ are an end and a beginning of $P${\rm )}
$H_k^+=P_2P^{n^+_k}P_1$, where $n_k, n_k^+\in {\mathbb Z}$;

\item[2)] For any $h_k\not\in {\mathcal P}$,  $H_k=H^+_k$;

\item[3)] For any base $\mu\not\in{\mathcal P}$, $H(\mu )=H^+(\mu
)$;

\item[4)] There exists a  cycle $c$ such that $h(c)\neq 1$ in
$F_{R(\Omega )}$ but $H^+(c)=1.$ \ei


\end{lemma}

Notice, that in the case described in the lemma, solution $H^+$
satisfies a proper equation. Solution $H^+$ is not necessarily a
solution of the generalized equation $\Omega $, but we will modify
$\Omega$ into a generalized equation $\Omega ({\mathcal P},BT)$.
This modification will be called the \emph{first minimal
replacement}. Equation  $\Omega ({\mathcal P},BT)$ will have the
following properties:

(1)
 $\Omega ({\mathcal P},BT)$ contains all the same parameter sections and
 closed sections which are not in $\mathcal P$, as $\Omega
 $;

(2) $H^+$ is a solution of $\Omega ({\mathcal P},BT)$;

(3) group $F_{R(\Omega ({\mathcal P},BT) )}$ is generated by the
same set of variables $h_1,\ldots ,h_{\delta}$;

(4)  $\Omega ({\mathcal P},BT)$ has the same set of bases as
$\Omega$ and possibly some new bases, but each new base is a
product of bases from $\Omega $;

(5) the mapping $h_i\rightarrow h_i$ is a proper homomorphism from
$F_{R(\Omega  )}$ onto $F_{R(\Omega ({\mathcal P},BT) )}$.

To obtain $\Omega ({\mathcal P},BT)$ we have to modify the closed
sections from $\mathcal P$.

The label of each cycle in $B\Gamma$ is a product of some bases
$\mu_1\cdots \mu _k$. Write a generalized equation
$\widetilde\Omega$ for the equations that say that $\mu _1\cdots
\mu _k=1$ for each cycle from $C^{(1)}$ representing the trivial
element  and for each cycle from $C^{(2)}$. Each $\mu _i$ is a
product $\mu_i=h_{i1}\cdots h_{it}$. Due to the first statement of
Lemma \ref{2.11}, in each product $H^+_{ij}H^+_{i,j+1}$ either
there is no cancellations between $H^+_{ij}$ and $H^+_{i,j+1}$, or
one of them is completely cancelled in the other. Therefore the
same can be said about each pair $H^+(\mu _i)H^+(\mu _{i+1})$, and
we can make a cancellation table without cutting items or bases of
$\Omega $.

Let $\widehat\Omega $ be a generalized equation obtained from
$\Omega$ by deleting bases from $\mathcal P\cup N\mathcal P$ and
items from $\mathcal P$ from the closed sections from $\mathcal
P$. Take a union of $\widetilde\Omega$ and $\widehat\Omega$ on the
disjoint set of variables, and add basic equations identifying in
$\widehat\Omega$ and $\widetilde\Omega$ the same bases that don't
belong to $\mathcal P$. This gives us $\Omega ({\mathcal P},BT).$

Suppose that
 $C^{(2)}$ for the equation $\Omega $ is either empty or contains one cycle.
 Suppose also that for each closed section from $\mathcal P$ in $\Omega $
 there exists a base $\mu$ such that the initial boundary of this section is
 $\alpha (\mu )$ and the terminal boundary is $\beta (\Delta (\mu ))$.

\begin{lemma}\label{4n} Suppose that the generalized equation $\Omega $
is periodized  with respect to a non-singular periodic structure
$\mathcal P$. Then for any periodic solution $H$ of $\Omega$ we
can choose a tree BT, some set of variables $S=\{h_{j_1},\ldots ,
h_{j_s}\}$ and a solution $H^+$ of $\Omega$ equivalent to $H$ with
respect to the group of canonical transformations $\bar A(\Omega
)$ in such a way that each of the bases $\lambda _i\in
BT\smallsetminus BT_0$ can be represented as $\lambda_i=\lambda
_{i1}h_{k_i}\lambda _{i2}$, where $h_{k_i}\in S$ and for any
$h_j\in S,$ $| H^+_{j}| < f_3| P|$, where $f_3$ is some
constructible function depending on $\Omega $. This representation
gives a new generalized equation $\Omega '$ periodic with respect
to a periodic structure $\mathcal P'$ with the same period $P$ and
all $h_{j}\in S$  considered as variables not from $\mathcal P'$.
The graph $B\Gamma '$ for the periodic structure $\mathcal P'$ has
the same set of vertices as $B\Gamma $, has empty set $C^{(2)}$
and $BT'=BT'_0$.

Let $c$ be a cycle from $C^{(1)}$ of minimal length, then
$H(c)=P^{n_c},$ where $|n_c|\leqslant 2\rho$. Using canonical
automorphisms from $A(\Omega )$ one can transform any solution $H$
of $\Omega$ into a solution $H^{+}$ such that for any $h_j\in S,$
$| H^{+}_j|\leqslant f_3d| c | .$ Let $\mathcal P '$ be a periodic
structure, in which all $h_i\in S$ are considered as variables not
from $\mathcal P '$, then $B\Gamma '$ has empty set $C^{(2)}$ and
$BT'=BT_0'$.
\end{lemma}

\begin{proof}   Suppose first that $C^{(2)}$ is empty. We prove
the statement of the lemma by induction on the number of edges in
$BT\smallsetminus BT_0$. It is true, when this set is empty.
Consider temporarily all   the edges in $BT\smallsetminus BT_0$
except one edge $e(\lambda)$ as edges corresponding to bases from
$N\mathcal P$. Then the difference between $BT_0$ and $BT$ is one
edge.

Changing $H(e(\lambda))$ by a transformation from $\bar A(\Omega
)$ we can  change only $H(e')$ for  $e'\in B\Gamma$ that could be
included into $BT\smallsetminus BT_0$ instead of $e$. For each
base $\mu\in N\mathcal P$, $H(\mu)=P_2(\mu)P^{n(\mu)}P_1(\mu)$,
for each base $\mu\in \mathcal P$,
$H(\mu)=P_2(\mu)P^{x(\mu)}P_1(\mu)$. For each cycle $c$ in
$C^{(1)}$ such that $h(c)$ represents the identity element we have
a linear equation in variables $x(\mu)$ with coefficients
depending on $n(\mu)$. We also know that this system has a
solution for arbitrary $x(\lambda)$ (where $\lambda\in
BT\smallsetminus BT_0$) and the other $x(\nu )$ are uniquely
determined by the value of $x(\lambda)$.

If we write for each variable $h_k\in {\mathcal P},$
$H_k=P_{2k}P^{y_k}P_{1k}$, then the positive unknowns $y_k$'s
satisfy the system of equations saying that $H(\mu)=H(\Delta (\mu
))$ for bases $\mu\in\mathcal P$ and equations saying that $\mu$
is a constant for bases $\mu\in N\mathcal P$. Fixing $x(\lambda)$
we automatically fix all the $y_k$'s. Therefore at least one of
the $y_k$ belonging to $\lambda$ can be taken arbitrary.   So
there exist some elements $y_k$ which can be taken as free
variables for the second system of linear equations. Using
elementary transformations over $\mathbb Z$ we can write the
system of equations for $y_k$'s in the form:
\begin{equation}\label{linear}
\begin{array}{cccc}
n_1y_1         & 0             & \cdots                  &=m_1y_k       +C_1             \\
             & n_2y_2         &   \cdots               &=m_2y_k     +C_2   \\

             &              & \ddots                      &    \\
\vdots   & \vdots   &                  & \ddots                  \\
             & \cdots   & n_{k-1}y_{k-1}       & =m_{k-1}y_k               +C_{k-1},         \\
\end{array}
\end{equation}
where $C_1,\ldots C_k$ are constants depending on parameters, we
can suppose that they are sufficiently large positive or negative
(small constants we can treat as constants not depending on
parameters). Notice that integers $n_1,m_1,\ldots
,n_{k-1},m_{k-1}$ in this system do not depend on parameters.  We
can always suppose that all $n_1,\ldots ,n_{k-1}$ are positive.
Notice that $m_i$ and $C_i$ cannot be simultaneously negative,
because in this case it would not be a positive solution of the
system. Changing the order of the equations we can write first all
equation with $m_i,C_i$ positive, then equations with negative
$m_i$ and positive $C_i$ and, finally,  equations with negative
$C_i$ and positive $m_i$.  The system will have  the form:
\begin{equation}\label{linear1}
\begin{array}{cccc}
n_1y_1         & 0             & \cdots                  &=|m_1|y_k       +|C_1|,             \\
               &              & \ddots       &                \\

& n_ty_t         &   \cdots               &=-|m_t|y_k     +|C_t|,   \\

             &              & \ddots       &                   \\
  &    &                  & \ddots                  \\
             &    & n_{s}y_{s}       & =|m_{s}|y_k               -|C_{s}|         \\
\end{array}
\end{equation}

If the last block (with negative $C_s$) is non-empty, we can take
a minimal $y_s$ of bounded value.
 Indeed, instead of $y_s$
we can always take a remainder of the division of $y_s$ by the
product $n_1\ldots n_{k-1}|m_1\ldots m_{k-1}|$, which is less than
this product (or by the product $n_1\ldots n_{k-1}|m_1\ldots
m_{k-1}|n_c$ if we wish to decrease $y_s$ by a multiple of $n_c$).
We respectively decrease $y_k$ and adjust $y_i$'s in the blocks
with positive $C_i$'s. If the third block is not present, we
decrease $y_k$ taking a remainder of the division of $y_k$ by
$n_1\ldots n_{k-1}$ (or by $n_1\ldots n_{k-1}n_c$) and adjust
$y_i$'s. Therefore for some $h_i$ belonging to a base which can be
included into $BT\smallsetminus BT_0$,  $\mid H^+(h_i)\mid <
f_3\mid P\mid .$ Suppose this base is $\lambda$, represent
$\lambda =\lambda _1h_i\lambda _2$. Suppose $e(\lambda
):v\rightarrow v_1$ in $B\Gamma$. Let $v_2,v_3$ be the vertices in
$B\Gamma$ corresponding to the initial and terminal boundary of
$h_k$. They would be the vertices in $\Gamma$, and $\Gamma$ and
$B\Gamma$ have the same set of vertices. To obtain the graph
$B\Gamma '$ from $B\Gamma$ we have to replace $e(\lambda)$ by
three edges $e(\lambda _1):v\rightarrow v_2$,
$e(h_k):v_2\rightarrow v_3$ and $e(\lambda _2):v_3\rightarrow
v_1$. There is no path in $BT_0$ from $v_2$ to $v_3$, because if
there were such  a path $p$, then we would have the equality
$h_k=h(c_1)h(p)h(c_2),$ in $F_{R(\Omega
 )},$ where $c_1$ and $c_2$ are cycles in $B\Gamma $ beginning in
vertices $v_2$ and $v_3$ respectively. Changing $H_k$ we do not
change $H(c_1),H(c_2)$ and $H(p)$, because all the cycles are
generated by cycles in $C^{(1)}$. Therefore there are paths
$r:v\rightarrow v_2$ and $r_1:v_3\rightarrow v_1$ in $BT_0,$ and
edges $e(\lambda _1), e(\lambda _2)$ cannot be included in
$BT'\smallsetminus BT'_0$ in $B\Gamma '$. Therefore $BT'=BT_0'$.
Now we can recall that all the edges except one in
$BT\smallsetminus BT_0$ were temporarily considered as edges in
$N\mathcal P$. We managed to decrease the number of such edges by
one. Induction finishes the proof.

If the set $C^{(2)}$ contains one cycle, we can temporarily
consider all the bases from $BT$ as parameters, and consider the
same system of linear equations for $y_i$'s. Similarly, as above,
at  least one $y_t$ can be bounded. We will bound as many  $y_i$'s
as we can. For the new periodic structure either $BT$ contains
less elements or the set $C^{(2)}$ is empty.

 The second part of the lemma follows from the remark that for
$\mu\in T$ left multiplication of $h(\mu )$ by $h(rcr^{-1})$,
where $r$ is the path in $T$ from $v_0$ to the initial vertex of
$\mu$, is an automorphism from $A(\Omega )$. \end{proof}

We call a solution $H^+$ constructed in Lemma \ref{4n} a {\em
solution equivalent to $H$ with  maximal number of short
variables}.

Consider now variables from $S$ as variables not from $\mathcal
P'$, so that for the equation $\Omega $ the sets $C^{(2)}$ and
$BT'\smallsetminus BT_0'$ are both empty. In this case we make the
{\em second minimal replacement}, which we will describe in the
lemma below.
\begin{definition} \label{overl} A pair of bases $\mu ,\Delta (\mu)$ is called an
overlapping pair if $\epsilon (\mu )=1$ and $\beta (\mu )> \alpha
(\Delta (\mu ))>  \alpha (\mu )$ or $\epsilon (\mu)=-1$ and $\beta
(\mu )< \beta (\Delta (\mu ))< \alpha (\mu )$. If a closed section
begins with $\alpha (\mu)$ and ends with $\beta (\Delta (\mu ))$
for an overlapping pair of bases we call such a pair of bases a
{\em principal overlapping pair} and say that a section is in {\em
overlapping form}.\end{definition}

Notice, that if $\lambda\in N\mathcal P$, then $H(\lambda )$ is
the same for any solution $H$, and we just write $\lambda$ instead
of $H(\lambda )$.

\begin{lemma} \label{3n'} Suppose that for the generalized equation $\Omega '$ obtained in
Lemma \ref{4n} the sets $C^{(2)}$ and $BT'\smallsetminus BT_0'$
are empty, $\mathcal P '$ is a non-empty periodic structure, and
each closed section from $\mathcal P '$ has a principal
overlapping pair. Then for each base $\mu\in\mathcal P '$ there is
a  fixed presentation for $h(\mu)=\prod (parameters)$ as a
product of elements $h(\lambda ), \lambda\in N\mathcal P$,
$h_k\not\in\mathcal P'$ corresponding to a path in $B\Gamma _0'$.
The maximal number of terms in this presentation is bounded by a
computable function of $\Omega$.
\end{lemma}

 \begin{proof}   Let $e$ be the edge in the graph $B\Gamma'$
corresponding to a base $\mu$ and suppose $e:v\rightarrow v'$.
There is a path $s$ in $BT'$ joining $v$ and $v'$ ,  and a cycle
$\bar c$ which is a product of cycles from $C^{(1)}$ such that
$h(\mu )=h(\bar c )h(s).$ For each cycle $c$ from $C^{(1)}$ either
$h(c)=1$ or $c$ can be written using only edges with labels not
from $\mathcal P'$; therefore, $\bar c$ contains only edges with
labels not from $\mathcal P '$. Therefore
\begin{equation}\label{cat}
h(\mu )=\prod (parameters)=h(\lambda _{i_1})\Pi _1\cdots h(\lambda
_{s_i})\Pi _s,\end{equation}
 where the doubles of all
$\lambda _i$ are parameters, and $\Pi _1,\ldots ,\Pi _s$ are
products of variables $h_{k_i}\not \in \mathcal P '$. \end{proof}

In the equality
\begin{equation}\label{cat1}
H(\mu )=H(\lambda _{i_1})\bar \Pi _1\cdots H(\lambda
_{s_i})\bar\Pi _s,\end{equation} where $\bar\Pi _1,\ldots ,\bar\Pi
_s$ are products of $H_{k_i}$ for variables $h_{k_i}\not \in
\mathcal P '$, the cancellations between two terms in the right
side are complete because the equality corresponds to a path in
$B\Gamma _0'$. Therefore the cancellation tree for the equality
(\ref{cat1}) can be situated on a horizontal axis with intervals
corresponding to $\lambda _i$'s directed either to the right or to
the left. This tree can be drawn on a $P$-scaled axis. We call
this one-dimensional tree a $\mu$-{\em tree}. Denote by
$I(\lambda)$ the interval corresponding to $\lambda$ in the
$\mu$-tree. If $I(\mu)\subseteq\bigcup _{\lambda _i\in N{\mathcal
P}}I(\lambda _i),$ then we say that $\mu$ is covered by
parameters. In this case a generalized equation corresponding to
(\ref{cat1}) can be situated on the intervals corresponding to
bases from $N{\mathcal P}$.

We can shift the whole $\mu$-tree to the left or to the right so
that in the new situation the uncovered part becomes covered by
the bases from $N\mathcal P$. Certainly, we have to make sure that
the shift is through the interval corresponding to a cycle in $
C^{(1)}$. Equivalently, we can shift any base belonging to the
$\mu$-tree through such an interval.

If $c$ is a  cycle from $C^{(1)}$ with shortest $H(c)$, then there
is a corresponding $c$-tree. Shifting this $c$-tree to the right
or to the left through the intervals corresponding to $H(c)$
bounded number of times we can cover every $H_i$, where $h_i\in S$
by a product $H(\lambda _{j_1})\bar\Pi_1\ldots H(\lambda
_{j_t})\bar \Pi _t$, where $\bar\Pi_ 1,\ldots ,\bar\Pi _t$ are
products of values of variables not from $\mathcal P$ and $\lambda
_{j_1},\ldots \lambda _{j_t}$ are bases from $N\mathcal P$.
Combining this covering together with the covering of $H(\mu )$ by
the product (\ref{cat1}), we obtain that
 $H([\alpha (\mu ),\beta (\Delta (\mu ))])$ is almost
covered by parameters, except for the short products $\bar \Pi$.
Let $h(\mu )$ be covered by \begin{equation}\label{cat2} h(\Lambda
_1)\Pi _1,\ldots ,h(\Lambda _s)\Pi _s,
\end{equation}
where $h(\Lambda _1),\ldots ,h(\Lambda _s)$ are parts completely
covered by parameters,  and $\Pi _1,\ldots ,\Pi _s$ are products
of variables not in $\mathcal P$. We also remove those bases from
$N\mathcal P$ from each $\Lambda _i$ which do not overlap with
$h(\mu )$. Denote by $f_4$ the maximal number of bases in
$N\mathcal P$ and $h_i\not \in\mathcal P$ in the covering
(\ref{cat2}).

If $\lambda _{i_1},\ldots ,\lambda _{i_s}$ are parametric bases,
then for any solution $H$ and any pair $\lambda _i,\lambda _j\in
\{\lambda _{i_1},\ldots ,\lambda _{i_s}\}$ we have either $\mid
H(\lambda _i)\mid <\mid H(\lambda _j)\mid$ or  $\mid H(\lambda
_i)\mid =\mid H(\lambda _j)\mid$ or $\mid H(\lambda _i)\mid >\mid
H(\lambda _j)\mid$. We call a {\em relationship between lengths of
parametric bases} a collection that consists of one such
inequality or equality for each pair of bases. There is only a
finite number of possible relationships between lengths of
parametric bases. Therefore we can talk about a parametric base
$\lambda$ of maximal length meaning that we consider the family of
solutions for which $H(\lambda )$ has maximal length.

\begin{lemma}\label{3n}Let $\lambda _{\mu}\in{N\mathcal P}$ be a base
of max length in the covering (\ref{cat2})
 for $\mu\in\mathcal P$. If  for a solution $H$ of $\Omega ,$ and for
each closed section $[\alpha (\mu ),\beta (\Delta (\mu)]$ in
$\mathcal P$,
 min$\ \mid H[\alpha (\nu),\alpha (\Delta (\nu))]\mid\leqslant \mid H(\lambda _{\mu })\mid,$
 where the minimum is taken for all pairs of overlapping bases for this section,
then one can transform $\Omega$ into one of the finite number
(depending on $\Omega $) of generalized equations $\Omega
(\mathcal P)$ which do not contain  closed sections from $\mathcal
P$ but contain the same other closed sections  except for
parametric sections. The content of  closed sections from
$\mathcal P$ is transferred using bases from $N\mathcal P$ to the
parametric part. This transformation is called the \emph{second
minimal replacement}.\end{lemma}

\begin{proof}   Suppose for a closed section $[\alpha
(\mu),\beta (\Delta (\mu ))]$ that there exists a base $\lambda$
in (\ref{cat2}) such that $\mid H(\lambda) \mid\geqslant \min \
(H(\alpha (\nu ),\alpha (\Delta (\nu ))),$ where the minimum is
taken for all pairs of overlapping bases for this section. We can
shift the cover  $H(\Lambda _1)\bar\Pi _1,\ldots ,H(\Lambda
_s)\bar\Pi _s$ through the distance $d_1=\mid H[\alpha (\mu
),\alpha (\Delta (\mu ))]\mid .$ Consider first the case when
$d_1\leqslant \mid H(\lambda )\mid$ for the largest base in
(\ref{cat2}). Suppose the part of $H(\mu )$ corresponding to
$\bar\Pi _i$ is not covered by parameters. Take the first base
$\lambda _j$ in (\ref{cat2}) to the right or to the left of $\bar
\Pi _i$ such that $\mid H(\lambda _j)\mid\geqslant d_1$. Suppose
$\lambda _j$ is situated to the left from $\bar\Pi _i$. Shifting
$\lambda _j$ to the right through a bounded by $f_4$ multiple of
$d_1$ we will cover $\bar\Pi _i$.

Consider  now the case when $d_1>  \mid H(\lambda )\mid$, but
there exists an overlapping pair $\nu ,\Delta (\nu)$ such that
$$d_2=\mid H[\alpha (\nu ),\alpha (\Delta (\nu ))]\mid\leqslant \mid H(\lambda ) \mid.$$
If the part of $H(\mu )$ corresponding to $\bar\Pi _i$ is not
covered by parameters, we take the first base $\lambda _j$ in
(\ref{cat2}) to the right or to the left of $\bar \Pi _i$ such
that $\mid H(\lambda _j)\mid\geqslant d_2$. Without loss of
generality we can suppose that $\lambda _j$ is situated to the
left of $\bar\Pi _i$. Shifting $\lambda _j$ to the right through a
bounded by $f_4$ multiple of $d_2$  we will cover $\bar\Pi _i$.

Therefore, if the first alternative in the lemma does not take
place, we can cover the whole section $[\alpha (\mu),\beta (\Delta
(\mu ))]$ by the bases from $N{\mathcal P}$, and transform $\Omega
$ into one of the finite number of generalized equations which do
not contain the closed section $[\alpha (\mu),\beta (\Delta (\mu
))]$ and have all the other non-parametric sections the same. All
the cancellations between two neighboring terms of any equality
that we have gotten are complete, therefore the coordinate groups
of new equations are quotients of $F_{R(\Omega )}$. \end{proof}

\subsection{Proof of Proposition \ref{pr:per}}\label{proofper}
The first statement follows from  Lemmas \ref{2.10''} and
\ref{2.10'}.   The second statement  follows from Lemmas \ref{2n},
\ref{2.11}, \ref{4n}, \ref{3n'}.

\section{Elimination process: splittings of coordinate groups}\label{se:5.3}
\subsection{Minimal solutions}
\label{5.5.1} Let $F = F(A\cup B)$ be a free group with basis
$A\cup B$, $\Omega$  be a generalized equation with constants from
$(A\cup B)^{\pm 1}$, and parameters $\Lambda$. Let $A(\Omega )$ be
an arbitrary  group of $(A\cup \Lambda )$-automorphisms of
$F_{R(\Omega )}$. For   solutions  $ H^{(1)}$ and $H^{(2)}$ of the
equation $\Omega$ in the group $F$ we write $H^{(1)}<  _{A(\Omega
)} H^{(2)}$ if there exists an endomorphism $\pi$ of the group $F$
which is an $(A,\Lambda )$-homomorphism, and an automorphism
$\sigma\in A(\Omega )$  such that the following conditions hold:
(1)  $\pi _{ H^{(2)}}=\sigma\pi _{ H^{(1)}}\pi $, (2) For all
active variables $d(H_k^{(1)})\leqslant d(H_k^{(2)})$ for all
$1\leqslant k\leqslant\rho$ and $d(H_k^{(1)})< d(H_k^{(2)})$ at
least for one such $k$.

We also define a relation $< _{cA(\Omega )}$ by the same way as $<
_{A(\Omega )}$ but with extra property: (3) for any $k,j$, if
$(H_k^{(2)})^{\epsilon}(H_j^{(2)})^{\delta}$ in non-cancellable,
then $(H_k^{(1)})^{\epsilon}(H_j^{(1)})^{\delta}$ in
non-cancellable ($\epsilon ,\delta =\pm 1$).  Obviously, both
relations are transitive.

A solution $\bar H$ of $\Omega$ is called {\em $A(\Omega
)$-minimal} if there is no any solution $ \bar H^\prime$ of the
equation $\Omega$ such that $\bar H^\prime < _{A(\Omega )} \bar
H.$
 Since the total length $\sum_{i = 1}^\rho l(H_i)$ of  a solution $\bar H$ is a
non-negative integer, every strictly decreasing chain of solutions
$\bar H >   \bar H^1 >   \ldots >   \bar H^k > _{A(\Omega )}
\ldots $ is finite. It follows that for every solution $\bar H$ of
$\Omega$ there exists a minimal solution $\bar H^0$ such that
$\bar H^0 < _{A(\Omega )}\bar H$.

\subsection{Splittings}\label{5.5.2} Assign to  some vertices $v$ of the
tree $T(\Omega)$ splittings of $F_{R(\Omega _v)}$ and groups of
canonical automorphisms of groups $F_{R(\Omega _v)}$ corresponding
to these splittings.

\begin{theorem} \label{th:spl} \

\bi \item[(1)] If $tp (v)=12,$ and the quadratic equation
corresponding to the quadratic section $\sigma$ is regular, then
either $F_{R(\Omega _v)}$ is a free product where one factor is a
closed surface group or there is a splitting of $F_{R(\Omega _v)}$
containing a QH-subgroup $Q$ corresponding to this quadratic
equation.

\item[(2)] If $tp (v)=2,$ and $\langle {\mathcal P}, R\rangle$ is
a singular periodic structure, then either $F_{R(\Omega _v)}$ is a
free product where one factor is a free abelian group or
$F_{R(\Omega _v)}$ has a splitting with an abelian vertex group
generated by the cycles $\bar c^{(1)}, \bar c^{(2)}$ and edge
group generated by $\bar c^{(1)}$ for this periodic structure.

\item[(3)] If $tp (v)=2,$ and $\langle {\mathcal P}, R\rangle$ is
a non-singular non-empty periodic structure, then $F_{R(\Omega
_v)}$ has an abelian splitting as an HNN extension with stable
letters corresponding to $h(e), e\in T\smallsetminus T_0.$

\item[(4)] Let $7\leqslant tp (v)\leqslant 10$.  There exists a
number $N=N(\Omega _v)$ such that if there is a branch of length
$N$ in $T(\Omega _v)$ beginning at $v$
 with each vertex $w$  of type  $7\leqslant tp (w)\leqslant  10,$
then $F_{R(\Omega _v )}$ splits as
 $F_{R(\Omega _v)}=F_{R({\rm Ker} \Omega )}\ast F(Z),$ where $F(Z)$ is a free group on $Z$.

\item[(5)] Let $tp (v)=15.$  There exists a number $N=N(\Omega
_v)$ such that if $T(\Omega_v)$ contains a branch with all
vertices of type $15$ beginning at $v$, then $F_{R(\Omega _v)}$
either has a splitting as in 3) or has a QH subgroup or both. \ei

Moreover, all the  splittings in Statements {\rm (1)--(5)} and the
number $N$ can be found effectively.
 \end{theorem}

 \begin{proof} Statement (1) follows from the argument in the description of Case 12.
 Statements (2) and (3) follow from Lemmas
 \ref{2.10''} and \ref{2n}.  Statement (4) follows from Lemma \ref{7-10}.
Statement (5) follows from Proposition \ref{3.4} to be proved
below. The statement about the effectiveness follows from the
effectiveness of the construction of a finite tree $T_0(\Omega )$
in the next two subsections. \end{proof}

 For each vertex $v$ such that $tp (v)=12$ the
group of automorphisms $A (\Omega _v)$ assigned to it is the
canonical  group of automorphisms of $F_{R(\Omega _v)}$
corresponding to $Q$ (and, therefore, identical on $\Lambda ).$

For each vertex $v$ such that $7\leqslant tp (v)\leqslant 10$ we
assign the group of automorphisms invariant with respect to the
kernel.

For each vertex $v$ such that $tp (v)=2,$ assign the group $\bar
A_v$ generated by the groups of
 automorphisms constructed in Lemma
\ref{2n} that applied to $\Omega _v$ and all possible non-singular
periodic structures of this equation.

Let $tp (v)=15$. Apply transformation $D_3$ and consider
$\Omega=\tilde \Omega _v$. Notice that the function $\gamma _i$ is
constant when $h_i$ belongs to some closed section of $\tilde
{\Omega _v}$. Applying $D_2$, we can suppose that the section
$[1,j+1]$ is covered exactly twice. We say now that this is a
quadratic section. Assign to the vertex $v$ the group of
automorphisms of $F_{R(\Omega )}$ acting identically on the
non-quadratic part.

\subsection{The finite subtree $T_0(\Omega )$ of
 $T(\Omega )$: cutting off long branches  }\label{5.5.3} For a generalized
equation $\Omega $ with parameters we construct a finite tree
$T_0(\Omega)$. Then we show that the subtree of $T(\Omega)$
obtained by tracing those paths in $T(\Omega)$ which actually can
happen for ``short'' solutions is a subtree of $T_0(\Omega )$.

According to Lemma \ref{3.2}, along an infinite path in
$T(\Omega)$ one can  either have $7\leqslant tp(v_k)\leqslant 10$
for all $k$ or $tp(v_k)=12$ for all $k$, or $tp(v_k)=15$ for all
$k$.
\begin{lemma}\label{3.3}{\rm[}Lemma 15 from \cite{KMIrc}{\rm]} Let $v_1\rightarrow v_2
\rightarrow\cdots\rightarrow v_k\rightarrow\cdots $ be an infinite
path in the tree $T(\Omega )$, and $7\leqslant tp(v_k)\leqslant
10$ for all $k$. Then among $\{\Omega _k\}$ some generalized
equation occurs infinitely many times. If $\Omega _{v_k}=\Omega
_{v_l}$, then $\pi (v_k,v_l)$ is an isomorphism invariant with
respect to the kernel.
\end{lemma}

\begin{lemma} Let $tp(v)=12.$ If a solution $\bar H$ of $\Omega _v$
is minimal with respect to the canonical group of automorphisms, then
there is a recursive function $f_0$ such that in the sequence
\begin{equation}\label{12'}
 (\Omega _v,\bar H)\rightarrow (\Omega _{v_1},\bar H^1)\rightarrow\cdots
 \rightarrow (\Omega _{v_i},\bar H ^i),\ldots,\end{equation}
corresponding to the path in $T(\Omega _v)$ and for the solution
$\bar H$, case 12 cannot be repeated more than $f_0$ times.
\end{lemma}
\begin{proof} If $\mu $ and $\Delta\mu$ both belong to the quadratic
section, then $\mu$ is called a {\em quadratic base}. Consider the
following set of generators for $F_{R(\Omega _v )}$: variables
from $\Lambda$ and quadratic bases from the active part. Relations
in this set of generators consist of the following three families.
\begin{enumerate}
\item [1)] Relations between variables in $\Lambda$.

\item [2)] If $\mu$ is an active base and $\Delta (\mu)$ is a
parametric base, and $\Delta (\mu )=h_{i}\cdots h_{i+t},$ then
there is a relation $\mu=h_i\cdots h_{i+t}$.

\item [3)] Since $\gamma _i$=2 for each $h_i$ in the active part
the product of $h_i\cdots h_j, $ where $[i,j+1]$  is a closed
active section, can be written in two different ways $w_1$ and
$w_2$ as a product of active bases. We write the relations
$w_1w_2^{-1}=1.$ These relations give a quadratic system of
equations with coefficients in the subgroup generated by $\Lambda
$.\end{enumerate}

When we apply the entire transformation in Case 12, the number of
variables is not increasing and the complexity of the generalized
equation is not increasing. Suppose the same generalized equation
is repeated twice in the sequence (\ref{12'}). for example,
$\Omega _j=\Omega _{j+k}$. Then $\pi (v_j,v_{j+k})$ is an
automorphism of $F_{R(\Omega _j)}$ induced by the automorphism of
the free product $\langle \Lambda\rangle  *B,$ where $B$ is a free
group generated by quadratic bases, identical on $\langle \Lambda
\rangle  $ and fixing all words $w_1w_2^{-1}$. Therefore, $\bar H
^j>  \bar H ^{j+k},$ which contradicts to the minimality of $\bar
H$. Therefore there is only a finite number (bounded by $f_0$) of
possible generalized equations that can appear in the sequence
(\ref{12'}).\end{proof}

Let $\bar H$ be a solution of the equation $\Omega$ with quadratic
part
 $[1,j+1]$.If $\mu $ belongs and $\Delta\mu$ does not  belong
to the quadratic section, then $\mu$ is called a {\em
quadratic-coefficient base}. Define the following numbers:
\begin{equation}\label{2.31}
d _1(\bar H)=\sum _{i=1}^{j}d(H_i),\end{equation}
\begin{equation}\label{2.32}
d _2(\bar H)=\sum _{\mu}d(H[\alpha (\mu),\beta
(\mu)]),\end{equation} where $\mu$ is a quadratic-coefficient
base.

\begin{lemma}
\label{2.8} Let $tp (v)=15$ For any solution $\bar H$ of $\Omega
_v$ there is a minimal solution $\bar H^+$, which is an
automorphic image of $\bar H$ with respect to the group of
automorphisms defined in the beginning of this section, such that
$$d_1(\bar H^+)\leqslant f_1(\Omega _v)\ \max \ \{d_2(\bar H^+),1\},$$
where $f_1(\Omega)$ is some recursive function.  \end{lemma}

\begin{proof} Consider instead of $\Omega _v$ equation $\Omega =(\tilde
\Omega _v)$ which does not have any boundary connections,
$F_{R(\Omega _v)}$ is isomorphic to $F_{R(\Omega )}.$ Consider a
presentation of $F_{R(\Omega _v)}$ in the set of generators
consisting of variables in the non-quadratic part and active
bases. Relations in this generating set consist of  the following
three families.
\begin{enumerate}
\item [1)] Relations between variables in the non-quadratic part.

\item [2)] If $\mu$ is a quadratic-coefficient base and  $\Delta
(\mu )=h_{i}\cdots h_{i+t}$ in the non-quadratic part, then there
is a relation $\mu=h_i\cdots h_{i+t}.$

\item [3)] Since $\gamma _i$=2 for each $h_i$ in the active part
the product $h_i\cdots h_j, $ where $[i,j+1]$  is a closed active
section, can be written in two different ways $w_1$ and $w_2$ as a
product of quadratic and quadratic-coefficient bases. We write the
relations $w_1w_2^{-1}=1.$
\end{enumerate}

Let $\bar H$ be a solution of $\Omega _v$ minimal with respect to
the canonical group of automorphisms of the
 free product $B_1*B,$ where $B$ is a free group generated by quadratic bases,
 and $B_1$ is a subgroup of $F_{R(\Omega _v)}$ generated by variables
 in the non-quadratic part, identical on $\langle \Lambda \rangle  $ and fixing all words $w_1w_2^{-1}$.

Consider the sequence
\begin{equation}\label{122}
 (\Omega ,\bar H)\rightarrow (\Omega _{v_1},\bar H^1)\rightarrow\ldots
 \rightarrow (\Omega _{v_i},\bar H^i),\ldots .\end{equation}

Apply now the entire transformations to the quadratic section of
$\Omega.$ As in the proof of the previous lemma, each time we
apply the entire transformation, we do not increase complexity
 and
do not increase the total number of items
 in the whole interval.
Since $\bar H$ is a solution of $\Omega _v$, if the same
generalized equation  appear in this sequence $2^{4^{j^2}}+1$
times then for some $j,j+k$ we have $\bar H^j > _c\bar H^{j+k}$,
therefore the same equation can only appear a bounded number of
times. Every quadratic base (except those that become matching
bases of length 1)
 in the quadratic part
can be transferred to the non-quadratic  part with the use of some
quadratic-coefficient base as a carrier base. This means that the
length of the transferred base is equal to the length of the part
of the quadratic-coefficient carrier base, which will then be
deleted. The double of the transferred base becomes a
quadratic-coefficient base. Because there are not more than $n_A$
 bases in the active part, this would give
$$d_1(\bar H')\leqslant n_A d_2(\bar H'),$$ for some solution $\bar H^+$
of the equation $\tilde\Omega _v$. But $\bar H^+$ is obtained from
the minimal solution $\bar H$ in a bounded number of steps.
\end{proof}

We call a path $v_1\rightarrow v_2 \rightarrow\ldots\rightarrow
v_k\rightarrow\ldots $ in $T(\Omega )$ for which $7\leqslant
tp(v_k)\leqslant 10$ for all $k$ or type 12  {\em prohibited} if
some generalized equation with $\rho$ variables occurs among
$\{\Omega _{v_i} \mid 1\leqslant i\leqslant \ell\}$ at least
$2^{(4\rho ^2)}+1$ times. We will define below also prohibited
paths in $T(\Omega )$, for which $tp(v_k)=15$ for all $k$. We will
need some auxiliary definitions.

Introduce a new parameter
$$\tau _v'=\tau _v+\rho -\rho _{v}',$$
where $\rho $ is the number of variables
 of the initial equation $\Omega $ and $\rho _{v}'$ the number of free
variables belonging to the non-active sections of the equation
$\Omega _v.$ We have $\rho _{v}'\leqslant \rho$ (see the proof of
Lemma \ref{3.2}), hence $\tau _v'\geqslant 0$. In addition if
$v_1\rightarrow v_2$ is an auxiliary edge, then $\tau _2'< \tau
_1'.$

Define by the joint induction on $\tau _v'$ a finite subtree
$T_0(\Omega _v)$ and a natural number $s(\Omega _v)$. The tree
$T_0(\Omega _v)$ will have $v$ as a root and consist of some
vertices and edges of $T(\Omega ) $ that lie higher than $v$. Let
$\tau _v'$=0; then in $T(\Omega )$ there can not be auxiliary
edges and vertices of type 15 higher than $v.$ Hence a subtree
$T_0(\Omega _v)$ consisting of vertices $v_1$ of $T(\Omega)$ that
are higher than $v$, and for which the path from $v$ to $v_1$ does
not contain prohibited subpaths, is finite.

Let now \begin{equation}\label{so} s(\Omega
_v)=\max_w\max_{\langle {\mathcal P},R\rangle  }\{\rho
_wf_2(\Omega _w,{\mathcal P},R),\ f_4 (\Omega _w',{\mathcal
P},R)\},\end{equation} where $w$ runs through all the vertices of
$T_0(v)$ for which $tp (w)=2$,
 $\Omega _w$ contains non-trivial non-parametric sections,
 $\langle {\mathcal P},R\rangle  $ is the set of non-singular periodic structures
of the equation $\tilde {\Omega} _w$,
  $f_2$ is a function appearing
in Lemma \ref{2.10'} ($f_2$ is present only if a periodic
structure has empty set $N{\mathcal P}$) and $\Omega _w'$ is
constructed as in Lemma \ref{4n}, where $f_4$ is a function
appearing in covering \ref{cat2}.

Suppose now that $\tau _v'>  0$ and that for all $v_1$ with $\tau
_{v_1}'< \tau _v'$ the tree $T_0(\Omega _{v_1})$ and the number
$s(\Omega _{v_1})$ are already defined. We begin with the
consideration of the paths
\begin{equation}\label{3.6}
r=v_1\rightarrow v_2\rightarrow \ldots\rightarrow
v_m,\end{equation} where $tp(v_i)=15\ (1\leqslant i\leqslant m)$.
We have $\tau _{v_i}'=\tau _v'.$

Denote by $\mu _i$ the carrier base of the equation $\Omega
_{v_i}$. The path (\ref{3.6}) will be called $\mu $-reducing if
$\mu _1=\mu$ and either there are no auxiliary edges from the
vertex $v_2$ and $\mu$ occurs in the sequence $\mu _1,\ldots ,\mu
_{m-1}$ at least twice, or there are auxiliary edges
$v_2\rightarrow w_1, v_2\rightarrow w_2\ldots ,v_2\rightarrow w_k$
from $v_2$ and $\mu$ occurs in the sequence $\mu _1,\ldots ,\mu
_{m-1}$ at least $\max _{1\leqslant i\leqslant k}s(\Omega _{w_i})$
times.

The path (\ref{3.6}) will be called {\em prohibited}, if it can be
represented in the form
\begin{equation}\label{3.7} r=r_1s_1\cdots r_ls_lr',\end{equation}
such that for some sequence of bases $\eta _1,\ldots ,\eta _l$ the
following three properties hold:
\begin{enumerate}
\item [1)] every base occurring at least once in the sequence $\mu
_1,\ldots ,\mu _{m-1}$ occurs at least $40n^2f_1(\Omega
_{v_2})+20n+1$ times in the sequence $\eta _1,\ldots ,\eta _l$,
where $n$ is the number of pairs of bases in equations $\Omega
_{v_i}$, \item [2)] the path $r_i$ is $\eta _i$-reducing; \item
[3)] every transfer base of some equation of path $r$ is a
transfer base of some equation of path $r'$.\end{enumerate} The
property of path (\ref{3.6}) of being prohibited is
algorithmically decidable. Every infinite path (\ref{3.6})
contains a prohibited subpath. Indeed, let $\omega$ be the set of
all bases occurring in the sequence $\mu _1,\ldots ,\mu _m,\ldots
$ infinitely many times, and $\tilde\omega$ the set of all bases,
that are transfer bases of infinitely many equations $\Omega
_{v_i}$. If one cuts out some finite part in the beginning of this
infinite path, one can suppose that all the bases in the sequence
$\mu _1,\ldots ,\mu _m,\ldots $ belong to $\omega$ and each base
that is a transfer base of at least one equation, belongs to
$\tilde\omega$. Such an infinite path for any $\mu\in\omega$
contains infinitely many non-intersecting $\mu$-reducing finite
subpaths. Hence it is possible to construct a subpath (\ref{3.7})
of this path satisfying the first two conditions in the definition
of a prohibited subpath. Making $r'$ longer, one obtains a
prohibited subpath.

Let $T'(\Omega _v)$ be a subtree of $T(\Omega _v)$ consisting of
the vertices $v_1$ for which the path from $v$ to $v_1$ in
$T(\Omega )$ contains neither prohibited subpaths nor vertices
$v_2$ with $\tau _{v_2}'< \tau _v',$ except perhaps $v_1$. So the
terminal vertices of $T'(\Omega _v)$ are either vertices $v_1$
such that $\tau _{v_1}'< \tau _v',$ or terminal vertices of
$T(\Omega _v)$. A subtree $T'(\Omega _v)$ can be effectively
constructed. $T_0(\Omega _v)$ is obtained by attaching of
$T_0(\Omega _{v_1})$ (already constructed by the induction
hypothesis) to those terminal vertices $v_1$ of $T'(\Omega _v)$
for which $\tau _{v_1}'< \tau _v'.$ The function $s(\Omega _v)$ is
defined by (\ref{so}). Let now $T_0(\Omega)=T_0(\Omega _{v_0}).$
This tree is finite by construction.

\subsection{Paths corresponding to minimal solutions of $\Omega$ are in $T_0(\Omega )$}\label{5.5.4}
Notice, that if $tp (v)\geqslant 6$ and $v\rightarrow w_1,\ldots
,v\rightarrow w_m$ is the list of principal outgoing edges from
$v$, then the generalized equations $\Omega _{w_1},\ldots ,\Omega
_{w_m}$ are obtained from $\Omega _v$ by the application of
several elementary transformations. Denote by $e$ a function that
assigns a pair $(\Omega _{w_i},\bar H ^{(i)})$ to the pair
$(\Omega _v,\bar H).$ For $tp (v)=4,5$ this function is identical.

If $tp (v)=15$ and there are auxiliary edges from the vertex $v$,
then the carrier base $\mu$ of the equation $\Omega _v$ intersects
$\Delta (\mu)$. For any solution $\bar H$ of the equation $\Omega
_v$ one can construct a solution $\bar H'$ of the equation $\Omega
_{v'}$ by $H'_{\rho _v+1}=H[1,\beta (\Delta (\mu))].$ Let
$e'(\Omega _v,\bar H)=e(\Omega _{v'},\bar H').$

In the beginning of this section we assigned to  vertices $v$ of
type 12, 15, 2 and such that $7\leqslant tp(v)\leqslant 10$ of $
T(\Omega)$
the groups of   automorphisms $A(\Omega _v)$.%
 Denote by $Aut (\Omega )$ the group of  automorphisms
of $F_{R(\Omega )}$ , generated by all groups \newline $\pi
(v_0,v)A(\Omega _v)\pi (v_0,v)^{-1}$, $v\in T_0(\Omega)$. (Here
$\pi (v_0,v)$ is an isomorphism, because $tp(v)\not = 1$.) We are
to formulate the main technical result of this section. The
following proposition states that every minimal solution of a
generalized equation $\Omega$ with respect to the group $A(\Omega
)$ either factors through one of the finite family of proper
quotients of the group $F_{R(\Omega )}$ or (in the case of a
non-empty parametric part) can be transferred to the parametric
part.
\begin{prop} \label{3.4} For any solution $\bar H$ of a generalized equation
$\Omega $ there exists a  terminal vertex $w$ of the tree
$T_0(\Omega )$ having type 1 or 2,  and a solution $\bar H^{(w)}$
of a generalized equation $\Omega _w$ such that
\begin {enumerate}
\item [(1)] $\pi _{\bar H}=\sigma\pi(v_0,w)\pi _{\bar H^{(w)}}\pi\
$ where $\pi$ is an endomorphism of a free group $F$ $\sigma\in
Aut (\Omega )$; \item [(2)] if $tp (w)=2$ and the equation $\Omega
_{w}$ contains  nontrivial non-parametric sections, then there
exists a primitive cyclically reduced word $P$ such that $\bar
H^{(w)}$ is periodic with respect to ${\mathcal P}$ and one of the
following conditions holds:
\begin{enumerate}\item the equation $\Omega _w$ is singular with
respect  to a periodic structure
\newline ${\mathcal P}(\bar H^{(w)},P) $ and the first minimal
replacement can be applied, \item it is possible to apply the
second minimal replacement and make the family of closed sections
in $\mathcal P$ empty.
\end{enumerate}
\end{enumerate}
\end{prop}

Construct a directed tree with paths from the initial vertex
\begin{equation}\label{3.8}
(\Omega ,\bar H)=(\Omega _{v_0},\bar H^{(0)})\rightarrow (\Omega
_{v_1}, \bar H^{(1)})\rightarrow\ldots \rightarrow(\Omega
_{v_u},\bar H^{(u)})\rightarrow\ldots\end{equation} in which the
$v_i$ are the vertices of the tree $T(\Omega )$ in the following
way. Let $v_1=v_0$ and let $\bar H^{(1)}$ be some solution of
 the equation $\Omega $, minimal with respect to the group of automorphisms $A(\Omega v_0)$
with the property $\bar H\geqslant \bar H^{(1)}.$

Let $i\geqslant 1$ and suppose the term $(\Omega _{v_i},\bar
H^{(i)})$ of the sequence (\ref{3.8}) has been already
constructed.  If $7\leqslant tp(v_i)\leqslant 10$ or $tp(v_i)=12$
and there exists a minimal solution $\bar H^+$ of $\Omega _{v_i}$
such that $\bar H^+< \bar H^{(i)}$ , then we set $v_{i+1}=v_i$,
$\bar H^{(i+1)}=\bar H^+.$

If $tp(v_i)=15,\  v_i\neq v_{i-1}$ and there are auxiliary edges
from vertex $v_i$: $v_i\rightarrow w_1,\ldots ,v_i\rightarrow w_k$
(the carrier base $\mu$ intersects with its double $\Delta
(\mu)$),  then there exists a primitive word $P$ such that
\begin{equation}\label{3.9}
H^{(i)}[1,\beta (\Delta (\mu))]\equiv P^rP_1, r\geqslant 2,\
P\equiv P_1P_2,\end{equation} where $\equiv$ denotes a graphical
equality. In this case the path (\ref{3.8}) can be continued along
several possible edges of $T(\Omega )$.

For each group of automorphisms assigned to vertices of type 2 in
the trees $T_0 (\Omega _{w_i})$, $i=1,\ldots ,k$ and non-singular
periodic structure including the closed section $[1,\beta (\Delta
(\mu)]$ of the equation $\Omega _{v_i}$ and corresponding to
solution $\bar H^{(i)}$ we replace $\bar H^{(i)}$ by a solution
$\bar H^{(i)+}$ with maximal number of short variables (see the
definition after Lemma \ref{4n}). This collection of short
variables can be different for different periodic structures.
Either all the variables in $\bar H^{(i)+}$ are short or there
exists a parametric base $\lambda _{\rm max}$ of maximal length in
the covering \ref{cat2}. Suppose there is a $\mu$-reducing path
(\ref{3.6}) beginning at $v_i$ and corresponding to $\bar
H^{(i)+}$. Let $\mu _1,\ldots ,\mu _m$ be the leading bases of
this path. Let ${\tilde H}^1=H^{(i)+},\ldots ,{\tilde H}^j$ be
solutions of the generalized equations corresponding to the
vertices of this path. If for some $\mu _i$ there is an inequality
$d({\tilde H}^j[\alpha (\mu _i),\alpha (\Delta (\mu
_i))])\leqslant d (\lambda _{\rm max}),$ we set $(\Omega
_{v_{i+1}},\bar H^{(i+1)})=e'(\Omega _{v_i},\bar H^{(i)})$ and
call the section $[1,\beta (\Delta (\mu ))]$ which becomes
non-active, {\em potentially transferable}.

If there is a singular periodic structure in a  vertex of type 2
of some tree $T_0(\Omega _{w_i}), i\in\{1,\ldots ,k\},$ including
the closed section $[1,\beta (\Delta (\mu)]$ of the equation
$\Omega _{v_i}$ and corresponding to the  solution $\bar H^{(i)}$,
we also include the possibility \newline $(\Omega _{v_{i+1}},\bar
H^{(i+1)})=e'(\Omega _{v_i},\bar H^{(i)}).$

In all of the other cases we set $(\Omega _{v_{i+1}},\bar
H^{(i+1)})=e(\Omega _{v_i},\bar H^{(i)+}),$ where
 $\bar H^{(i)+}$ is a solution with maximal number of short variables and minimal solution of
$\Omega _{v_{i}}$ with respect to the canonical group of
automorphisms $P_{v_i}$ (if it exists). The path (\ref{3.8}) ends
if $tp (v_i)\leqslant 2.$

We will show that in the path (\ref{3.8}) $v_i\in T_0(\Omega)$. We
use induction on $\tau '$.
Suppose $v_i\not\in T_0(\Omega ),$ and let $i_0$ be the first of
such numbers. It follows from the construction of $T_0(\Omega )$
that there exists $i_1< i_0$ such that the path from $v_{i_1}$
into $v_{i_0}$ contains a subpath prohibited in the construction
of $T_2(\Omega _{v_{i_1}}).$ From the minimality of $i_0$ it
follows that this subpath goes from $v_{i_2}\ \ (i_1\leqslant i_2<
i_0)$ to $v_{i_0}$. It cannot be that $7\leqslant tp(v_i)\leqslant
10$ or $tp(v_i)=12$ for all $i_2\leqslant i\leqslant i_1$, because
there will be two indices $p< q$ between $i_2$ and $i_0$ such that
$\bar H^{(p)}=\bar H^{(q)}$, and this gives a contradiction,
because in this case it must be by construction $v_{p+1}=v_p$. So
$tp(v_i)=15$ ($i_2\leqslant i\leqslant i_0$).

Suppose we have a subpath (\ref{3.6}) corresponding to the
fragment
\begin{equation}\label{3.11}
(\Omega _{v_1},\bar H^{(1)})\rightarrow (\Omega _{v_2}, \bar
H^{(2)})\rightarrow\ldots \rightarrow(\Omega _{v_m},\bar
H^{(m)})\rightarrow\ldots\end{equation} of the sequence
(\ref{3.8}). Here $v_1,v_2,\ldots,v_{m-1}$ are vertices of the
tree $T_0(\Omega)$, and for all vertices $v_i$ the edge
$v_i\rightarrow v_{i+1}$ is principal.

As before, let $\mu _i$ denote the carrier base of $\Omega
_{v_i}$, and $\omega =\{\mu _1,\ldots ,\mu _{m-1}\},$ and
$\tilde\omega $ denote the set of such bases which are transfer
bases for at least one equation in (\ref{3.11}). By $\omega _1$
denote the set of such bases $\mu $ for which either $\mu$ or
$\Delta (\mu)$ belongs to $\omega\cup\tilde\omega $; by $\omega
_2$ denote the set of all the other bases. Let
$$\alpha (\omega)=\min(\min _{\mu\in\omega _2}\alpha (\mu),j),$$
where $j$ is the boundary between active and non-active sections.
Let $$X_{\mu}\circeq H[\alpha (\mu),\beta (\mu)].$$ If $(\Omega
,\bar H)$ is a member of sequence (\ref{3.11}), then denote
\begin{equation}\label{3.12}
d _{\omega}(\bar H)=\sum _{i=1}^{\alpha (\omega )-1}d(H_i),
\end{equation}

\begin{equation}\label{3.13}
\psi _{\omega}(\bar H)=\sum _{\mu\in\omega
_1}d(X_{\mu})-2d_{\omega}(\bar H).\end{equation}

Every item $h_i$ of the section $[1,\alpha (\omega)]$ belongs to
at least two bases, and both bases are in $\omega _1$, hence $\psi
_{\omega}(\bar H)\geqslant 0.$

Consider the quadratic part of $\tilde\Omega _{v_1}$ which is
situated to the left of $\alpha(\omega )$. The  solution $\bar
H^{(1)}$ is minimal with respect to the canonical group of
automorphisms corresponding to this vertex. By Lemma \ref{2.8} we
have
\begin{equation}
\label{3.14} d_1(\bar H^{(1)})\leqslant f_1(\Omega _{v_1})d_2(\bar
H^{(1)}).
\end{equation}

 Using this inequality we estimate the length of the interval
participating in the process $d_{\omega}(\bar H^{(1)})$ from above
by a product of $\psi _{\omega }$ and some function depending on
$f_1$. This will be inequality \ref{3.19}. Then we will show that
for a prohibited subpath the length of the participating interval
must be reduced by more than this figure (equalities \ref{3.27},
\ref{3.28}). This will imply that there is no prohibited subpath
in the path \ref{3.11}.

Denote by $\gamma _i(\omega)$ the number of bases $\mu\in\omega
_1$ containing $h_i$. Then
\begin{equation}\label{3.15}
\sum _{\mu\in\omega _1}d(X_{\mu}^{(1)})=\sum
_{i=1}^{\rho}d(H_i^{(1)}) \gamma _i(\omega),\end{equation} where
$\rho =\rho (\Omega _{v_1}).$ Let $$I=\{i\mid 1\leqslant
i\leqslant\alpha (\omega)-1 \&\gamma _i=2\}$$ and $$J=\{i\mid
1\leqslant i\leqslant\alpha (\omega)-1 \&\gamma _i>  2\}.$$ By
(\ref{3.12})
\begin{equation}\label{3.16}
d_{\omega}(\bar H^{(1)})=\sum _{i\in I}d(H_i^{(1)})+ \sum _{i\in
J}d(H_i^{(1)})= d_1(\bar H^{(1)})+\sum _{i\in
J}d(H_i^{(1)}).\end{equation} Let  $(\lambda ,\Delta(\lambda))$ be
a pair of quadratic-coefficient bases of the equation $\tilde
\Omega _{v_1}$, where $\lambda$ belongs to the nonquadratic part.
This pair can appear only from the bases $\mu\in\omega _1$. There
are two types of quadratic-coefficient bases.

{\em Type} 1. $\lambda$ is situated to the left of the boundary
$\alpha (\omega)$. Then $\lambda$ is formed by items $\{h_i\mid
i\in J\}$ and hence $d(X_{\lambda})\leqslant\sum _{i\in
J}d(H_i^{(1)}).$ Thus the sum of the lengths
$d(X_{\lambda})+d(X_{\Delta (\lambda)})$ for quadratic-coefficient
bases of this type is not more than $2n\sum _{i\in
J}d(H_i^{(1)}).$

{\em Type} 2. $\lambda$ is
 situated to the right of the boundary $\alpha
(\omega)$. The sum of length of the quadratic-coefficient bases of
the second type is not more than $$2\sum _{i=\alpha
(\omega)}^{\rho}d(H_i^{(1)})\gamma _i(\omega).$$

 We have \begin{equation}\label{3.17}
 d_2(\bar H^{(1)})\leqslant 2n\sum _{i\in J}d(H_i^{(1)})+2\sum _{i=\alpha (\omega)}^{\rho}d(H_i^{(1)})\gamma
_i(\omega).\end{equation} Now (\ref{3.13}) and (\ref{3.15}) imply
\begin{equation}\label{3.18}
\psi _{\omega}(\bar H^{(1)}_i)\geqslant \sum _{i\in
J}d(H_i^{(1)})+\sum _{i=\alpha (\omega)}^{\rho}d(H_i^{(1)})\gamma
_i(\omega).\end{equation} Inequalities (\ref{3.14}),
(\ref{3.16}),(\ref{3.17}),(\ref{3.18}) imply
\begin{equation}\label{3.19}
d_{\omega}(\bar H^{(1)})\leqslant \max \{\psi _{\omega}(\bar
H^{(1)})(2nf_1(\Omega _{v_1})+1),\ f_1(\Omega _{v_1})\}.
\end{equation}

 From the definition of Case 15 it follows that all the words
 $H^{(i)}[1,\rho _i+1]$ are the ends of the word $H^{(1)}[1,\rho _1+1]$,
 that is
 \begin{equation}\label{3.20}
 H^{(1)}[1,\rho _1+1]\doteq U_iH^{(i)}[1,\rho _i+1].\end{equation}
 On the other hand bases $\mu\in\omega _2$  participate in these
 transformations neither as carrier bases nor as transfer bases; hence
 $H^{(1)}[\alpha (\omega ),\rho _1+1]$ is the end of the word $H^{(i)}[
 1,\rho _i+1]$, that is
 \begin{equation} \label{3.21}
 H^{(i)}[1,\rho _i+1]\doteq V_iH^{(1)}[\alpha
(\omega ),\rho _1 +1].\end{equation} So we have
\begin{equation}\label{3.22} d_{\omega}(\bar H
^{(i)})-d_{\omega}(\bar H ^{(i+1)})=
d(V_i)-d(V_{i+1})=d(U_{i+1})-d(U_{i})= d(X_{\mu
_i}^{(i)})-d(X_{\mu _{i}}^{(i+1)}).\end{equation} In particular
(\ref{3.13}),(\ref{3.22}) imply that $\psi _{\omega }(\bar
H^{(1)})=\psi _{\omega }(\bar H^{(2)})=\ldots =\psi _{\omega
}(\bar H^{(m)})=\psi _{\omega }.$ Denote  the number (\ref{3.22})
 by $\delta _i$.

Let the path (\ref{3.6}) be $\mu$-reducing, that is either $\mu
_1=\mu$ and $v_2$ does not have auxiliary edges and $\mu$ occurs
in the sequence $\mu _1,\ldots ,\mu _{m-1}$ at least twice, or
$v_2$ does have auxiliary edges $v_2\rightarrow w_1,\ldots
v_2\rightarrow w_k $ and the base $\mu$ occurs in the sequence
$\mu _1,\ldots ,\mu _{m-1}$ at least $\max _{1\leqslant i\leqslant
k}s(\Omega _{w_i})$ times. Estimate $d(U_m)=\sum
_{i=1}^{m-1}\delta _i$ from below. First notice that if $\mu
_{i_1}=\mu _{i_2}=\mu (i_1< i_2) $ and $\mu _i\not =\mu$ for $i_1<
i< i_2$, then
\begin{equation}
\label{3.23} \sum _{i=i_1}^{i_2-1}\delta _i\geqslant
d(H^{i_1+1}[1,\alpha (\Delta (\mu_{i_1+1}))]).
\end{equation}
Indeed, if $i_2=i_1+1,$ then $\delta _{i_1}=d(H^{(i_1)}[1,\alpha
(\Delta (\mu))]=d(H^{(i_1+1)}[1,\alpha (\Delta (\mu))].$ If $i_2 >
i_1+1,$ then $\mu _{i_1+1}\not = \mu$ and $\mu$ is a transfer base
in the equation $\Omega _{v_{i_1+1}}.$ Hence $\delta
_{i_1+1}+d(H^{(i_1+2)}[1,\alpha (\mu)])=d(H^{(i_1+1)}[1,\alpha
(\mu _{i_1+1})]).$ Now (\ref{3.23}) follows from
$$\sum _{i=i_1+2}^{i_2-1}\delta _i\geqslant d(H^{(i_1+2)}[1,\alpha (\mu)]).$$
So if $v_2$ does not have outgoing auxiliary edges, that is the
bases $\mu _2$ and $\Delta(\mu _2)$ do not intersect in the
equation $\Omega _{v_2}$; then (\ref{3.23}) implies that

$$\sum _{i=1}^{m-1}\delta _i\geqslant d(H^{(2)}[1,\alpha (\Delta\mu _2)])\geqslant d(
X_{\mu _2}^{(2)})\geqslant
d(X_{\mu}^{(2)})=d(X_{\mu}^{(1)})-\delta _1,$$ which implies that
\begin{equation}\label{3.24}
\sum _{i=1}^{m-1}\delta
_i\geqslant\frac{1}{2}d(X_{\mu}^{(1)}).\end{equation}

Suppose now there are outgoing auxiliary edges from the vertex
$v_2$: $v_2\rightarrow w_1,\ldots ,v_2\rightarrow w_k$. The
equation $\Omega _{v_1}$ has some solution. Let $H^{(2)}[1,\alpha
(\Delta (\mu _2))]\doteq Q$, and $P$ a  word (in the final $h$'s)
such that $Q\doteq P^{d}$, then $X_{\mu _2}^{(2)}$ and $X_{\mu
}^{(2)}$ are beginnings of the word $H^{(2)}[1,\beta (\Delta (\mu
_2))],$ which is a beginning of $P^{\infty}$. Denote $M=\max
_{1\leqslant j\leqslant k} s(\Omega _{w_j}).$

By the construction of (\ref{3.8})we either have
\begin{equation}\label{3.25}
X_{\mu}^{(2)}\doteq P^rP_1, P\doteq P_1P_2,r< M.
\end{equation}
or for each base  $\mu _i,\ i\geqslant 2,$ there is an inequality
$d(H^{(i)}(\alpha (\mu _i),\alpha (\Delta (\mu _i))))\geqslant
d(\lambda )$ and therefore
\begin{equation}\label{korova}
d(X_{\mu}^{(2)})<  M d( H^{(i)}[\alpha (\mu _i),\alpha (\Delta
(\mu _i))]).\end{equation}

Let $\mu _{i_1}=\mu _{i_2}=\mu; i_1< i_2; \mu _i\not = \mu$ for
$i_1< i< i_2.$ If
\begin{equation} \label{3.26a}
d(X_{\mu _{i_1+1}}^{(i_1+1)})\geqslant 2 d(P)\end{equation} and
$H^{(i_1+1)}[1,\rho _{i_1+1}+1]$ begins with a cyclic permutation
of $P^3$, then \newline $d(H^{(i_1+1)}[1,\alpha (\Delta (\mu
_{i_1+1}))])> d(X^{(2)}_{\mu})/M$. Together with (\ref{3.23}) this
gives $$\sum _{i=i_1}^{i_2-1}\delta _i> d(X_{\mu }^{(2)})/M.$$ The
base $\mu$ occurs in the sequence $\mu _1,\ldots ,\mu _{m-1}$ at
least $M$ times, so either (\ref{3.26a}) fails for some
$i_1\leqslant m-1$ or $\sum _{i=1}^{m-1}\delta _i\
(M-3)d(X^{(2)}_{\mu})/M.$

If (\ref{3.26a}) fails, then the inequality $d(X_{\mu
_i}^{(i+1)})\leqslant d(X_{\mu _{i+1}}^{(i+1)}) ,$ and the
definition (\ref{3.22})  imply that
$$\sum _{i=1}^{i_1}\delta _i\geqslant d(X_{\mu}^{(1)})-d(X_{\mu _{i_1+1}}^
{(i_1+1)})\geqslant (M-2)d(X_{\mu}^{(2)})/M;$$ so everything is
reduced to the second case.

Let $$\sum _{i=1}^{m-1}\delta _i\geqslant
(M-3)d(X_{\mu}^{(1)})/M.$$ Notice that (\ref{3.23}) implies for
$i_1=1$ $\sum _{i=1}^{m-1}\delta _i\geqslant d(Q)\geqslant d(P)$;
so $$\sum _{i=1}^{m-1}\delta _i\geqslant
max\{1,M-3\}d(X^{(2)}_{\mu})/M.$$ Together with (\ref{3.25}) this
implies $\sum _{i=1}^{m-1}\delta _i\geqslant
\frac{1}{5}d(X_{\mu}^{(2)})=\frac{1}{5} (d(X_{\mu}^{(1)})-\delta
_1).$ Finally,
\begin{equation}\label{3.26}
\sum _{i=1}^{m-1}\delta _i\geqslant\frac{1}{10} d(X_{\mu}^{(1)}).
\end{equation}
Comparing (\ref{3.24}) and (\ref{3.26}) we can see that for the
$\mu$-reducing path (\ref{3.6}) inequality (\ref{3.26}) always
holds.

Suppose now that the path (\ref{3.6}) is prohibited; hence it can
be represented in the form (\ref{3.7}). From  definition
(\ref{3.13}) we have $\sum _{\mu \in\omega
_1}d(X_{\mu}^{(m)})\geqslant \psi _{\omega}$; so at least for one
base $\mu\in\omega _1$ the inequality
$d(X_{\mu}^{(m)})\geqslant\frac{1}{2n}\psi _{\omega}$ holds.
Because $X_{\mu}^{(m)}\doteq (X_{\Delta (\mu)}^{(m)})^{\pm 1},$ we
can suppose that $\mu\in\omega\cup\tilde{\omega}.$ Let $m_1$ be
the length of the path $r_1s_1\cdots r_ls_l$ in (\ref{3.7}). If
$\mu\in\tilde{\omega}$ then  by the third part of the definition
of a prohibited path there exists $m_1\leqslant i\leqslant m$ such
that $\mu$ is a transfer base of $\Omega _{v_i}$. Hence, $d(X_{\mu
_i}^{(m_1)})\geqslant d(X_{\mu _i}^{(i)})\geqslant
d(X_{\mu}^{(i)})\geqslant
d(X_{\mu}^{(m)})\geqslant\frac{1}{2n}\psi _{\omega}.$ If
$\mu\in\omega$, then take $\mu$ instead of $\mu _i$. We proved the
existence of a base $\mu\in\omega$ such that
\begin{equation}\label{3.27}
d(X_{\mu }^{(m_1)})\geqslant\frac{1}{2n}\psi
_{\omega}.\end{equation}
 By the definition of a prohibited path, the inequality
$d(X_{\mu}^{(i)})\geqslant d(X_{\mu}^{(m_1)})$ $(1\leqslant
i\leqslant m_1),$ (\ref{3.26}), and (\ref{3.27}) we obtain
\begin{equation}\label{3.28}
\sum _{i=1}^{m_1-1}\delta _i\geqslant \max\{\frac{1}{20n}\psi
_{\omega},1\}(40n^2f_1+20n+1).
\end{equation}

By (\ref{3.22}) the sum in the left part of the inequality
(\ref{3.28}) equals $d_{\omega}(\bar H^{(1)})-d_{\omega}(\bar
H^{(m_1)});$ hence
$$d_{\omega}(\bar H^{(1)})\geqslant \max \{\frac{1}{20n}\psi _{\omega},1\}
(40n^2f_1 +20n +1) ,$$ which contradicts (\ref{3.19}).

This contradiction was obtained from the supposition that there
are prohibited paths (\ref{3.11}) in the path (\ref{3.8}). Hence
(\ref{3.8}) does not contain prohibited paths. This implies that
$v_i\in T_0(\Omega )$ for all $v_i$ in (\ref{3.8}).
 For all $i$ $v_i\rightarrow v_{i+1}$ is an edge of a finite tree. Hence
the path (\ref{3.8}) is finite. Let $(\Omega _{w},\bar H^{w})$ be
the final term of this sequence. We show that $(\Omega _{w},\bar
H^{w})$ satisfies all the  properties formulated in the lemma.

The first property is obvious.

Let $tp (w)=2$ and let $\Omega _w$ have non-trivial non-parametric
part. It follows from the construction of (\ref{3.8}) that if
$[j,k]$ is a non-active section for $\Omega _{v_i}$ then
$H^{(i)}[j,k]\doteq H^{(i+1)}[j,k]\doteq\ldots \doteq
H^{(w)}[j,k]$. Hence (\ref{3.9}) and the definition of $s(\Omega
_v)$ imply that the word $h_1\cdots h_{\rho _w}$ can be subdivided
into subwords $h[i_1,i_2],\ldots ,h[i_{k-1},i_k],$ such that for
any $a$ either $H^{(w)}$ has length 1, or $h[i_a,i_{a+1}]$ does
not participate in basic and coefficient equations, or
$H^{(w)}[i_a,i_{a+1}]$ can be written as
\begin{equation}
\label{3.29} H^{(w)}[i_a,i_{a+1}]\doteq P_a^rP_a';\ P_a\doteq
P_a'P_a''; r\geqslant \max_{\langle \mathcal P,R\rangle  }\ \max
\{\rho _{w}f_2(\Omega _w,P,R), f_4(\Omega _w')\},\end{equation}
where $P_a$ is a primitive word, and $\langle \mathcal P,R\rangle
$ runs
 through all the periodic
structures of ${\tilde \Omega} _w$  such that either one of them
is singular or for a  solution with maximal number of short
variables with respect to the group of extended automorphisms all
the closed sections are potentially transferable. The proof of
Proposition \ref{3.4} will be completed after we prove the
following statement.
\begin{lemma}\label{33} If $tp (w)=2$ and every closed section belonging to a
periodic structure $\mathcal P$ is potentially transferable
{\rm(}the definition is given in the construction of $T_0$ in case
15{\rm)}, one can apply the second minimal replacement and get a
finite number {\rm(}depending on periodic structures containing
this section in the vertices of type 2 in the trees $T_0(w_i),
i=1,\ldots ,m$ {\rm)} of possible generalized equations containing
the same closed sections not from $\mathcal P$ and not containing
closed sections from $\mathcal P$.\end{lemma}

\begin{proof} From the definition of a potentially transferable
section it follows that after finite number of transformations
depending on $f_4(\Omega '_u,\mathcal P)$, where $u$ runs through
the vertices of type 2 in the trees $T_0(w_i), i=1,\ldots ,m$, we
obtain a cycle that is shorter than or equal to $d(\lambda _{\rm
max}).$ This cycle is exactly $h[\alpha (\mu _i),\alpha (\Delta
(\mu _i)]$ for the base $\mu _i$ in the $\mu$-reducing subpath.
The rest of the proof of Lemma \ref{33} is a repetition of the
proof of Lemma \ref{3n}. \end{proof}

\subsection{Construction of the decomposition tree $T_{{\rm dec}}(\Omega)$ from $T_0(\Omega )$}
\label{5.5.5}

We can define now a decomposition tree $T_{{\rm dec}}(\Omega ).$
To obtain $T_{{\rm dec}}(\Omega )$ we add some edges to the
terminal vertices of type 2 of $T_0(\Omega )$. Let $v$ be a vertex
of type 2 in $T_0(\Omega )$. If there is no periodic structures in
$\Omega _v$ then this is a terminal vertex of $T_{{\rm
dec}}(\Omega )$. Suppose there exists a finite number of
combinations of different periodic structures ${\mathcal
P}_1,\ldots ,{\mathcal P}_s$ in $\Omega _v$. If some ${\mathcal
P}_i$ is singular, we consider a generalized equation $\Omega
_{u({\mathcal P}_1,\ldots ,{\mathcal P}_s)}$ obtained from $\Omega
_v({\mathcal P}_1,\ldots ,{\mathcal P}_s)$ by the first minimal
replacement corresponding to ${\mathcal P}_i$. We also draw the
edge $v\rightarrow u=u({\mathcal P}_1,\ldots ,{\mathcal P}_s)$.
This vertex $u$ is a terminal vertex of $T_{{\rm dec}}(\Omega )$.
If all ${\mathcal P}_1,\ldots ,{\mathcal P}_s$ in $\Omega _v$ are
not singular, we can suppose that for each periodic structure
${\mathcal P}_i$ with period $P_i$ some values of variables in
${\mathcal P}_i$ are shorter than $2|P_i|$ and values of some
other variables are shorter than $f_3(\Omega _v)|P_i|,$ where
$f_3$ is the function from Lemma \ref{4n}. Then we apply the
second minimal replacement. The resulting generalized equations
$\Omega _{u_1},\ldots ,\Omega _{u_t}$ will have empty
non-parametric part. We draw the edges $v\rightarrow u_1,\ldots ,
v\rightarrow u_t$ in $T_{{\rm dec}}(\Omega )$. Vertices
$u_1,\ldots ,u_t$ are terminal vertices of $T_{{\rm dec}}(\Omega
)$.

Combining this construction with Theorem \ref{th:spl} and
Proposition \ref{3.4} we obtain the following theorem.
\begin{theorem}\label{th:onelevel} Let $\Omega$ be a generalized
equation without parameters. To each branch of $T_{{\rm
dec}}(\Omega)$ with terminal vertex $w$ one can assign a splitting
{\rm(}maybe degenerate{\rm)} of $F_{R(\Omega )}$ such that every
solution $H$ of $\Omega$ corresponding to this branch can be
transformed by a canonical automorphism corresponding to this
splitting into a solution $H^+$ of $\Omega _w$.  The group
$F_{R(\Omega _w)}$ is a proper quotient of $F_{R(\Omega
)}.$\end{theorem}

\section{Structure of solutions, the solution tree
$T_{{\rm sol}}(\Omega,\Lambda)$}\label{se:5.5}

Let $\Omega = \Omega(H)$ be a generalized  equation in variables
$H$ with the set of bases $B_\Omega = B \cup \Lambda$.  Let
$T_{{\rm dec}}(\Omega)$ be  the tree constructed in Subsection
\ref{5.5.5} for a generalized equation $\Omega$ with parameters
$\Lambda$.

  Recall that in a leaf-vertex $v$  of $T_{{\rm dec}}(\Omega)$  we have the
coordinate group $F_{R(\Omega_v)}$ which is a proper homomorphic
image of $F_{R(\Omega)}$.  We define  a new  transformation $R_v$
(we call it {\it leaf-extension}) of the tree $T_{{\rm
dec}}(\Omega)$ at the leaf vertex $v$. We take the union of two
trees $T_{{\rm dec}}(\Omega)$ and $T_{{\rm dec}}(\Omega_v)$ and
identify the vertices $v$ in both trees (i.e., we extend the tree
$T_{{\rm dec}}(\Omega)$ by gluing the tree $T_{{\rm
dec}}(\Omega_v)$ to the vertex $v$).  Observe that if the equation
$\Omega_v$ has non-parametric  non-constant sections  (in this
event we call $v$ {\em a terminal vertex}), then $T_{{\rm
dec}}(\Omega_v)$ consists of a single vertex, namely $v$.

  Now we construct a solution tree $T_{{\rm sol}}(\Omega)$ by induction
starting at $T_{{\rm dec}}(\Omega)$.  Let $v$ be  a leaf
non-terminal vertex of $T^{(0)} = T_{{\rm dec}}(\Omega)$. Then we
apply the transformation $R_v$ and obtain a new tree $T^{(1)} =
R_v(T_{{\rm dec}}(\Omega))$. If there exists a leaf  non-terminal
vertex $v_1$ of $T^{(1)}$, then we apply the transformation
$R_{v_1}$, and so on. By induction we construct a strictly
increasing sequence of trees

\begin{equation}
\label{eq:5.3.1} T^{(0)} \subset  T^{(1)} \subset \ldots \subset
T^{(i)}  \subset \ldots .\end{equation} This sequence is finite.
Indeed, suppose to the contrary that the sequence is infinite and
hence  the union $T^{(\infty)}$ of this sequence is an infinite
tree in which every vertex has a finite degree. By  Konig's lemma
there is  an infinite branch $B$  in $T^{(\infty)}$. Observe that
along any infinite  branch in $T^{(\infty)}$  one has to encounter
infinitely many proper epimorphisms. This contradicts the fact
that $F$ is equationally Noetherian.

Denote the union of the sequence of the trees (\ref{eq:5.3.1}) by
$T_{{\rm sol}}(\Omega,\Lambda)$. We call $T_{{\rm
sol}}(\Omega,\Lambda)$ the {\em solution tree of $\Omega$ with
parameters $\Lambda$}.
 Recall that with every edge $e$ in $T_{{\rm dec}}(\Omega)$ (as well as in
$T_{{\rm sol}}(\Omega,\Lambda)$)
  with the initial vertex $v$ and
 the terminal  vertex $w$ we associate an epimorphism
 $$\pi_e: F_{R(\Omega_v)} \rightarrow  F_{R(\Omega_v)}.$$
 It follows that every connected (directed) path $p$ in the graph gives
rise to a
 composition of homomorphisms which we denote by $\pi_p$. Since
$T_{{\rm sol}}(\Omega,\Lambda)$
 is a tree the path $p$ is completely defined by its initial and
terminal vertices $u, v$;
 in this case we  sometimes write $\pi_{u,v}$ instead of $\pi_p$.
  Let $\pi_v$ be the homomorphism
 corresponding to the path from the initial vertex $v_0$ to a given
vertex $v$,
  we call it the {\it canonical epimorphism} from $F_{R(\Omega)}$ onto
 $F_{R(\Omega_v)}$.

 Also,
 with some vertices $v$ in the tree
$T_{{\rm dec}}(\Omega)$, as well as in the tree $T_{{\rm
sol}}(\Omega,\Lambda)$, we associate groups of canonical
automorphisms $A(\Omega _v)$ or extended automorphisms  $\bar
A(\Omega _v)$ of the coordinate group $F_{R(\Omega_v)}$ which, in
particular, fix  all variables  in the non-active part of
$\Omega_v$.  We can suppose that the group  $\bar A(\Omega _v)$ is
associated to every vertex, but for some vertices it is trivial.
Observe also, that canonical epimorphisms map parametric parts
into parametric parts (i.e., subgroups generated by variables in
parametric parts).

Recall that  writing $(\Omega,U)$ means that $U$ is a solution of
$\Omega$. If $(\Omega,U)$ and $\mu \in B_{\Omega}$, then by
$\mu_U$ we denote the element
\begin{equation}
\label{eq:5.3.2} \mu_U = [u_{\alpha(\mu)} \cdots u_{\beta(\mu)-1}
]^{\varepsilon(\mu)}. \end{equation}
 Let $B_U = \{\mu_U \mid \mu \in B\}$ and
$\Lambda_U = \{\mu_U \mid \mu \in \Lambda\}$. We refer to  these
sets as the set of values of bases from $B$ and the set of values
of parameters from $\Lambda$ with respect to the solution $U$.
Notice, that the value  $\mu_U$ is given in (\ref{eq:5.3.2}) as a
value of one fixed word mapping  $$P_\mu(H) = [h_{\alpha(\mu)}
\cdots h_{\beta(\mu)-1} ]^{\varepsilon(\mu)}. $$ In vector
notation we can write that
  $$B_U =  P_B(U), \ \ \  \Lambda_U = P_{\Lambda}(U),$$
  where $P_B(H)$ and $P_{\Lambda}(H)$ are corresponding word mappings.

  The following result explains the name of the
tree $T_{{\rm sol}}(\Omega,\Lambda)$.

\begin{theorem}
\label{th:5.3.1} Let $\Omega = \Omega(H,\Lambda)$ be a generalized
equation in variables $H$ with parameters $\Lambda .$  Let
$T_{{\rm sol}}(\Omega,\Lambda)$ be the solution tree for $\Omega$
with parameters.
  Then the following conditions hold.
\begin{enumerate}
\item For any solution $U$ in $F$  of the generalized equation
$\Omega$ there exist:  a path \newline $v_0 \rightarrow  v_1
\rightarrow \ldots \rightarrow v_n = v$ in $T_{{\rm
sol}}(\Omega,\Lambda)$ from the root vertex $v_0$ to a terminal
vertex $v$, an  abelian splitting of each of $F_{R(\Omega
_{v_0})},\ldots ,F_{R(\Omega _{v_n})},$ a sequence of canonical
automorphisms $\sigma = (\sigma_0, \ldots, \sigma_n ), \sigma_i
\in  A(\Omega _{v_i}),$ and a solution $U_v$ in $F$ of the
generalized equation $\Omega _v$ such that the solution $U$
(viewed as a homomorphism $F_{R(\Omega)} \rightarrow F$) is equal
to the following composition of homomorphisms
\begin{equation}
 \label{eq:5.3.3}
 U=\Phi_{\sigma,U_v} =
\sigma_0 \pi_{v_0,v_1}\sigma_{1} \cdots  \pi_{v_{n-1},v_n}\sigma
_n U_v.\end{equation}

\item For any path $v_0 \rightarrow  v_1 \rightarrow \ldots
\rightarrow v_n = v$ in $T_{{\rm sol}}(\Omega,\Lambda)$ from the
root vertex $v_0$ to a terminal vertex $v$, a pair consisting of a
sequence of canonical automorphisms $\sigma = (\sigma_0, \ldots,
\sigma_n ), \sigma_i \in  A(\Omega _{v_i}),$ and a solution $U_v$
 in $F$ of the generalized equation $\Omega_v$, the composition $\Phi_
{\sigma,U_v}$,  gives a solution in $F$ of the corresponding group
equation $\Omega^* = 1$; moreover, every solution  in $F$ of
$\Omega^* =1$ can be obtained this way.

 \item  For each  terminal vertex $v$ in $T_{{\rm sol}}(\Omega,\Lambda)$
there exists a word mapping $Q_v(H_v)$ such that for any solution
$U_v$ of $\Omega_v$ and any solution $U = \Phi_{\sigma,U_v}$ from
(\ref{eq:5.3.3}) the values of the parameters $\Lambda$ with
respect to $U$ can be written as $\Lambda_U = Q_v(U_v)$ (i.e.,
these values do not depend on $\sigma$)  and the word $Q_v(U_v)$
is reduced as written.

\end{enumerate}
\end{theorem}
\begin{proof} Statements (1) and (2) follow from the construction of
the tree $T_{{\rm sol}}(\Omega,\Lambda)$. To verify  (3) we need
to invoke the argument above this theorem which  claims that the
canonical automorphisms associated with generalized equations in
$T_{{\rm sol}}(\Omega,\Lambda)$ fix  all variables  in the
parametric part and, also, that  the  canonical epimorphisms map
variables from the parametric part into themselves.
\end{proof}

\begin{theorem} \label{Ase} For any finite system  $S(X)=1$ over
a free group $F$, one  can find effectively a
finite family of nondegenerate triangular quasi-quadratic systems
$U_1,\ldots, U_k$ and word mappings $p_i: V_F(U_i) \rightarrow
V_F(S)$ $(i = 1, \ldots,k)$ such that for every $b \in V_F(S)$
there exists $i$ and $c \in V_F(U_i)$ for which $b = p_i(c)$, i.e.
$$
V_F(S) = p_1(V_F(U_1)) \cup \cdots \cup p_k(V_F(U_k))
$$
and all sets $p_i(V_F(U_i))$ are irreducible; moreover, every
irreducible component of $V_F(S)$ can be obtained as a closure of
some $p_i(V_F(U_i))$ in the Zariski topology.
\end{theorem}

\begin{proof}  Each
solution of the system $S(X)=1$ can be obtained as $X=p_i(Y_i),$
where $Y_i$ are variables of $\Omega =\Omega _i$ for a finite
number of generalized equations. We have to show that all
solutions of $\Omega ^*$ are solutions of some NTQ system. We can
use Theorem \ref{th:5.3.1} without parameters and Theorem
\ref{th:spl}. In this case $\Omega _v$ is an empty equation with
non-empty set of variables. In other words $F_{R(\Omega
_v)}=F*F(h_1,\ldots ,h_{\rho}).$ To each of the branches of
$T_{{\rm sol}}$ we assign an NTQ system from the formulation of
the theorem. Let $\Omega _w$ be a leaf vertex in $T_{{\rm dec}}.$
Then $F_{R(\Omega _w)}$ is a proper quotient of $F_{R(\Omega )}$.
Consider the path $v_0, v_1, \ldots, v_n = w$ in $T_{{\rm
dec}}(\Omega)$ from the root vertex $v_0$ to a terminal vertex
$w$. All the groups $F_{R(\Omega _{v_i})}$ are isomorphic. There
are the following five possibilities.

1. $tp (v_{n-1})= 2$. In this case there is a singular periodic
structure on $\Omega _{v_{n-1}}$.  By Lemma \ref{2.10''},
$F_{R(\Omega _{v_{n-1}})}$ is a fundamental group of a graph of
groups with  one vertex group $K$,  some free abelian vertex
groups, and some edges defining HNN extensions of $K$. Recall that
making the first minimal replacement we first replaced
$F_{R(\Omega _{v_{n-1}})}$  by a finite number of proper quotients
in which the edge groups corresponding to abelian vertex groups
and HNN extensions are maximal cyclic in $K$. Extend the
centralizers of the edge groups of $\Omega _{v_{n-1}}$
corresponding to HNN extensions by stable letters $t_1,\ldots
,t_k$. This new group that we denote by $N$ is the coordinate
group of a quadratic equation  over $F_{R(\Omega _w)}$ which has a
solution in $F_{R(\Omega _w)}$.

In all the other cases  $tp (v_{n-1})\neq 2$.

2. There were no auxiliary edges from vertices $v_0, v_1, \dots,
v_n = w$, and if one of the Cases 7--10 appeared at one of these
vertices, then it only appeared a bounded (the boundary depends on
$\Omega _{v_0}$) number of times in the sequence. In this case we
replace $F_{R(\Omega )}$ by $F_{R(\Omega _w)}$

3. $F_{R(\Omega _w)}$ is a term in a free decomposition of
$F_{R(\Omega _{v_{n-1}})}$ ($\Omega _w$ is a kernel of a
generalized equation $\Omega _{v_{n-1}}$). In this case we also
consider $F_{R(\Omega _w)}$ instead of $F_{R(\Omega)}$.

4. For some $i$ $tp(v_i)=12$ and the path $v_i,\ldots ,v_n=w$ does
not contain vertices of type $ 7-10, 12$ or $15$. In this case
$F_{R(\Omega)}$ is the coordinate group of a quadratic equation.

5. The path $v_0, v_1, \dots, v_n = w$ contains  vertices of type
15.  Suppose\linebreak $v_{i_j},\ldots ,v_{i_j+k_j},$ $j=1,\ldots
,l$ are all blocks of  consecutive vertices of type 15 This means
that $tp(v_{i_j+k_j+1})\neq 15$ and $i_j+k_j+1<i_{j+1}$. Suppose
also that none of the previous cases takes place. To each
$v_{i_j}$ we assigned a quadratic equation and a group of
canonical automorphisms corresponding to this equation. Going
along the path $v_{i_j},\ldots ,v_{i_j+k_j},$ we take minimal
solutions corresponding to some non-singular periodic structures.
Each such structure corresponds to a representation of
$F_{R(\Omega _{v_{i_j}})}$ as an HNN extension. As in the case of
a singular periodic structure, we can suppose that the edge groups
corresponding to  HNN extensions are maximal cyclic and not
conjugated in $K$. Extend the centralizers of the edge groups
corresponding to HNN extensions by stable letters $t_1,\ldots
,t_k$.  Let $N$ be the new group. Then $N$ is the coordinate group
of a quadratic system of equations over $F_{R(\Omega
_{v_{i_j+k_j+1}})}$. Repeating this construction for each
$j=1,\ldots ,l$, we construct NTQ system over $F_{R(\Omega _w)}$.

Since $F_{R(\Omega _w)}$ is a proper quotient of $F_{R(\Omega )}$,
the theorem can  now be proved by induction.

\end{proof}

\section{Maximal standard quotients and canonical embeddings of $\mathcal F$-groups}
\label{se:canonical}

In this section we construct  standard quotients and canonical
embeddings of freely indecomposable   finitely generated fully
residually free groups into NTQ groups.

\subsection{JSJ-decompositions are sufficient splittings}
\label{se:max-standard}

In this section we prove that a non-degenerate JSJ decomposition
$D$  of a finitely generated freely indecomposable fully
residually free group $G$ is a sufficient splitting of $G$, i.e.,
the standard maximal quotient $G/R_D$ is a proper quotient. Let,
as usual,  $F$ be a fixed free non-abelian  group.

\begin{theorem}\label{th:suff} Let $G$ be a finitely generated freely
indecomposable {\rm[}freely indecomposable modulo $F${\rm]}
 fully residually free group $G$. Then a non-degenerate JSJ
 $\mathbb Z$-
 decomposition {\rm[}modulo $F${\rm]} of $G$  is a sufficient splitting of $G$.
\end{theorem}

\begin{proof}  \label{3.2'}
Let $D$ be a non-degenerate $JSJ$ decomposition of the group
 $G=F_{R(U)}$.  Recall, that by ${\mathcal GE}(U)$ we denote the
 finite set of the initial generalized equations associated with $U = 1$
  (see Section \ref{se:4-2}). For each $\Omega \in {\mathcal GE}(U)$
  we are going to construct a tree $T_{mq}(\Omega)$ as follows.
 Let $T_{{\rm dec}}(\Omega)$ be the decomposition tree which has
 been  constructed  in Subsection \ref{5.5.5}. Suppose $v$ is a leaf vertex of
 $T_{{\rm dec}}(\Omega)$ and  $\Omega _v$ is the generalized
equation associated with $v$ in  $T_{{\rm dec}}(\Omega).$ By the
construction the group $F_{R(\Omega _v)}$ is a proper quotient of
$F_{R(\Omega )}$. If the image of $G$ in $F_{R(\Omega _v)}$ under
the quotient epimorphism $F_{R(\Omega)} \rightarrow F_{R(\Omega
_v)}$ is not a proper quotient of $G$, then we extend the tree
$T_{{\rm dec}}(\Omega)$ by gluing up the tree $T_{{\rm
dec}}(\Omega _v)$ to the vertex $v$ (by identifying the vertex $v$
in $T_{{\rm dec}}(\Omega)$ with the root vertex in $T_{{\rm
dec}}(\Omega _v)$). Denote the resulting tree by $T_{{\rm
dec}}(\Omega)^\prime$. Again, if there is a leaf vertex $v^\prime
\in T_{{\rm dec}}(\Omega)^\prime$ such that the canonical image of
$G$ in the group associated with $v^\prime$ is not a proper
quotient of $G$ then we extend $T_{{\rm dec}}(\Omega)^\prime$ by
gluing up the tree $T_{{\rm dec}}(\Omega _v^\prime)$  to the
vertex $v^\prime$ (we refer to this operation as to {\em leaf
extension}). More formally, let $T_0 = T_{{\rm dec}}(\Omega)$.
Suppose a tree $T_i$ with associated groups and homomorphisms is
already constructed. If there is a vertex $w \in T_i$ such that
the image of $G$ in the group associated with $w$ is not a proper
quotient of $G$ then we glue up the tree $T_{{\rm dec}}(\Omega
_w)$  to the vertex $w$ (by identifying $w$ with the root vertex
in $T_{{\rm dec}}(\Omega _w)$) and denote the resulting tree by
$T_{i+1}$.

\medskip
{\bf Claim. } {\em The sequence
 \begin{equation} \label{eq:Ti}
 T_0 \subset T_1 \subset T_2 \subset \ldots
 \end{equation}
 is finite.}

 \medskip
 Indeed, suppose the sequence (\ref{eq:Ti}) is infinite. Then the
 graph $T = \bigcup_{i=1}^{\infty} T_i$ is infinite. Then by Konig's lemma
 there is an infinite branch in $T$, say,
  \begin{equation}
  \label{eq:brunch}
  v_0 \rightarrow v_1 \rightarrow v_2 \rightarrow \ldots .
  \end{equation}

Therefore there exists an infinite subsequence, say,
 $$v_{i_0} \rightarrow v_{i_1} \rightarrow v_{i_2} \rightarrow \ldots $$
such that for every vertex $v_{i_j}$ in this sequence at some step
the tree $T_{{\rm dec}}(\Omega _{v_{i_j}})$ was glued to the
vertex $v_{i_j}$. It follows that the homomorphism corresponding
to the edge $v_{i_j -1} \rightarrow v_{i_j}$ is a proper
epimorphism. Hence there are infinitely many proper epimorphisms
associated with the edges of the branch (\ref{eq:brunch}). This
implies that there exists an infinite sequence of epimorphisms of
finitely generated residually free groups which contains
infinitely many proper epimorphisms - contradiction with the
equationally Noetherian property of free groups. This proves the
claim.

 Denote the union tree $T$ of the finite sequence (\ref{eq:Ti}) by $T_{mq}(\Omega )$.
Fix some terminal vertex $w$ of $ T_{mq}(U,\Omega)$.
Each homomorphism from the  image of $G$ in
$F_{R(\Omega _w)}$ into $F$ factors through a finite number of
fully residually free quotients. Denote one such quotient by $\bar
G=F_{R(\bar U)}*F(t_1,\ldots ,t_{k})$, where $\bar U$ is
irreducible, and $t_1,\ldots ,t_k$ are free variables from
quadratic equations.

Let $v_0$ be the initial vertex of $ T_{{\rm dec}}(\Omega )$,
$v_1$ be the terminal vertex of $T_{{\rm dec}}(\Omega)$, $v_2$ be
the terminal vertex of $T_{{\rm dec}}(\Omega _{v_1})$ and so on,
$v_n=w$. Each solution of $U=1$ corresponding to this branch of $
T_{mq}(U,\Omega)$ can be transformed by  canonical automorphisms
corresponding to vertices of $T_{{\rm dec}}(\Omega)$  (which are
canonical automorphisms from $A_{D_0}$, where $D_0$ is a
decomposition of $F_{R(\Omega)}$, see Subsection \ref{5.5.2}) and,
subsequently, by canonical automorphisms corresponding to vertices
of $T_{{\rm dec}}(\Omega _{v_i}), i=1,\ldots ,n$ (which are
canonical automorphisms from $A_{D_{i}}$, where $D_i$ is a
decomposition  of $F_{R(\Omega _{v_i})})$, into a solution
satisfying a proper equation.


By Lemma \ref{le:1.6}, for each QH subgroup $Q$ in $F_{R(\Omega
_{v_i})}$, either the intersection of $G$ with some conjugate
$Q^g$ is of finite index in $Q^g$ and is a QH subgroup of $G$,
 or $G$ is conjugated into the fundamental
group of a connected component of the graph of groups obtained
from $D$ by removing the vertex and edges corresponding to $Q$.
Case 1) in the lemma cannot take place because $G$ is freely
indecomposable. In the case when $G$ is conjugated into the
fundamental group of the graph of groups obtained from $D$ by
removing the vertex and edges corresponding to $Q$,  each solution
of $U=1$ can be obtained from a solution of $\Omega _{v_i}$
minimal with respect to the canonical group of automorphisms
corresponding to $Q.$ If $G$ has trivial intersection with all
conjugates of some edge group of $D_i$, then again each solution
of $U=1$ can be obtained from a solution of $\Omega _{v_i}$
minimal with respect to the group generated by the canonical Dehn
twists along this edge. If the intersection $G\cap Q^g$ has a
finite index, then it is a $QH$ subgroup $Q_1$ of $G$ and the
group of automorphisms of $Q_1$ induced by canonical automorphisms
of $Q^g$ has finite index in the group of canonical automorphisms
of $Q^g$. Let $\bar F_{R(\Omega _{v_i})}$ be the factor that
contains $G$ in some Grushko decomposition of $F_{R(\Omega
_{v_i})}.$ Let $A(\bar F_{R(\Omega _{v_i})})$ be the group of
outer canonical automorphisms of $\bar F_{R(\Omega _{v_i})}$
induced from canonical automorphisms of $F_{R(\Omega _{v_i})}$. By
Lemma \ref{le:1.6}, the group of canonical automorphisms of $\bar
F_{R(\Omega _{v_i})}$ corresponding to $Q^g$ (where $Q$ is a QH
subgroup of $\bar F_{R(\Omega _{v_i})}$) that induce canonical
automorphisms of $G$ has finite index in the subgroup of $A(\bar
F_{R(\Omega _{v_i})})$ corresponding to $Q$. Denote by $A_G(\bar
F_{R(\Omega _{v_i})})$ the subgroup of $A(\bar
 F_{R(\Omega _{v_i})})$  generated by canonical Dehn twists corresponding to
the splittings of $\bar F_{R(\Omega _{v_i})}$ that induce
non-trivial splittings of $G$. The subgroup $A_{ind}$ of $A_G(\bar
F_{R(\Omega _{v_i})})$ consisting of automorphisms  that induce
automorphisms of $G$, has finite index in $A_G(\bar F_{R(\Omega
_{v_i})}).$

This implies that using canonical automorphisms corresponding to
the JSJ decomposition of $G$, we can transform each solution
$\phi$ of $U=1$ into a  solution that is not longer than $\phi$
and satisfies one of the finite number of proper equations (which
can be effectively found) or, equivalently, satisfies a proper
equation which represents this disjunction of equations. Since a
minimal solution in the equivalence class of $\phi$ factors
through some branch of $T_{mq}(\Omega )$ for one of the finite
number of generalized equations $\Omega _w$, it must satisfy this
proper equation.

\end{proof}

\subsection{Quasi-quadratic systems associated with $\mathbb{Z}$-splittings}

Let $G$ be a finitely generated freely indecomposable fully
residually free group and $D$ a  non-degenerate JSJ
$\mathbb{Z}$-decomposition of $G$. Combining foldings and
slidings, we can transform $D$ into an abelian decomposition in
which each vertex with non-cyclic abelian subgroup that is
connected to some rigid vertex, is connected to only one vertex
which is rigid. We suppose from the beginning that $D$ has this
property. Let $L_i, i \in I$ be a (finite) collection of all
standard maximal fully residually free quotients of $G$ relative
to $D$ and $\pi_i :G\rightarrow L_i$ the corresponding canonical
epimorphisms. By Corollary \ref{cy:maxstq} one can enumerate these
quotients as
 $$L_1, \ldots,  L_k,  \dots,  L_{k+s}, \ \ k \geqslant 1, s \geqslant 0,$$
in such a way that  for each  $i = 1, \ldots, k$ all the
restrictions of $\pi_i$ onto  rigid subgroups of $D$, onto edge
subgroups in $D$, and onto the subgroups of the abelian vertex
groups $A$ generated by the images of all the edge groups of edges
adjacent to $A$, are monomorphisms. Let $H_1, \ldots, H_q$ be
rigid subgroups of $G$ in $D$. Consider  the Grushko's
decomposition of $L_i$ compatible with the collection of subgroups
$\pi_i(H_1), \ldots, \pi_i(H_q)$:
 $$L_i = L_{i,1} \ast \cdots \ast L_{i,\alpha_i} \ast F(t_1, \dots,
 t_{\beta_i}),$$
where $F(t_1, \ldots, t_{\beta_i})$ is the maximal free factor
that does not contain conjugates of any of the subgroups
$\pi_i(H_1), \dots, \pi_i(H_q)$.
 Denote
 $${\bar G_i} = L_{i,1} \ast \cdots \ast L_{i,\alpha_i}.$$

 In this section we show that
$G$ embeds into a  group $G_i$ which is a particular extension of
the group ${\bar G_i}$, $i = 1, \ldots, k$. More precisely, $G_i$
is the fundamental group of a graph of groups
 $\Gamma_i$ obtained from a single vertex $v$ with the associated
 vertex group $G_v = {\bar G_i}$ by
 adding  finitely many edges corresponding to  extensions of centralizers (viewed as
 amalgamated products) and finitely many QH vertices connected only to $v$ (see Figure
 15).
\begin{figure}[here]
\centering{\mbox{\psfig{figure=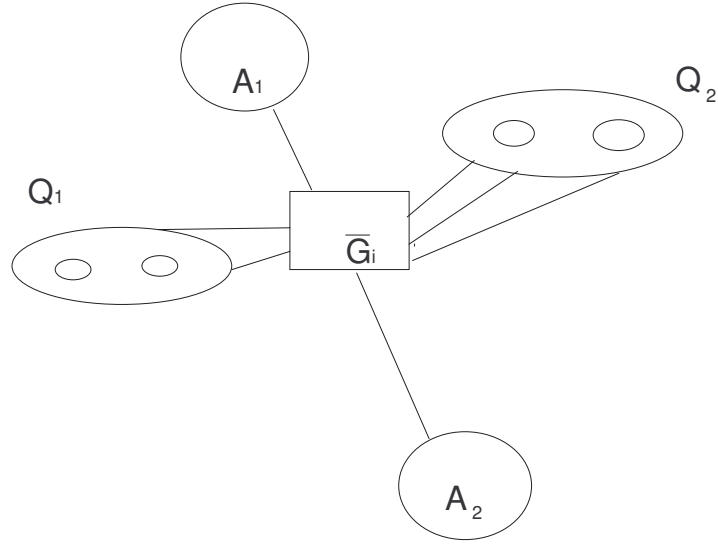,height=3in}}} \caption{ The
graph of groups $\Gamma_i$ of the coordinate group of the
quasi-quadratic system over $\bar G_i$.} \label{NTQ}
\end{figure}

 We start with description of the extensions $G_i$ of ${\bar G_i}$.
  Fix an integer $i$, $i=1,\ldots ,k$ and, to simplify notation,
denote ${\bar G_i}$ by ${\bar G}$, $\pi_i$ by $\pi$, and $G_i$ by
$G^\ast$. Let $\bar G=P_1\ast\cdots\ast P_{\alpha}$ be the Grushko
decomposition of $\bar G$. Then by  construction of $\bar G$, each
factor in this decomposition contains a conjugate of the image of
some rigid subgroup or an edge group in $D$. Let $g_1,\ldots,g_l$
be a fixed finite generating set of ${\bar G}$. For an edge $e \in
D$ we fix  a generator  $d_e$ of the cyclic edge group $G_e$ (or
generators of an abelian edge group connecting a non-cyclic
abelian vertex group to a rigid subgroup).
 The required
extension $G^\ast$ of ${\bar G}$ is constructed  in three steps.
On each step we extend the centralizers $C_{\bar G}(\pi(d_e))$ of
some edges $e$ in $D$ or add a QH subgroup. Simultaneously, for
every edge $e \in D$ we associate an element $s_e \in
C_{G^\ast}(\pi(d_e))$.

Step 1. Let $E_{rig}$ be the set of all edges between rigid
subgroups in $D$.   For each edge $e\in E_{rig}$ we do the
following. If
 $$
 G=A\ast _{\langle d_e\rangle}B  \  \mbox{and} \; {\bar G}= \pi (A)\ast _{\langle
\pi(d_e)\rangle}\pi (B)
$$
 or
$$
G=A\ast _{\langle d_e\rangle} \ \mbox{and} \; {\bar G}= \pi(A)\ast
_{\langle \pi(d_e)\rangle}
 $$
 then we delete the edge $e$ from $E_{rig}$.  In this case we associate the trivial element $s_e = 1$
  with $e$. Denote the resulting set of edges by $E^\prime$.
 One can define an equivalence relation $\sim$ on $E^\prime$ assuming for
$e, f \in E^\prime$ that
 $$e \sim f \Longleftrightarrow \left(\exists g_{ef} \in {\bar G}\right) \left ( g_{ef}^{-1}C_{\bar G}(\pi(e))g_{ef} =
 C_{\bar G}(\pi(f)) \right).
  $$
Let $E$ be a set of representatives of equivalence classes of
$E^\prime$ modulo $\sim$.

Now we construct a group $G^{(1)}$  by extending every centralizer
$C_{\bar G}(\pi(d_e))$ of ${\bar G}$, $e \in E$ as follows. Let
$$[e] = \{e = e_1, \ldots, e_{q_e}\}$$
and $y_e^{(1)}, \ldots, y_e^{(q_e)}$  be new letters corresponding
to the elements in $[e]$. Then put
 $$G^{(1)} = \left< {\bar G}, y_e^{(1)},
\ldots, y_e^{(q_e)} (e \in E) \mid [C(\pi(d_e)),y_e^{(j)}]
 = 1, [y_e^{(i)},y_e^{(j)}] = 1 (i,j = 1, \dots, q_e)\right>.$$

One can  associate with $G^{(1)}$  the following system of
equations over $\bar G$:

\begin{equation}\label{eq:edges}
[{\bar g_{es}},y_e^{(j)}]=1,\ [y_e^{(i)},y_e^{(j)}] = 1, \ \ i,j =
1, \ldots, q_e, \ s = 1, \ldots, p_e, \ e\in E,
\end{equation}
 where $y_e^{(j)}$ are new variables and the elements ${\bar g_{e1}}, \ldots, {\bar g_{ep_e}}$
 are constants from ${\bar G}$  which generate  the centralizer $C(\pi(d_e))$.
 We assume that the constants ${\bar g_{ej}}$
 are given as words in the generators $g_1, \ldots, g_l$ of ${\bar
 G}$. We associate the element $s_{e_i}=y_e^{(i)}$ with the edge $e=e_i$.

Step 2. Let $A$ be a non-cyclic abelian vertex group in $D$ and
$A_e$ the subgroup of $A$ generated by the images in $A$ of the
edge groups of edges adjacent to $A$. Then $A = Is(A_e) \times
A_0$ where $Is(A_e)$ is the isolator of $A_e$ in $A$ (the minimal
direct factor containing $A_e$) and $A_0$ a  direct complement of
$Is(A_e)$ in $A$. Notice, that the restriction of $\pi$ on
$Is(A_e)$ is a monomorphism (since $\pi$ is injective on $A_e$ and
$A_e$ is of finite index in $Is(A_e)$). For each non-cyclic
abelian vertex group $A$ in $D$ we extend the centralizer of
$\pi(Is(A_e))$ in $G^{(1)}$ by the abelian group $A_0$ and denote
the resulting group by $G^{(2)}$. Observe, that since
$\pi(Is(A_e)) \leqslant {\bar G}$ the group $G^{(2)}$ is obtained
from ${\bar G}$ by extending finitely many centralizers of
elements from ${\bar G}$.

 If the  abelian group $A_0$ has rank
$r$ then the system of equations associated with the abelian
vertex group $A$ has the following form
\begin{equation}\label{eq:barS} [y_p, y_q]=1, [y_p, {\bar d_{ej}}]=1, \ \
 \  p,q=1,\ldots ,r, j = 1, \ldots, p_e,\end{equation}
 where $y_p, y_q$ are new variables and the elements ${\bar d_{e1}}, \ldots, {\bar
d_{ep_e}}$
 are constants from ${\bar G}$  which generate  the subgroup $\pi(Is(A_e))$.
 We assume that the constants ${\bar d_{ej}}$
 are given as words in the generators $g_1, \ldots, g_l$ of ${\bar
 G}$.

 Denote by
   $$
   \bar S_i(y_1,\ldots ,y_t,g_1, \ldots, g_l) =1
   $$
the   union of the system (\ref{eq:edges}) and all the systems
(\ref{eq:barS})  for every abelian non-cyclic vertex group $A$ in
$D$. Here $y_1,\ldots ,y_t$ are all the new  variables that occur
in the equations (\ref{eq:edges}), (\ref{eq:barS}).

Step 3.  A QH subgroup $Q$ such that $\pi (Q)$ is a QH subgroup of
$\bar G$ of the same size (see Section 1.6) is called a {\em
stable QH subgroup}. Let $Q$ be a non-stable QH subgroup in $D$.
Suppose $Q$ is given by a presentation
$$\prod _{i=1}^{n}[x_i, y_i]p_1\cdots p_{m}=1.$$
 where  there are exactly $m$
outgoing edges $e_1,\ldots ,e_m$ from $Q$  and $\sigma(G_{e_i}) =
\langle p_i\rangle $,  $\tau(G_{e_i}) = \langle c_i\rangle $ for
each edge $e_i$. We add a QH vertex $Q$ to $G^{(2)}$ by
introducing new generators and the following quadratic relation

\begin{equation}\label{QQQC}\prod _{i=1}^{n}[ x_i,
y_i](c_1^{\pi_i})^{z_1}\cdots
(c_{m-1}^{\pi_i})^{z_{m-1}}c_m^{\pi_i}=1\end{equation}
 to the presentation of $G^{(2)}$.
Observe, that in the relations (\ref{QQQC}) the coefficients in
the original quadratic relations for $Q$ in $D$ are
 replaced by their images in $\bar G$.

Similarly,  one introduces QH vertices for non-orientable QH
subgroups in $D$. We refer to Section \ref{QHquad} for details on
quadratic relations associated with non-orientable QH subgroups.

The resulting group is denoted by $G^\ast = G^{(3)}$.

 Let
 $$
 S=S(z_1,\ldots ,z_m,g_1,\ldots,g_l)=1
  $$
  be union of the quadratic systems of the type (\ref{QQQC})
  for orientable and non-orientable
non-stable QH subgroups in $D$.   Here the set of variables
$\{z_1,\ldots ,z_m\}$ is  union of all variables in quadratic
equations of the type (\ref{QQQC}). Let $D^*$ be a splitting of
$G^*$ with the following vertex groups: factors $P_1,\ldots
,P_{\alpha}$ in the compatible reduced free decomposition of $\bar
G$, abelian vertex groups corresponding to relations
(\ref{eq:edges}), (\ref{eq:barS}), and QH subgroups added to $\bar
G$. Let $T^*$ be a maximal subtree in the underlying graph $\Gamma
^*$. So far we defined elements $s_e$ only for edges between two
rigid subgroups of $D$.

Now we define an $F$-homomorphism $\psi:G \rightarrow G^\ast$ and
show that $\psi$ is injective. In the process we will associate
elements $s_e$ to the remaining edges of $D$. We need the
following lemma.

 \begin{lemma}
 \label{le:extend-hom-amalgam}\
 \bi
 \item[(1)] Let $H = A\ast_{\langle d\rangle}B$ and $\pi:H \rightarrow {\bar
 H}$ be a homomorphism such that the restrictions of $\pi$ on $A$
 and $B$ are injective. Put
  $$H^\ast = \langle {\bar H}, y \mid [C_{{\bar H}}(\pi(d)),y] = 1
  \rangle.$$
   Then for every $u \in C_{H^\ast}((\pi(d))$,  $u \not \in C_{{\bar
   H}}(\pi(d))$, a  map
    $$\psi(x) = \left \{ \begin{array}{lr}  \pi(x), & \  x \in A,
    \\
                                    \pi(x)^u, & \ x \in B.
                         \end{array}
                         \right. $$
    gives rise to a monomorphism $\psi:H \rightarrow H^\ast$.

    \item[(2)] Let $H = \langle A, t \mid d^t = c\rangle$  and $\pi:H \rightarrow {\bar
 H}$ be a homomorphism such that the restriction of $\pi$ on $A$
 is injective. Put
  $$H^\ast = \langle {\bar H}, y \mid [C_{{\bar H}}(\pi(d)),y] = 1
  \rangle.$$
   Then for every $u \in C_{H^\ast}((\pi(d))$,  $u \not \in C_{{\bar
   H}}(\pi(d))$, a map
    $$\psi(x) = \left \{ \begin{array}{lr}  \pi(x), & \  x \in A,
    \\
                                    u\pi(x), & \ x = t.
                         \end{array}
                         \right. $$
    gives rise to a monomorphism $\psi:H \rightarrow H^\ast$.
    \ei
\end{lemma}

\begin{proof}
We prove (1), a similar argument proves  (2). Clearly, the map
$\psi$ agrees on the intersection $A \cap B$, so it defines a
homomorphism $\psi:H_{rig} \rightarrow G$. Straightforward
verification shows that $\psi$ maps reduced forms of elements in
$H_{rig}$ into reduced forms of elements in $G$. This  proves that
$\psi$ is injective.
\end{proof}

We  define the homomorphism $\psi:G \rightarrow G^\ast$ with
respect to the splitting $D$ of $G$. Let $T$ be the maximal
subtree of $D$. First, we define $\psi$ on the fundamental group
of the graph of groups induced from $D$ on $T$.
 Notice that the subgroup  $F$ is elliptic in $D$, so there is a
rigid vertex $v_0 \in T$ such that $F \leqslant G_{v_0}$. Mapping
$\pi$ embeds $G_{v_0}$ into ${\bar G}$, hence into $G^\ast$.

 Let $P$ be a path $v_0 \rightarrow v_1 \rightarrow \ldots
\rightarrow v_n$ in $T$ that starts at $v_0$. With each edge $e_i
= (v_{i-1} \rightarrow v_i)$ between two rigid vertex groups  we
have already associated the  element $s_{e_i}$. Let us associate
elements to other edges of $P$:

a) if $v_{i-1}$ is a rigid vertex, and $v_i$ is either abelian or
QH, then $s_{e_i}=1$;

b) if $v_{i-1}$ is a QH vetrex, $v_i$ is rigid or abelian, and the
image of $e_i$ in the decomposition  $D^*$ of $G^*$ does not
belong to $T^*$, then $s_{e_i}$ is the stable letter corresponding
to the image of $e_i$;

c) if $v_{i-1}$ is a QH vertex  and $v_i$ is rigid or abelian, and
the image of $e_i$ in the decomposition of $G^*$ belongs to $T^*$,
then $s_{e_i}=1$.

d) if $v_{i-1}$ is an abelian vertex with $G_{v_{i-1}}=A$   and
$v_i$ is a QH vertex, then $s_{e_i}$ is an element from $A$ that
belongs to $A_0$.

Since two abelian vertices  cannot be connected by an edge in
$\Gamma$, and we can suppose that  two QH vertices are not
connected by an edge, these are all possible cases.

We now define the embedding $\psi$ on the fundamental group
corresponding to the path $P$ as follows:
 $$\psi(x) = \pi(x)^{s_{e_i}\cdots s_{e_1}}\ \mbox{for } \ x\in G_{v_i}.$$
This map is a monomorphism by Lemma \ref{le:extend-hom-amalgam}.
Similarly  we define $\psi$ on the fundamental group of the graph
of groups induced from $D$ on $T$. We  extend it to $G$ using the
second statement of  Lemma \ref{le:extend-hom-amalgam}.

We  constructed a finite number of  proper fully residually free
quotients of $G$: $\bar G,\ldots ,\bar G_k$ where the systems
$S=1\ \wedge\ \bar S_i=1$ have solutions. Let  $\bar G_i=F_{R(\bar
U_i)}.$  Now $U_{D,i}(z_1,\dots ,z_m,y_1,\dots ,y_t,g_1,\dots
,g_l)=1$ is the system
 $$
\begin{array}{rl}
S(z_1,\ldots ,z_m,g_{1},\ldots,g_{l}) & = 1 \\
\bar S_i(y_1,\ldots ,y_t,g_{1},\ldots,g_{l}) &  = 1\\
 \bar U_i(g_{1},\ldots ,g_l) &  = 1
\end{array}
$$

 If $G=F_{R(U)}$ is a surface group, we take the system $S=1$  the same as
$U=1$, $\bar S_i=1$ and $\bar U_i=1$  trivial. If $G$ is a free
abelian group, we take $\bar S=1$ the same as $U=1$, and $S=1,\bar
U=1$ trivial. If $G$ is cyclic, we do not construct $ U_{D,i}=1$.

We summarize the construction in the following theorem.

\begin{theorem}  \label{qq} Let $U(X) =1$, $X = \{x_1, \ldots,
x_n\},$ be a finite irreducible system of equations over $F$ such
that $G=F_{R(U)}$ is freely indecomposable {\rm (}modulo $F${\rm
)} if
 $U(X) = 1$ has non-trivial
coefficients. Let $D$ be a JSJ $\mathbb Z$-decomposition of $G$.
Then one can effectively construct finitely many equations
$$V_1(X) = 1, \ldots, V_s(X) = 1$$
with coefficients in $F$  and finitely many systems
 $$U_{D,1} = 1, \ldots, U_{D,k} = 1$$
   over $F$ with  embeddings
 $$\phi_1: F_{R(U)} \rightarrow F_{R(U_{D,1})}, \ \ldots, \ \phi_k: F_{R(U)} \rightarrow
 F_{R(U_{D,k})} $$
  such that:
 \begin{enumerate}
 \item  $V_i \not \in R(U), i = 1, \ldots, s.$
 \item $U_{D,i}(z_1,\ldots ,z_m,y_1,\ldots ,y_t,
 g_1,\ldots ,g_l)=1$ is a
system of the following type
 $$
\begin{array}{rr}
S(z_1,\ldots ,z_m,g_1,\ldots ,
g_{l}) & = 1 \\
\bar S_i(y_1,\ldots ,y_t,g_{1},\ldots ,g_{l}) &  = 1\\
 \bar U_{i}(g_{1},\ldots ,g_{l}) &  = 1
\end{array}
$$
where $\bar U_{i}=1$ is irreducible,  $S=1$ is a non-degenerate
quadratic system in variables $z_1,\ldots ,z_m$ over $\bar
G_i=F_{R(\bar U_i)}$, variables $z_1,\ldots ,z_m$ are  standard
generators of $MQH$ subgroups of $G$ and stable letters
corresponding to non-tree edges adjacent to QH subgroups, $\bar
S_i=1$ corresponds to extensions of centralizers of elements from
$F_{R(\bar U_i)}.$ There is a splitting ${D_i}^*$ of the
fundamental group of the system $U_{D,i}$ modulo free factors of
$\bar G_i$; MQH vertex groups in ${D_i}^*$ correspond to equations
from the system $S=1$, abelian vertex groups correspond to
equations from the system $\bar S_i=1$.
  \item
$F_{R(\bar U_i)}\ast F_(t_1,\ldots ,t_{\beta _i})$ is a proper
quotient of $F_{R(U)}$ with the canonical  epimorphism   $\pi _i:
F_{R(U)}\rightarrow F_{R(\bar U_i)}\ast F_(t_1,\ldots ,t_{\beta
_i})$, $i = 1, \ldots, k$, which can be found effectively.
  \item
$V(U)=p_1(V(U_{D,1}))\cup\cdots\cup p_k(V(U_{D,k}))\cup V(U\
\wedge\ V_1)\cup  \cdots \cup V(U\ \wedge\ V_s),$
 where the word
mappings $p_i= (p_{i,1}, \ldots, p_{i,n}), i = 1, \ldots, k$,
correspond to the embeddings $\phi _1,\ldots ,\phi _k$, i.e., for
every $i$ and every $x_j \in X$
  $$
  \phi_i(x_j) =p_{i,j}(z_1,\ldots ,z_m,y_1,\ldots ,y_t,g_{1},\ldots
,g_{l}). $$
  \item
   Let $\Phi _i=\{\phi_i \sigma\pi_i \psi \mid \sigma
\in A_{{D_i}^*}, \psi \in V({\bar U}_i)\}.$
   Then $\Phi _i \subseteq  p_i(V(U_{D,i}))$ and
  $$V(U)=\Phi _1\cup\cdots\cup \Phi _k\cup V(U\ \wedge\ V_1)\cup  \cdots \cup V(U\ \wedge\
  V_s).$$\end{enumerate}\end{theorem}

\begin{proof} We have already constructed the systems
$$U_{D,1}=1,\ldots ,U_{D,k}=1$$ and embeddings
 $$\phi_1: F_{R(U)} \rightarrow F_{R(U_{D,1})}, \ \dots, \ \phi_k: F_{R(U)} \rightarrow
 F_{R(U_{D,k})}. $$

Recall that $L_1,\dots ,L_{k+s}$ is a family of all standard
maximal fully residually free quotients for $G$ such  that for
each  $i = 1, \ldots, k$ all the restrictions of $\pi_i$ onto
rigid subgroups of $D$, onto edge subgroups in $D$, and onto the
subgroups of the abelian vertex groups $A$ generated by the images
of all the edge groups of edges adjacent to $A$, are
monomorphisms. Consider now fully residually free quotients
$L_{k+1},\ldots , L_{k+s}.$ For each quotient $L_{k+i}$ there
exists an element $V_i$ from a rigid subgroup or from $A_e$ for
some abelian vertex group $A$ such that $\pi_i(V_i)=1$. Since all
automorphisms from $A_D$ only conjugate rigid subgroups, all
solutions of $U=1$ that can be represented as a composition
$\sigma\pi _i$, where $\sigma\in A_D$, satisfy the equation
$V_i=1$.

In order to prove statements (4) and (5) we just notice that the
set $\phi _i\sigma,$ where $\sigma\in A_{D^*}$ is the set of
canonical automorphisms from $A_D$ except the automorphisms
corresponding to stable QH subgroups and edges between rigid
subgroups for which we did not extend their centralizers. Using
these automorphisms one  can transform an arbitrary solution into
a solution in one of the maximal fully residually free quotients.
\end{proof}

\section{Effective free decompositions}\label{efd}

In this section we prove the following result, which is one of the
corner stones in constructing JSJ decompositions affectively.

\begin{theorem}
\label{th:1.7}  There is an algorithm which for every  finitely
generated fully residually free group $G$ and its finitely
generated subgroup $H$  determines whether or not $G$ has a
nontrivial  free decomposition modulo $H$. Moreover, if $G$ does
have such a decomposition, the algorithm finds one (by giving
finite generating sets of the factors).
\end{theorem}

\begin{proof}

 Let $G$ be a finitely generated fully residually free group, $G=\langle Z\mid U(Z)\rangle$,
  where $U(Z) = 1$ is an irreducible system over $F$. Let $H$
  be a finitely generated subgroup of $G$ given by a finite set of generators $H =
  \langle h_1(Z), \ldots, h_n(Z) \rangle$, where $h_i(Z)$ is a word in generators
  $Z$. For every $i = 1,
\ldots, n$ add the equation $y_i = h_i(X)$ to the system $U(Z) = 1$. We denote
the resulting system  by $U_H(Z,Y) = 1$.

   Let
 \begin{equation}
 \label{eq:1.6}
  G \simeq G_1 \ast \cdots \ast G_k
  \end{equation}
   be a free decomposition of $G$ such that
 $H \leqslant G_1$, $G_1$ is freely indecomposable modulo $H$, and $G_i$ is a freely
 indecomposable non-trivial group for $i \neq 1$. We deduce from (\ref{eq:1.6}) that
 $G = F_{R(S)}$ for some irreducible system $S(X) = 1$ over $F$ which splits as
 $S(X) = S_1(X_1) \cup \cdots \cup S_k(X_k)$; here $S_i(X_i)$ is a finite set of defining
 relations of $G_i$.  Notice that the system $S_1(X_1) = 1$ may contain constants from $F$,
 but $S_i(X_i) = 1$
 are coefficient free. Obviously, the system $S(X) = 1$ is rationally equivalent to
 $U(Z) = 1$ as well as to $U_H(Z,Y) = 1$.
Let $F= F(A)$, $F_i$ be an isomorphic copy of $F$ and $\alpha_i: F \rightarrow F_i$
   be a fixed isomorphism ($i = 2, \ldots,k$). Denote by $F^*$ the  free product
  $$F^* = F_1 \ast F_2 \ast\cdots \ast F_k$$
   in which $F_1=F$. For every $\mu\in Hom(G,F)$ there exists
   $\lambda _{\mu}\in Hom(G, F^*)$ defined for $x\in X_i$ as
   $\lambda _{\mu}(x)= \alpha_i(\mu |_{G_i}(x))$. Every solution
   of $S(X)=1$ in $F$ can be obtained by a substitution from a
   solution $\lambda$ such that
$\lambda (x)\in F_i$ for $x\in X_i$. Since $G$ is discriminated by
$F$ we can chose $\lambda$ such that
   $\lambda(x) \neq 1$ for any $x \in X$ (otherwise we can replace $X$ by a smaller set
   of generators for $G$).

  Denote by $\phi_{\lambda}$ the solution of the system $U_H = 1$
which is obtained from $\lambda$ by changing coordinates. Varying
the solution $\lambda$  we will get a set of solutions $$\Phi =
\{\phi_{\lambda} \mid \lambda \},$$ which discriminates $F_{R(S)}$
into $F^*$. We may assume that all solutions in $\Phi$ have
syllable length in $F^*$ not greater than a constant $K$ depending
on the expressions of generators $Z$ in terms of generators $X$.

Since the set ${\mathcal GE}(U_H)$ of generalized equations
corresponding to the system $U_H(Z,Y) = 1$ is finite, there exists
a generalized equation $\Omega_{\Phi} \in  {\mathcal GE}(U_H)$ and
the canonical  embedding $\gamma: F_{R(U_H)} \rightarrow
F_{R(\Omega _{\Phi})}$ such that some infinite  subset of $\Phi$
factors through $\gamma$ and discriminates $F_{R(U_H)}$.   For
simplicity we denote this subset again by $\Phi$ and also identify
$G \simeq F_{R(U_H)}$ with its image $\gamma(F_{R(U_H)})$ in
$F_{R(\Omega _{\Phi})}$.

 We claim that given $\Omega _{\Phi}$ one can effectively construct
  a generalized equation $\Omega$
 with the following properties:
  \begin{itemize}
  \item [a)]  there is an embedding $\phi: F_{R(\Omega _{\Phi})} \rightarrow F_{R(\Omega )}$;
  \item [b)] the set of  all closed sections of $\Omega$ is partitioned into a disjoint union of
 subsets $C_1, \ldots, C_k$ such that any base from $C_i$ has
 its dual base in $C_i$, which implies that
  $$F_{R(\Omega )} \simeq  F_{R(\Omega_1 )} \ast \cdots \ast F_{R(\Omega_k )}$$
  where $\Omega_i$ is the generalized equation situated over the section $C_i$;
  \item [c)]    $\phi(\gamma ( F _{R(S_i)}))\leqslant F_{R(\Omega _i)}$.
  \end{itemize}

Let $\hat T_{{\rm dec}}(U_H,\Omega_{\Phi})$ be the tree
constructed in Subsection \ref{3.2'}.  Since the tree $\hat
T_{{\rm dec}}(U_H,\Omega_{\Phi})$ is finite  there is a branch $B$
of this tree from the root $v_0$ to a leaf, say $w$, such that an
infinite discriminating subset of homomorphisms from $\Phi$ goes
along this branch. Again we denote this subset by $\Phi$. Let $e =
(v,u)$ be an edge in $B$, $\pi_{B,u}$ the canonical homomorphism
from $\Omega_{\Phi}$ into $F_{R(\Omega_u)}$ along the branch $B$,
$\pi^*_{B,u}$ the restriction of $\pi_{B,u}$ onto $F_{R(U_H)}$,
and $G(u) = \pi^*_{B,u}(F_{R(U_H)})$. When doing the Elimination
process along the branch $B$ the following cases may occur at the
vertex $v$.
 \begin{enumerate}
\item [1)]  Neither  a  quadratic section nor a periodic structure
 occurs in constructing the generalized equation $\Omega_v$ at $v$.
 \item [2)]   A  quadratic section or a periodic structure
 occurs in constructing the generalized equation $\Omega_v$ at $v$ and the
 splitting induced  on $G(v)$ from $F_{R(\Omega_v)}$ is insufficient,
 i.e., the homomorphism $\pi^*_{B,u} : F_{R(U_H)} \rightarrow G(v)$ is an isomorphism.
 \item [3)] A  quadratic section or a periodic structure
 occurs in constructing the generalized equation $\Omega_v$ at $v$ and the
 splitting induced  on $G(v)$ from $F_{R(\Omega_v)}$ is sufficient,
 i.e., the homomorphism $\pi^*_{B,u} : F_{R(U_H)} \rightarrow G(v)$ is not
 an  isomorphism.  Notice that this case can occur
   only if $u = w$.
   \end{enumerate}

Observe that given a branch $B$ one can  effectively construct the
homomorphism $\pi_{B,u}$ because the construction of the tree
$\hat T_{{\rm dec}}(U_H,\Omega_{\Phi})$ is effective. Also,
effectiveness of the homomorphism $\gamma$ implies that one can
effectively find the subgroup $F_{R(U_H)}$ in $\Omega_{\Phi}$
hence the homomorphism $\pi^*_{B,u}$. It follows that one can
effectively decide  which of the cases above occurs at a given
vertex $v \in B$. Indeed, if
 the case 1) occurs at $v$ then it is obvious from the process. If the case
1) does not occur then to distinguish the case 2) from the case 3)
one needs to check whether the homomorphism $\pi^*_{B,u}$
 is an isomorphism  or not. This
can be done effectively by Theorem \ref{le:0.2.}.

It follows that one can effectively find the set $V(3)$ of all the
vertices $v \in \hat T_{{\rm dec}}(U_H,\Omega_{\Phi})$ for which
the third case holds. There are two cases here. Considering these
cases we identify $G$ with its image $G(v)$.

Case 1. Suppose there is no vertex from $V(3)$  in the branch $B$.
In this event either the vertex $v=w$ is a leaf and it has type 2
(in the tree $T(\Omega _{\Phi})$) or $\Omega _w$ is obtained from
$\Omega _v$ by a second minimal relacement and $w$ has type 2.
Therefore, in the equation $\Omega_w$ all the bases are located in
the non-active part. Observe that by the construction of the tree
$\hat T_{{\rm dec}}(U_H,\Omega_{\Phi})$ the non-active part of
$\Omega _w$  does not have periodic structures. Hence it coincides
with the parametric part. There are several cases to consider.
Considering these cases we identify $G$ with its image $G(w)$.

Case 1.1. All the pieces of the bases corresponding to the
generators $Z$ of $G$ in $\Omega _{\Phi}$ are eventually
transferred to the parametric part of $\Omega_w$. In this case the
group $G$ is freely indecomposable modulo $H$. Indeed, observe
first that each homomorphism  $\phi$ from $\Phi$ maps all the
items from the parametric part of $\Omega_w$ to the subgroup $F_1$
of $F^*$ . This follows from the fact that $\phi(H) \leqslant F_1$
and the property that  solutions of generalized equations are
reduced as written, so the image under $\phi$ of a piece of a base
corresponding to a generator, say $h$ of $H$ occurs in
 $\phi(h)$ as a subword, so it again belongs to $F_1$. This
 implies that $\phi(G) \leqslant F_1$ for every $\phi \in \Phi$ - hence
 $G = G_1$, thus $G$ is freely indecomposable.

Case 1.2. Not all pieces of the bases corresponding to the
generators $Z$ of $G$ in $\Omega _{\Phi}$ are eventually
transferred to the parametric part of $\Omega_w$ (some of them
were deleted as matched bases). In this case $F_{R(\Omega_w)}$ is
a non-trivial free product of a free group $F_w$ generated by free
variables of $\Omega_w$ and the subgroup $H_w$ generated by the
items in the parametric part of $\Omega_w$.

 Case 1.2a. If $G$ is not elliptic in the splitting
 $F_{R(\Omega_w)} = F_w \ast H_w$ then it has a non-trivial
 induced decomposition and one can find it effectively by Theorem
 \ref{le:0.2.new}. This finishes the proof of the theorem in this
 case.

 Case 1.2b. If $G$ is conjugated into $F_w$ then by Theorem \ref{int} one can
effectively verify this and find a conjugator. In this event one
has an effective embedding of $G$ into a  free group $F_w$, so $G$
is free and using the  standard techniques  from free groups it
possible to find a basis of $G$, thus proving the theorem.

Case 1.2c. If $G$ is conjugated into $H_w$ (in this case $G
\leqslant H_w$) then the coordinate group of the generalized
equation $\Omega_{w,H}$  which is located on the parametric part
of $\Omega_w$ is a proper quotient of the group $F_{R(\Omega_w)}$.
In this case we have an effective embedding of $G$ into
$F_{R(\Omega_{w,H})}$ and induction finishes the proof.

 Case 2. Suppose now that there is a vertex $v$ as in the case 3) in the branch $B$.
The group $F_{R(\Omega_v)}$  admits a sufficient splitting $D$.
Recall that the set of the canonical automorphisms $A_D$ of the
splitting $D$ consists of canonical Dehn twists related to some
elementary ${\mathbb Z}$-splittings or abelian splittings from
$D$. Considering this case we again identify $G$ with its image
$G(v)$. Let $\sigma$ be a quadratic section or a section
corresponding to a periodic structure which occurs in $\Omega_v$.
If $\sigma$ is quadratic, we add  a non-active ``duplicate"
$\sigma^*$  of $\sigma$  to the generalized equation $\Omega_v$.
Namely, we  add a copy $\sigma^*$ of the section $\sigma$  to the
end of the equation $\Omega_v$, make $\sigma^*$  non-active, and
for each non-variable base  $b$ in $\sigma$ identify the duplicate
$b^*$ of $b$ in $\sigma^*$ with $b$ by remembering a graphical
equality $b=b^*$. If $\sigma$ corresponds to a periodic structure,
we take as $\sigma ^*$ a section defining the corresponding
extension of a centralizer, $\sigma ^*$ in this case is an abelian
section. Denote this equation  $\Omega _{v,\sigma}$. For a minimal
solution with respect to $A_D$ it will be transformed into some
$\Omega _{w,\sigma}$, and the equality $b=b^*$ will be transformed
into some expression without cancellation of $b^*$ as a product of
bases of $\Omega _w$. Construct a new ``bookkeeping" tree
$T_{book}(U_H,\Omega _{\Phi})$ the following way. First we
construct an auxiliary tree $T_{0,{\rm book}}(U_H,\Omega
_{\Phi})$. The initial vertex of this tree is the same as the
initial vertex of $\hat T_{{\rm dec}}(U_H,\Omega _{\Phi})$. All
the other vertices will be terminal vertices. Let $w$ be a
terminal vertex of $T_{{\rm dec}}(U_H,\Omega _{\Phi})$. Let $v$ be
a vertex that precedes $w$. To each such vertex $v$ we assign some
quadratic section  $\sigma$ for which a minimal solution with
respect to canonical automorphisms corresponding to this section
satisfies a proper equation or a section corresponding to a
periodic structure in $\Omega _v$.
   If $\sigma$ is a quadratic section, let $\Omega _v'$ be a generalized equation
  obtained from $\Omega _v$ by deleting
  $\sigma$ from $\Omega _v$. Notice that the complexity of $\Omega
_v'$ is strictly less than the complexity of $\Omega _v$. Let
$\Omega _{v,\sigma }'$ be obtained from $\Omega _{v,\sigma }$ by
deleting the quadratic section $\sigma $.
 If $\sigma $ is a section corresponding to a periodic structure, then we assign to the
terminal vertex of $T_{0,{\rm book}}(U_H,\Omega _{\Phi})$   the
generalized equation  $\Omega _{w,\sigma }$. If $\sigma $ is
quadratic, then we assign to the terminal vertex of $T_{0,{\rm
book}}(U_H,\Omega _{\Phi})$ the generalized equation $\Omega
_{v,\sigma }'$. To obtain $T_{\rm book}(U_H, \Omega _{\Phi})$ we
first construct $T_{0,{\rm book}}(U_H, \Omega _{\Phi})$. If for a
terminal vertex of $T_{0,{\rm book}}(U_H,\Omega _{\Phi} )$ $\sigma
$ corresponds to a periodic structure, then we paste to this
vertex $T_{0,{\rm book}}(U_H,\Omega _{w,\sigma})$. If a terminal
vertex of $T_{0,{\rm book}}(U_H,\Omega _{\Phi})$ has quadratic
section $\sigma $, then we paste to this vertex $T_{0,{\rm
book}}(U_H,\Omega _{v,\sigma }')$. We iterate this process. The
tree $T_{\rm book}(U_H, \Omega _{\Phi})$ is finite because on each
step of the construction we either decrease the complexity of the
generalized equation from the previous step or we do not increase
the complexity and the group generated by all the variables except
the ones on the sections $\sigma ^*$ is a proper quotient of such
group on the previous step.

If $ \Phi$ is constructed using infinitely many Dehn twists
related to an edge group $\langle u\rangle \leqslant
F_{R(\Omega_v)}$, then $u \in G_v$. Moreover, because of our
choice of
 the set of homomorphisms $\Phi$, the element $u$ has to belong to one of the free
 factors $\pi_B(G_i)$ of  $G_v$, which are homomorphic images of the initial free
 factors of $G$.
 Indeed, suppose
 that the cyclic  syllable length of $u$ in $G_v$ is more than 1, then, applying sufficiently
 many Dehn twists related to $u$, one can obtain an element of arbitrary large
 syllable  length. In this event, applying automorphisms from $A_D$, we can obtain
 a set of solutions with arbitrary syllable  length in $F$, which contradicts to
 the choice of $\Phi$.

Let $\sigma$ be a quadratic section for $F_{R(\Omega _v)}$  and
let $Q$ be the corresponding QH subgroup of $F_{R(\Omega _v)}$.
Denote by $\Omega _v'$ the generalized equation obtained  from
$\Omega _v$ by removing $\sigma$. Then the complexity of $\Omega
_v'$ is strictly less than the complexity of $\Omega _v$. By
induction we can decide if $F_{R(\Omega _v')}$ can be represented
as in b). Represent $F_{R(\Omega _v)}$ as a fundamental group of a
graph of groups with one vertex group $Q$ and the others subgroups
of $F_{R(\Omega _v')}$. Denote this decomposition by $D_1$. By
Theorem \ref{le:0.2.new} we can find effectively the induced
decomposition $D_1(G)$ of $G$.

The following cases occur.

Case 2.1. $G$ is conjugated into $F_{R(\Omega _v')}$. In this case
we apply induction.

Case 2.2.  $G$ is conjugated into $Q$, then $G$ is free.

Case 2.3. $G$ is not conjugated into $F_{R(\Omega _v')}$ and all
the conjugates of $Q^q$ have the intersection $Q^q\cap G$ of
infinite index in $Q^q$. Then by the corollary from Lemma
\ref{int},  $G$ is freely decomposable modulo $H$, and the
decomposition $D_1(G)$ which can be found effectively by Theorem
\ref{le:0.2.new}, is a free decomposition.

Case 2.4. Some conjugate of $Q^q$ has the intersection $Q^q\cap G$
of finite index in $Q^q$. Without loss of generality we can
suppose $Q_G=Q\cap G$ is of finite index in $Q$. As a group of a
surface  with boundaries, $Q_G$ is freely indecomposable modulo
the edge groups. Let $D_2(G)$ be a free decomposition of $G$
modulo $H$. Either $Q_G$ is a QH subgroup of one of the factors in
the decomposition $D_2(G)$ , or one of the edge groups in $Q_G$ is
hyperbolic in $D_2(G)$.  In the second case all  other edge groups
can be all together included into a basis of the free group $Q_G$,
they generate a free factor in the decomposition of $Q_G$. The
other factor in the free decomposition of $Q_G$ is a free factor
in a free decomposition of $G$. In this case a minimal solution of
$U_H=1$ with respect to canonical automorphisms of $Q_G$ does not
satisfy a proper equation and we can pass to this minimal
solution.

Obviously, all abelian subgroups of $G$ are elliptic in any free
decomposition of $G$.

To finish case 3) apply induction to the generalized equation in
the terminal vertex of $T_{0,{\rm book}}(U_H,\Omega _{\Phi})$.
Going along some branch of the tree $T_{{\rm book}}(U_H,\Omega
_{\Phi})$ we either obtain a free decomposition of $G$ as one of
the induced decompositions considered in the process or arrive to
a generalized equation, in which the set
 of closed active and parametric sections is partitioned into disjoint subsets of independent sections
 described above. Adjoining the  sections from $\sigma ^*$, we suppose that they all have
 syllable length one, so if the expressions for constant bases
 from $\sigma ^*$ include variables from
 several independent close sections on the next level of the tree $T_{{\rm book}}(U_H,\Omega _{\Phi})$,
 we join these sections into the same component.
Going from the leaf vertices of this tree to the initial vertex,
we obtain the generalized equation $\Omega$ with the properties
a)--c). If in a)--c) we have $k=1$ for all the branches of
$T_{{\rm book}}(U_H,\Omega _{\Phi})$, then $G$ is freely
indecomposable modulo $H$.
\end{proof}

\begin{cy}
\label{cy:Grushko}
 There is an algorithm which for every  finitely
generated fully residually free group $G$ finds its Grushko
decomposition.
\end{cy}

\begin{cy}   There is an algorithm which for every  finitely
generated fully residually free group $G$ and its finitely
generated subgroup $H$ finds a free decomposition of $G$ modulo
$H$ $$ G \simeq G_1 \ast \cdots \ast G_n$$ such that $H \leqslant
G_1$ and $G_1$ is freely indecomposable modulo $H$, $G_i$ is
freely indecomposable non-trivial group   for every $i = 2,
\ldots, n$.
\end{cy}

\begin{cy}\label{compat}
There is an algorithm which for every  finitely generated fully
residually free group $G$ and a family of finitely generated
subgroups $H_1,\ldots ,H_m$  finds a reduced compatible free
decomposition of $G$ modulo $\{H_1,\ldots ,H_m\}$ $$ G \simeq G_1
\ast \cdots \ast G_n$$ such that $H_1\leqslant G_1$.
 This decomposition of $G$ may be
trivial.
\end{cy}
\begin{proof} Notice that since conjugates of different factors in a free decomposition
intersect trivially, a subgroup can be conjugated only into one
factor in a given decomposition. We try all possible combinations
for subgroups that can be conjugated into the same factor as
$H_1$. Suppose they are $H_1,\ldots ,H_s$. Let $H_i$ be generated
by $h_{i1},\ldots ,h_{ij_i}$. We represent $h_{ik}=g_{i}^{-1}\bar
h_{ik}g_i$ for $k=2,\ldots ,j_i$ and put all the elements $h_{1k},
k=1,\ldots ,i_1',\ \bar h_{ik}, i=1,\ldots s,\ k=1,\ldots ,j_i$
into the parametric part of the initial generalized equations. We
include the case when some of $g_i$'s are trivial. As in the proof
of the theorem, we can find a free decomposition of $G$ such that
all $H_1,\ldots ,H_s$ are conjugated into the same factor, and
$H_1$ belongs to this factor. If the obtained decomposition is
non-trivial, apply induction.\end{proof}

\begin{theorem}
\label{th:free-group} There  is an algorithm which for every
finitely generated fully residually free group $G$ determines
whether $G$ is free or not, and if it is the algorithm finds a
free basis of $G$.
\end{theorem}
 \begin{proof} Let $G$ be a finitely generating fully residually
 free group. By Corollary \ref{cy:Grushko} one can effectively
 find a Grushko decomposition of $G$:
  $$G = G_1 \ast \cdots \ast G_m$$
so the factors $G_i$, which are given by finite sets of generators
$G_i = \langle Y_i \rangle$, are freely indecomposable. By Theorem
\ref{le:0.2.2} one can find finite presentations of the groups
$G_i$. Since the word problem is decidable in every group $G_i$
one can effectively verify whether each one of the groups $G_i$ is
abelian or not. If one of the groups $G_i$ is  non-abelian then
$G$ is not free. Otherwise, one has to check if all the abelian
groups $G_i$ are cyclic or not. This can be done algorithmically
since the finitely generated abelian groups $G_i$ are given by
their finite presentations. If one of them is not cyclic then $G$
is not free, and if all of them are cyclic then $G$ is a free
product of infinite cyclic groups, hence free. This proves the
theorem.
\end{proof}

\section{Homomorphisms of finitely generated  groups into fully residually free groups}
\label{se:homs}

\subsection{Hom-diagrams}
\label{subsec:hom-diag}

Let $G$ and $H$ be finitely generated groups. Recall that by
$Hom(G,H)$ we denote the set of all homomorphisms from $G$ into
$H$. If $G$ contains $H$ as a fixed distinguished subgroup then by
$Hom_H(G,H)$ we denote the set of $H$-homomorphisms from $G$ onto
$H$ (the set of retracts from $G$ into $H$). In this section  for
a free non-abelian group $F$ we describe the set $Hom_F(G,F)$ as
well as $Hom(G,F)$ in terms of hom-diagrams.

Let $T$ be a rooted directed tree with the set of vertices $V =
V(T)$ and the set of edges $E = E(T)$. Suppose that for every
vertex $v \in V$  there is a  group $G_v$ and a group of
automorphisms $A_v \leqslant {\rm Aut}(G_v)$ associated with $v$,
and with every directed edge $e = (u,v) \in E$ there is a
homomorphism $\pi_{u,v} : G_u \rightarrow G_v$ associated with
$e$. We  call this structure a  {\em hom-diagram}. Suppose a group
$G = G_{v_0}$ is assigned to the root vertex $v_0$ in $T$. Then
every branch (directed path) $B$ in $T$ from $v_0$ to a leaf
vertex $v_n$
 $$B: \ v_0 \rightarrow v_1 \rightarrow \ldots \rightarrow v_n$$
 determines  a set of homomorphisms $\Phi_B \subseteq Hom(G,G_v)$
 as follows.  Put
 $$A_B = A_{v_0} \times \cdots \times A_{v_n}.$$
If $\sigma = (\sigma_0, \ldots, \sigma_n) \in A_B$ then
 $$\phi_{B,\sigma} = \sigma_0  \pi_{v_0,v_1}  \sigma_2  \pi_{v_1,v_2}
  \ldots  \pi_{v_{n-1},v_n} \sigma_n$$
is a homomorphism from $Hom(G,G_{v_n})$. Set
 $$\Phi_B = \{\phi_{B,\sigma} \mid \sigma \in A_B \}.$$
Now assume that  a  group $H$ is given and for every leaf vertex
$v$ in $T$ the group $G_v$ associated with $v$ is a free product
$G_v = H_{v,1}\ast\ldots H_{v,s} \ast F_v$ of  (non-free)
subgroups $H_{v,i}$ of $H$ and a free group $F_v = F_v(X_v)$ with
basis $X_v$. We denote such diagrams by ${\mathcal C} = {\mathcal
C}(T,G,H)$ and refer to them as to $Hom(G,H)$-{\em diagrams}. If
$v$ is  a leaf vertex of $T$ then embeddings $\lambda _i: H_{v,i}
\rightarrow H$ together with a map $\xi : X_v \rightarrow H$ give
rise to the canonical homomorphism
$$\psi_{\lambda, \xi} : G_v \rightarrow H.$$
 It follows that the branch $B$, and a triple  $\sigma$,  $\lambda =\{\lambda _1,\ldots ,\lambda _s\},$ $ \xi$
 determine a homomorphism $\phi_{B,\sigma, \lambda, \xi}  = \phi_{B,\sigma}  \psi_{\lambda, \xi}$
which is composition of $\phi_{B,\sigma}$ and $  \psi_{\lambda,
 \xi}$:
 $$\phi_{B,\sigma, \lambda, \xi} : G \rightarrow  H. $$
The set $\Phi_{\mathcal{C},B}$ of all such homomorphisms is called
the {\it fundamental sequence } of the hom-diagram ${\mathcal C}$
determined by the branch $B$. Union of all fundamental sequences
of the hom-diagram ${\mathcal C}$ is denoted by
 $${\mathcal C}Hom(G,H) = \bigcup \{\Phi_{\mathcal{C},B} \mid B \ \mbox{is  a  branch  of} \ T\}, $$
 we  refer to it as to the set of all homomorphisms
from $Hom(G,H)$ determined by ${\mathcal C}$. If ${\mathcal
C}Hom(G,H) = Hom(G,H)$ then we say that the hom-diagram ${\mathcal
C}(T,G,H)$ is {\em complete}.

We say that a hom-diagram ${\mathcal C}$ is {\em canonical} if the
following conditions are satisfied:
 \begin{enumerate}
  \item [1)] all groups $G_v$ associated with vertices in $T$ are {\em
  fully residually $H$}, except, perhaps, for the group $G$ associated
  with the root vertex $v_0$;
  \item [2)] all homomorphisms $\pi_e$ associated with edges of $T$ are
  {\em surjective} and {\em proper} (i.e., with a non-trivial kernel);
  \item [3)] the groups of automorphisms $A_v$ are canonical with respect to
  some JSJ decomposition of free quotients in a Grushko decomposition of $G_v$ (see Section \ref{se:split}).
  \end{enumerate}

In the case where $H$ is a fixed subgroup of $G$ one can consider
only the set of retracts $Hom_H(G,H)$. In this event we assume
that in the corresponding hom-diagram ${\mathcal C}$  all
homomorphisms assigned to edges in $T$ are $H$-homomorphisms, the
vertex groups $G_v$ associated with  leaves $v$ in $T$ are
isomorphic to $H \ast  H_{v,1}\ast\ldots \ast H_{v,s} \ast F_v$,
where $H_{v_i}$ are (non-free) subgroups of $H$, and the mappings
$\lambda _i: H\ast H_{v_i} \rightarrow H$ are identical on $H$ and
embeddings on $H_v$.
 In this case
we say that ${\mathcal C}$ is $H$-{\em complete} or
$Hom_H(G,H)$-{\em complete} if $Hom_H(G,H)= {\mathcal
C}Hom_H(G,H)$.

We say that a hom-diagram can be found effectively if all groups
$G_v$ are given by some  finite presentations, for all groups
$A_v$  particular finite generated sets are given, all
homomorphisms $\pi_e$ and the automorphisms generating $A_v$ are
given by their images on the corresponding finite generating sets.

\subsection{Description of homomorphisms from a  finitely generated fully
residually free groups into a
free non-abelian group}
\label{3.6'}

Let $G$ be a finitely generated fully residually free group. In
this section we describe effectively the sets $Hom(G,F)$ and
$Hom_F(G,F)$ in terms of hom-diagrams.

\begin{theorem}
\label{th:hom-machine-frf} Let $G$ be a  group from $\mathcal F$
and $F$ a free group. Then:
 \begin{enumerate}
  \item [(1)]  there exists a complete canonical $Hom(G,F)$-diagram ${\mathcal C}$.
 \item [(2)] if $F$ is a fixed subgroup of $G$ then there exists an $F$-complete canonical
 $Hom_F(G,F)$-diagram ${\mathcal C}$.
 \end{enumerate}
Moreover,    the Hom-diagrams from {\rm (1)} and {\rm (2)} can be
found effectively.
\end{theorem}
\begin{proof}

Since $G$ is fully residually free, it can we represented as $G =
F_{R(U)}$ for the irreducible system $U = 1$ over $F$. We
construct an embedding tree $T_{emb}(G)$ by induction.
 Observe from the beginning that $T_{emb}(G)$  is a  rooted tree oriented from the root
 $v_0$.  We  associate the
group $G$ with $v_0.$  Let $\bar U_i, i=1,\ldots, t$ be the
systems corresponding to $U = 1$ from Theorem \ref{qq}. We add new
vertices $v_{11},\ldots ,v_{1t}$ and connect $v_0$ with each
$v_{1i}$ by  an edge $e_{1i}$. We associate  the group $F_{R(\bar
U_i)} \ast F(t_1,\ldots ,t_{k_i})$ with the vertex $v_{1i}$, and
the canonical epimorphism $$\pi_{1i}: F_{R(U)} \rightarrow
F_{R(\bar U_i)} \ast F(t_1,\ldots ,t_{k_i})$$
 with the edge $e_{1i}$,  $i=1,\ldots, t.$
 By Corollary \ref{cy:Grushko} one  can find a Grushko's
decomposition of $F_{R(\bar U_i)}$  effectively.  Notice, that
freely indecomposable non-abelian factors in the Grushko's
decomposition of $F_{R(\bar U_i)}$  again satisfy Theorem
\ref{qq}. This allows one to continue the construction of a tree
by induction.  The resulting rooted tree is denoted by $T_{\rm
emb}(G)$, we refer to it as to   the {\em embedding tree}. By
construction, every vertex in $T_{\rm emb}(G)$ has only finitely
many outgoing edges. We claim that this tree is finite. Indeed,
since the group $F$ is equationally Noetherian it follows that
there are no infinite sequences of proper surjective homomorphisms
of irreducible coordinate groups. Now the claim follows from
Koenig's lemma.

By construction $G$ is embedded into the NTQ group corresponding
to each branch of $T_{\rm emb}(G).$ The solution tree $T=T_{{\rm
sol}}(U)$ which in an underlying tree of the Hom-diagram for
$Hom(G,F)$ or $Hom _F(G,F)$ can be obtained from $T_{\rm emb}(G)$
by adding edges connecting $v_0$ to the initial vertices of
$T_{\rm emb}(F_{R(U\wedge V_i)})$ for each equation $V_i=1$ from
Theorem \ref{qq} and adding similar edges on each level of the
construction. The groups of automorphisms $A_v$ assigned to
vertices of $T_{{\rm sol}}(U)$ are canonical groups of
automorphisms corresponding to JSJ decompositions of these groups.
By Theorem \ref{qq} the Hom-diagram corresponding to $T_{{\rm
sol}}(U)$ is complete and canonical. It can be found effectively
by Theorem \ref{modulo}.\end{proof}

\subsection{Hom-diagrams of finitely generated
groups}

In 1987 Razborov proved the following remarkable result
\cite{Razborov}: for every finitely generated group $G$ and a free
 non-abelian group $F$  one can
effectively find a complete constructive hom-diagram ${\mathcal C}
= \mathcal{C}(T,G,F)$.
 In his diagram the associated groups were
not, in general, fully residually free (only residually free), the
homomorphisms were not necessary onto, and the groups of
automorphisms were not canonical. In the paper \cite{KMIrc} the
authors rewrote  the Razborov's elimination process which allowed
them to show that all the associated groups $G_v$ are subgroups of
the coordinate groups of NTQ-systems, hence fully residually free.
Still the diagram was far from being canonical.

The following is a rephrasing of Theorem 5 from the preprint
\cite{KM7}.

\begin{theorem}
\label{th:hom-machine} Let $G$ be a finitely generated group G and
$F$ a free group. Then:
 \begin{enumerate}
  \item [(1)]  there exists a complete canonical
 $Hom(G,F)$-diagram ${\mathcal C}$.
 \item [(2)] if $F$ is a fixed subgroup of $G$ then there exists an $F$-complete canonical
 $Hom_F(G,F)$-diagram ${\mathcal C}$.
 \end{enumerate}
Moreover,  if the group $G$ is finitely presented then  the
Hom-diagrams from {\rm (1)} and {\rm (2)} can be found
effectively.
\end{theorem}
\begin{proof} It follows directly from Theorem \ref{qq} and the
construction of $T_{{\rm sol}}$  in Subsection \ref{3.6'}.
\end{proof}

\subsection{Solutions of systems of equations in free groups and
hom-diagrams}

In this section we describe effectively all solutions of a finite
system of equations in a free non-abelian group $F$.

\begin{theorem} \label{th:12.3}
Let $S(X)=1$ be a finite system of equations and \[G=(F*F[X])/
ncl(S).\] Then one can effectively construct an $F$-complete
canonical $Hom(G,F)$-diagram ${\mathcal C}$ such that all
solutions of $S(X)=1$ in $F$ are exactly all $F$-homomorphisms in
${\mathcal C}$.\end{theorem}

\section{Free Lyndon length functions on NTQ groups.}\label{lef}

 \subsection{Lyndon length functions}

In paper \cite{Ly3} Lyndon introduced an abstract  notion of a
length function on a group, which gives an axiomatic approach to
Nielsen cancellation argument. In this section we briefly discuss
length functions on groups.

 Let $G$ be a group and $A$  an
ordered abelian group. Then a function $l: G \rightarrow A$ is
called a {\it (Lyndon) length function} on $G$ if the following
conditions hold:
\begin{enumerate}
\item [(L1)] $\forall\ x \in G:\ l(x) \geqslant 0$ and $l(1)=0$;
\item [(L2)] $\forall\ x \in G:\ l(x) = l(x^{-1})$; \item [(L3)]
$\forall\ x, y, z \in G:\ c(x,y) > c(x,z) \rightarrow c(x,z) =
c(y,z)$,

\noindent where $c(x,y) = \frac{1}{2}(l(x)+l(y)-l(x^{-1}y))$.
\end{enumerate}
In general $c(x,y) \notin A$, but $c(x,y) \in A_{\mathbb{Q}}$,
where $A_{\mathbb{Q}}$ is a $\mathbb{Q}$-completion of $A$, so in
the axiom (L3) we view $A$ as a subgroup of $A_{\mathbb{Q}}$.

Below we list several extra axioms which describe some special
classes of Lyndon length functions. Sometimes we write $x_1 \circ
\cdots \circ x_n$ instead of $x_1 \ldots x_n$ if $l(x_1 \cdots
x_n) = l(x_1) + \cdots + l(x_n)$.
\begin{enumerate}
\item [(L4)] $\forall\ x \in G:\ c(x,y) \in A.$
\end{enumerate}

\noindent A length function $l: G \rightarrow A$ is called {\it
free}, if it satisfies
\begin{enumerate}
\item [(L5)] $\forall\ x \in G:\ x \neq 1 \rightarrow l(x^2) >
l(x).$
\end{enumerate}
A group $G$ has a Lyndon length function $l: G \rightarrow A$,
which satisfies axioms (L4)-(L5) if and only if $G$ acts freely on
some $A$-tree. This is a remarkable result due to  Chiswell (see
\cite{Ch}).

A length function $l: G \rightarrow A$ is  called {\it regular} if
it satisfies the following {\it regularity} axiom:
\begin{enumerate}
\item [(L6)] $\forall x, y \in G, \exists u, x_1, y_1 \in G:$
$$x = u \circ x_1 \ \& \  y = u \circ y_1 \ \& \ l(u) = c(x,y).$$
\end{enumerate}
The regularity condition is crucial for Nielsen's cancellation
method, which is the base for generalized equations and
elimination process from Sections \ref{se:ge} and \ref{se:5} for
solving equations over $F$.

\subsection{Free regular Lyndon length functions on HNN-extensions}
\label{frllf}

 In this section we will
show how to extend a free regular length function given on a group
   $H \in {\mathcal F}$ onto an HNN extension of $H$ of a particular type in such a way that the extended length
 function is again free and regular.

Let $H$ be a group from $\mathcal F$, $\Lambda$ an ordered abelian
group, and $\ell_H :H \rightarrow \Lambda$ a free regular length
function on $H$.

Suppose that
 \begin{equation}
 \label{eq:HNNcyclic}
 G=\langle H,x\mid x^{-1}vx=u\rangle,
 \end{equation}
 where elements $u, v \in H$  satisfy the following conditions:
 \begin{enumerate}
 \item [1)]  $\ell_H(u) = \ell_H(v)$;
 \item [2)]  $u,v$ are cyclically reduced, i.e., $u^2 = u \circ u, v^2 = v \circ v$, and  not proper powers;
 \item [3)]  the subgroups $\langle u \rangle $ and $\langle v \rangle $ are conjugately separated in $H$, i.e.,
  $\langle u \rangle^h \cap \langle v \rangle = 1$ for every $h \in H$.
 \end{enumerate}

Define  a  function
 $$\ell_G :G \rightarrow \Lambda _0=\Lambda \oplus {\mathbb Z}$$
 as follows. Suppose $g \in G$ is given in its
reduced form (as an element of the  HNN extension $G$)
$$g=g_1x^{\alpha _1}g_2\cdots g_mx^{\alpha _m}g_{m+1}.$$
  For $n \in \mathbb{N}$  define
$$\ell _1(g,n)=\ell _H(g_1(v^nu^n)^{\varepsilon(\alpha_1)}g_2\cdots g_m(v^nu^n)^{\varepsilon (\alpha_m)}
g_{m+1})-m\ell _H(v^nu^n),$$
 where $\varepsilon(\alpha _k)= \left
\{ \begin{array}{rr}
1, & if \ \alpha _k >  0, \\
-1, & if \ \alpha _k <  0
\end{array} \right. $

It has been shown in  \cite{MRL,MRS} that for every $g\in G$ there
exists a positive integer $n_0=n_0(g)$ such that $\forall n
\geqslant n_0\ \ell _1(g,n)=\ell _1(g,n_0)$. Put  $\ell _1(g) =
\ell _1(g,n_0(g))$ and define
 $$\ell_G(g)=\left(\ell _1(g),\sum _{i=1}^m|\alpha _i|\right).$$
Then the function $\ell_G :G\rightarrow \Lambda \oplus {\mathbb
Z}$ is well-defined.

The following result follows from \cite{MRL} and  \cite{MRS}, a
more detailed version is contained in  \cite{KMS1}.

\begin{theorem}
\label{th:freelf} Let $H$ be a group from $\mathcal F$  with a
free  regular length function $\ell_H:H \rightarrow \Lambda$ and
 $$G=\langle H,x\mid x^{-1}vx=u\rangle$$
 an HNN extension of $H$ of the type (\ref{eq:HNNcyclic}). Then
 the function
  $$\ell_G :G\rightarrow \Lambda \oplus {\mathbb Z},$$
  constructed above,  is a free regular length function on $G$.
\end{theorem}

\subsection{Coordinate groups of regular quadratic equations}

Let $G$ be a finitely generated freely indecomposable fully
residually free group. Let $U=1$ be an irreducible system over
$F$, such that $G=F_{R(U)}$. Let $Q$ be a MQH subgroup of $G$.
Consider a relation
\begin{equation}\label{or} \prod _{i=1}^n[ x_i, y_i]\prod _{j=1}^k
c_j^{z_j}d^{-1}=1\end{equation} or
\begin{equation}\label{nonor}
\prod _{r=1}^nx_r^2\prod _{l=1}^k c_l^{z_l}d^{-1}=1\end{equation}
corresponding to $Q$. We define the basic sequence
$$\Gamma = (\gamma_1, \gamma_2, \ldots, \gamma_{K(m,n)})$$
 of canonical automorphisms of the group  $G$ corresponding to the
 vertex group $Q$ the same way as we did for the corresponding
 quadratic equation.  A family of solutions $\Psi=\{\phi\beta \mid \phi\in \Gamma _P\}$
 for $U=1$, where $\beta$ is an arbitrary  solution of $U=1$  in general position for $Q$
 and $\Gamma _P$ is a
 positive unbounded family of canonical automorphisms of $H$ is called
 {\em positive unbounded with respect to $Q$}.

Take a solution $\beta :G\rightarrow F$ in general position with
respect to $Q$. Denote by $S_{\beta}=1$ an equation
\begin{equation} \prod _{i=1}^n[ x_i, y_i]\prod _{j=1}^k
c_j^{\beta z_j}d^{-\beta}=1\end{equation} corresponding to the
subgroup $Q$. Take a positive unbounded family of solutions for
the equation $S_{\beta}=1$ and consider a generalized equation
$\Omega =\Omega _{S_{\beta}}$  for this family of solutions.

\begin{lemma} The group $F_{R(\Omega )}$ is a free product of a group
isomorphic to $F_{R(S_{\beta })}$ and a free group.\end{lemma}
\begin{proof}
By the implicit function theorem (\cite{KMLift}, Theorem A) there
is a surjective $F_{R(S_{\beta})}$-homomorphism $\theta:
F_{R(\Omega )}\rightarrow F_{R(S_{\beta})}$. We can represent
$F_{R(\Omega )}=F_{R(S)}*F_1$ for some quadratic equation $S=1$
and free group $F_1$, therefore $F_{R(S_{\beta})}$ is isomorphic
to a subgroup of finite index of $F_{R(S)}$, and the retract is an
isomorphism between $F_{R(S)}$ and $F_{R(S_{\beta})}$. Hence
$F_{R(\Omega )}$ in the variables of $\Omega $ is a free product
of a group isomorphic $F_{R(S_{\beta })}$ and a free group with
some generators $h_{i_1},\ldots ,h_{i_p}$. \end{proof}

Assign to the free generators values corresponding to a minimal
solution of $\Omega$. Denote by $\Omega '$ the new generalized
equation. Then the coordinate group of $F_{R(\Omega ')}$ is
isomorphic to $F_{R(S_{\beta})}.$
\begin{lemma} For any positive number $N$ there is a solution of
$\Omega '$ such that one pair of variable bases in longer than $N$
times the maximal length of other variable bases.\end{lemma}
\begin{proof}
Fix some solution of $\Omega '$. We can make one of the pairs of
variable bases of $\Omega '$ much longer than all the others by
pre-composing this solution with a power of a suitable
automorphism of $F_{R(\Omega ')}$.\end{proof}

Now we apply the entire transformation and bring $\Omega '$ into
the overlapping form $\Omega ''$ (see Definition \ref{overl}).
\begin{lemma}
Generalized equation $\Omega ''$ does not have open
boundaries.\end{lemma}
\begin{proof} If $\Omega ''$ had an open boundary, the coordinate
group would split as a free product modulo $F$. This is not the
case.\end{proof}

Equation $\Omega ''$ consists of one closed section and this
section is partitioned as $\mu\lambda _1\cdots\lambda _n=\Delta
(\lambda _{i_1})\cdots \Delta (\lambda _{i_s})\Delta (\mu )$ where
some of the bases are quadratic-coefficient bases. Since $\Omega
''$ does not have open boundaries, the group $F_{R(\Omega
'')}\approx F_{R(S_{\beta})}$ can be generated by $F$ and the
bases of $\Omega$ and has relations:
$$\mu\lambda _1\cdots\lambda _n=\Delta (\lambda
_{i_1})\cdots \Delta (\lambda _{i_s})\Delta (\mu ),\ \mu =\Delta
(\mu ),\lambda _i=\Delta (\lambda _i)$$ for each $\lambda _i$. It
is, therefore, generated by $F$ and the set of representatives of
dual bases and has in these generators just one defining relation
$$\mu\lambda _1\cdots\lambda _n=\lambda
_{i_1}\cdots \lambda _{i_s}\mu .$$ Denote by $F[\lambda ]$ a free
group generated by $F$ and all the bases $\lambda _i$. We proved
the following

\begin{lemma}
The group $F_{R(S_{\beta })}$ is an HNN extension of $F[\lambda ]$
with  stable letter $\mu$ and relation
$$\mu\lambda _1\cdots\lambda _n\mu ^{-1}=\lambda
_{i_1}\cdots \lambda _{i_s} .$$ \end{lemma}

\subsection{Free regular Lyndon length functions on coordinate groups of
regular quadratic equations}

In this section we will construct  a  free regular Lyndon length
function $\ell _{\beta}$  from $F_{R(S_{\beta})}$ into $Z^2$
corresponding to $\Omega ''.$ Here $Z^2$ is ordered
lexicographically from the right.

In the case of the group $F_{R(S_{\beta })}$ the function $\ell
_1$ from Section \ref{frllf} is just the ordinary length of a
reduced word in the free group $F[\lambda ]$. In particular for
each variable base $\lambda _i$ we have $\ell _{\beta} (\lambda
_i)=\ell _1(\lambda _i)=(1,0).$ The base $\mu$ is infinitely
longer than $\lambda _i$, namely $\ell _{\beta}(\mu )=(0,1).$
Futhermore, $\ell _{\beta}(\lambda _{i_1}\cdots \lambda
_{i_s})=\ell _{\beta}(\lambda _{1}\cdots \lambda _{n})$ and $\
c(\mu , (\lambda _{i_1}\cdots \lambda _{i_s})^k)= k\ell
_{\beta}(\lambda _{i_1}\cdots \lambda _{i_s}),$ $c(\mu ^{-1},
(\lambda _{1}\cdots \lambda _{n})^{-k})= k\ell _{\beta}(\lambda
_{1}\cdots \lambda _{n})$ for any positive number $k$.

In particular, we just proved the following result.

\begin{theorem}
\label{th:13.6} If $G$ is the coordinate group of a quadratic
equation, then there exists a free regular Lyndon length function
$\ell :G\rightarrow {\mathbb Z}^2,$ where ${\mathbb Z}^2$ is
ordered lexicographically.
\end{theorem}

Two elements $x$ and $y$ in a group $G$ with the length function
$\ell$ are said to be {\em commensurable} if there exist natural
numbers $m$ and $n$ such that $m\ell (x)\geqslant \ell (y)$ and
$n\ell (y)\geqslant \ell (x)$.

Notice, that  for all possible solutions $\beta$ in general
position with respect to $Q$ and coordinate groups
$F_{R(S_{\beta})}$ there is a finite number of possible
generalized equations in the overlapping form up to the values of
constant bases.

Let now $\mathcal Q=\{Q_1,\ldots ,Q_k\}$ be a family of all the
MQH vertex groups in the JSJ decomposition of $G$ (modulo
subgroups $K_1,\dots, K_s$). For each solution $\beta :
G\rightarrow F$ we can construct a group $F_{R(S_{\beta})}$  which
is the coordinate group of a quadratic system corresponding to
$\mathcal Q$ over $F$. Denote this system ${\mathcal Q}'$.
Starting with a positive unbounded family of solutions
corresponding to ${\mathcal Q}'$ we can similarly construct a
length function $\ell _{\beta}: F_{R(S _{\beta})}\rightarrow Z^2.$

Let $\mathcal A=\{A_1,\ldots ,A_l\}$ be a family of all the
abelian vertex groups in the JSJ decomposition of $G$ (modulo
subgroups $K_1,\dots, K_s$). For each group $F_{R(S_{\beta})}$
where the images of the edge groups adjacent to subgroups from
$\mathcal A$ are non-trivial, we can assign a group
$F_{R(P_{\beta})}$ which is obtained by an extension of
centralizers of these images by the groups from $\mathcal A$.
Denote this family in $F_{R(P_{\beta})}$ by ${\mathcal A}'$. If
$m$ is the maximum of ranks of the groups from ${\mathcal A}'$,
then one can construct a free regular Lyndon's length function $L
_{\beta}: F_{R(P _{\beta})}\rightarrow Z^m.$

\subsection{Free regular Lyndon length functions on NTQ groups}

In this section we construct   regular free length functions with
values in $\mathbb{Z}^n$  on a coordinate group of an NTQ system.

\begin {theorem}
\label{th:13.7} If $G$ is the coordinate group of an NTQ system,
then there exists a free regular Lyndon length function $\ell
:G\rightarrow {\mathbb Z}^t,$ where ${\mathbb Z}^t$ is ordered
lexicographically.
\end{theorem}
\begin{proof}

  Consider  a NTQ system
${\mathcal S}$ :
$$S_1(X_1,\ldots ,X_n)=1,$$
$$\vdots $$
$$S_n(X_n)=1.$$
Construct a solution $X=X_1\cup\cdots \cup X_n$ inductively.
 If $S_n=1$ is a regular quadratic equation, we
denote $\Theta _n=\phi \beta$, where $\beta$ is some solution in a
general position and $\phi$ belongs to a positive unbounded family
of automorphisms for $F_{R(S_n)}$. If $S_n=1$ is a system of the
type $U_{com}$ (see Definition 39 in \cite{Imp}), and
$U_{com}:\{[y_p,y_s]=1, [y_p,u]=1, p,s=1,\ldots ,r\}$, we take a
solution $y_1=u^N, \ldots ,y_r=u^{N^r}$ and also denote it by
$\theta _n$. Now, if
 $S_i=1$ is a regular quadratic equation, we denote $\Theta _i=\phi \Theta _{i+1}$, where
 $\phi$ belongs to a positive unbounded family
of automorphisms. If $S_i=1$ is a system of the type
$U_{com}:\{[y_p,y_s]=1, [y_p,u^{\theta _{i+1}}]=1, p,s=1,\ldots
,r\}$, we take a solution $y_1=u^{N\theta _{i+1}}, \ldots
,y_r=u^{N^r\theta _{i+1}}$ and also denote it by $\theta _i$.
Finally, $\theta _1$ is a solution for $X.$ A family of solutions
 $\Theta=\Theta _1$  is
called {\em positive unbounded} family of solutions for the NTQ
system. One can construct a generalized equation for a positive
unbounded family of solutions on $F_{R(S)}$ and then by induction
 construct a regular length function $\ell :F_{R(S)}\rightarrow
Z^t$.
\end{proof}

Every positive unbounded family  of solutions $\Theta$ is a
generic family of solutions for ${\mathcal S}.$ We recall the
definition of a generic family of solutions of a NTQ system.
\begin{definition}\label{generic}
A family of solutions $\Theta$ of a regular  quadratic equation
$S(X)=1$ over a group $G$  is called {\em generic} if for any
equation $V(X,Y,A)=1$ with coefficients in $G$ the following is
true: if for any solution $\theta\in\Theta$ there exists a
solution of $V(X^{\theta},Y,A)=1$ in $G$, then $V=1$ admits a
complete $S$-lift. If the equation $S(X)=1$ is empty
($G_{R(S)}=G*F(X)$) we always take as a generic family a sequence
of growing different Merzlyakov's words (defined in \cite{Imp},
Section 4.4).

Let $W(X,A)=1$ be a NTQ system that consists of equations
$S_1(X_1,\ldots ,X_n)=1,\ldots ,S_n(X_n)=1$. A family of solutions
$\Psi$ of  $W(X,A)=1$ is called {\em generic} if $\Psi=\Psi
_1\ldots \Psi _n,$ where $\Psi _i$ is a generic family of
solutions of $S_i=1$ over $G_{i+1}$ if $S_i=1$ is a regular
quadratic system, and $\Psi _i$ is a discriminating family for
$S_i=1$ if it is  a system of the type $U_{\rm com}$.
\end{definition}

\section{Effective construction of JSJ decompositions of groups from $\mathcal{F}$.}
\label{se:eff-JSJ}

In this section we prove the following result.

\begin{theorem} \label{modulo} There exists an algorithm to obtain a cyclic [abelian] $JSJ$
decomposition for a f.g. fully residually free group $G$ modulo a
given finite family of finitely generated  subgroups. The
algorithm constructs a presentation of $G$ as a fundamental group
of a JSJ decomposition.
\end{theorem}

\subsection{Construction of $T_{JSJ}$}\label{const}

In this subsection we will extend the tree $T(\Omega )$ to a tree
$T_{ext}(\Omega )$ which will serve as an auxiliary tree for the
construction of the tree $T_{JSJ}$. This tree will be used to
construct a JSJ decomposition of the group $G$.

We are looking for a JSJ decomposition modulo  subgroups
$F\leqslant K_1,\ldots ,K_s$ of  $G$. Let $K_i$ be generated by
$h_{i1},\ldots ,h_{ij_i}$. We represent $h_{ik}=g_{i}^{-1}\bar
h_{ik}g_i$ for $k=2,\ldots ,j_i$ and put all the elements $h_{1k},
k=1,\ldots ,i_1',\ \bar h_{ik}, i=1,\ldots s,\ k=1,\ldots ,j_i$
into the parametric parts of the initial generalized equations. We
include the case when some of $g_i$'s are trivial.

We do not consider in the construction of $T_{ext}(\Omega )$
vertices of type 1 as final vertices. In other words, at each
vertex $v$ of type 1 we glue a tree $T(\Omega _v)$ and iterate
this process. We introduce also a new case 12.1, which can have
place together with case 12 and produces an auxiliary edge.

{\it Case 12.1} Announce all the quadratic closed sections to be
non-active sections. The corresponding edge is auxiliary.

\begin{prop} Let ${\mathcal Q} =\{Q_1,\ldots ,Q_s\}$ and $\mathcal A=\{A_1,\ldots ,A_l\}$
be respectively all MQH subgroups and all abelian vertex groups in
an abelian JSJ decomposition of a freely indecomposable fully
residually free group $G$ {\rm (}modulo subgroups $K_1,\ldots
,K_s${\rm )}. There exists a generalized equation $\Omega$ for $G$
such that a discriminating family of F-solutions of $U=1$ onto
different groups $F_{R(P_{\beta} )}$ that take subgroups from
$\mathcal Q$ and $\mathcal A$ onto their images in ${\mathcal Q}'$
and ${\mathcal A}'$  factors through $\Omega$. Moreover, one can
effectively construct a finite subtree $T_{JSJ}(\Omega )$ of
$T_{ext}(\Omega)$ with the property that such a family of
solutions factors through $\Omega _v$ for a leaf vertex $v$
{\rm(}of type $2${\rm)} of $T_{JSJ}(\Omega )$.
\end{prop}

\begin{proof} When we consider a solution $\bar H$ of the
generalized equation $\Omega$, it is a solution in some group
$F_{R(P_{\beta})}$ which has a regular free length function $L
_{\beta}$. Since we consider here an arbitrary solution from the
family, we denote by $\ell$ the length function $L _{\beta}$ which
corresponds to the given solution. All these length functions are
similar.

If we have an infinite path in $T(\Omega )$
$$v_1\rightarrow v_2\rightarrow\dots\rightarrow v_r\rightarrow\ldots ,$$ then
by Lemma \ref{3.2} we have three possibilities: 1) a linear case,
2) a quadratic case, 3) case 15.

1) In this case $G$ must be conjugated into the coordinate group
of the kernel of  $F_{R(\Omega _{v_n})},$ and we can work with the
kernel of $\Omega _{v_n}$ instead of it.

2) In this case there is an auxiliary edge. Transformation of Case
12.1 and 15.1 decreases the complexity $\tau '$.

3) Since there is only a finite family of generalized equations
for $U=1$, there exists at least one of them such that  a
discriminating family of solutions onto groups $F_{R(P_{\beta })}$
that are retracts on $F$ and map $\mathcal Q$, $\mathcal A$ onto
${\mathcal Q}'$, ${\mathcal A}'$ factors through $\Omega$. Denote
this family by $\Psi$. Let ${\mathcal Q} _{\rm inf}$ be a family
of $QH$ vertex groups $Q$ of $F_{R(\Omega )}$ such that $Q^g\cap
H$ has infinite index in $Q^g$ for any $g\in F_{R(\Omega )}$. We
take for each solution in $\Psi$ an equivalent solution minimal
with respect to groups of automorphisms of ${\mathcal Q}_{\rm
inf}$ among solutions that extend to a solution of $\Omega $.

Denote by ${\mathcal Q}_{\rm fin}$ the set of QH vertex groups $Q$
of $F_{R(\Omega )}$ such that  for some $g\in F_{R(\Omega )}$, the
intersection $Q^g\cap G$ has a finite index in $Q^g$. Let a
quadratic section $[j_1,j_2]$ of $\tilde\Omega$ (see
transformation D3) correspond to a QH subgroup conjugated into one
of the subgroups in ${\mathcal Q}_{\rm fin}$.

Let $\bar H$ be a solution of $\Omega$. Denote by $d_1(\bar H)$
the length of $\bar H[j_1,j_2]$ and by
 $d_2(\bar H)$ the sum of the lengths of the quadratic-coefficient bases on
$[j_1,j_2]$. If   $d_2(\bar H)\geqslant (r,1)$ for some $r\in Z$,
then there exists a solution $\bar H^+$ of $\Omega$ equivalent to
$\bar H$ with respect to the group of automorphisms of $Q\cap H$
such that
\begin{equation}\label{normal}
d_1(\bar H)\leqslant f_1(\Omega)d_2(\bar H),
\end{equation}
 where $f_1(\Omega)$ is a function from Lemma 32. If $d_2(\bar
H)=(r,0)$ for some $r\in Z$, then there exists a solution $\bar
H^+$ of $\Omega$ equivalent to $\bar H$ with respect to the group
of automorphisms of $Q\cap H$ such that
\begin{equation}\label{ubnormal}
d_1(\bar H)\leqslant f_1(\Omega)(0,1),\end{equation} where
$f_1(\Omega)$ is a function similar to the one from Lemma 32,
maybe slightly increased.

For quadratic sections corresponding to ${\mathcal Q}_{fin}$ we
replace $\bar H$ by $\bar H^+$.
 Denote the new family of solutions of $U=1$ by $\Psi _{\rm min}$.
  We will show that every path in
$T(\Omega )$ for a family $\Psi _{\rm min}$ after finite number of
steps either does not have active bases or comes to a vertex $v$
such that $\Omega _v$ begins with an overlapping pair of bases,
and, therefore, has an outgoing auxiliary edge.

\begin{lemma}\label{15} Suppose for the subinterval $[1,\delta]$
of the interval $I$ of a generalized equation $\Omega $ and a
solution $\bar H$ of $\Omega $ we have $d_1(\bar
H[1,\delta])\leqslant f_1(\Omega )d_2(\bar H[1,\delta]).$ Let
$\Omega =\Omega _{v_1}$ and there is a path
\begin{equation}\label{44}
r=v_1\rightarrow v_2\rightarrow \ldots\rightarrow
v_m,\end{equation} where $tp(v_i)=15\ (1\leqslant i\leqslant m)$
in the tree $T(\Omega ),$ the leading base of\/ $\Omega _{v_i} \
(1\leqslant i\leqslant m) $ does not overlap with its double, and
the initial boundary of\/ $\Omega _{v_m}$ is to the left of the
boundary $\delta$ of\/ $\Omega$. Then this path cannot contain a
prohibited subpath.
\end{lemma}
\begin{proof} The proof repeats the proof of Proposition \ref{3.4} for a
generalized equation in a group with non-archimedian length
function. We consider a subpath (\ref{3.6}) corresponding to the
fragment
\begin{equation}\label{3.11'}
(\Omega _{v_1},\bar H^{(1)})\rightarrow (\Omega _{v_2}, \bar
H^{(2)})\rightarrow\ldots \rightarrow(\Omega _{v_m},\bar
H^{(m)})\rightarrow\ldots\end{equation} of the sequence
(\ref{3.8}). Here $v_1,v_2,\ldots,v_{m-1}$ are vertices of the
tree $T_0(\Omega)$, and for all vertices $v_i$ the edge
$v_i\rightarrow v_{i+1}$ is principal.

Let $\mu _i$ denote the carrier base of $\Omega _{v_i}$, and
$\omega =\{\mu _1,\ldots ,\mu _{m-1}\},$ and $\tilde\omega $
denote the set of such bases which are transfer bases for at least
one equation in (\ref{3.11'}). By $\omega _1$ denote the set of
such bases $\mu $ for which either $\mu$ or $\Delta (\mu)$ belongs
to $\omega\cup\tilde\omega $; by $\omega _2$ denote the set of all
the other bases. Let
$$\alpha (\omega)=\min(\min _{\mu\in\omega _2}\alpha (\mu),\delta).$$
 Let $X_{\mu}\doteq H[\alpha (\mu),\beta (\mu)].$ If $(\Omega
,\bar H)$ is a member of sequence (\ref{3.11'}), then denote
\begin{equation}\label{3.12'}
d _{\omega}(\bar H)=\sum _{i=1}^{\alpha (\omega )-1}\ell (H_i),
\end{equation}

\begin{equation}\label{3.13'}
\psi _{\omega}(\bar H)=\sum _{\mu\in\omega _1}\ell
(X_{\mu})-2d_{\omega}(\bar H).\end{equation}

Every item $h_i$ of the section $[1,\alpha (\omega)]$ belongs to
at least two bases, and both bases are in $\omega _1$, hence $\psi
_{\omega}(\bar H)\geqslant 0.$

Consider the quadratic part of $\tilde\Omega _{v_1}$ which is
situated to the left of $\alpha(\omega )$. The  solution $\bar
H^{(1)}$ is minimal with respect to the canonical group of
automorphisms corresponding to this vertex. We have
\begin{equation}
\label{3.14'} d_1(\bar H^{(1)})\leqslant f_1(\Omega
_{v_1})d_2(\bar H^{(1)}).
\end{equation}
Using this inequality we estimate the length of the interval
participating in the process $d_{\omega}(\bar H^{(1)})$ from above
by a product of $\psi _{\omega }$ and some function depending on
$f_1$. This will be inequality (\ref{3.19'}). Then we will show
that for a prohibited subpath the length of the participating
interval must be reduced by more than this figure (equalities
(\ref{3.27'}), (\ref{3.28'})). This will imply that there is no
prohibited subpath in the path (\ref{3.11'}).

Denote by $\gamma _i(\omega)$ the number of bases $\mu\in\omega
_1$ containing $h_i$. Then
\begin{equation}\label{3.15'}
\sum _{\mu\in\omega _1}\ell (X_{\mu}^{(1)})=\sum _{i=1}^{\rho}\ell
(H_i^{(1)}) \gamma _i(\omega),\end{equation} where $\rho =\rho
(\Omega _{v_1}).$ Let $$I=\{i\mid 1\leqslant i\leqslant\alpha
(\omega)-1 \&\gamma _i=2\}$$ and $$J=\{i\mid 1\leqslant
i\leqslant\alpha (\omega)-1 \&\gamma _i> 2\}.$$ By (\ref{3.12'})
\begin{equation}\label{3.16'}
d_{\omega}(\bar H^{(1)})=\sum _{i\in I}\ell (H_i^{(1)})+ \sum
_{i\in J}\ell (H_i^{(1)})= d_1(\bar H^{(1)})+\sum _{i\in J}\ell
(H_i^{(1)}).\end{equation} Let  $(\lambda ,\Delta(\lambda))$ be a
pair of quadratic-coefficient bases of the equation $\tilde \Omega
_{v_1}$, where $\lambda$ belongs to the nonquadratic part. This
pair can appear only from the bases $\mu\in\omega _1$. There are
two types of quadratic-coefficient bases.

{\em Type} 1. $\lambda$ is situated to the left of the boundary
$\alpha (\omega)$. Then $\lambda$ is formed by items $\{h_i\mid
i\in J\}$ and hence $\ell (X_{\lambda})\leqslant\sum _{i\in J}\ell
(H_i^{(1)}).$ Thus the sum of the lengths $\ell (X_{\lambda})+\ell
(X_{\Delta (\lambda)})$ for quadratic-coefficient bases of this
type is not more than $2n\sum _{i\in J}\ell (H_i^{(1)}).$

{\em Type} 2. $\lambda$ is
 situated to the right of the boundary $\alpha
(\omega)$. The sum of length of the quadratic-coefficient bases of
the second type is not more than $2\sum _{i=\alpha
(\omega)}^{\rho}\ell (H_i^{(1)})\gamma _i(\omega).$

 We have \begin{equation}\label{3.17'}
 d_2(\bar H^{(1)})\leqslant 2n\sum _{i\in J}\ell
 (H_i^{(1)})+2\sum _{i=\alpha (\omega)}^{\rho}\ell (H_i^{(1)})\gamma
_i(\omega).\end{equation} Now (\ref{3.13'}) and (\ref{3.15'})
imply
\begin{equation}\label{3.18'}
\psi _{\omega}(\bar H^{(1)}_i)\geqslant \sum _{i\in J}\ell
(H_i^{(1)})+\sum _{i=\alpha (\omega)}^{\rho}\ell (H_i^{(1)})\gamma
_i(\omega).\end{equation} Inequalities (\ref{3.14'}),
(\ref{3.16'}), (\ref{3.17'}), (\ref{3.18'}) imply
\begin{equation}\label{3.19'}
d_{\omega}(\bar H^{(1)})\leqslant \max \{\psi _{\omega}(\bar
H^{(1)})(2nf_1(\Omega _{v_1})+1),\ f_1(\Omega _{v_1})\}.
\end{equation}

 From the definition of Case 15 it follows that all the words
 $H^{(i)}[1,\rho _i+1]$ are the ends of the word $H^{(1)}[1,\rho _1+1]$,
 that is
 \begin{equation}\label{3.20'}
 H^{(1)}[1,\rho _1+1]\doteq U_iH^{(i)}[1,\rho _i+1].\end{equation}
 On the other hand bases $\mu\in\omega _2$  participate in these
 transformations neither as carrier bases nor as transfer bases; hence
 $H^{(1)}[\alpha (\omega ),\rho _1+1]$ is the end of the word $H^{(i)}[
 1,\rho _i+1]$, that is
 \begin{equation} \label{3.21'}
 H^{(i)}[1,\rho _i+1]\doteq V_iH^{(1)}[\alpha
(\omega ),\rho _1 +1].\end{equation} So we have
\begin{equation}\label{3.22'} d_{\omega}(\bar H
^{(i)})-d_{\omega}(\bar H ^{(i+1)})= \ell (V_i)-\ell
(V_{i+1})=\ell (U_{i+1})-\ell (U_{i})= \ell (X_{\mu
_i}^{(i)})-\ell (X_{\mu _{i}}^{(i+1)}).\end{equation} In
particular (\ref{3.13}),(\ref{3.22'}) imply that $\psi _{\omega
}(\bar H^{(1)})=\psi _{\omega }(\bar H^{(2)})=\dots =\psi _{\omega
}(\bar H^{(m)})=\psi _{\omega }.$ Denote  the number (\ref{3.22'})
 by $\delta _i$.

Let the path (\ref{3.6}) be $\mu$-reducing, that is either $\mu
_1=\mu$ and $v_2$ does not have auxiliary edges and $\mu$ occurs
in the sequence $\mu _1,\dots ,\mu _{m-1}$ at least twice, or
$v_2$ does have auxiliary edges $v_2\rightarrow w_1,\dots,
v_2\rightarrow w_k $ and the base $\mu$ occurs in the sequence
$\mu _1,\dots ,\mu _{m-1}$ at least $\max _{1\leqslant i\leqslant
k}s(\Omega _{w_i})$ times. Estimate $\ell (U_m)=\sum
_{i=1}^{m-1}\delta _i$ from below. First notice that if $\mu
_{i_1}=\mu _{i_2}=\mu (i_1< i_2) $ and $\mu _i\not =\mu$ for $i_1<
i< i_2$, then
\begin{equation}
\label{3.23'} \sum _{i=i_1}^{i_2-1}\delta _i\geqslant \ell
(H^{i_1+1}[1,\alpha (\Delta (\mu_{i_1+1}))]).
\end{equation}
Indeed, if $i_2=i_1+1,$ then $\delta _{i_1}=\ell
(H^{(i_1)}[1,\alpha (\Delta (\mu))]=\ell (H^{(i_1+1)}[1,\alpha
(\Delta (\mu))].$ If $i_2> i_1+1,$ then $\mu _{i_1+1}\not = \mu$
and $\mu$ is a transfer base in the equation $\Omega
_{v_{i_1+1}}.$ Hence $\delta _{i_1+1}+\ell (H^{(i_1+2)}[1,\alpha
(\mu)])=\ell (H^{(i_1+1)}[1,\alpha (\mu _{i_1+1})]).$ Now
(\ref{3.23'}) follows from
$$\sum _{i=i_1+2}^{i_2-1}\delta _i\geqslant \ell (H^{(i_1+2)}[1,\alpha (\mu)]).$$
So if $v_2$ does not have outgoing auxiliary edges, that is the
bases $\mu _2$ and $\Delta(\mu _2)$ do not intersect in the
equation $\Omega _{v_2}$; then (\ref{3.23'}) implies that

$$\sum _{i=1}^{m-1}\delta _i\geqslant \ell (H^{(2)}[1,\alpha (\Delta\mu _2)])\geqslant \ell (
X_{\mu _2}^{(2)})\geqslant \ell (X_{\mu}^{(2)})=\ell
(X_{\mu}^{(1)})-\delta _1,$$ which implies that
\begin{equation}\label{3.24'}
\sum _{i=1}^{m-1}\delta _i\geqslant\frac{1}{2}\ell
(X_{\mu}^{(1)}).\end{equation}

Suppose now that the path (\ref{3.6}) is prohibited; hence it can
be represented in the form (\ref{3.7}). From  definition
(\ref{3.13'}) we have $\sum _{\mu \in\omega _1}\ell
(X_{\mu}^{(m)})\geqslant \psi _{\omega}$; so at least for one base
$\mu\in\omega _1$ the inequality $\ell
(X_{\mu}^{(m)})\geqslant\frac{1}{2n}\psi _{\omega}$ holds. Because
$X_{\mu}^{(m)}\doteq (X_{\Delta (\mu)}^{(m)})^{\pm 1},$ we can
suppose that $\mu\in\omega\cup\tilde{\omega}.$ Let $m_1$ be the
length of the path $r_1s_1\ldots r_ls_l$ in (\ref{3.7}). If
$\mu\in\tilde{\omega}$ then  by the third part of the definition
of a prohibited path there exists $m_1\leqslant i\leqslant m$ such
that $\mu$ is a transfer base of $\Omega _{v_i}$. Hence, $\ell
(X_{\mu _i}^{(m_1)})\geqslant \ell (X_{\mu _i}^{(i)})\geqslant
\ell (X_{\mu}^{(i)})\geqslant \ell
(X_{\mu}^{(m)})\geqslant\frac{1}{2n}\psi _{\omega}.$ If
$\mu\in\omega$, then take $\mu$ instead of $\mu _i$. We proved the
existence of a base $\mu\in\omega$ such that
\begin{equation}\label{3.27'}
\ell (X_{\mu }^{(m_1)})\geqslant\frac{1}{2n}\psi
_{\omega}.\end{equation}
 By the definition of a prohibited path, the inequality
$\ell (X_{\mu}^{(i)})\geqslant \ell (X_{\mu}^{(m_1)})$
$(1\leqslant i\leqslant m_1),$ (\ref{3.24'}), and (\ref{3.27'}) we
obtain
\begin{equation}\label{3.28'}
\sum _{i=1}^{m_1-1}\delta _i\geqslant \max\{\frac{1}{20n}\psi
_{\omega},1\}(40n^2f_1+20n+1).
\end{equation}

By (\ref{3.22'}) the sum in the left part of the inequality
(\ref{3.28'}) equals $d_{\omega}(\bar H^{(1)})-d_{\omega}(\bar
H^{(m_1)});$ hence
$$d_{\omega}(\bar H^{(1)})\geqslant \max \left\{\frac{1}{20n}\psi _{\omega},1\right\}
(40n^2f_1 +20n +1) ,$$ which contradicts  to (\ref{3.19'}).

This contradiction was obtained from the supposition that there
are prohibited paths (\ref{3.11'}) in the path (\ref{3.8}). Hence
(\ref{3.8}) does not contain prohibited paths.

By Lemma \ref{15} in the case of inequality (\ref{normal}) we
cannot have a prohibited subpath.

Suppose we have inequality (\ref{ubnormal}) and the total length
of the section $[1,\alpha (\omega )]$ is greater than $(s,1)$ for
some $s\in Z$.  The only infinite path we can have is path
(\ref{44}).
 Instead
of inequality (\ref{3.19}) we now have
\begin{equation}\label{56}
d_{\omega}(\bar H^{(1)})\leqslant (0,1)2n(f_1(\Omega _{v_1})+1).
\end{equation}
 For a
$\mu$-reducing path we still have inequality (\ref{3.24}), and,
automatically (\ref{3.26}). In the absence of overlapping pair
bases there exists a base $\mu\in\omega$ such that
\begin{equation}
\label{66} \ell (X_{\mu}^{(m_1)})\geqslant (r,1)\end{equation} for
some $r\in Z$. By the definition of the prohibited subpath, we
obtain
\begin{equation}\label{67}
d_{\omega }(\bar H^{(1)})-d_{\omega}(\bar H^{(m_1)})=\sum
_{i=1}^{{m_1}-1}\delta _i\geqslant (1/10)(r,1) (40n^2f_1+20n+1)
\end{equation} which contradicts to (\ref{56}).
The Lemma is proved. \end{proof} The proposition is
proved.\end{proof}

\subsection{Periodic structures and overlapping pair sections}
\label{const1} Leaf vertices of $T_{JSJ}$ have type 2. We now
replace each length function
$L_{\beta}:F_{R(P_{\beta})}\rightarrow {\mathbb Z}^m$ by the
function $\ell _{\beta}:F_{R(S_{\beta})}\rightarrow {\mathbb Z}^2$
by specializing elements from $\mathcal A$. Namely, for each
abelian group $A_i$ instead of a homomorphism onto $A_i'$ we take
a positive unbounded family of homomorphisms onto $F$. Therefore,
solutions of generalized equations at final vertices of $T_{JSJ}$
will be now considered as solutions in $F_{R(S_{\beta})}.$ Like in
the previous subsection, if $H$ is a solution of a generalized
equation in $F_{R(S_{\beta})}$, then by $\ell (H_i)$ we denote
$\ell _{\beta }(H_i)$.
\begin{definition} A triple of elements $(A,T,A')$ of
$F_{R(S_{\beta})}$ forms an overlapping pair section for
$S_{\beta}=1$ if a quadratic generalized equation corresponding to
a quadratic equation  $S_{\beta}=1$ forms a closed section $A\circ
T=T\circ A'$,  two occurrences of $T$ in this equation overlap,
and $A$ is not a proper power.\end{definition}

 A generalized equation
$\Omega $ at a vertex of type 2 can be written as a union of
closed sections each having a principal overlapping pair of bases.

\begin{lemma} If $(B,T,B')$ is an overlapping pair section for a (regular) quadratic equation
$S_{\beta}=1$, then
\begin{enumerate}
\item for any $U,V\in F_{R(S_{\beta})}$ such that $T=U\circ V$,
either $V$ infinitely smaller {\rm(}$\ll${\rm)} than $T$ or $U\ll
T$,

\item In the case $U\ll T$, $U=B^n.$\end{enumerate}
\end{lemma}

\begin{definition}\label{11}
A solution ${H}$ of a generalized equation $\Omega$ is called {\em
periodic with respect to a period $P$,}  where $P$ is not a proper
power, if for every non-parametric closed section $[i,j]$
containing at least one base the word $H[i,j]$ can be represented
in one of the following forms:

\bi \item[1)]
$$H[i,j] =A^rA_1$$

\centerline{$(r \geqslant 2, A=A_1A_2,$ $A$ is a primitive
element.) }

\item[2)] $H[i,j]=A^rTA_1=TA'^rA_1$, where  $(A,T,A')$ forms an
overlapping pair section for some $S_{\beta}=1$, with $T$
infinitely larger than $A$ and $c (T, A^n)=\ell (A^n)$, $r\in
{\mathbb Z}$, and $A'=A_1\circ A_2.$ \ei
 In addition,  it is
possible to define a function $\ell _1$ in the definition of $\ell
_{\beta}$ such a way that $\ell _{\beta}(A) \leqslant \ell
_{\beta}(P)$ and for at least one such closed section $[i,j]$ the
element
 $A$
must be a cyclic shift of the element $P^{\pm 1}$.\end{definition}

In our new situation Lemmas \ref{2.9b}--\ref{2.11} still hold.
Instead of Lemma \ref{4n} we have the following

\begin{lemma}\label{4n'} Suppose that the generalized equation $\Omega $
is periodized with respect to a non-singular periodic structure
$\mathcal P$. Then for any periodic solution $H$ of $\Omega$ we
can choose a tree BT, some set of variables $S=\{h_{j_1},\ldots ,
h_{j_s}\}$ and a solution $H^+$ of $\Omega$ equivalent to $H$ with
respect to the group of canonical transformations $\bar A(\Omega )
$ such a way that each of the bases $\lambda _i\in
BT\smallsetminus BT_0$ can be represented as $\lambda_i=\lambda
_{i1}h_{k_i}\lambda _{i2}$, where $h_{k_i}\in S$ and for any
$h_j\in S,$ $\ell (H^+_{j})\leqslant f_3\ell (P)$, where $f_3$ is
some constructible function depending on $\Omega $. This
representation gives a new generalized equation $\Omega '$
periodic with respect to a periodic structure $\mathcal P'$ with
the same period $P$ and all $h_{j}\in S$  considered as variables
not from $\mathcal P'$. The graph $B\Gamma '$ for the periodic
structure $\mathcal P'$ has the same set of vertices as $B\Gamma
$, has empty set $C^{(2)}$ and either $BT'=BT'_0$ or the
difference between $BT'$ and $BT'_0$ consists of one infinitely
large edge.

Let $c$ be a cycle from $C^{(1)}$ of minimal length, $\ell
(c)=n_c\ell (P)$. Using canonical automorphisms from $A(\Omega )$
one can transform any solution $H$ of $\Omega$ into a solution
$H^{+}$ such that for any $h_j\in S,$ $\ell (H^{+}_j)\leqslant
f_3\ell (c).$ Let $\mathcal P '$ be a periodic structure, in which
all $h_i\in S$ are considered as variables not from $\mathcal P
'$, then $B\Gamma '$ has empty set $C^{(2)}$ and  either
$BT'=BT'_0$ or the difference  between $BT'$ and $BT'_0$ consists
of one infinitely large edge.
\end{lemma}

The proof is the same as the proof of Lemma \ref{4n}.

Notice that if a solution $H$ of a generalized equation is
periodic with respect to some period $P$, and $(B,T,B')$ forms an
overlapping pair section for some $S_{\beta}=1$, with $T$
infinitely larger than $B$, then either $B$ is a cyclic shift  of
$P$ or the length function $\ell _1$ in the definition of $\ell
_{\beta}$ can be redefined such a way that $\ell
_{\beta}(P)\leqslant\ell _{\beta}(B)$. Therefore this lemma can be
applied to each closed section that contains an infinitely long
overlapping pair of bases. Notice also that for any periodic
structure $\mathcal P$ the statements of  Lemmas \ref{2.9b}--\ref
{2n} are true because all the other periodic structures can be
temporary considered as parameters.

\subsection{Recognizing  QH  and abelian subgroups using
$T_{JSJ}$}\label{recog}

Let $v$ be a leaf vertex of $T_{JSJ}$ such that a discriminating
family of solutions of $U=1$ onto groups $F_{R(S_{\beta})}$
factors through $\Omega _v$.  A closed section from $\Omega _v$
either belongs to some periodic structure or can be considered as
an overlapping form for some quadratic generalized equation
$\Omega _{\rm quadr}$ with some extra bases situated on it.
  For each such
overlapping pair section
\begin{equation}\label{nu}
\nu _1\cdots\nu _s\Delta (\mu)=\mu\theta _1\cdots\theta
_k\end{equation} replace this section  by
 two closed sections
$\lambda _2\Delta (\mu _1)=\mu _1\Delta (\lambda _2)$ and $\lambda
_3\Delta (\mu _2)=\mu _2\Delta (\lambda _3)$. Put all the other
bases from the left part of the  overlapping pair section
(\ref{nu}) to $\lambda _2\Delta (\mu _1)=\mu _1\Delta (\lambda
_2)$  and all the other bases from the right part of the
overlapping pair section (\ref{nu}) to $\lambda _3\Delta (\mu
_2)=\mu _2\Delta (\lambda _3)$. Keep also the basic equations $\nu
_i =\Delta (\nu _i)$, where $\alpha (\nu _1)=\alpha (\lambda _2),\
\beta (\nu _s)=\beta (\lambda _2)$ and $\theta _i=\Delta (\theta
_i)$, where $\alpha (\theta  _1)=\alpha (\lambda _3),\  \beta
(\theta _k)=\beta (\lambda _3).$ Let $X_{\mu _1}$ and $X_{\mu _2}$
have minimal necessary lengths (defined by Lemma \ref{4n'}) to
accommodate all the extra bases for a minimal solution. Apply the
entire transformation to this new generalized equation. After a
bounded number of steps all the extra bases will be situated on
the bases that where quadratic-coefficient bases of generalized
equations $\Omega _{\rm quadr}$. Denote these new generalized
equations by $\Omega _{\rm quadr}'$. Replace in the generalized
equation $\Omega _v$ the closed sections that we have changed and
which were denoted by $\Omega _{\rm quadr}$ by the generalized
equations $\Omega _{\rm quadr}'$. Denote this generalized equation
by ${\Omega _v'}$.

\begin{prop} There is an algorithm to find all MQH and abelian
vertex groups in a JSJ decomposition of $G$ {\rm(}modulo subgroups
$K_1,\ldots ,K_s${\rm)}.\end{prop}

\begin{proof} If $F_{R(\Omega _v')}$ is a free product, then $G$ is
conjugated into a factor in this free decomposition.  We can
consider all possible periodic structures corresponding to each
vertex $v$ of type 2. All MQH subgroups of $G$ can be obtained as
$G\cap Q^g$, where $Q$ is a QH subgroup of $F_{R(\Omega _v')}$ for
some $v$ and $\Omega _v'$ and $G\cap Q^g$ has finite index in
$Q^g$. Then this subgroup $Q^g$ belongs to the factor in a free
decomposition of $F_{R(\Omega _v')}$ and we will consider this
factor $K=F_{R({\mathcal T})}$ instead of $F_{R(\Omega _v')}$. By
Lemma \ref{KMRS1} all intersections of conjugates of two subgroups
can be found effectively.

\begin{lemma}\label{QHS} Let $\mathcal Q$ be a family of  MQH vertex subgroups of a
freely indecomposable fully residually free group $G$ {\rm(}modulo
subgroups $K_1,\ldots ,K_s${\rm)} Then each $Q\in\mathcal Q$ is
conjugated to one of the QH subgroups for $K.$\end{lemma}

\begin{proof}
 For a discriminating family of
solutions $\Psi$ of $U=1$ which is positive unbounded with respect
to $\mathcal Q$ equation $\mathcal T$ has a solution. A
discriminating subfamily $\Psi _1$ of $\Psi$ which is also
positive unbounded with respect to $\mathcal Q$ factors through
one of the systems $U_i=1$ from Theorem \ref{qq}.  By the Implicit
function theorem (\cite{Imp}, Theorem 12, see also \cite{KMLift},
Theorem 10) $\mathcal T$ has a solution which takes subgroups from
$\mathcal Q$  to subgroups from $\mathcal Q'$ in some group
obtained from $F_{R(U_i)}$ by adding roots of a finite number of
elements from abelian vertex groups and by extending $F_{R(\bar
U_i)}$ by some group without sufficient splitting. Therefore
elements of QH subgroups of $K$ that have QH subgroups of $G$ as
finite index subgroups actually coincide with these QH subgroups
of $G$, and $Q=Q_i^g$.\end{proof}

 We constructed
$\Omega _v$ for a discriminated family of solutions of $U=1$ and
length functions $L_{\beta}$, the infinitely long variables cannot
disappear in the process of construction of $T_{JSJ},$ therefore
each abelian vertex group of $G$ will be seen in some periodic
structure of $F_{R(\Omega _v')}.$ Suppose an abelian group is
given by relations $[y_i,y_j]=[y_i,u]=1,\ i,j=1,\ldots ,r$. After
we replaced length function $L_{\beta}$ by $\ell _{\beta}$ we have
a solution $y_1=u^N,\ldots ,y_r=u^{N^r}$ for an increasing
sequence $\{N\}$. Therefore each abelian vertex group of $G$ is
conjugated into an abelian vertex group of $F_{R(\Omega _v')}$
that corresponds to a periodic structure. If for a periodic
structure the set of cycles $\bar c^{(2)}$ contains at least one
cycle, then by Lemma \ref{2.10''}, $\bar c^{(1)}$ and $\bar
c^{(2)}$ generate an abelian vertex group of $F_{R(\Omega _v')}$,
and each abelian vertex group of $G$ is conjugated into one of the
abelian groups obtained this way. The proposition is proved.
\end{proof}
\subsection{Recognizing rigid subgroups}

First we discuss how to construct a JSJ decomposition of $G$
modulo $F$.
 Let $D$ be a
cyclic [abelian] such JSJ decomposition of $G=F_{R(U)}$. Consider
a discriminating family of solutions for $U=1$ which for each edge
group $G_e$ of $D$  contains infinitely increasing powers of the
corresponding canonical Dehn' twists. If the edge group is
non-cyclic abelian with generators $y_1,\ldots ,y_k$, then we take
infinitely increasing powers of the  product of Dehn' twists
corresponding to $y_1^N,\ldots ,y_k^{N^k}.$ Each solution from
this family we pre-compose with a positive unbounded family of
solutions for each MQH subgroup of $G$ and a family of solutions
of the form $y_1=u^N,\ldots ,y_r=u^{N^r}$ for an increasing
sequence $\{N\}$ for each abelian vertex group of $G$ given by
relations $[y_i,y_j]=[y_i,u]=1,\ i,j=1,\ldots ,r$. All Dehn's
twists that we consider fix elements from $F$. For different edges
of $D$ corresponding canonical Dehn's twists commute, they also
commute with canonical automorphisms corresponding to MQH and
abelian vertex groups. Since all solutions of $U=1$ factor through
a finite number of generalized equations, there is a generalized
equation $\Omega$ such that a family of solutions with the
properties described above factors through this generalized
equation. In this case we will have a HNN splitting of
$F_{R(\Omega _v')}$ and, therefore, a splitting of $K$ that
induces the splitting of $G$ along $e$.
 We are able now to construct
effectively an abelian JSJ decomposition for $G$.  Considering
different possibilities for periodic structures we have different
decompositions for $F_{R(\Omega _v')}$. Each such decomposition
consists of several abelian and QH vertex groups, one non-abelian
non-QH vertex group, and some edges corresponding to HNN
extensions. A JSJ decomposition $D$ of $G$ is induced by one of
the decompositions of $F_{R(\Omega _v')}$, which we denote by $D
(\Omega _v')$. In addition we already know all QH and abelian
vertex groups for $D$. Denote by $S$ the  non-abelian non-QH
vertex group  in the decomposition  $D(\Omega _v')$, by $w_S$ the
corresponding vertex and by $t_1,\ldots ,t_k$ the stable letters
corresponding to edges with initial and terminal vertex $w_S$.  We
 know that each rigid subgroup in the JSJ decomposition of $G$ is the intersection of $G$
 with $S^{t_{i_1}\ldots t_{i_j}}$ where the product ${t_{i_1}\ldots t_{i_j}}$
 depends on the choice of the maximal subforest with rigid vertices and maximal subtree in
 the graph corresponding to the JSJ decomposition $D$ of $G$.
  Let $T$ be the Bass-Serre tree corresponding to
$D(\Omega _v')$. We can effectively construct the induced
decomposition of $G$ by Theorem \ref{le:0.2.new}. After we
obtained induced abelian decompositions for different leaf
vertices of $T_{JSJ}$ and for different periodic structures
associated with the same vertex, we can unfold them and we have to
decide which of the induced  decompositions is a JSJ decomposition
of $G$. We can compare QH and abelian vertex groups for different
induced decompositions and choose the ones with maximal  family of
such subgroups. We can transform all these decompositions into the
form
 with the same QH and abelian vertex groups.  If all the
 rigid vertex groups of
 such a decomposition are elliptic in all the other decompositions,
 then it is an abelian  JSJ decomposition. Denote it by $D$.
 This decomposition is unique up to slidings, conjugations, and
modifying boundary monomorphisms by conjugation. We can now
collapse all edges with non-cyclic edge groups and obtain a cyclic
JSJ decomposition.

Similarly we obtain an abelian JSJ decomposition modulo  subgroups
$K_1,\ldots ,K_s$.

To construct a cyclic [abelian] JSJ decomposition $D$ of $G$ (not
modulo $F$) we fix one of the rigid subgroups of $D$ (if exists)
and begin with a discriminating family of solutions that contains
increasing powers of Dehn's twists for all the edges of $D$, where
all Dehn's twists fix elements of  this subgroup. The rest of the
argument is the same. In the case when $D$ does not have rigid
subgroups, we begin  with a discriminating family of solutions for
$U=1$. Indeed, in this case all the edges have an abelian group as
one of the vertex groups, and, therefore, they ``will be seen'' in
some decomposition of some of the $F_{R(\Omega v')}.$

\section{Homomorphisms into NTQ groups}
\label{sec:NTQ}

The following result was proved in \cite{Kh}. It can be  obtained
by applying the elimination process to equations over an NTQ group
(with regular free Lyndon length function) as described in
Subsections \ref{const} and \ref{const1}.

\begin{theorem}
\label{th:hom-machineNTQ} Let $G$ be a finitely generated group
and $N$ an NTQ group. Then:
 \begin{enumerate}
  \item   there exists a complete canonical
 $Hom(G,N)$-diagram ${\mathcal C}$;
 \item  if $N$ is a fixed subgroup of $G$ then there exists an $N$-complete canonical
 $Hom_N(G,N)$-diagram ${\mathcal C}$.
 \end{enumerate}
Moreover,  if the group $G$ is finitely presented then  the
Hom-diagrams from {\rm (1)} and {\rm (2)} can be found
effectively.
\end{theorem}

This theorem implies the following result

\begin{cy}
\label{cor:15.2} Let $G$ be a finitely generated group and $H$ an
$\mathcal F$-group. Then:
 \begin{enumerate}
  \item   there exists a complete canonical
 $Hom(G,H)$-diagram;
 \item  if $H$ is a fixed subgroup of $G$ then there exists an $H$-complete canonical
 $Hom_H(G,H)$-diagram.
 \end{enumerate}
\end{cy}

\section{Some applications to equations in
$\mathcal F$-groups} \label{sec:mathcal F}

Consider monomorphisms from $\mathcal F$-group $G$ to $\mathcal
F$-group $H$. One can define an equivalence relation on the set of
all such monomorphisms:  monomorphisms $\phi$ and $\psi$ are
equivalent if $\psi$ is a composition of $\phi$ and conjugation by
an element from $H$.

\begin{theorem} \label{monom} Let $G$ {\rm[}respectively, $H${\rm]} be a non-abelian $\mathcal{F}$-group, and let $\mathcal A=\{A_1,\ldots
,A_n\}$ {\rm[}respectively, $\mathcal B=\{B_1,\ldots ,B_n\}${\rm]}
be a finite set of non-conjugated maximal abelian subgroups of $G$
{\rm [}respectively, $H${\rm ]} such that the abelian
decomposition of $G$ modulo $\mathcal A$ is trivial. The number of
equivalence classes of monomorphisms from $G$ to $H$ that map
subgroups from $\mathcal A$ onto conjugates of the corresponding
subgroups from $\mathcal B$ is finite. A set of representatives of
the equivalence classes can be effectively found.\end{theorem}

\begin{cy}\label{lem:BPproperty} Let $G$ be a non-abelian
$\mathcal{F}$-group, and let $\mathcal A=\{A_1,\ldots ,A_n\}$ be a
finite set of maximal abelian subgroups of $G$. Denote by
$Out(G;{\mathcal A})$ the set of those outer automorphisms of $G$
which map each $A_i\in {\mathcal A}$ onto a conjugate of it. If
$Out(G;{\mathcal A})$ is infinite, then $G$ has a non-trivial
abelian splitting, where each subgroup in $\mathcal A$ is
elliptic. There is an algorithm to decide whether $Out(G;{\mathcal
A})$ is finite or infinite. If $Out(G;{\mathcal A})$ is infinite,
the algorithm finds the splitting. If $Out(G;{\mathcal A})$ is
finite, the algorithm finds all its elements.
\end{cy}

Let $G,H,\mathcal A, \mathcal B$ be as in the formulation of the
theorem. Let \[G=\langle x_1,\ldots ,x_k\mid r_1(x_1,\dots
,x_k)=1, \ldots ,r_{\ell}(x_1,\ldots ,x_k)=1\rangle .\]  Let
$A_i=\langle a_{i1},\ldots , a_{ik_i}\rangle,$ $B_i=\langle
b_{i1},\ldots , b_{ik_i}\rangle $, $i=1,\ldots ,k_i.$ Let $C_i$ be
a cyclic group generated by $c_{i}$ if $A_i$ is cyclic, and free
abelian group generated by $c_{i1},\ldots ,c_{i,2k_i}$ otherwise.

\begin{lemma}\label{lem:JSJdouble} Let $G,H,\mathcal A, \mathcal B$ be as in the formulation of the
theorem, $C_1,\ldots,C_n$ be free abelian groups as above. Let
$$K=\langle G,H,C_1,\ldots, C_n,t_1,\dots,t_{n}\mid t_n=1,
a_{i1}^{t_i}=c_{i}=b_{i1}, \mbox{ if } \ k_i=1,\ $$ $$
a_{ij}^{t_i}=c_{i,k_i+j}, b_{ij}=c_{ij}, \mbox{ if } \ k_i\neq 1,\
i=1,\ldots ,n,\  j=1,\ldots, k_i \rangle .$$
 Then the reduced graph of groups
$\Delta$ that corresponds to the given presentation of $K$ is an
abelian JSJ decomposition of $K$ modulo  $H$.
\end{lemma}

\begin{proof} Let $\tilde{\Delta}$ be the Bass-Serre tree
corresponding to $\Delta$, and let $Y$ be a minimal $K$-tree.
Assume that $H\subset K$ is elliptic when acting on $Y$. It
follows that ${B}_i$ is conjugated to an elliptic subgroup,
therefore is an elliptic subgroup itself, for each $i$. Therefore
each $C_i$ and $A_i$ is an elliptic subgroup. Now, suppose that
$G$ is not elliptic when acting on $Y$. In the non-trivial
splitting that $G$ inherits from its action on $Y$, each subgroup
$A_i$ is elliptic, a contradiction to our assumption that the
decomposition of $G$ relative to $\mathcal{A}$ is trivial. Hence,
$G$ is elliptic when acting on $Y$. Therefore, the vertex $v$ (or
$w_i$) fixed by $G$ [or $C_i$] in $\tilde{\Delta}$ can be mapped
to a vertex fixed by $G$ [or $C_i$] in $Y$. If $y\in Y$ (or
${y_i}\in Y$) is a vertex fixed by $G$ [or ${C_i}$], then the edge
joining $v$ and $t_i.w_i$ can be mapped to the path joining ${y}$
and the vertex fixed by $t_iC_it_i^{-1}$. Thus, we have defined a
simplicial map from the extended fundamental domain of
$\tilde{\Delta}$ to $Y$. Extend this map equivariantly to obtain a
simplicial map from $\tilde{\Delta}$ into $Y$. Since $Y$ is a
minimal $G$-tree, this latter map is onto; the claim follows.
\end{proof}

\textsc{Proof of the theorem.}
 Consider the following system $W=1$ of equations in variables
$x_1,\ldots ,x_k,
 y_1,\ldots ,y_{n}$ over $H$:
 $$ r_1(x_1,\ldots ,x_k)=1,\ldots ,r_{\ell}(x_1,\ldots ,x_k)=1,
 a_{i1}^{y_i}=b_{i1},\verb if \  k_i=1,\ $$ $$
 [a_{ij}^{y_i},b_{is}]=1, \  \verb if \  k_i\neq 1,\ i=1,\ldots, n, y_n=1,$$
where $a_{i1},\ldots ,a_{ik_i}$ are generators of the  subgroup
$A_i$
 of $G$ written in variables $x_1,\ldots ,x_k$, and
  $b_{i1},\ldots ,b_{ik_i}$ are generators of the
 subgroup $B_i$
 of $H$ written in variables $h_1,\ldots ,h_m$. If the system $W=1$ is irreducible,
 then $K=H_{R(W)}$. If $W=1$ is not irreducible, then $V(W)$ is a union of irreducible
 subvarieties over $H$.
 Denote by $W_1=1,\ldots
 ,W_t=1$ the systems that define these subvarieties. Suppose that
 there is $W_j=1$ such that the canonical $H$-homomorphism from
 $H_{R(W)}$ onto $H_{R(W_i)}$ is a monomorphism on $G$ (we can verify this effectively).
 Since we
 will consider systems $W_j=1$ independently, we skip the index
 now and will write $W=1$ instead of $W_j=1$.
 Let the subgroups $C_i$ be mapped into free abelian subgroups
 $\bar C_i=\langle c_{i1},\ldots ,c_{it_i}\rangle$, and let
 $Z_i=\langle z_{i1},\ldots ,z_{ik_i}\rangle $ be the image of
 $A_i$ in $\bar C_i$. Notice that $Z_i$ and $B_i$ generate $\bar
 C_i$.
Then $H_{R(W)}$ has the following presentation:
 $$\langle G,H,\bar C_1,\ldots, \bar C_n,t_1,\dots,t_{n}\mid t_n=1,
a_{ij}^{t_i}=z_{ij}, b_{ij}=c_{ij}, i=1,\ldots ,n, j=1,\ldots,
k_i\rangle .$$

Indeed, this group is embedded into a series of extensions of
centralizers of the image of $K$ in $H_{R(W_j)}$, therefore it is
fully residually free and is itself $H_{R(W_j)}$.

Similarly to Lemma \ref{lem:JSJdouble} one can show that the
reduced graph of groups $\Delta$ that corresponds to the given
presentation of $H_{R(W)}$ is a JSJ decomposition of $H_{R(W)}$
modulo  $H$.

The group $H$ can be canonically embedded into an NTQ group $N$
(see  Section \ref{se:canonical}). One can define a free regular
Lyndon's length function from $N$ onto $\Lambda=Z^m$ which is also
a free length function on $H$ (Section \ref{lef}). Denote this
length function by $\ell$.

 Consider a parametric generalized equation $\Omega$ for
 $W=1$ with $H$-bases as parameters.
 We consider only those generalized equations $\Omega$ for which
$G\leqslant F_{R(\Omega)}\leqslant N_{R(\Omega )}$.
 Construct the tree $T(\Omega)$
  as described in Subsection \ref{se:5.2} with respect to the
  length function $\ell$. The tree will be the same as in
 Section \ref{const}.
 It is enough to trace possible paths in $T(\Omega )$
 for a minimal solution of $W=1$ such that the values of
 $y_i$ are commensurable with the value of  maximal elements in
 $B_i$. Indeed, when we take a minimal solution we only change values of $y_i$ and conjugate $G$.
  Denote the subtree corresponding to these paths by
 $T'(\Omega)$.
 In this case, constructing the subtree  $T'(\Omega)$, we do not have to consider infinite
 branches corresponding to Cases 15 and 12. Indeed, even if
 $F_{R(\Omega)}$ contains a QH subgroup $Q$, by Lemma \ref{le:1.6}, the group $H_{R(W)}$ is conjugated into some subgroup
 $P_i$ from Lemma 2.11, therefore we can consider only minimal
 solutions of $\Omega$ with respect to $Q$.
If, reducing $\Omega$ to the terminal equations we  obtain some
free variables because the boundary between the variables $h_i$
and $h_{i+1}$ does not touch any base, we can consider $ker
(\Omega)$ instead of $\Omega$. Indeed, in this case $\bar
H_{R(W)}$ is embedded into $F_{R(ker (\Omega ))}$. Therefore Cases
7--10 ( Subsection \ref{se:5.2}) can only appear a bounded number
of times in the subtree  $T'(\Omega)$ corresponding to the actual
paths for a minimal solution. Hence this subtree is finite.

In the leaf vertices $v$ of $T'(\Omega)$
 we obtain generalized equations of four types:
\begin{enumerate}
 \item [1)] generalized equations $\Omega _v$ with intervals labelled by generators $h_1,\ldots ,h_m$
 of $H$;

\item [2)] generalized equations $\Omega _v$ such that
$F_{R(\Omega _v)}$ is a proper quotient of $F_{R(\Omega )}$, and
the image of $G$ in $ F_{R(\Omega _v)}$ is a proper quotient of
$G$;

\item [3)] generalized equations $\Omega _v$ such that
$F_{R(\Omega _v)}$ is a proper quotient of $F_{R(\Omega )}$, and
the image of $G$ in $ F_{R(\Omega _v)}$ is isomorphic to $G$;

\item [4)]  generalized equations $\Omega _v$ with vertex $v$ of
type 2, such that the solution can be taken minimal with respect
to all the periodic structures.
\end{enumerate}
Therefore in the terminal vertices of $T'_{{\rm dec}}(\Omega )$ (
Subsection \ref{5.5.5}) we only have Cases (1)--(3). In Case 1) we
solved the  generalized equation $\Omega $ in $N$, and we just
have to check whether the image of $G$ belongs to $H$ and is
isomorphic to $G$.

In  Case 2) we cannot have a monomorphism from $G$ to $H$, so we
do not continue. In Case 3) we apply the leaf-extension
transformation ( Subsection \ref{se:max-standard}) at this leaf
vertex. Let $\bar K$ be the image of $H_{R(W)}$ in $F_{R(\Omega
)}.$ The group $\bar K$ is freely indecomposable modulo $H$.
Moreover, a JSJ decomposition of $\bar K$ modulo $H$ either has
two vertices stabilized by $G$ and $H$ and vertices with abelian
vertex groups that are images of $\bar C_1,\ldots ,\bar C_n$ or
has one vertex with the stabilizer containing $H$ and the image of
$G$ and some vertices with abelian vertex groups.  Therefore, we
again only consider minimal (and commensurable with elements of
edge groups in the decomposition of $\bar K$) solutions of $\Omega
_v$ with respect to the groups of canonical automorphisms.  The
number of times when we have case 3) is finite. Therefore after a
finite number of steps we either show that there are no solutions
of $W=1$ which are monomorphisms on $G$ or end up with Case 1). We
now check which solutions of $W=1$ in $N$ belong to $H$. This can
be done effectively by Theorem \ref{th:0.4}. \hfill $\Box$

Applying the elimination process to a finite system of equations
and inequalities (with constants) in an NTQ group  one can also
prove the following result.
\begin{theorem} There is an algorithm to solve finite systems of
equations and inequalities in an NTQ group. The universal theory
(with constants) of an NTQ group is decidable.\end{theorem}

\end{document}